\documentclass[pdflatex,sn-mathphys-ay]{sn-jnl}
\usepackage[T1]{fontenc}
\usepackage[utf8]{inputenc}
\usepackage{graphicx,amsmath,amssymb,booktabs,longtable,array,xcolor}
\usepackage{listings}
\usepackage{algorithm,algpseudocode}
\usepackage{float,placeins}
\usepackage{tikz}
\usetikzlibrary{arrows.meta,calc}

\newcommand{\R}{\mathbb{R}}
\newcommand{\catlabel}[1]{\ifcsname catlabel@#1\endcsname\csname catlabel@#1\endcsname\else\texttt{#1}\fi}
\expandafter\def\csname catlabel@variablePlayers\endcsname{\texttt{variable-}$N$}
\newcommand{\cat}[1]{\hyperref[cat:#1]{\catlabel{#1}}}
\newcommand{\caseref}[1]{\hyperref[case:#1]{\texttt{#1}}}

\newcommand{\code}[1]{\texttt{\detokenize{#1}}}
\newcommand{\dt}{\,\mathrm{d}t}
\DeclareMathOperator*{\argmin}{arg\,min}
\newcommand{\myheader}[1]{\par\smallskip\noindent\textbf{#1}\;}
\newcommand{\casefigheight}{.66\textheight}
\allowdisplaybreaks

\hypersetup{hypertexnames=false,bookmarksdepth=2}

\lstdefinelanguage{Julia}{morekeywords={function,end,using,include,for,in,if,else,elseif,return,struct,const,module,begin,let,do,true,false,nothing},
 morekeywords=[2]{Game,Options,StateConstraint,run_game,main,create_synthetic_data,collect,zeros},
 sensitive=true,morecomment=[l]{\#},morestring=[b]"}
\definecolor{lstkeyword}{RGB}{0,90,160}
\definecolor{lstfunction}{RGB}{120,40,140}
\definecolor{lstcomment}{RGB}{90,110,90}
\definecolor{lststring}{RGB}{170,60,20}
\definecolor{lstbackground}{RGB}{247,247,247}
\lstdefinestyle{numbered}{numbers=left,numberstyle=\tiny\color{black!55},numbersep=5pt,xleftmargin=2em,
 extendedchars=true,literate={₁}{{$_1$}}1 {₂}{{$_2$}}1}

\begin{document}

\title[CorleoneGame: A library of biological dynamic games]{CorleoneGame: A library of biological dynamic games with a companion Julia solver}

\author*[1,2]{\fnm{Sebastian} \sur{Sager}}\email{sager@ovgu.de}
\author[1,2]{\fnm{Christoph} \sur{Plate}}
\author[1]{\fnm{Julius} \sur{Martensen}}
\author[1,3]{\fnm{Frank W.} \sur{Ohl}}
\affil[1]{\orgname{Otto von Guericke University}, \orgaddress{\city{Magdeburg}, \country{Germany}}}
\affil[2]{\orgname{Max Planck Institute for Dynamics of Complex Technical Systems}, \orgaddress{\city{Magdeburg}, \country{Germany}}}
\affil[3]{\orgname{Leibniz Institute for Neurobiology}, \orgaddress{\city{Magdeburg}, \country{Germany}}}

\abstract{
Biological agents, from regulatory modules and microbial strains to organisms and populations, act on a common environment while each responds to its own costs, benefits, and limits. A dynamic game makes such interactions explicit by assigning who controls which decision, what each participant values, and which restrictions constrain unilateral deviations.
 
We present a library of 25 biological dynamic games with 2 to 10 players spanning molecular, cellular, organismal, and population scales and competitive, antagonistic, mutualistic, and regulatory interactions. Instances cover single seasons, cycles that return to their initial state, or episodic tasks; seven include duration choice and some players' feasible sets are coupled.
For each game we provide biological motivation, an executable Julia specification with player-owned controls, objectives and restrictions, a parameter instance, and a reference open-loop generalized Nash equilibrium candidate with interpreted trajectories.

The solver CorleoneGame.jl computes candidates by sequential optimal-control best responses and audits them for feasibility and profitable unilateral deviations on a coarse and a refined control grid. All 25 reference computations satisfy the audit tolerances on both grids.
Using a plant--herbivore game as a running example, we show how an executable specification yields an audited reference candidate that qualitatively reproduces a published prediction, how biological hypotheses become instance variants with different candidates, and how the candidates generate synthetic data for inverse games.
}

\maketitle

\section{Introduction}\label{sec:introduction}
Biological systems adapt to their environment, but also change it. A mitochondrial DNA variant that expands within a cell draws on the energy shared with the wild type; a bacterial strain that secretes an amino acid feeds a partner strain that competes with it for carbon; a person moving a handle together with a partner changes the forces the partner feels; and a tree that defends its leaves changes the food available to hatching leafminer larvae. 
Describing any of these actions in isolation can miss the strategic interaction. A dynamic game in continuous time, also called a differential game \citep{BasarOlsder1998}, makes the alternatives explicit: it specifies the interacting states, the decisions assigned to each participant, the consequences each participant values, the restrictions on its decisions, and when the interaction ends or repeats.

\paragraph*{Biological and mathematical context.}
Optimality models connect mechanisms with proposed criteria and feasible alternatives \citep{Parker1990-zj}. Examples include foraging decisions \citep{StephensKrebs1986} and, at the molecular scale, temporal enzyme allocation, a classical optimal-control problem \citep{Klipp2002-sg,Oyarzun2009-ef}; dynamic optimality models are established in behavioral biology \citep{mcnamara2001optimality}. When one organism changes another's returns, game theory provides a language for interaction \citep{MaynardSmithPrice1973,Smith1982}, and game-theoretic reasoning has also been applied to biochemical and cellular processes \citep{Schuster2008,hummert2014evolutionary}. Biological dynamic games already describe seasonal plant defense and herbivore hatching \citep{Low2013ColdWar}, renewable-stock harvesting \citep{Dockner1989Fishery,Suri2008Harvesting}, mechanically coupled human movements \citep{Braun2009Motor,Chackochan2019JointAction}, and territorial space use \citep{HamelinLewis2010Territoriality}. Other game instances in the library adapt dynamical mechanisms, such as Wilson--Cowan neural populations \citep{WilsonCowan1972}, by introducing separate controls and objectives. Appendix~\ref{app:cases} records this provenance for each game instance.

The games cited above use different notation, control ownership, feasibility conventions, and levels of numerical documentation. A common library makes them accessible to biological interpretation and computational comparison.
General-purpose resources such as GAMUT \citep{Nudelman2004} and OpenSpiel \citep{Lanctot2019} illustrate the value of shared game families and evaluation protocols; they address normal-form and extensive-form games, mostly without biological interpretation. SpaceGym collects discrete and differential games from noncooperative space operations, such as pursuit--evasion, for evaluating autonomy algorithms \citep{Allen2023SpaceGym}. Optimal-control collections provide a complementary foundation for a single decision maker \citep{Sager2012-ts,OptimalControlProblemsRepository}. We are not aware of a collection of dynamic games providing biological interpretations, control ownership, individual objectives and restrictions, and open-source code.

\paragraph*{Game-theoretic setting.}
The library uses an open-loop information structure: each player commits to a time-dependent schedule of its decisions at the start of the interaction \citep{BasarOlsder1998}. Shared resources, joint floors, and periodic return bands put opponents' decisions into each player's feasible set; the equilibrium concept is therefore a generalized Nash equilibrium \citep{FacchineiKanzow2010}. Such equilibria are typically not isolated (Section~\ref{sec:equilibrium}). Feedback equilibria, in which players respond to the current state \citep{BasarOlsder1998,FridovichKeil2020}, evolutionarily stable strategies \citep{MaynardSmithPrice1973,Smith1982}, and Stackelberg equilibria are different concepts not covered by our library (Section~\ref{sec:assumptions}).

\paragraph*{Forward and inverse methods.}
Constrained games can be approached through best-response decomposition, in which each player's optimal reply to fixed opponents, its best response, is computed in turn, coupled optimality systems, or variational inequalities \citep{FacchineiKanzow2010}. Augmented Lagrangian methods \citep{KanzowSteck2016}, ALGAMES \citep{LeCleach2022ALGAMES}, and DG-SQP \citep{ZhuBorrelli2023} offer alternatives to the sequential best-response baseline used here. This baseline reuses any optimal-control solver as a black box but has no general convergence guarantee (Section~\ref{sec:methods}); its results are therefore reported with independent feasibility and deviation audits. Feedback-game methods such as iterative linear--quadratic games address a different information structure \citep{FridovichKeil2020}. A variety of different optimal control algorithms exists \citep{Andersson2018,Verschueren2022,PattersonRao2014}.

For inverse analysis, objectives and constraints deserve the same attention as differential equations. Low investment may reflect a small benefit, a large cost, or a binding budget. Generalized inverse optimal control treats these alternatives as model components to be inferred \citep{Banga2025}. \citet{BangaVillaverde2026Hierarchical} order biological inverse problems from parameter identification and model discovery to inverse optimal control and inverse differential games; for the latter, inverse problems have been studied for open-loop equilibria with known dynamics and constraints \citep{Molloy2020-vg,molloy2022inverse}. The library supplies forward instances for the two upper levels (Section~\ref{sec:tupelo-inverse}); inverse recovery is left to future work.

\paragraph*{Contributions and organization of the paper.}
Section~\ref{sec:formulation} defines the open-loop generalized Nash equilibrium problem and a game specification consisting of mechanism, preferences, and restrictions. Section~\ref{sec:assumptions} states general modeling assumptions. 

Section~\ref{sec:cases} presents 25 game instances with biological provenance, classifications (Tables~\ref{tab:categories} and~\ref{tab:overview}), and dimensions (Table~\ref{tab:reference-results}). They span four biological scales, four interaction types, and three temporal classes, and they have 2 to 10 players, 2 to 12 controls including owned durations, and up to 420 scalar constraint functions. 
Section~\ref{sec:glance} summarizes every instance by biological scale, and Appendix~\ref{app:cases} contains all 25 full specifications in alphabetical order. 

Section~\ref{sec:methods} presents the novel modular open-source solver \code{CorleoneGame.jl} that extends the optimal control solver \code{Corleone.jl} to dynamic games. For every instance, it provides a reference open-loop generalized Nash equilibrium candidate that passes the feasibility and unilateral-deviation audits on a coarse and a refined control grid.

Section~\ref{sec:uses} works through the tupelo--leafminer game, from the executable specification through an audited reference candidate and alternative hypotheses, formulated as instance variants, to synthetic observations for inverse games. Along the way, it addresses teachers and contributors of new instances, algorithm developers, biologists, and researchers in inverse games. Section~\ref{sec:outlook} discusses limitations and next steps.

\section{Game formulation and biological assumptions}\label{sec:formulation}
\subsection{Players, decisions, and dynamics}
The $N$ players are labeled $\mathcal P=\{A,B,\ldots\}$ with $|\mathcal P|=N$. On a time horizon $[0,T]$, the biological state $x(t)\in\R^{n_x}$ satisfies the initial value problem
\begin{equation}\label{eq:dynamics}
 \dot x(t)=f(t,x(t),u(t);p),\qquad x(0)=x^0(p).
\end{equation}
for a given function $f$ describing the dynamics and initial values $x^0(p)$.
Here $p$ denotes model parameters and fixed data, and the explicit dependence on $t$ represents prescribed environmental inputs, such as seasonal or periodic supply, light, or activity drives. The vector $u=(u_i)_{i\in\mathcal P}$ collects player-owned decisions. Ownership identifies who may change a decision. Time-varying controls are measurable schedules with componentwise bounds
\begin{equation}\label{eq:control-set}
 U=\prod_{i\in\mathcal P}[\underline u_i,\overline u_i],\qquad u(t)\in U\text{ for almost every }t\in[0,T].
\end{equation}
In a slight abuse of notation, $u_i$ also contains bounded time-independent decisions owned by player $i$; for these, the bounds in \eqref{eq:control-set} apply to a constant value. The library uses such a scalar decision only for the duration $T$, owned by player A in seven instances (Section~\ref{sec:owned-duration}). Other parameters or initial values could be assigned to players in the same manner. The decision space $\mathcal U_i$ of player $i$ consists of all $u_i$ whose time-varying components are measurable and satisfy their bounds almost everywhere and whose scalar components are constant and satisfy their bounds; $\mathcal U=\prod_{i\in\mathcal P}\mathcal U_i$.

\begin{figure}[tb]
\centering
\definecolor{playerA}{RGB}{0,114,178}
\definecolor{playerB}{RGB}{230,159,0}
\definecolor{leafgreen}{RGB}{86,150,70}
\definecolor{leaflight}{RGB}{190,220,160}
\definecolor{bark}{RGB}{120,85,55}
\begin{tikzpicture}[x=1cm,y=1cm,font=\small,
  state/.style={draw=black!60,rounded corners=3pt,fill=black!4,minimum width=2.1cm,minimum height=.75cm,align=center},
  flow/.style={-{Stealth[length=2.2mm]},thick,black!70},
  inhibit/.style={-{Bar[width=2.6mm]},very thick,playerA},
  control/.style={-{Stealth[length=2.2mm]},very thick},
  goal/.style={draw,rounded corners=3pt,thick,align=left,inner sep=4pt,text width=5.9cm,font=\footnotesize}]
\begin{scope}[shift={(-5.3,0.35)}]
 \fill[bark] (-.16,-1.35) -- (.16,-1.35) -- (.12,.1) -- (-.12,.1) -- cycle;
 \draw[bark,thick] (0,-.2) -- (.45,.25) (0,-.45) -- (-.4,-.05);
 \foreach \x/\y/\r in {0/.75/.62,-.5/.35/.5,.52/.4/.5,-.25/1.15/.45,.3/1.1/.45,0/.3/.5}
  \fill[leafgreen!85!black,opacity=.9] (\x,\y) circle (\r);
 \foreach \x/\y in {-.3/.8,.25/.6,-.1/.25,.45/.95,-.55/.45,.05/1.2}
  \fill[leaflight] (\x,\y) circle (.07);
 \fill[black!25] (0,-1.38) ellipse (.8 and .08);
 \node[text=playerA,font=\small\bfseries,align=center] at (0,-1.75) {A: tree (tupelo)};
\end{scope}
\begin{scope}[shift={(5.3,.35)}]
 \fill[leafgreen] (-1.15,-.55) .. controls (-.6,.75) and (.55,.85) .. (1.2,.25)
                  .. controls (.55,-.55) and (-.5,-1.0) .. cycle;
 \draw[leafgreen!60!black,thick] (-1.15,-.55) .. controls (-.2,-.05) and (.5,.15) .. (1.2,.25);
 \draw[leafgreen!60!black] (-.55,-.27) -- (-.35,.2) (-.1,-.1) -- (.05,.35) (.35,.02) -- (.5,.4)
                           (-.55,-.27) -- (-.45,-.62) (-.1,-.1) -- (0,-.55) (.35,.02) -- (.45,-.35);
 \draw[leaflight,line width=1pt] (-.75,.12) .. controls (-.55,.45) and (-.3,-.05) .. (-.15,.3);
 \draw[leaflight,line width=2.2pt] (-.15,.3) .. controls (0,.55) and (.2,.1) .. (.35,.33);
 \foreach \x/\y in {-.85,-.72,-.95}{}
 \foreach \x/\y in {-.8/-.35,-.65/-.45,-.92/-.2}
  \fill[white,draw=black!40] (\x,\y) ellipse (.06 and .04);
 \foreach \i in {0,...,4}
  \fill[playerB!85!black] ($(.38,.34)+(\i*.075,.02*\i)$) circle (.055);
 \fill[black!70] (.72,.43) circle (.035);
 \node[text=playerB!80!black,font=\small\bfseries,align=center] at (0,-1.4) {B: leafminer};
\end{scope}
\node[state] (L) at (-.6,1.35) {leaf area $L$};
\node[state,minimum width=1.9cm] (H) at (-.6,-.3) {larvae $H$};
\node[state,minimum width=1.6cm] (E) at (2.7,-.3) {eggs $E$};
\node[font=\footnotesize] (P) at (-.6,-1.25) {pupation $\kappa\rho(L)H$};
\draw[flow] (E) -- node[below,font=\footnotesize] {hatching} (H);
\draw[flow,dashed] (H) -- node[left,font=\footnotesize] {feeding $b\rho(L)H$} (L);
\draw[flow] (H) -- (P);
\draw[inhibit] (-4.2,.15) .. controls (-2.8,.05) .. (H.west);
\draw[inhibit] (-4.2,.15) .. controls (-2.2,-2.75) and (1.9,-2.7) .. (E.south);
\node[text=playerA,font=\footnotesize,align=left,anchor=south west] at (-4.42,.75) {defense $u_A(t)$:\\kills eggs and larvae};
\draw[control,playerB] (4.1,.45) .. controls (3.4,.75) .. (E.north);
\draw[control,playerB] (4.1,.45) .. controls (2.6,1.25) and (1.1,.95) .. (H.north east);
\node[text=playerB!80!black,font=\footnotesize,anchor=south] at (3.55,.8) {hatching $u_B(t)$};
\node[goal,draw=playerA,fill=playerA!6] at (-3.3,-3.05) {%
 \textbf{\color{playerA}A minimizes} $J_A$: keep much valuable leaf area at a low carbon cost of defense\\[2pt]
 \emph{optional:} effort and carbon limits $U_A$, $C_A$};
\node[goal,draw=playerB,fill=playerB!8] at (3.3,-3.05) {%
 \textbf{\color{playerB!80!black}B minimizes} $J_B$: let many larvae pupate early at a low hatching cost\\[2pt]
 \emph{optional:} effort limit $U_B$};
\end{tikzpicture}
\caption{The tupelo--leafminer game. The tree (A) schedules its defense $u_A(t)$, which kills eggs and larvae; the leafminer (B) schedules the hatching hazard $u_B(t)$ of its eggs. Larvae feed on the leaves and leave the game by pupation. Black solid arrows denote transitions of individuals (hatching from eggs to larvae, pupation out of the larval stage); the dashed arrow denotes the loss of leaf area through larval feeding. Colored arrows denote controls that drive a transition (hatching changes eggs and larvae); colored lines ending in a bar denote controls that remove individuals (defense kills eggs and larvae). Each player minimizes its own objective; the allowances are optional restrictions, disabled in the reference instance. Section~\ref{sec:tupelo-setup} defines all terms.}\label{fig:tupelo-schematic}
\end{figure}
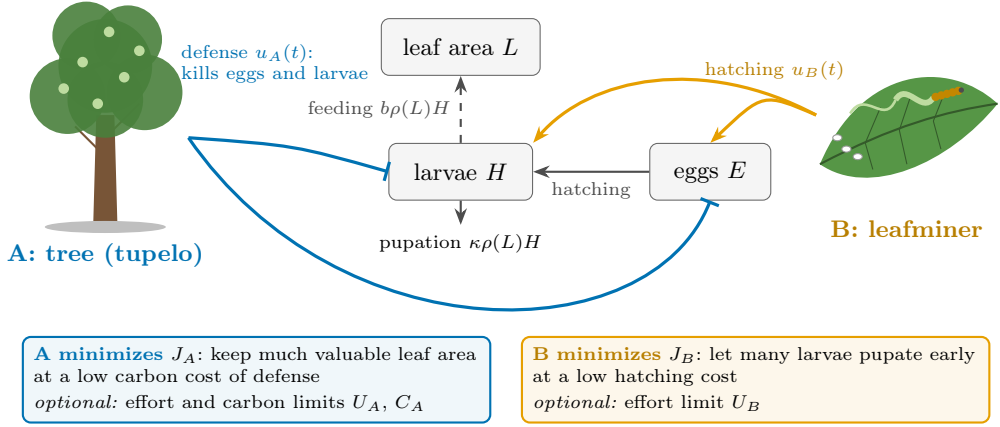

Figure~\ref{fig:tupelo-schematic} illustrates these objects for the tupelo--leafminer game, the worked example of Section~\ref{sec:uses}: each of the two players owns one control $u_A$ and $u_B$, the controls act on three biological states $x = (L, E, H)$, and each player minimizes its own objective $J_A$ and $J_B$, optionally subject to its own restrictions, as introduced in Section~\ref{sec:objectives}.

\subsection{Objectives and restrictions}\label{sec:objectives}
All objectives are minimized; benefits to be maximized enter with a negative sign. Player $i$ minimizes a weighted sum of $m_i$ features,
\begin{equation}\label{eq:objective}
\begin{aligned}
 J_i(u)&=\sum_{k=1}^{m_i}w_i^kJ_{ik}(u),\\
 J_{ik}(u)&=\Phi_{ik}(x(T),u)
 +\int_0^T L_{ik}(t,x(t),u(t);p)\dt
 +\sum_\ell Q_{ik\ell}(x(t_{ik\ell});p),
\end{aligned}
\end{equation}
where terms may be absent, $w_i^k\geq0$, and $\sum_kw_i^k=1$; in total, an instance has $n_w=\sum_{i\in\mathcal P}m_i$ features. The terminal term $\Phi_{ik}$ may depend on the scalar components of $u$, which allows a direct duration cost. Interior-time terms $Q_{ik\ell}$ at fixed times $t_{ik\ell}$ occur only in \caseref{motor}, the two other terms are standard Bolza-type terms. Each feature is divided by a fixed positive scale, listed with the parameters, which makes the features of one player comparable in magnitude; we absorb these scales into $J_{ik}$.

Restrictions beyond control bounds are written as nonnegative slacks,
\begin{equation}\label{eq:constraints}g_i[x,u;p]\geq0,\end{equation}
where the brackets indicate functionals of the whole trajectory and of the decisions. They are evaluated on the numerical state $y$, which augments $x$ by quadratures of objective terms and budgets and, in \caseref{ciliate} and \caseref{phage}, uses logarithmic coordinates for some biological components. Quadratures are excluded from $n_x$. Every scalar restriction $c$ assigned to player $i$ has the bound form
\begin{equation}\label{eq:bound-form}
 \ell_c\leq y_{q_c}(t)\leq r_c\qquad\text{for all }t\in[0,T]\text{ (path) or for }t=T\text{ (terminal)},
\end{equation}
with constrained component $q_c$, bounds $\ell_c\in\R\cup\{-\infty\}$ and $r_c\in\R\cup\{+\infty\}$, and an equality for $\ell_c=r_c$. The library uses four kinds of restrictions: path floors and ceilings on biological states; integral budgets, represented by a quadrature $\dot z=b(t,x,u;p)$, $z(0)=0$, and the terminal bound $z(T)\leq B$; terminal targets; and, in periodic instances, return bands
\begin{equation}\label{eq:return-band}
 x_j(0)-\delta\leq x_j(T)\leq x_j(0)+\delta,\qquad j=1,\ldots,n_x,
\end{equation}
imposed in every player's problem with an instance tolerance $\delta>0$; quadratures are excluded.
A \emph{joint} restriction is an identical restriction assigned to several players, i.e., a copy in each of their $g_i$. Optional restrictions carry a switch $\sigma\in\{0,1\}$; $\sigma=0$ removes the restriction from the instance, and $\sigma=1$ enables the stated finite bound.

\subsection{Open-loop generalized Nash equilibrium}\label{sec:equilibrium}
In the open-loop information structure, decisions are functions of time alone, and each player commits to its schedule at $t=0$. Let $u^{i,v_i}$ replace only player $i$'s decisions by $v_i$. Define $\mathcal F_i(u_{-i})$ as the set of $v_i\in\mathcal U_i$ for which the absolutely continuous solution of \eqref{eq:dynamics} under $u^{i,v_i}$ satisfies player $i$'s restrictions \eqref{eq:constraints}. Only the restrictions assigned to player $i$ enter $\mathcal F_i$: restrictions are assigned to the players who are bound by them. For example, the herbivore's hatching allowance in the allowance variant of \caseref{tupelo} (Section~\ref{sec:tupelo-variants}) restricts the herbivore but not the plant. The best-response problem is
\begin{equation}\label{eq:best-response-ocp}
 \min_{v_i\in\mathcal F_i(u_{-i})}J_i(u^{i,v_i}),
\end{equation}
with the set of global minimizers
\begin{equation}\label{eq:best-response-set}
 \mathcal K_i(u_{-i})=\argmin_{v_i\in\mathcal F_i(u_{-i})}J_i(u^{i,v_i}).
\end{equation}
An open-loop generalized Nash equilibrium satisfies
\begin{equation}\label{eq:gne}
 u_i^*\in\mathcal K_i(u_{-i}^*)\qquad\text{for every }i\in\mathcal P.
\end{equation}
Every deviation recomputes the complete coupled state. Generalized feasibility arises when $\mathcal F_i$ depends on opponents' decisions after eliminating the states. Private control-effort budgets can retain an ordinary Nash structure. Shared biological return conditions or resource floors generally couple feasibility.

When the same restriction is imposed in several players' problems, generalized Nash equilibria are typically not isolated \citep{FacchineiKanzow2010}. A frequently used refinement is the variational equilibrium, in which all players share one multiplier per shared restriction; the equilibrium definition \eqref{eq:gne} does not impose it. Section~\ref{sec:methods} computes candidates for \eqref{eq:gne} by solving discretized versions of \eqref{eq:best-response-ocp} and audits them for feasibility and profitable unilateral deviations (Section~\ref{sec:audit}). Different forward methods can therefore return different, equally valid equilibria of the same instance.

\subsection{Owned duration}\label{sec:owned-duration}
In many biological interactions, when an episode ends or how long a cycle lasts is itself decided by one participant: a seed cohort determines by its germination speed when the shared water-uptake phase ends, a host determines when an infection is cleared, and the fleet that sets the opening schedule of a fishing ground determines how long a rotation lasts. Assigning the duration to this player lets the game express a trade-off that a fixed horizon hides, between ending early, e.g., to limit exposure or to repeat the cycle more often, and the effort needed to reach the terminal state in time, which the other players shape through the common trajectory.
In seven instances of the library, player A owns the horizon, i.e., $T$ is a scalar component of $u_A$ with bounds $\underline T\leq T\leq\overline T$. We write all schedules on normalized time $s\in[0,1]$, with $t=Ts$, $\widetilde x(s)=x(Ts)$, and $\widetilde u(s)=u(Ts)$, so that
\begin{equation}\label{eq:time}
 \widetilde x'(s)=T f(Ts,\widetilde x(s),\widetilde u(s);p).
\end{equation}
We adopt the convention that, when A changes $T$, the opponents keep their schedules as functions of $s$, so that all physical schedules are rescaled. This convention is part of the game definition; keeping physical-time schedules fixed would define a different game. An accumulated running cost becomes $T\int_0^1L_{ik}(Ts,\widetilde x(s),\widetilde u(s);p)\,ds$, whereas a time average becomes $\int_0^1L_{ik}(Ts,\widetilde x(s),\widetilde u(s);p)\,ds$. This distinction matters when the duration is variable and applies equally to total and mean-effort budgets. In \caseref{fishery} and \caseref{predatorPrey}, the objectives are cycle averages. Non-owners influence the duration only through the trajectory and the restrictions that the owner must respect.

\subsection{Specification components and instance variants}\label{sec:specification}
The formulation separates three components of a game specification:
the \emph{mechanism}, i.e., the dynamics $f$, the data $p$, and the initial state $x^0$ in \eqref{eq:dynamics};
the \emph{preferences}, i.e., the features $J_{ik}$ and weights $w_i$ in \eqref{eq:objective}; and
the \emph{restrictions}, i.e., the control bounds \eqref{eq:control-set} and the restrictions $g_i$ with their bounds, tolerances, and switches in \eqref{eq:constraints}--\eqref{eq:return-band}.
Together with the player set, the control ownership, and the horizon, they define a specification
\[
 \mathcal G=\bigl(f,p,x^0;\ (J_{ik},w_i)_{i\in\mathcal P};\ (g_i,\sigma)_{i\in\mathcal P};\ \text{ownership},\,T\bigr).
\]
Bounds of restrictions belong to the restrictions even where the implementation stores them with the parameters. An instance fixes all entries of $\mathcal G$ (Section~\ref{sec:curation}). An \emph{instance variant} changes one or more entries while keeping the players, the ownership, and the information structure. Accordingly, we speak of preference, restriction, and mechanism variants. They correspond to three kinds of biological hypotheses: what a player values, e.g., leaf area relative to the cost of defense; what a player can do, e.g., defend only within a carbon budget; and how the interaction works, e.g., whether feeding saturates. Alternative explanations of the same observed behavior often differ in only one of these components, so comparing such variants, and recovering them from data, is a central use of the library. Section~\ref{sec:tupelo-variants} uses such variants to express biological hypotheses, and Section~\ref{sec:tupelo-inverse} treats entries of $\mathcal G$ as unknowns of inverse problems.

\subsection{Assumptions and scope}\label{sec:assumptions}

The following assumptions apply throughout this paper.

\paragraph*{Players, objectives, and constraints.}
Each instance describes deterministic aggregate dynamics conditional on prescribed initial conditions, parameters, and environmental inputs; spatial, stochastic, and individual-level mechanisms enter only through the stated reductions. Players represent regulatory modules, microbial strains, organisms, or populations. Assigning independent objectives to modules within one organism is an additional modeling hypothesis. Biomass, production, and effort terms are proxies; a link to reproductive success is needed before a population objective receives a selection interpretation. Many coefficients, weights, and restrictions are synthetic, selected to expose nonlinear dynamics, different player behavior, and binding restrictions. Abstract traits, occurrence associations, and hypothetical module objectives are not species calibrations, causal interaction estimates, or measured goals.

\paragraph*{Information structure.}
Open-loop schedules may approximate inherited developmental programs or repeatable allocation schedules in predictable environments. This interpretation requires no conscious foresight; \citet{BangaVillaverde2026Hierarchical} stress that assuming a game rather than evolved regulatory programs needs independent justification in each application. Open-loop schedules do not describe responses to unexpected changes in state or environment: open-loop Nash equilibria are in general not subgame perfect, i.e., if the state is perturbed during the episode, the continuation of an equilibrium schedule need not be an equilibrium of the continuation game. Feedback (Markov) equilibria describe state-dependent responses and generally differ from the open-loop equilibria reported here. Within-episode physiology and intergenerational strategy change operate on different modeled timescales. The games are simultaneous and noncooperative; the sequential best-response updates of Section~\ref{sec:methods} can be read as a learning process in which each player adapts to the others' latest schedules, but they serve here only to compute an equilibrium, so their number of rounds and their intermediate profiles are not a model of biological learning or evolution. A biological leader in migration or germination is not a Stackelberg leader. A duration owner specifies who can alter the common schedule; for regulatory modules and predator--prey cycles, this is a synthetic allocation of decision authority. The additional degree of freedom is balanced by the return band \eqref{eq:return-band} in the three periodic instances and by a terminal target in the four transfer tasks.
Accordingly, the open-loop Nash conditions established here should not by themselves be interpreted as demonstrating evolutionary stability or resistance to invasion; establishing such properties would require an explicit evolutionary framework and the corresponding stability conditions \citep{MaynardSmithPrice1973,Smith1982}.

\paragraph*{Temporal classes.}
A seasonal episode can connect to the next season through an unmodeled transition $x_{k+1}(0)=R(x_k(T),e_k)$ that incorporates offseason conditions $e_k$, mortality, reproduction, dormancy, or replenishment. The game interface can represent such transitions through dynamics or endpoint restrictions, but the library currently implements none. Periodic instances describe interactions that repeat cycle after cycle, such as a daily activity cycle or a recurring supply of a culture: biologically, the state at the end of one cycle is the start of the next, and the environment repeats with the same period. Mathematically, the return band \eqref{eq:return-band} imposes this repetition: every biological state must end close to its initial value. For fixed initial values, this may lead to modeling artefacts, i.e., end-of-cycle corrections in the controls. The shared return band also couples all players (Section~\ref{sec:equilibrium}). Episodic instances describe individual tasks, recovery, or experiments with resets. Only periodic instances impose a biological return.

\section{The library}\label{sec:cases}
\subsection{Curation}\label{sec:curation}
A \emph{model family} specifies a mechanism with configurable preferences and restrictions (Section~\ref{sec:specification}). An \emph{instance} fixes all entries of the specification $\mathcal G$ of Section~\ref{sec:specification}; changing weights or enabling a restriction therefore creates a new instance and can change feasibility or the existence of an equilibrium. A \emph{numerical record} additionally identifies grids, solver settings, and the resulting trajectories and audits.

Every instance of the curated library has a complete, executable specification with distinct control ownership, documented biological motivation and temporal interpretation, and a reference computation. In addition, the players must be plausible decision units, the central mechanism must be supported by a primary source, every biological state must influence the dynamics, an objective, or a restriction, and the specification must couple the players' responses. Models whose game content rested on decorative states or on a single player's optimal control were not included. Repeated patterns, such as squared-effort budgets, permit controlled comparisons across different coupling mechanisms. Other contrasts include private versus opponent-dependent feasibility, state and control dimensions, variable player counts, and fixed versus owned durations.

\subsection{Classifications and attributes}\label{sec:classifications}

Table~\ref{tab:categories} defines three exclusive classifications and two attributes. Each instance has exactly one biological scale, one temporal class, and one interaction type. The attributes \cat{variablePlayers} and \cat{freeTime} are nonexclusive. 

\begin{table}[!htbp]
\caption{Category identifiers and their description.}\label{tab:categories}
\small
\newlength{\categorylabelwidth}
\settowidth{\categorylabelwidth}{\texttt{competition}}
\begin{tabular}{@{}l@{\hspace{.6em}}p{\dimexpr\textwidth-\categorylabelwidth-.6em\relax}@{}}
\toprule Category & Definition\\\midrule
\multicolumn{2}{@{}l}{\emph{Biological scale (exclusive)}}\\
\phantomsection\label{cat:molecular}\texttt{molecular} & Biochemical species, molecular pools, or subcellular processes. \\
\phantomsection\label{cat:cellular}\texttt{cellular} & Cells or microbes as players or states. \\
\phantomsection\label{cat:organismal}\texttt{organismal} & Organisms, their movements, or their neural control. \\
\phantomsection\label{cat:population}\texttt{population} & Populations, communities, or colonies of multicellular organisms. \\\midrule
\multicolumn{2}{@{}l}{\emph{Temporal class (exclusive)}}\\
\phantomsection\label{cat:seasonal}\texttt{seasonal} & Finite seasonal episode; any offseason transition is outside the model. \\
\phantomsection\label{cat:periodical}\texttt{periodical} & All biological states return within a stated tolerance. \\
\phantomsection\label{cat:episodical}\texttt{episodical} & Single task, recovery, experiment, or startup transient. \\\midrule
\multicolumn{2}{@{}l}{\emph{Interaction type (exclusive)}}\\
\phantomsection\label{cat:competition}\texttt{competition} & Competition for a shared resource or reproductive opportunity. \\
\phantomsection\label{cat:antagonism}\texttt{antagonism} & Predation, infection, or herbivory as a central interaction. \\
\phantomsection\label{cat:mutualism}\texttt{mutualism} & Reciprocal benefits from exchange or service; competition may coexist. \\
\phantomsection\label{cat:regulation}\texttt{regulation} & Players regulate one coupled system toward different targets, e.g., regulatory modules of one organism or partners in a joint task. \\\midrule
\multicolumn{2}{@{}l}{\emph{Attributes (nonexclusive)}}\\
\phantomsection\label{cat:variablePlayers}\catlabel{variablePlayers} & Player count $N$ adjustable through model parameters. \\
\phantomsection\label{cat:freeTime}\texttt{freeTime} & Horizon owned as a bounded scalar decision parameter. \\
\bottomrule
\end{tabular}
\end{table}

Seasonal and episodic instances impose no return condition. Periodic instances impose the return band \eqref{eq:return-band} with the tolerance $\delta$ given by the instance parameter \code{periodicity_tolerance}. The library contains seven seasonal, nine periodic, and nine episodic instances. There are seven free-duration instances: four transfer tasks and three periodic instances. The other six periodic instances have fixed external clocks. The biological rationale differs across scales. Externally forced cycles represent recurring supply or activity regimes, such as periodic feeding of maintained cultures (\caseref{algae}, \caseref{metabolicPathway}, \caseref{microbial}) or daily activity (\caseref{mitochondria}); metabolic and organellar processes are known to be organized in time \citep{Tu2005MetabolicCycle,MitochondrialClock2020}. Consumer--resource cycles, in which evolution and population dynamics interact on similar time scales \citep{Yoshida2007Cryptic}, motivate \caseref{predatorPrey}, whose cycle is autonomous, and \caseref{ciliate}, whose cycle is driven by seasonal productivity. The specific forcing, controls, and individual costs remain part of each documented adaptation. Ten model families define their player count through parameters (\cat{variablePlayers}); for periodic families, a changed player count or changed kinetics may require new initial states compatible with the cycle.

\subsection{The game instances at a glance}\label{sec:glance}
This subsection summarizes the 25 instances by biological scale, following the order of Table~\ref{tab:categories}. Each summary names the players, their decisions, and what couples them; Appendix~\ref{app:cases} gives the complete specifications.

\paragraph*{Molecular scale.}
In \caseref{circadianCycle}, the ten-state oscillator of the clock proteins PER and TIM of \citet{LeloupGoldbeter1998PerTim} keeps its published kinetics; player A modulates \emph{per} transcription and owns the cycle duration, while player B, the light-input pathway, controls TIM degradation. Both players prefer balanced PER/TIM levels at low effort, and return bands on all ten concentrations couple their deviation sets.
In \caseref{metabolicPathway}, two expression modules each pay for their own enzyme in a branched pathway whose common product both value, a public-goods situation in which free-riding is expected. Substrate is supplied periodically, and module-specific byproducts inhibit each module's own reaction, so enzyme investment must be timed relative to supply and byproduct load.
\caseref{mitochondria} treats a wild-type and a variant subpopulation of mitochondrial DNA (mtDNA) as players that invest in biogenesis while drawing on one shared energetic pool over a daily activity cycle. The variant has a larger replication allowance but contributes less energy and is lost faster.
In \caseref{proteinFolding}, two chaperone systems, i.e., helper proteins that assist folding, with different efficiencies share a consumable chaperone capacity during recovery from a stress pulse. Both value the same native protein pool and have private effort budgets, so the game isolates free-riding on assistance effort.
All molecular instances assign independent objectives to regulatory modules or subpopulations within one cell, which is an additional modeling hypothesis (Section~\ref{sec:assumptions}); three of the four are periodic.

\paragraph*{Cellular scale.}
In \caseref{algae}, four algal species share dissolved CO$_2$ and bicarbonate pools and allocate their uptake machinery between the two carbon forms over one supply cycle of a maintained culture. The uptake traits follow the four abstract types of \citet{Zhou2022Algae}, and joint return bands on all biomasses and resources couple all players.
In \caseref{batchFermentation}, an upstream strain converts substrate into an intermediate that two downstream strains turn into product; the upstream strain also chooses when the batch ends. Production costs a growth burden, time is costly to the upstream strain, and a product target shared by all strains closes the episode.
In \caseref{biofilm}, three strains on one surface divide effort between matrix production, which protects their own cells, and dispersal, while they compete for substrate and oxygen. The strains differ in how they value attached biomass relative to dispersal output, which creates a timing problem between protection, growth, and departure.
\caseref{ciliate} couples a bacterium that expresses grazing defense with a ciliate, a single-celled grazer, that adjusts its feeding activity over an annual cycle; both phenotypes follow the commanded levels with a delay. Return bands on densities and phenotypes couple the two deviation sets.
In \caseref{immuneClearance}, the host activates cytotoxic responses, which kill infected cells, and antibody responses and ends the infection episode, whereas the virus invests in immune evasion to maximize its accumulated viral load. The clearance target restricts only the host; the virus is bound by its own evasion budget alone.
\caseref{microbial} extends the two-strain cross-feeding game of \citet{BangaVillaverde2026Hierarchical}, which these authors use as a common example for all four levels of inverse problems, to ten strains in reciprocal pairs, each secreting one metabolite that its partner requires while all compete for carbon under a periodic feed. With 21 biological states and 420 scalar constraint functions, it is the largest instance of the library.
In \caseref{phage}, a bacterial population invests in resistance at a cost in growth, and a phage population adjusts the average infectivity of its progeny within one growth episode. Expressed resistance and infectivity follow the commanded levels with a delay.
\caseref{woundRepair} assigns inflammatory recruitment, fibroblast recruitment, and oxygen support to three players who weigh effort, excess inflammation, and collagen quality differently. Inflammatory recruitment owns the completion time, and the episode ends when the wound deficit has fallen to a prescribed residual size; a free duration thus arises without periodicity.
The cellular instances cover all four interaction types; five of the eight are episodic.

\paragraph*{Organismal scale.}
In \caseref{avoidance}, two decentralized control modules time their effort in a signaled avoidance trial: acting builds readiness but depletes a shared, slowly recovering motivational resource. Each module benefits in part from the other's readiness, so the players are partially aligned.
\caseref{motor} follows the joint-action experiments of \citet{Chackochan2019JointAction}: two people move mechanically coupled handles through different private via-points toward a common target. Its linear dynamics and quadratic costs make it a candidate for an analytical cross-check of numerical equilibria.
In \caseref{neural}, two regulators drive the excitatory and inhibitory populations of a Wilson--Cowan circuit \citep{WilsonCowan1972} toward different activity targets under a periodic background drive. Both activities must return to their initial values after one cycle.
All three organismal instances are regulation games between modules or partners; \caseref{motor} retains published mechanical coefficients.

\paragraph*{Population scale.}
In \caseref{antColonies}, three neighboring colonies allocate workers among food foraging, water collection, and territorial aggression during one season. Stored food becomes brood and later workers, so allocations have delayed consequences; effort budgets and a food-store floor restrict each colony.
\caseref{fishery} models one rotation of a rotationally managed fishing ground: two fleets harvest a logistic stock, and fleet A also chooses the rotation length. Catch quotas and an approximate stock return connect consecutive openings; the stock and cost structure follow \citet{Suri2008Harvesting}.
In \caseref{foraging}, three foragers deplete a shared renewable food patch and convert food into reserves and a retained harvest. A management floor keeps the patch above 70\% of its carrying capacity and thereby couples their feasible sets.
In \caseref{plantNitrogen}, three plants compete for an early-spring nitrogen pulse that builds photosynthetic capacity while they shade one another. A shared floor on soil nitrogen represents the concentration at which net uptake ceases.
\caseref{plantWater} lets four plants allocate root investment across three soil layers during a season with a wet first and a dry second half. The plants specialize on different layers, and a seasonal water-uptake allowance limits each of them.
In \caseref{pollination}, two plant species invest in nectar and two pollinator species allocate their visits among the plants during one flowering season. Effort budgets and population floors restrict all four players.
\caseref{predatorPrey} describes an autonomously oscillating community of three prey and two predator species; prey invest in defense or growth, predators in hunting or maintenance. Prey species A chooses the cycle length, which no external clock imposes.
In \caseref{seedGermination}, three seed cohorts draw on one water pool during imbibition, the water uptake that precedes germination; fast activation secures early emergence but consumes water that the neighbors need. The window ends when the leading cohort reaches radicle emergence, i.e., its embryonic root breaks through the seed coat, and cohort A owns this duration.
\caseref{sexRatio} adapts local mate competition \citep{Hamilton1967SexRatios}, in which brothers compete for mates and mothers therefore produce few sons: three foundresses, females that found a local brood, choose the male fraction of their offspring and provision their brood from a shared, seasonally replenished resource. A synthetic floor on this resource couples their feasible sets.
\caseref{tupelo} is the plant--herbivore game discussed in Figure~\ref{fig:tupelo-schematic} and in Section~\ref{sec:uses}.
Seven of the ten population instances are seasonal, and seven describe competition.

\subsection{Common instance schema}
All instances are documented in the same schema. The appendix records biological context and provenance before the equations, followed by the parameter defaults and an interpretation of the saved trajectory. Published equations or coefficients, reductions, and introduced game components are identified separately. Table~\ref{tab:overview} lists the classifications and attributes of every instance. Its dimensions appear in Table~\ref{tab:reference-results} and, per player, with each specification.
\FloatBarrier
\begingroup\small\setlength{\tabcolsep}{3pt}
\begin{longtable}{llllll}
\caption{Library overview with classifications of Table~\ref{tab:categories}.}\label{tab:overview}\\
\toprule Case & Scale & Temporal & Type & \texttt{variable-}$N$ & \texttt{freeTime}\\\midrule\endfirsthead
\toprule Case & Scale & Temporal & Type & \texttt{variable-}$N$ & \texttt{freeTime}\\\midrule\endhead
\hyperref[case:algae]{\texttt{algae}} & cellular & periodical & competition & yes & no \\
\hyperref[case:antColonies]{\texttt{antColonies}} & population & seasonal & competition & yes & no \\
\hyperref[case:avoidance]{\texttt{avoidance}} & organismal & episodical & regulation & no & no \\
\hyperref[case:batchFermentation]{\texttt{batchFermentation}} & cellular & episodical & mutualism & no & yes \\
\hyperref[case:biofilm]{\texttt{biofilm}} & cellular & episodical & competition & yes & no \\
\hyperref[case:ciliate]{\texttt{ciliate}} & cellular & periodical & antagonism & no & no \\
\hyperref[case:circadianCycle]{\texttt{circadianCycle}} & molecular & periodical & regulation & no & yes \\
\hyperref[case:fishery]{\texttt{fishery}} & population & periodical & competition & no & yes \\
\hyperref[case:foraging]{\texttt{foraging}} & population & episodical & competition & yes & no \\
\hyperref[case:immuneClearance]{\texttt{immuneClearance}} & cellular & episodical & antagonism & no & yes \\
\hyperref[case:metabolicPathway]{\texttt{metabolicPathway}} & molecular & periodical & regulation & no & no \\
\hyperref[case:microbial]{\texttt{microbial}} & cellular & periodical & mutualism & yes & no \\
\hyperref[case:mitochondria]{\texttt{mitochondria}} & molecular & periodical & regulation & no & no \\
\hyperref[case:motor]{\texttt{motor}} & organismal & episodical & regulation & no & no \\
\hyperref[case:neural]{\texttt{neural}} & organismal & periodical & regulation & no & no \\
\hyperref[case:phage]{\texttt{phage}} & cellular & episodical & antagonism & no & no \\
\hyperref[case:plantNitrogen]{\texttt{plantNitrogen}} & population & seasonal & competition & yes & no \\
\hyperref[case:plantWater]{\texttt{plantWater}} & population & seasonal & competition & yes & no \\
\hyperref[case:pollination]{\texttt{pollination}} & population & seasonal & mutualism & yes & no \\
\hyperref[case:predatorPrey]{\texttt{predatorPrey}} & population & periodical & antagonism & yes & yes \\
\hyperref[case:proteinFolding]{\texttt{proteinFolding}} & molecular & episodical & regulation & no & no \\
\hyperref[case:seedGermination]{\texttt{seedGermination}} & population & seasonal & competition & no & yes \\
\hyperref[case:sexRatio]{\texttt{sexRatio}} & population & seasonal & competition & yes & no \\
\hyperref[case:tupelo]{\texttt{tupelo}} & population & seasonal & antagonism & no & no \\
\hyperref[case:woundRepair]{\texttt{woundRepair}} & cellular & episodical & regulation & no & yes \\
\bottomrule\end{longtable}\endgroup

\section{Solving the games with CorleoneGame.jl}\label{sec:methods}
\subsection{Interface and optimal-control subproblems}
The solver \code{CorleoneGame.jl} builds on the optimal-control package \code{Corleone.jl} \citep{CorleoneRepository}. The game extension, the case files, and this manuscript were developed with an AI coding agent under the authors' supervision; its role is described in the statement on the use of generative AI at the end of the article. Each case file implements the specification $\mathcal G$ of Section~\ref{sec:specification} together with numerical options. The solver fixes the opponents' decisions and assembles only the active player's objective and restrictions, i.e., a discretization of the best-response problem \eqref{eq:best-response-ocp}. Controls are piecewise constant on $M$ equal intervals, owned scalar decisions such as the duration are optimization parameters, and path restrictions are imposed at sample points within each interval. The reference records use single shooting, which eliminates the states by numerical integration, and the interior-point solver Ipopt. Appendix~\ref{app:numerics} lists the integrator, tolerances, derivatives, and further options.

\subsection{Sequential response iteration and refinement}\label{sec:iteration}
Algorithm~\ref{alg:baseline} states the implemented baseline. Updates are sequential, so later players see earlier updates from the same round. Each response is damped, $u_i\leftarrow(1-\alpha)u_i+\alpha\widehat u_i$, including owned scalar decisions. With $\alpha=1$, the sweep is the nonlinear Gauss--Seidel-type method described by \citet{FacchineiKanzow2010}; damping with $\alpha<1$ is a heuristic modification. The same damped best-response iteration with a subsequent unilateral-deviation check is used in the Level-4 case study of \citet{BangaVillaverde2026Hierarchical}. If the whole sequence of profiles of such a method converges, its limit is an equilibrium, but conditions that guarantee convergence are typically not known or very restrictive, and generalized problems require additional properties of the constraints \citep{FacchineiKanzow2010}. The iteration therefore has no general convergence guarantee for this library. In \caseref{mitochondria}, for example, tight return bands let successive best responses drift along a family of generalized equilibria, which are not isolated; the looser band of the reference instance avoids this drift.

\begin{algorithm}[t]
\caption{Sequential response baseline with independent audits}\label{alg:baseline}
\begin{algorithmic}[1]
\Require Game specification, grids $M_0<M_1$, initial profile, tolerances, damping $\alpha$
\For{$M\in\{M_0,M_1\}$}
 \If{$M=M_1$} \State Transfer the coarse controls by interval averaging \EndIf
 \For{$k=1,\ldots,k_{\max}$}
  \For{$i\in\mathcal P$ in the order $A,B,\ldots$}
   \State Solve player $i$'s shooting problem with current opponents fixed
   \State Independently check the response (stop on failure)
   \State $u_i\leftarrow(1-\alpha)u_i+\alpha\widehat u_i$ for owned controls and scalar decisions
  \EndFor
  \State Reintegrate the completed profile with tighter tolerances
  \State Audit all unilateral deviations against that same profile
  \State Save objectives, gains, violations, and a trajectory snapshot
  \If{$\max_iV_i\leq\varepsilon_g$ and $\max_i\widehat\Delta_i\leq\varepsilon_J$}
   \State Mark this grid converged and leave the round loop
  \EndIf
 \EndFor
 \State If $k_{\max}$ is reached, report this grid as not converged
\EndFor
\State \Return Profile, per-grid convergence flags, histories, and saved outputs
\end{algorithmic}
\end{algorithm}

The defaults are $M_0=10$ and $M_1=40$. The damping ranges from $\alpha=0.15$ to $\alpha=1$; instance-specific values are listed in Appendix~\ref{app:numerics}.

\subsection{Audits and acceptance}\label{sec:audit}
An audit asks, for each player separately, whether it could have done better by changing only its own schedule while all other players keep theirs. It evaluates every player's deviation against the same completed profile and thus tests the equilibrium condition \eqref{eq:gne}. With a positive reporting scale $\gamma_i$, the observed normalized gain is
\begin{equation}\label{eq:gain}
 \widehat\Delta_i=\frac{\max\{0,J_i(u)-J_i(u^{i,\widehat v_i})\}}{\gamma_i},
\end{equation}
where $\widehat v_i$ is the lowest-cost feasible response found by a local solve from a single start per player and round; all reference records use $\gamma_i=1$. The violation $V_i$ of player~$i$ sums, over its controls and its restrictions in the bound form \eqref{eq:bound-form}, the largest excess over the bounds at the audit sample times \eqref{eq:violation}. A grid is accepted when $\max_iV_i\leq\varepsilon_g=10^{-6}$ and $\max_i\widehat\Delta_i\leq\varepsilon_J=2\cdot10^{-5}$. Each round uses one response sweep and one full deviation audit, hence $2N$ optimal-control calls. In addition, every response is checked by an independent reintegration. %

\subsection{Reference records}\label{sec:records}

Table~\ref{tab:reference-results} summarizes the saved \code{iterations.csv} files. 

All 25 records contain 10- and 40-interval stages that satisfy the acceptance test. Each instance provides the saved trajectories, controls, iteration histories including the controls after every outer round, and plots of its reference computation. Across the instances, the largest final normalized gain is $1.959\cdot10^{-5}$ (\caseref{foraging}), and the largest violation is $9.070\cdot10^{-7}$ (\caseref{antColonies}). The coarse stage requires 2--41 rounds, and the refined stage 1--35. The round counts are unrelated to the problem dimensions (rank correlations with $N$, $n_x$, $n_u$, $n_w$, and $n_g$ between $-0.27$ and $0.29$): the small games \caseref{tupelo} and \caseref{neural} need 40 and 31 coarse rounds, whereas the largest instance, \caseref{microbial}, needs 3 and 1. The median CPU time of one optimal-control solve ranges from 0.025\,s to 4.4\,s on the coarse grid and from 0.16\,s to 76\,s on the refined grid, where it is 5 to 27 times larger. It grows with the state and control dimensions (rank correlations of about 0.7 with $n_x$ and $n_u$, at most 0.37 with $N$ and $n_g$), but on the refined grid the three instances with exact Hessians (Appendix~\ref{app:numerics}) are the most expensive per solve, ahead of \caseref{microbial} with 21 states and 420 scalar constraint functions.
\begingroup\small\setlength{\tabcolsep}{3pt}
\begin{longtable}{lrrrrrrrrrr}
\caption{Dimensions and results for 10 (coarse) and 40 (refined) control intervals. %
CPU is the median time [s] of one optimal-control solve.}\label{tab:reference-results}\\
\toprule &&&&&&\multicolumn{2}{c}{Coarse grid}&\multicolumn{2}{c}{Refined grid}&\\\cmidrule(lr){7-8}\cmidrule(lr){9-10}
Case & $N$ & $n_x$ & $n_u$ & $n_w$ & $n_g$ & Rounds & CPU & Rounds & CPU & $T$\\\midrule\endfirsthead
\toprule &&&&&&\multicolumn{2}{c}{Coarse grid}&\multicolumn{2}{c}{Refined grid}&\\\cmidrule(lr){7-8}\cmidrule(lr){9-10}
Case & $N$ & $n_x$ & $n_u$ & $n_w$ & $n_g$ & Rounds & CPU & Rounds & CPU & $T$\\\midrule\endhead
\hyperref[case:algae]{\texttt{algae}} & 4 & 6 & 4 & 12 & 48 & 41 & $0.215$ & 13 & $1.350$ & $60$ \\
\hyperref[case:antColonies]{\texttt{antColonies}} & 3 & 11 & 9 & 9 & 6 & 20 & $0.902$ & 35 & $7.999$ & $6$ \\
\hyperref[case:avoidance]{\texttt{avoidance}} & 2 & 4 & 2 & 8 & 0 & 5 & $0.025$ & 2 & $0.157$ & $12$ \\
\hyperref[case:batchFermentation]{\texttt{batchFermentation}} & 3 & 6 & 4 & 9 & 9 & 3 & $1.371$ & 1 & $37.316$ & $4.07393$ \\
\hyperref[case:biofilm]{\texttt{biofilm}} & 3 & 8 & 6 & 9 & 9 & 12 & $0.601$ & 5 & $5.775$ & $12$ \\
\hyperref[case:ciliate]{\texttt{ciliate}} & 2 & 4 & 2 & 6 & 18 & 13 & $0.158$ & 10 & $1.057$ & $20$ \\
\hyperref[case:circadianCycle]{\texttt{circadianCycle}} & 2 & 10 & 3 & 6 & 44 & 6 & $4.448$ & 6 & $75.646$ & $24.0425$ \\
\hyperref[case:fishery]{\texttt{fishery}} & 2 & 1 & 3 & 4 & 8 & 2 & $0.049$ & 2 & $0.649$ & $5.39474$ \\
\hyperref[case:foraging]{\texttt{foraging}} & 3 & 7 & 3 & 6 & 6 & 7 & $0.168$ & 10 & $1.862$ & $10$ \\
\hyperref[case:immuneClearance]{\texttt{immuneClearance}} & 2 & 6 & 4 & 5 & 5 & 4 & $0.804$ & 2 & $20.724$ & $2.93398$ \\
\hyperref[case:metabolicPathway]{\texttt{metabolicPathway}} & 2 & 6 & 2 & 4 & 26 & 9 & $0.055$ & 4 & $0.404$ & $12$ \\
\hyperref[case:microbial]{\texttt{microbial}} & 10 & 21 & 10 & 40 & 420 & 3 & $2.231$ & 1 & $15.370$ & $80$ \\
\hyperref[case:mitochondria]{\texttt{mitochondria}} & 2 & 5 & 2 & 6 & 22 & 10 & $0.205$ & 11 & $1.169$ & $40$ \\
\hyperref[case:motor]{\texttt{motor}} & 2 & 16 & 4 & 10 & 0 & 22 & $1.121$ & 13 & $11.735$ & $0.5$ \\
\hyperref[case:neural]{\texttt{neural}} & 2 & 2 & 2 & 6 & 8 & 31 & $0.126$ & 20 & $0.752$ & $100$ \\
\hyperref[case:phage]{\texttt{phage}} & 2 & 5 & 2 & 6 & 2 & 15 & $0.086$ & 8 & $0.464$ & $8$ \\
\hyperref[case:plantNitrogen]{\texttt{plantNitrogen}} & 3 & 7 & 3 & 9 & 6 & 11 & $0.207$ & 11 & $1.976$ & $12$ \\
\hyperref[case:plantWater]{\texttt{plantWater}} & 4 & 19 & 12 & 12 & 4 & 27 & $0.955$ & 4 & $8.882$ & $8$ \\
\hyperref[case:pollination]{\texttt{pollination}} & 4 & 6 & 6 & 12 & 8 & 12 & $0.428$ & 4 & $4.444$ & $8$ \\
\hyperref[case:predatorPrey]{\texttt{predatorPrey}} & 5 & 5 & 11 & 15 & 55 & 22 & $0.769$ & 2 & $5.288$ & $17.4166$ \\
\hyperref[case:proteinFolding]{\texttt{proteinFolding}} & 2 & 4 & 2 & 4 & 2 & 9 & $0.057$ & 3 & $0.413$ & $8$ \\
\hyperref[case:seedGermination]{\texttt{seedGermination}} & 3 & 7 & 4 & 9 & 9 & 2 & $0.078$ & 1 & $1.336$ & $3.07319$ \\
\hyperref[case:sexRatio]{\texttt{sexRatio}} & 3 & 7 & 3 & 9 & 6 & 10 & $0.104$ & 3 & $1.358$ & $6$ \\
\hyperref[case:tupelo]{\texttt{tupelo}} & 2 & 3 & 2 & 7 & 0 & 40 & $0.064$ & 17 & $0.740$ & $8$ \\
\hyperref[case:woundRepair]{\texttt{woundRepair}} & 3 & 5 & 4 & 9 & 9 & 2 & $0.090$ & 1 & $1.361$ & $7.33042$ \\
\bottomrule\end{longtable}\endgroup

\FloatBarrier

\section{Working with the library: the tupelo game}\label{sec:uses}
This section follows one instance, the tupelo--leafminer game, through the steps a user of the library might take. Section~\ref{sec:tupelo-setup} translates biological assumptions into an executable specification, and Section~\ref{sec:tupelo-reference} computes and audits a reference candidate. Section~\ref{sec:tupelo-variants} changes the specification to express an alternative hypothesis, namely that defense is limited by a carbon budget rather than by the declining value of the leaves, and shows that it leads to a different candidate. Section~\ref{sec:tupelo-inverse} uses the resulting candidates to set up synthetic inverse problems. Each step relies on the objects introduced in the previous one.
This path is also a teaching exercise.

\subsection{From biological assumptions to an executable game}\label{sec:tupelo-setup}
The tupelo--leafminer interaction illustrates strategic timing within a season \citep{Low2013ColdWar}. Player A, the tree, invests in leaf defense $u_A$; player B, the leafminer, an insect whose larvae feed inside the leaves, chooses a hatching hazard $u_B$ for its eggs. Figure~\ref{fig:tupelo-schematic} in Section~\ref{sec:formulation} summarizes the game. Leaf area $L$, eggs $E$, and larvae $H$ satisfy
\begin{equation}\label{eq:tupelo-main}
\begin{aligned}
 \dot L&=-b\rho(L)H,\\
 \dot E&=-(u_B+\mu_E+d_Eu_A)E,\\
 \dot H&=u_BE-(\mu_H+d_Hu_A+\kappa\rho(L))H.
\end{aligned}
\end{equation}
Defense reduces egg and larval survival; hatching controls when larvae encounter defense and available leaf tissue. The reference instance follows the model of \citet{Low2013ColdWar}, in which consumption is proportional to larval density and larvae pupate, i.e., enter the pupal stage and leave the leaf, at a constant rate, i.e., $\rho(L)=1$. The alternative of Section~\ref{sec:tupelo-variants} uses the saturation $\rho(L)=L/(K+L)$ instead.
The plant values leaf area weighted by $p(t)$ and pays tissue-weighted defense costs; the herbivore values maturation weighted by $\omega(t)$ and pays an egg-weighted hatching cost:
\[
\begin{aligned}
 J_A&=-w_A^1\!\int_0^T\! p(t)L\dt+w_A^2\!\int_0^T\! u_AL\dt+w_A^3\!\int_0^T\! u_A^2L\dt-w_A^4L(T),\\
 J_B&=-w_B^1\!\int_0^T\! \omega(t)\kappa\rho(L)H\dt+w_B^2\!\int_0^T\! u_B^2E\dt-w_B^3H(T),
\end{aligned}
\]
with the photoperiod (day-length) curve $p(t)=e^{-(t/t^*)^z}$ and $\omega(t)=e^{-\delta_Bt}$ of \citet{Low2013ColdWar}. With $w_A^1=w_A^2$, the first two terms of $J_A$ are proportional to their net carbon gain $\int(p(t)-u_A)L\dt$, and the first term of $J_B$ is their pupation payoff. The reference weights are $w_A=(0.49,0.49,0.02,0)$ and $w_B=(0.98,0.02,0)$. The small quadratic terms are not part of the published payoffs; they regularize the best responses, whose optimal controls are otherwise singular or bang-bang \citep{Low2013ColdWar}.

The rates follow coevolution scenario~11 of \citet{Low2013ColdWar} (defense lethality $0.5$ per day, half of the larvae lost to background mortality, larval life span 25 days), mapped to one time unit of 25 days, so that $T=8$ corresponds to a season of 200 days:
\begin{center}\small
\begin{tabular}{@{}lcccccccc@{}}\toprule
 & $\kappa$ & $\mu_H$ & $d_H$ & $d_E$ & $\mu_E$ & $b$ & $\delta_B$ & $z$\\\midrule
 per day \citep{Low2013ColdWar} & $0.04$ & $0.04$ & $0.5$ & $0.01$ & $0$ & $0.06$ & $0.005$ & $9$\\
 per time unit (reference) & $1$ & $1$ & $12.5$ & $0.25$ & $0$ & $1.5$ & $0.125$ & $9$\\\bottomrule
\end{tabular}
\end{center}
Egg mortality is neglected, as in their model. \citet{Low2013ColdWar} state no value for the photoperiod scale $t^*$; we assume $t^*=6$ (150 days), so that the leaf value halves near the end of the growing season. The season length of 200 days is inferred from the time axes of their figures.

Optional allowances restrict the players' effort and the tree's carbon use:
\[
 \sigma_{U_A}\Bigl(U_A-\int_0^T u_A\dt\Bigr)\geq0,\quad
 \sigma_{U_B}\Bigl(U_B-\int_0^T u_B\dt\Bigr)\geq0,\quad
 \sigma_{C_A}\Bigl(C_A-\int_0^T u_AL\dt\Bigr)\geq0.
\]
The carbon allowance $C_A$ depends on the opponent through leaf area; the two effort allowances do not. As in \citet{Low2013ColdWar}, all switches are $\sigma=0$ in the reference instance; the alternative of Section~\ref{sec:tupelo-variants} enables them. 

The implementation loads \code{CorleoneGame.jl}, which defines the module \code{CorleoneGame}, and an instance:
\begin{lstlisting}[language=Julia]
include("tupelo/tupelo.jl")
using .CorleoneGame
a = Tupelo.Parameters()
\end{lstlisting}
With the coefficients of Appendix~\ref{case:tupelo}, the dynamics callback of \code{tupelo/tupelo.jl} implements mechanism and quadratures for objectives and constraints:
\begin{lstlisting}[language=Julia,style=numbered,firstline=47,lastline=64,firstnumber=47]
leaf_value(t,a)=a.photoperiod ? exp(-(t/a.t_star)^a.z) : exp(-a.p_decay*t)

function dynamics!(dx, x, controls, t, a)
    L, E, H = x[1:3]
    u, v = controls
    g = a.saturation ? L / (a.K + L) : one(L)
    dx[1] = -a.b * g * H
    dx[2] = -(v + a.muE + a.dE * u) * E
    dx[3] = v * E - (a.muH + a.dH * u + a.kappa * g) * H
    wp = a.w_A ./ a.scales_P
    wh = a.w_B ./ a.scales_H
    dx[4] = -(wp[1] * leaf_value(t, a) - wp[2] * u - wp[3] * u^2) * L
    dx[5] = -wh[1] * exp(-a.q_decay * t) * a.kappa * g * H + wh[2]*v^2*E
    dx[6] = u
    dx[7] = v
    dx[8] = u * L
    return nothing
end
\end{lstlisting}
Line~47 defines the leaf value $p(t)$, with an exponential alternative $e^{-\delta_At}$. Line~52 selects $\rho(L)$, and lines~53--55 implement \eqref{eq:tupelo-main}. Lines~56--59 accumulate the running parts of $J_A$ and $J_B$, including the fixed feature scales, and lines~60--62 accumulate the three allowance integrals.\footnote{Code names map to the notation of \eqref{eq:tupelo-main} as \code{u}$=u_A$, \code{v}$=u_B$, \code{g}$=\rho(L)$, \code{muE}$=\mu_E$, \code{muH}$=\mu_H$, \code{dE}$=d_E$, \code{dH}$=d_H$, \code{kappa}$=\kappa$, \code{p_decay}$=\delta_A$, and \code{q_decay}$=\delta_B$; \code{nothing} corresponds to $\sigma=0$.} The objective callback adds the terminal terms, and the constraint callback turns every enabled allowance into a terminal bound on the corresponding integral, owned by the restricted player:
\begin{lstlisting}[language=Julia,style=numbered,firstline=66,lastline=75,firstnumber=66]
function objective(x,player,a)
    return player==1 ? x[4]-a.w_A[4]/a.scales_P[4]*x[1] :
        x[5]-a.w_B[3]/a.scales_H[3]*x[3]
end

function constraints(a)
    bounds = StateConstraint[]
    for (player,index,bound,label) in ((1,6,a.U_A,"budget u₁"),(2,7,a.U_B,"budget u₂"),(1,8,a.C_A,"carbon budget"))
        bound===nothing || push!(bounds,StateConstraint(;player,index,label,upper=bound))
    end
\end{lstlisting}
The remaining lines of \code{constraints} only number the plot labels. The \code{game} constructor (lines~107--121 of the case file) collects these pieces with state and control labels, the ownership of $u_A$ by A and $u_B$ by B, the control bounds \eqref{eq:control-set}, the initial state, and the indices of the quadratures and running objectives.

A new game instance follows the same construction: define its parameters and callbacks, then provide a \code{Game}. A contributed instance consists of one case file with a \code{Parameters} structure, the dynamics, objective, and constraint callbacks, validation functions, and a \code{game} constructor. 

\subsection{Computing and auditing a reference candidate}\label{sec:tupelo-reference}
\begin{figure}[tb]
\centering
\includegraphics[width=\textwidth]{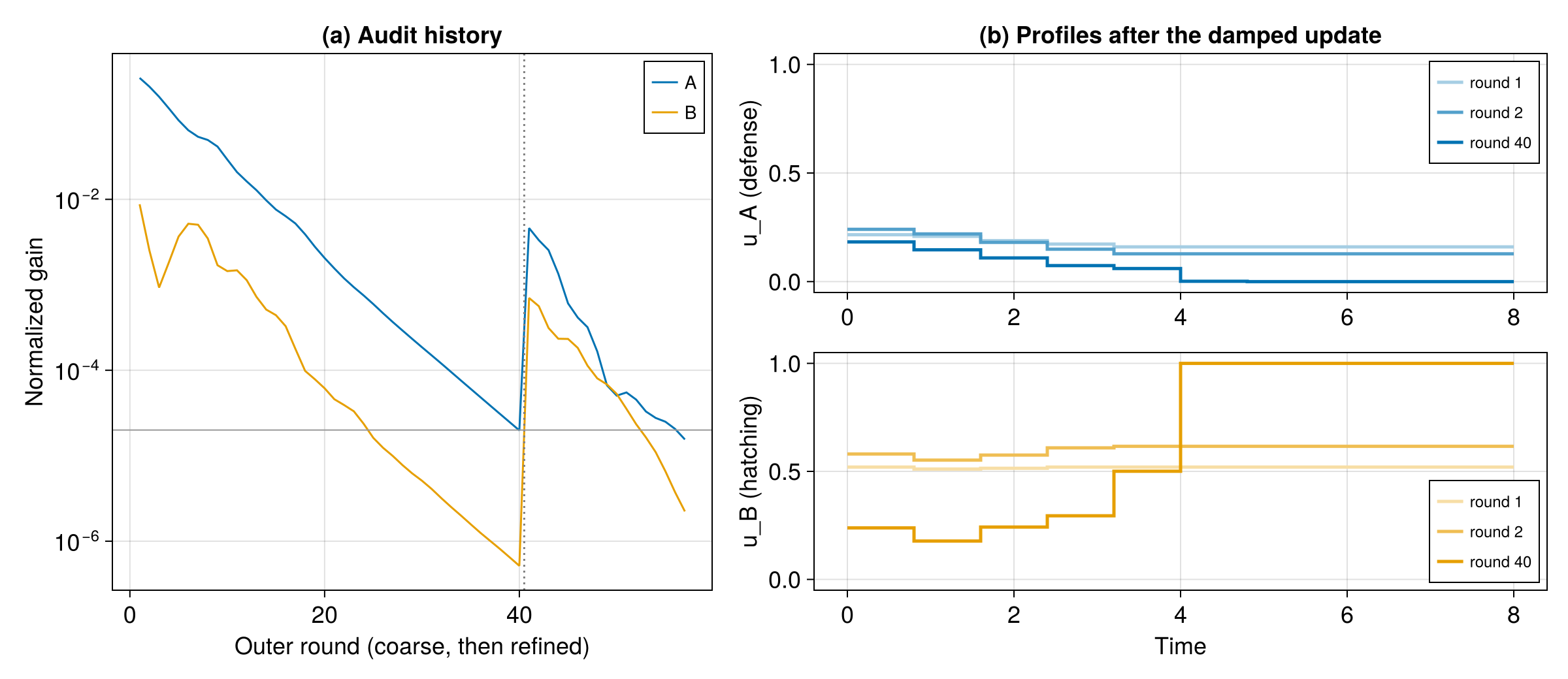}
\caption{Results of Algorithm~\ref{alg:baseline} applied to the tupelo reference instance. (a)~Normalized gains of the tree (Player A, blue) and the leafminer (B, ochre) per outer round on the coarse grid (rounds 1--40) and the refined grid (after the dotted line); the gray line marks the tolerance $\varepsilon_J$. Violations are not shown, since the instance has no restrictions. (b)~Profiles of defense $u_A$ (blue) and hatching $u_B$ (ochre) after the damped update of coarse rounds 1, 2, and 40, from light to dark.}\label{fig:tupelo-iterations}
\end{figure}

The specification of Section~\ref{sec:tupelo-setup} is solved by Algorithm~\ref{alg:baseline} (Section~\ref{sec:iteration}) and accepted by the audit of Section~\ref{sec:audit}:
\begin{lstlisting}[language=Julia]
g = Tupelo.game(a)
options = Tupelo.algorithm_options(intervals=10, refined_intervals=40)
result = run_game(g; options, folder="tupelo/teaching_output")
\end{lstlisting}
The instance uses the damping $\alpha=0.2$. With $\alpha=0.5$, the normalized gains of the tree oscillated between about $10^{-2}$ and $5\cdot10^{-2}$ over 200 rounds without converging; \citet{Low2013ColdWar} also report that simple best-response iteration did not converge for their model and use damping. Figure~\ref{fig:tupelo-iterations}(a) shows the audit history of the reference record. The coarse grid converges after 40 rounds and the refined grid after 17. No restriction is enabled, so all violations are zero, and the normalized gain of the tree determines the round count. Figure~\ref{fig:tupelo-iterations}(b) shows how the profiles change over the outer rounds on the coarse grid. The first rounds keep defense and hatching nearly constant in time. By round 40, defense on the coarse grid decreases to zero at $t=4$, and hatching is low early and maximal after $t=4$. The refined grid then resolves the early hatching pulse and the end of defense.

The reference candidate combines low, early defense with delayed mass hatching. The tree defends with at most $0.30$, mostly between $0.1$ and $0.2$, and stops defending at $t\approx4.4$. The leafminer releases a small hatching pulse at the start, hatches at a reduced hazard of about $0.2$--$0.3$ while defense lasts, and hatches at maximal hazard from $t\approx4.2$ until shortly before $T$. 
The corresponding states are shown in Figure~\ref{fig:tupelo-observations}. Larvae remain at low density, below $0.08$, with a maximum at $t\approx4.5$; eggs decline to $0.26$ by $t=4$ and are nearly exhausted by $t=6$. Leaf area falls to $0.50$; about half of the loss (46\%) occurs after $t=4$.

These features agree qualitatively with the coevolutionary results of \citet{Low2013ColdWar}: defense early in the season with a low larval density permitted, relaxation of defense before the end of the season, an early pulse of ``sneakers'', and mass hatching after defense has ended, which shifts consumption to later in the season. The onset of maximal hatching at $t\approx4.2$ corresponds to day 105 of the mapped season, close to the mass emergence on days 106--127 reported by \citet{Low2013ColdWar}. One difference remains: the leafminer does not avoid defense completely but hatches at a reduced rate while defense lasts, and defense is used rather than only threatened. We interpret this as a consequence of the open-loop information structure. 

The saved iteration history is also the unit of comparison for algorithm developers. Algorithm~\ref{alg:baseline} can be replaced by another forward dynamic-game solver while retaining the same dynamics, objectives, ownership, and restrictions. Alternatively, one may keep the response iteration and replace its control transcription or nonlinear-programming backend. Comparable records contain feasibility and unilateral-deviation evidence on the same instance, as in Figure~\ref{fig:tupelo-iterations}, together with runtime environment and discretization choices. The reference candidate is the baseline for the variants of Section~\ref{sec:tupelo-variants} and the data source of Section~\ref{sec:tupelo-inverse}.

\subsection{Alternative hypotheses as instance variants}\label{sec:tupelo-variants}
The reference candidate of Section~\ref{sec:tupelo-reference} results from specific modeling choices. A biological hypothesis about the interaction corresponds to an instance variant (Section~\ref{sec:specification}), i.e., to changes in one or more components of the specification in Section~\ref{sec:tupelo-setup}. For tupelo, \emph{preference} variants change the weights $w_A,w_B$ or the value functions $p,\omega$; \emph{restriction} variants change the switches $\sigma$ or the bounds $U_A,U_B,C_A$; and \emph{mechanism} variants change the dynamics \eqref{eq:tupelo-main}, e.g., the saturation $\rho$, the rates, or a fixed hatching schedule.
Many variants are possible. A positive terminal weight $w_B^3$ rewards larvae that mature after $T$ (preference variant); \code{C_A=1.0} enables the opponent-dependent carbon allowance, whereas \code{U_A=2.0} enables only an effort allowance (restriction variants); fixing $u_B$ to a prescribed schedule reduces the game to A's optimal-control problem (mechanism variant). A preference variant alters what a player prefers, whereas a restriction variant alters what a player may do; students can identify why enabling $C_A$ also changes how B's decisions affect A's feasible set.

As an example, the alternative instance \code{tupelo/tupelo_allowance.jl} expresses a different biological hypothesis: the tree's defense is limited by a carbon budget rather than by the declining value of its leaves, in a regime of mild defense lethality. It differs from the reference in all three components.
\begin{itemize}
\item \emph{Restrictions:} it enables the carbon allowance $C_A=1$ and the effort allowances $U_A=2$ and $U_B=3.5$.
\item \emph{Mechanism:} it uses the saturation $\rho(L)$ with $K=0.15$, egg mortality $\mu_E=0.025$, mild lethality $d_E=0.05, d_H=1.2$, and rates $b=0.55$, $\mu_H=0.1$, $\kappa=0.65$.
\item \emph{Preferences:} it uses the exponentially declining leaf value $e^{-0.2t}$, $\omega(t)=e^{-0.06t}$, and the weights $w_A=(100,8,15,0)/123$ and $w_B=(40,1,0)/41$.
\end{itemize}

The alternative converges in 19 rounds on 10 intervals and in 5 rounds on 40 intervals, with damping $\alpha=0.5$. On the coarse grid, the normalized gains fall below $\varepsilon_J$ from round~10 on, whereas the violation of the carbon allowance delays convergence until round~19; here feasibility, not optimality, governs the round count. The tree raises its defense to about $0.44$ almost immediately and lowers it to zero by $t\approx4.4$, when its carbon allowance is exhausted ($\int_0^Tu_AL\dt=C_A$). The effort allowance is not binding ($\int_0^Tu_A\dt\approx1.15<U_A$). The leafminer releases a first pulse, keeps hatching low while defense is strong, and hatches at maximal hazard during $t\approx3.4$--$5$, exhausting its hatching allowance. Larvae stay near $0.17$--$0.19$ during the first half of the season, and leaf area levels off at $0.54$.

The two candidates support different biological readings. In the reference, defense ends because the remaining value of the leaves no longer justifies it, larval density stays low while defense lasts, and 46\% of the leaf loss occurs after $t=4$. In the alternative, defense ends because the carbon budget is exhausted, larvae feed at high density from the start, and only 33\% of the leaf loss occurs after $t=4$. The alternative thus describes an ongoing conflict in the early part of the season, whereas the reference shifts herbivory to later in the season, as in \citet{Low2013ColdWar}. Figure~\ref{fig:tupelo-observations} compares both candidates with synthetic observations.

\subsection{Synthetic observations and inverse games}\label{sec:tupelo-inverse}
\begin{figure}[tb]
\centering
\includegraphics[width=.85\textwidth]{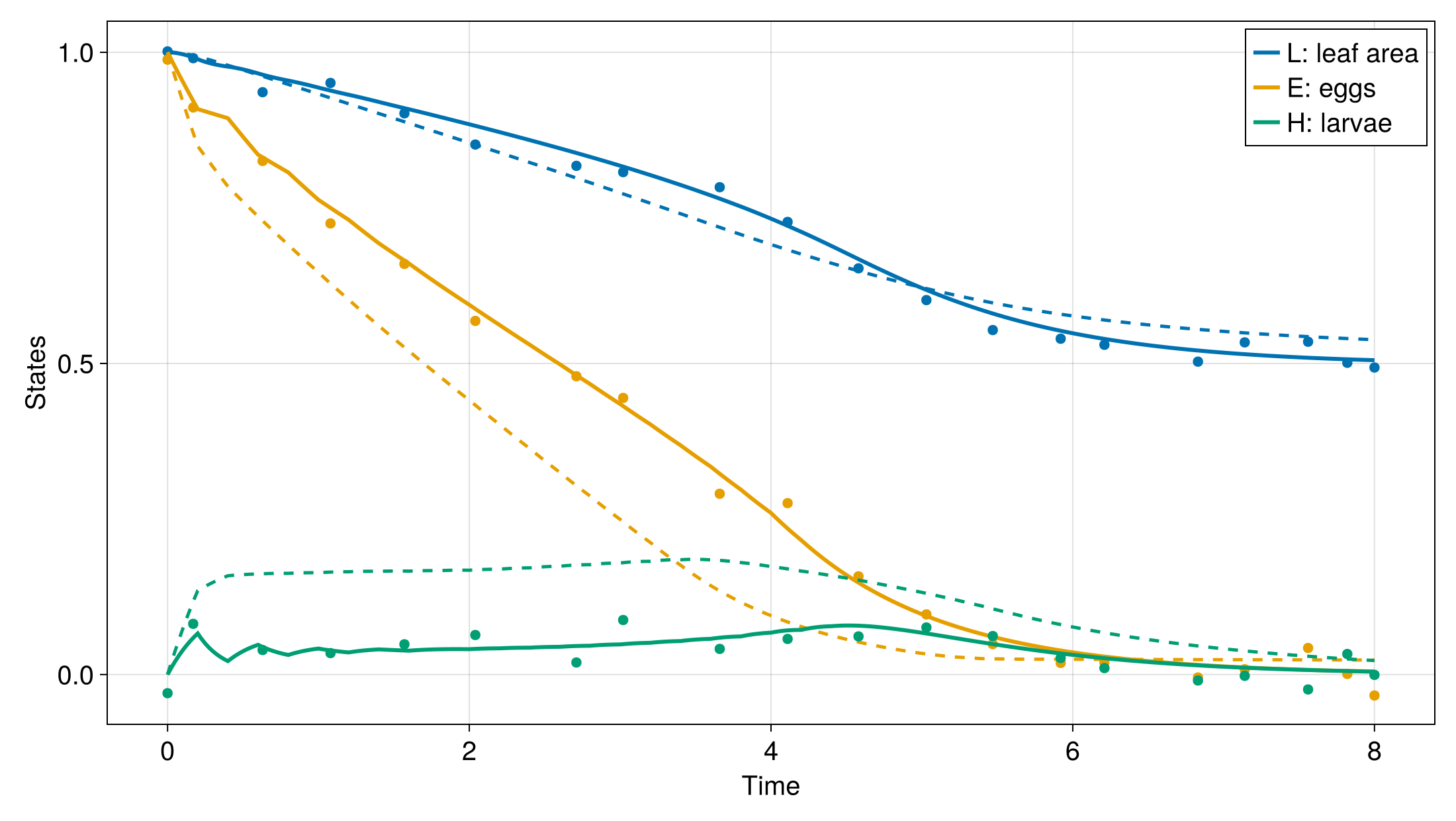}
\caption{States of the tupelo reference candidate (solid lines) and of the allowance alternative of Section~\ref{sec:tupelo-variants} (dashed lines), together with the synthetic observations generated from the reference candidate (points). The observation-generation script specifies times, noise level, and seed, so this experiment can be repeated from the saved trajectory.}\label{fig:tupelo-observations}
\end{figure}

The saved reference candidate of Section~\ref{sec:tupelo-reference} can generate observations without another optimization. The supplied \code{tupelo_synthetic.jl} linearly interpolates all three biological states at 20 irregular times and adds independent Gaussian noise with standard deviation $0.02$ and seed $20260908$:
\begin{lstlisting}[language=Julia]
include("tupelo/tupelo_synthetic.jl")
main()  # writes tupelo/output/tupelo_synthetic.csv
\end{lstlisting}
Figure~\ref{fig:tupelo-observations} shows these observations with the reference candidate. The general \code{create_synthetic_data} utility additionally accepts observation times, a random seed, and state-specific noise levels.

The variants of Section~\ref{sec:tupelo-variants} define the hypothesis space of an inverse problem. \citet{BangaVillaverde2026Hierarchical} arrange biological inverse problems in four levels of increasing ambition: (1) parametric identification of a known model, (2) model discovery, in which structure and parameters are inferred, (3) inverse optimal control, which infers the objective of a single decision maker, and (4) inverse differential games, which infer the individual objectives of several interacting agents. In the first two levels the controls are known inputs; in the upper two they are the behavior to be explained. The tupelo variants map onto this hierarchy. Mechanism variants with changed rates are Level-1 problems when the schedules $u_A$ and $u_B$ are observed, and a changed saturation $\rho$ is a Level-2 question. Fixing $u_B$ to its observed schedule reduces the game to the tree's optimal-control problem, whose unknown weights $w_A$ form a Level-3 problem. Preference variants with unknown $w_A$ and $w_B$ are the Level-4 problem, studied for open-loop equilibria by \citet{Molloy2020-vg,molloy2022inverse}. Restriction variants add a dimension that the four levels do not separate: generalized inverse optimal control treats restrictions as model components to be inferred \citep{Banga2025}, in Levels 3 and 4 alike. Figure~\ref{fig:tupelo-observations} also shows the candidate of the allowance alternative, which differs from the reference in all three components. Differences in the data can therefore not be attributed to one component; \citet{BangaVillaverde2026Hierarchical} call such trade-offs between objective and mechanism cross-level compensation. Leaf area separates the candidates only weakly: the trajectories differ by at most $0.046$. Eggs and larvae separate them clearly, with differences of up to $0.19$ and $0.14$. The root-mean-square residuals of the alternative with respect to the observations are $0.10$ for eggs and $0.09$ for larvae, compared with $0.02$ for the reference and the noise standard deviation $0.02$. Whether such data suffice to identify individual components is an identifiability question that this forward comparison does not answer.

A controlled inverse benchmark on this instance can proceed as follows. The dynamics and the observation model, i.e., observed states, times, and noise distribution, are given, and a subset of the specification is hidden, for example $w_A$ and $w_B$ (preferences) or the switches and bounds of the allowances (restrictions); in terms of Section~\ref{sec:specification}, the hidden entries of $\mathcal G$ are the unknowns. A method succeeds if it recovers the hidden components within a stated tolerance and if the recovered game, solved and audited as in Sections~\ref{sec:tupelo-reference} and~\ref{sec:audit}, reproduces the observations. Because generalized Nash equilibria are typically not isolated (Section~\ref{sec:equilibrium}), a recovered specification can reproduce the data through another equilibrium of the same game, and a correct specification can fail to reproduce them through equilibrium selection. \citet{BangaVillaverde2026Hierarchical} name further Level-4 difficulties that the library allows to study in a controlled way: permutation symmetry of players with similar roles, cooperative mimicry of a single aggregate optimizer, and an unknown number of players, which the \cat{variablePlayers} instances can vary. Reports should therefore state the forward solver, its initializations, and whether several observation scenarios or initializations were used to assess ambiguities between weights, restrictions, and equilibrium selection. No inverse experiment is reported here.

\FloatBarrier

\section{Discussion and outlook}\label{sec:outlook}
The library provides 25 executable game specifications from the molecular to the population scale, each with a parameter instance, an audited reference candidate, and an interpretation of the candidate. The common specification makes preferences, restrictions, and decision ownership explicit, so that alternative biological hypotheses become instance variants that can be solved and compared with the same method.

The tupelo example in Section~\ref{sec:uses} shows this use in a chain of steps: an executable specification yields an audited reference candidate; changes of preferences, restrictions, or mechanisms turn biological hypotheses into instance variants; and the candidates generate synthetic observations whose hypothesis space for inverse games is defined by these variants. The reference instance, mapped from a published scenario of \citet{Low2013ColdWar}, reproduces their delayed mass hatching qualitatively, but not their complete conflict avoidance, which we attribute to the open-loop information structure: in their model, the plant can respond to the current larval density with induced defense, so that the threat of defense alone deters early hatching, whereas a schedule fixed at the start of the season cannot respond to larval density, and defense must therefore be used rather than only threatened.

The reference candidates have the limitations stated in Sections~\ref{sec:assumptions} and~\ref{sec:audit}: they are open-loop, need not be unique, and are accepted on two control grids by local solves; they carry no certificate of an exact continuous-time equilibrium. Most parameters are synthetic. The sequential iteration has no general convergence guarantee, and in some periodic instances only a looser return band yields convergence.

These limitations define the next steps. Inverse games can start from the hypothesis spaces of Section~\ref{sec:uses}, ideally with real-world data. Feedback equilibria for selected instances would show whether conclusions such as the tupelo comparison depend on the information structure, and hence which hypotheses remain plausible for given synthetic or real-world data. Simultaneous or extrapolated updates, multiple shooting, and methods for equilibrium selection can be compared on the same audits. Calibrated instances, multiseason models with an explicit offseason map, and scalability studies with a variable number of players extend the library in further directions.
The code, the saved candidates, and the scripts that generate all tables and figures are available (see the code availability statement). Each instance is also documented in the category of dynamic games of the benchmark library \url{https://mintoc.de} \citep{Sager2012-ts}, whose category system includes the classifications of Table~\ref{tab:categories}, so that new instances and methods can be added and compared with the reference results.

\backmatter
\section*{Statements and Declarations}
\paragraph*{Funding.}
The project has been funded by
 the European Regional Development Fund under the European Union's Horizon Europe Research and Innovation Program via grants
 intelAlgen ZS/2023/12/182064,
 timingMatters ZS/2023/12/182063,
and Center for Dynamic Systems ZS/2023/12/182075,
which we gratefully acknowledge.

\paragraph*{Competing interests.} None.

\paragraph*{Ethics approval and consent to participate.} Not applicable.

\paragraph*{Consent for publication.} Not applicable.

\paragraph*{Data availability.}
No experimental data were used. The saved reference candidates (controls, trajectories, iteration histories, and solver logs) and the synthetic observations of Section~\ref{sec:uses} are part of the code repository.

\paragraph*{Materials availability.} Not applicable.

\paragraph*{Code availability.}
The solver, the implementations of all 25 game instances, and the scripts that generate all tables and figures of this paper are available at \url{https://github.com/SciML/Corleone.jl} with fixed package versions.

\paragraph*{Author contributions.}
Conceptualization: SS; Methodology: SS; Software: JM, CP, SS; Validation: SS; Writing -- original draft: SS; Writing -- review and editing: all authors; Supervision: SS; Funding acquisition: SS; Biological interpretation: FO, KR, SS.

\paragraph*{Use of generative AI.}
The project was developed with extensive use of Claude Code 2.1.285 using Opus 5.5. The agent's behavior was constrained by instructions that prohibited invented references, data, and results and required unverified statements to be marked. The authors specified all tasks in approximately 100 human-in-the-loop prompts. The agent implemented parts of \code{CorleoneGame.jl}; formulated and implemented the game instances from biological settings chosen by the authors; ran all computations and generated tables and figures; searched the literature; and drafted and revised parts of the text. The authors checked all model equations, provenance statements, references, and the overall text. They take full responsibility for the contents of this article.

\appendix
\clearpage
\section{Game specifications}\label{app:cases}
Each game instance uses the notation from Section~\ref{sec:formulation}.
Players are indexed by $\mathcal P=\{A,B,\ldots\}$ with $|\mathcal P|=N$;
player-indexed tuples and parameter rows follow this alphabetical order.
Numeric subscripts on the full state vector $x$ or numerical state $y$ denote
component positions, not player names. The position map is
$\iota(A)=1,\iota(B)=2,\ldots$.
Local biochemical and resource indices retain their stated meanings.
Unless a game instance defines a smaller group, sums over player indices run over $\mathcal P$.

Parameter listings reproduce the Julia parameter fields, with tuples and matrices
in code order. The correspondence to mathematical notation should be clear from the context.
The fixed horizon is \code{T}; for free-time instances,
\code{duration} initializes the optimized horizon reported in the numerical record.
Initial physical values and bounds are specified by these fields; objective and
budget quadratures start at zero. The field \code{initial} is the initial control guess of the solver, and \code{numerical_x0} is the initial numerical state including quadratures. A \code{nothing} entry disables an optional
restriction. Enabled-bound tables use numerical-state coordinates $y$, including
quadratures and any log transformation. Return bands are written in physical biological coordinates $x$; control boxes are listed separately.
A local $g_i\geq0$ slack convention is equivalent to \eqref{eq:constraints}
with the same sign convention. Categories describe the enabled instance.

\myheader{Reading the figures.}
The first panel shows the displayed differential states. The second shows piecewise-constant controls; parentheses identify their owner. The third, when present, shows constraint slacks. A lower bound has slack $x_j-\underline x_j$ and an upper bound has slack $\overline x_j-x_j$. A shared restriction is drawn once and lists all owners. Slack labels $g_1,g_2,\ldots$ follow the row order of the instance's enabled-bounds table; return bands are omitted from the numbering. Solid curves denote path restrictions; dotted curves and terminal markers denote endpoint restrictions. Values along a dotted curve before $T$ are contextual, not additional path requirements.
The label is \emph{violated} below $-10^{-4}$, \emph{active} up to $10^{-3}$, and otherwise \emph{inactive}, based on the minimum path slack or terminal slack. Equalities use the absolute residual with a negative sign for this classification. These display thresholds differ from the solver's acceptance tolerance. Return bands are hidden in the plot. Axes show physical time.

\myheader{A common template for transfer tasks.}
Three free-duration transfer tasks (\caseref{batchFermentation}, \caseref{seedGermination}, \caseref{woundRepair}) share one template: a duration cost for the owner, squared effort, an instance-specific running cost, a terminal target shared by all players, private effort budgets, and one path restriction per player. These instances are deliberately weakly coupled: an equilibrium candidate is found within a few outer rounds, the path restrictions are synthetic and inactive in the reference records, and the timing is governed by the shared target and the owner's time cost. They serve as easy benchmark instances and as a controlled comparison of coupling mechanisms under a common structure. The fourth transfer task, \caseref{immuneClearance}, differs: the clearance target restricts only the host.

\clearpage
\renewcommand\floatpagefraction{.5}
\begingroup
\def\caseidentifier{algae}%
\subsection{algae: Carbon-uptake allocation in an algae community}\label{case:algae}
\myheader{Biological interpretation.} Microalgae in lakes and cultures take up inorganic carbon either as dissolved CO$_2$, which
diffuses passively but is scarce at high pH, or as bicarbonate, which requires active transport
through carbon-concentrating mechanisms. Species differ in their affinities for the two carbon
forms, and cells can shift the allocation of their uptake machinery between them. When the carbon
supply fluctuates, for example through a recurring feeding or mixing regime in a maintained
culture, the community goes through cycles of bloom and decline in which the relative
availability of CO$_2$ and bicarbonate changes. We consider one full supply cycle of a maintained
culture: biomasses and resources must return to their initial levels, and each strain maximizes
its average gross production while avoiding extreme specialization.

\cite{Zhou2022Algae} study monocultures and pairwise competition among
\emph{Scenedesmus quadricauda}, \emph{Chlorella vulgaris},
\emph{Synedra ulna} and \emph{Phormidium} sp. (a cyanobacterium),
under different atmospheric CO$_2$ levels.
Their model includes carbon uptake, quotas, carbonate chemistry and light.

Proposed player $i$: one species, allocating uptake machinery
between dissolved CO$_2$ and bicarbonate using $u_i$.
All species share the resource pools. This provides a community
benchmark for either pairs or $N=4$ players.
The source has fixed uptake traits and does not solve a differential
game. The controlled reduction and objectives below are our hypotheses.

The four model affinity types in the source are abstract types; they must
not be assigned to the four named species without calibration. Players A--D use the types
HCHH, HCLH, LCHH, and LCLH.
\begin{figure}[p]
\centering
\includegraphics[width=\textwidth,height=\casefigheight,keepaspectratio]{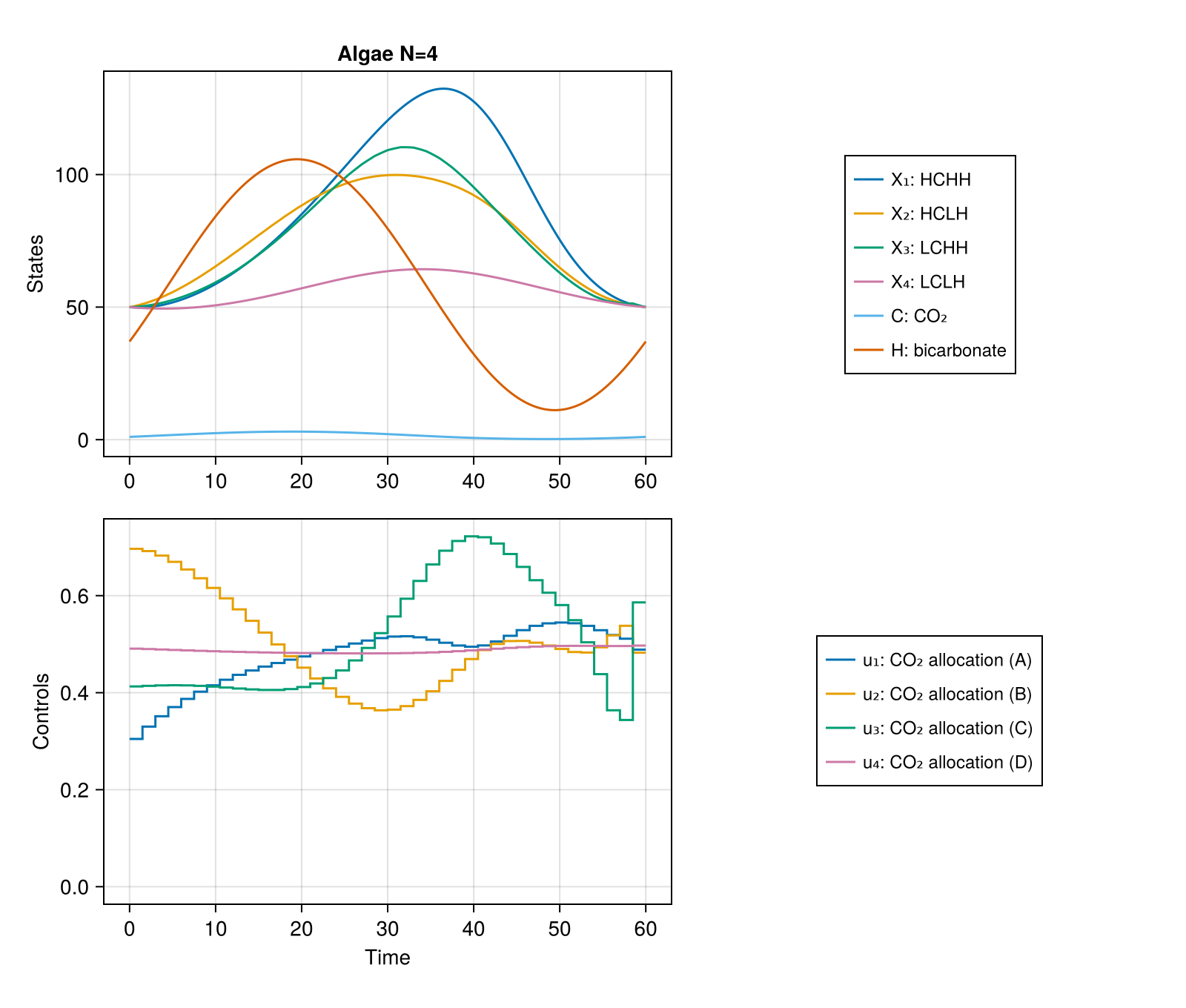}
\caption[Reference equilibrium candidate for \texttt{\caseidentifier}]{Reference equilibrium candidate for \texttt{\caseidentifier}. \textit{Interpretation.} The carbon supply peaks in the first half of the cycle, so CO$_2$ and bicarbonate rise to maxima
near $t\approx20$ ($3.0$ and $106$) and fall to minima near $t\approx48$--$50$. The four strains
respond with a delayed bloom: biomasses increase up to a factor of $2.65$ (strain A) and peak
between $t\approx31$ and $t\approx37$, after which resource depletion and declining supply cause
a synchronous decline back to the initial level. Strain D, with low affinities for both carbon
forms, grows least and keeps a nearly balanced allocation. The allocation strategies differ
according to the uptake traits: B, with low bicarbonate affinity, starts strongly CO$_2$-oriented
and shifts toward bicarbonate while bicarbonate is abundant, whereas A and C increase their CO$_2$
allocation during and after the bloom, when bicarbonate becomes scarce; C reaches the strongest
CO$_2$ preference ($0.72$) near $t\approx39$. No additional path or budget restrictions are enabled, and the
return conditions close all biomasses and resource pools within the tolerance.}\label{fig:\caseidentifier}
\end{figure}
\myheader{Literature and provenance.} The source's four abstract uptake types have the following numerical values:
\begin{center}
\begin{tabular}{lrrrr}
\hline
Type & $K_C$ & $K_H$ & $V_C$ & $V_H$\\\hline
HCHH & 1 & 30 & 0.2 & 0.2\\
HCLH & 1 & 1200 & 0.2 & 0.4\\
LCHH & 40 & 30 & 0.4 & 0.2\\
LCLH & 40 & 1200 & 0.4 & 0.4\\\hline
\end{tabular}
\end{center}
For this \emph{dimensionless} dynamic game, we use the same numerical
contrasts, not interpreting them as a calibrated conversion of source units.
Rate units of the source are not needed for this dimensionless reduction.

Its carbonate/quota model is richer than our reduction. Fitting the proposed
ODE to named-species observations requires new parameter estimation and
checking whether the reduced dynamics are adequate.
\myheader{Dimensions.} $N=4$, $n_x=6$, $n_{\rm ODE}=10$, $(n_{u,i})_{i\in\mathcal P}=(4\times 1)$, $(n_{g,i})_{i\in\mathcal P}=(4\times 12)$.
\myheader{Categories.} \cat{cellular}, \cat{periodical}, \cat{competition}, \cat{variablePlayers}.
\myheader{Game specification.} Let $X_i$ be biomass and $C,H$ the CO$_2$ and bicarbonate carbon pools.
Use consistent dimensionless units for this initial benchmark:
\[
q_{C,i}=u_iV_{C,i}\frac{C}{K_{C,i}+C},\quad
q_{H,i}=(1-u_i)V_{H,i}\frac{H}{K_{H,i}+H},\quad
\mu_i=Y_i(q_{C,i}+q_{H,i}).
\]
\[
\begin{aligned}
\dot X_i&=(\mu_i-d_i)X_i,\qquad i\in\mathcal P,\\
\dot C&=\lambda(C_{\rm in}(t)-C)-k_+C+k_-H-\sum_{i\in\mathcal P}q_{C,i}X_i,\\
\dot H&=\lambda(H_{\rm in}(t)-H)+k_+C-k_-H-\sum_{i\in\mathcal P}q_{H,i}X_i.
\end{aligned}
\]
This $N+2$ state reduction omits cell quotas, variable pH and self-shading;
$k_+,k_-$ approximate conversion at buffered conditions.
Each uptake vanishes at depletion. Conversion cancels in $\dot C+\dot H$.
Biomass loss leaves the modeled live-carbon pools.

Let $B_0>0$ be a fixed biomass scale. Player $i$ minimizes
\[
J_i=-(w_i^1+w_i^3)\frac1T\int_0^T\frac{\mu_iX_i}{B_0}\,\mathrm dt
+w_i^2\frac1T\int_0^T[u_i^2+(1-u_i)^2]\,\mathrm dt.
\]
\[
w_i\geq0,\quad\sum_{k=1}^3w_i^k=1.
\]
The first term is the average gross biomass production. In the periodical setting the
terminal biomass is fixed by the return band, so a terminal-biomass term would be nearly
constant; its weight $w_i^3$ is therefore added to the production weight. The quadratic term penalizes specialization of
uptake investment; it is not an observed metabolic expenditure. Production benefits create
incentives to deplete resources needed by competitors.

\emph{Periodical setting.} The inflow concentrations are modulated as
$C_{\rm in}(t)=C_{\rm in}[1+a\sin(2\pi t/T_{\rm per})]$ and $H_{\rm in}(t)=H_{\rm in}[1+a\sin(2\pi t/T_{\rm per})]$
with forcing period $T_{\rm per}=T$. The fixed initial state lies on a reference cycle obtained with constant
allocations: phase-zero biomasses are fixed, and loss rates $d_i$ and initial resource pools
are calibrated so that this reference cycle closes. Every player's problem contains the return
conditions $|x_j(T)-x_j(0)|\leq\delta$ for all biomasses and resource pools. The tolerance
$\delta=10^{-4}$ is absolute; since bicarbonate reaches $H\approx106$ and CO$_2$ about $3$,
the band is relatively much tighter for $H$ than for the biomasses.

Enabled bounds:
\begin{center}\small\begin{tabular}{lll}\toprule Owners & Bound & Location\\\midrule
A, B, C, D & $|x_j(T)-x_j(0)|\leq 0.0001$, $j=1,\ldots,6$ & terminal (return) \\
\bottomrule\end{tabular}\end{center}
\myheader{Parameter instance.} %
\begin{lstlisting}
types = [:HCHH, :HCLH, :LCHH, :LCLH]
forcing_period = 60.0
T = 60.0
forcing_amplitude = 0.9
periodicity_tolerance = 0.0001
B0 = 1.0
K_C = [1.0, 1.0, 40.0, 40.0]
K_H = [30.0, 1200.0, 30.0, 1200.0]
V_C = [0.2, 0.2, 0.4, 0.4]
V_H = [0.2, 0.4, 0.2, 0.4]
Y = [1.0, 1.0, 1.0, 1.0]
d = [0.1139088931, 0.06327295407, 0.06733780299, 0.01670186393]
lambda = 0.2
k_plus = 0.3
k_minus = 0.01
C_in = 2.0
H_in = 60.0
x0 = [1.0, 1.0, 1.0, 1.0, 1.046509686, 37.02538762]
weights = [0.45 0.1 0.45; 0.45 0.1 0.45; 0.45 0.1 0.45; 0.45 0.1 0.45]
initial = [0.5, 0.5, 0.5, 0.5]
control_lower = [0.0, 0.0, 0.0, 0.0]
control_upper = [1.0, 1.0, 1.0, 1.0]
numerical_x0 = [1.0, 1.0, 1.0, 1.0, 1.046509686, 37.02538762, 0.0, 0.0, 0.0, 0.0]
\end{lstlisting}

\endgroup

\begingroup
\def\caseidentifier{antColonies}%
\subsection{antColonies: Competition for food, water, territory}\label{case:antColonies}
\myheader{Biological interpretation.} Neighboring ant colonies in arid habitats compete for food and water during the active season:
workers leave the nest to forage, and colonies defend foraging areas against neighbors through
territorial aggression, which reduces the rivals' harvest but costs workers. Food brought home is
stored and converted into brood, which later matures into new workers, so a colony's allocation of
its workforce has delayed consequences. We model one foraging season of three colonies that share
food and water sources; each colony allocates its workers between food foraging, water collection,
and aggression, subject to an effort budget and a minimum food store that must be maintained throughout the
season.

Players are $N\geq2$ colonies. Shared states $F,W$ are available food and water;
$A_i,B_i,C_i$ are workers, brood, and stored food.
$x=(F,W,(A_i)_{i\in\mathcal P},(B_i)_{i\in\mathcal P},(C_i)_{i\in\mathcal P})$.
Each colony chooses nonnegative intensities $(f_i,w_i,a_i)$ for food foraging,
water collection, and territorial activity. Define worker shares
\[
(s_{fi},s_{wi},s_{ai})=\frac{(f_i,w_i,a_i)}{1+f_i+w_i+a_i}.
\]
The residual share performs nest tasks; total external allocation is below one.
All controls depend on time alone.
\begin{figure}[p]
\centering
\includegraphics[width=\textwidth,height=\casefigheight,keepaspectratio]{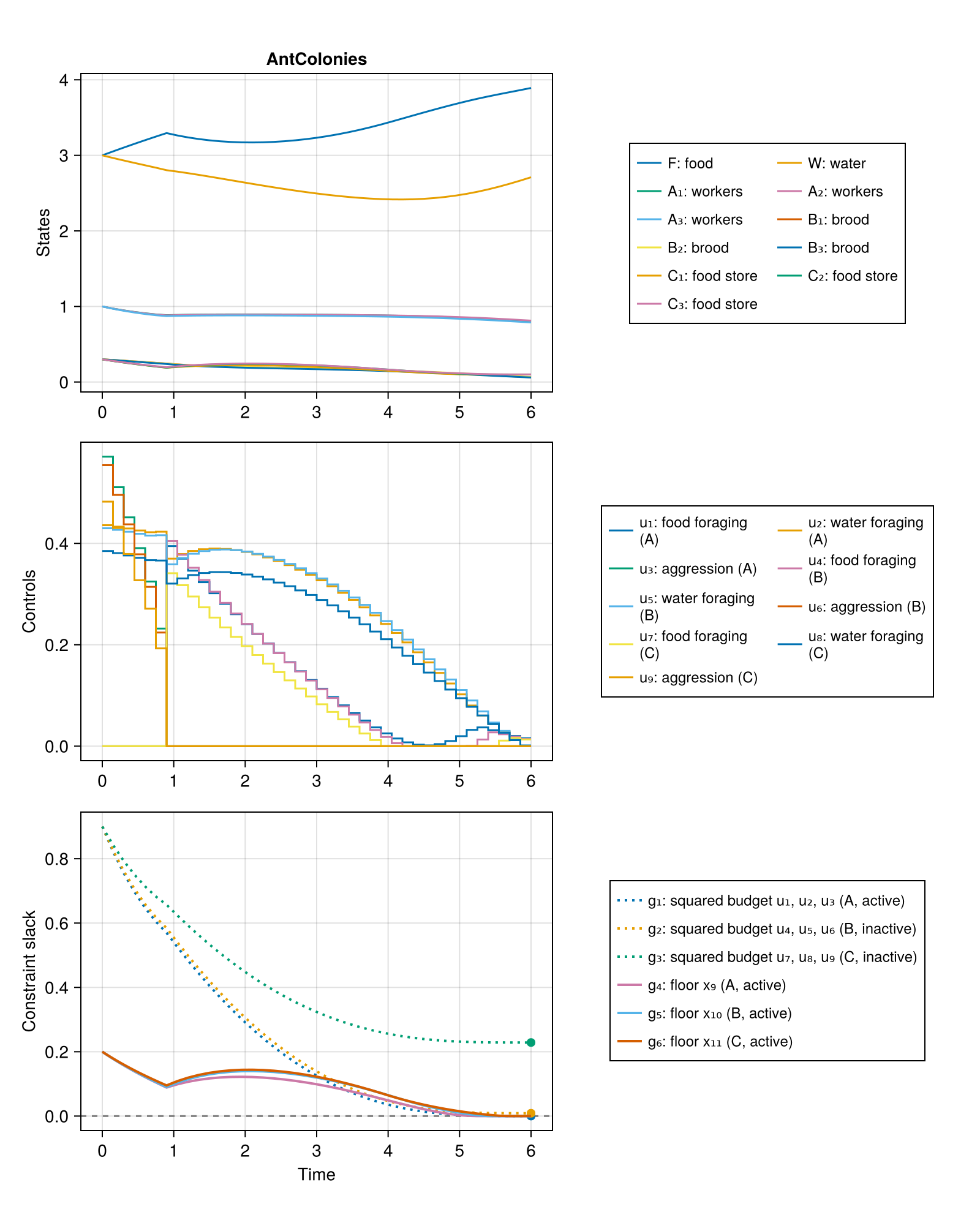}
\caption[Reference equilibrium candidate for \texttt{\caseidentifier}]{Reference equilibrium candidate for \texttt{\caseidentifier}. \textit{Interpretation.} The season opens with a territorial contest: during the first time unit all three colonies invest
heavily in aggression (up to $0.57$) together with water collection, while food foraging is zero.
Because conflict mortality of a colony decreases with its own aggression ($\chi Q_i/(1+s_{ai}A_i)$), early aggression is partly defensive: each colony protects its workforce against the others' aggression, a mutual arms race that ends once worker losses outweigh the protection.
Aggression declines quickly and stops at $t\approx0.9$, after which food foraging starts abruptly
and water collection continues. Food foraging then tapers off first, ending near $t\approx4$, and
water collection declines more slowly until the end of the season; colony C, which weights effort
costs most, forages least. Because nobody harvests food during the contest, the shared food
source rises at the start, dips slightly while foraging is intense, and then accumulates again;
water is drawn down to a minimum near $t\approx4.2$ and recovers once collection slows. Conflict
and foraging risk reduce the worker forces by about $20\%$, brood declines from $0.3$ to about
$0.06$ as it matures into workers faster than stored food and water produce new brood, and the food stores are run down to their
protective floor at the end of the season ($g_4$--$g_6$ active). Colony A also exhausts its effort
budget ($g_1$ active), B nearly so, whereas C retains unused effort.}\label{fig:\caseidentifier}
\end{figure}
\myheader{Literature and provenance.} The source \citep{AdlerGordon2003Ants} models resource exploitation and interference between neighboring ant colonies. We retain these mechanisms in an aggregate dynamic game, adding water limitation, brood maturation, and food stores as explicit proposed components. We do not reproduce the source spatial model or claim its fitted parameters.

Optimization and conflict among foraging ranges of neighboring colonies are analyzed by \citet{AdlerGordon2003Ants}. A differential-game analysis of territorial space use, for wolf packs, is given by \citet{HamelinLewis2010Territoriality}.

All numerical defaults are synthetic and dimensionless; no fitted biological values or source equilibrium results are asserted.

\citep{BenchmarkantWater} reports water collection under experimental water stress. This supports treating water acquisition as a distinct colony task.
\citep{BenchmarkantResource} links food collection, internal storage, and redistribution. Our single food-store compartment aggregates these processes.
\citet{Rao2020NutritionalAnts} provide a model linking nutrition, brood rearing, and colony dynamics. Our maturation and production terms are a simpler proposed representation.

These components come from different species or model settings. Their combination
is a generic benchmark, not a species-specific reconstruction; worker castes
and individual communication are not resolved.
\myheader{Dimensions.} $N=3$, $n_x=11$, $n_{\rm ODE}=17$, $(n_{u,i})_{i\in\mathcal P}=(3\times 3)$, $(n_{g,i})_{i\in\mathcal P}=(3\times 2)$.
\myheader{Categories.} \cat{population}, \cat{seasonal}, \cat{competition}, \cat{variablePlayers}.
\myheader{Game specification.} Rival territorial pressure is $Q_i=\sum_{j\in\mathcal P\setminus\{i\}}s_{aj}A_j$.
Food intake $H_i$, water collection $V_i$, and brood production $G_i$ are
\[
H_i=\frac{\alpha_i s_{fi}A_iF}{(h+F)(1+\gamma Q_i)},\qquad
V_i=\frac{\beta_i s_{wi}A_iW}{h+W},
\]
\[
G_i=bA_i\frac{C_i}{h+C_i}\frac{V_i}{h+V_i}.
\]
$h=\texttt{half}$ is the common normalized saturation scale;
$\alpha_i,\beta_i$ are collection rates; $\gamma$ is interference strength.
Brood production requires both stored food and current water acquisition.
\[
\dot F=I_F-\delta F-\sum_{i\in\mathcal P}H_i,\qquad
\dot W=I_W-\delta W-\sum_{i\in\mathcal P}V_i.
\]
$I_F,I_W$: external supplies; $\delta$: resource loss.

\[
\dot A_i=\eta B_i-
\left(d_A+r_f s_{fi}+\frac{\chi Q_i}{1+s_{ai}A_i}\right)A_i,
\]
\[
\dot B_i=G_i-(\eta+d_B)B_i,\qquad
\dot C_i=H_i-d_CC_i-c_BG_i.
\]
$\eta$: maturation; $d_A,d_B,d_C$: worker, brood, and storage losses;
$r_f$: foraging risk; $\chi$: conflict mortality; $c_B$: food cost per brood
(parameter fields \texttt{maturation, worker\_loss, brood\_loss, store\_loss,
foraging\_risk, conflict, brood\_cost}).
Each loss flux vanishes when the corresponding stock is zero;
$G_i=0$ at $C_i=0$. Thus nonnegative states are invariant.

Each player minimizes
\[
J_i=\frac1T\int_0^T\sum_{k=1}^3w_i^k\ell_{ik}\,dt,
\qquad w_i^k\geq0,\quad\sum_kw_i^k=1.
\]
There is no terminal objective. Each weight row is independently configurable.
With $E_i=f_i^2+w_i^2+a_i^2$,
\[
(\ell_{i1},\ell_{i2},\ell_{i3})=(-B_i,\ E_i,\ -A_i).
\]
The colony trades brood production, costly recruitment intensity, and
worker retention. Squared intensity is a behavioral cost proxy, not a literal
worker count or a calibrated metabolic expenditure.

$\dot Z_i=E_i$, $Z_i(0)=0$, $Z_i(T)\leq\mathrm{budgets}_i$.
Food reserve path restriction: $C_i(t)\geq\mathrm{floors}_i$.
The default three-colony game has three recruitment budgets and three path
reserve floors. Competition affects reserve feasibility.

For $q$ biological states, append $N$ cost integrals then $N$ budget integrals,
all initially zero; the cost integrals include the factor $1/T$, the ODE dimension is $q+2N$, and $J_i=x_{q+\iota(i)}(T)$.

Enabled bounds:
\begin{center}\small\begin{tabular}{lll}\toprule Owners & Bound & Location\\\midrule
A & $y_{15}(T)\leq 0.9$ & terminal \\
B & $y_{16}(T)\leq 0.9$ & terminal \\
C & $y_{17}(T)\leq 0.9$ & terminal \\
A & $0.1\leq y_{9}(t)$ & path \\
B & $0.1\leq y_{10}(t)$ & path \\
C & $0.1\leq y_{11}(t)$ & path \\
\bottomrule\end{tabular}\end{center}
\myheader{Parameter instance.} %
\begin{lstlisting}
N = 3
T = 6.0
food_supply = 0.8
water_supply = 0.7
resource_loss = 0.15
food_rate = [0.9, 1.0, 1.1]
water_rate = [0.8833333333, 0.8666666667, 0.85]
half = 0.4
interference = 6.0
birth = 1.2
maturation = 0.6
brood_loss = 0.08
worker_loss = 0.12
foraging_risk = 0.08
conflict = 0.6
store_loss = 0.1
brood_cost = 0.8
x0 = [3.0, 3.0, 1.0, 1.0, 1.0, 0.3, 0.3, 0.3, 0.3, 0.3, 0.3]
weights = [0.15 0.15 0.7; 0.15 0.2 0.65; 0.15 0.25 0.6]
budgets = [0.9, 0.9, 0.9]
floors = [0.1, 0.1, 0.1]
initial = [0.2, 0.2, 0.2, 0.2, 0.2, 0.2, 0.2, 0.2, 0.2]
control_lower = [0.0, 0.0, 0.0, 0.0, 0.0, 0.0, 0.0, 0.0, 0.0]
control_upper = [1.0, 1.0, 1.0, 1.0, 1.0, 1.0, 1.0, 1.0, 1.0]
numerical_x0 = [3.0, 3.0, 1.0, 1.0, 1.0, 0.3, 0.3, 0.3, 0.3, 0.3, 0.3, 0.0, 0.0, 0.0, 0.0, 0.0, 0.0]
\end{lstlisting}

\endgroup

\begingroup
\def\caseidentifier{avoidance}%
\subsection{avoidance: Avoidance and decentralized control}\label{case:avoidance}
\myheader{Biological interpretation.} In signaled active avoidance, an animal learns that a warning stimulus predicts an aversive
event and that a timely action prevents it. The episode considered here is a single trial block:
a warning cue is presented for the first half of the horizon and then withdrawn. Two
decentralized control modules (for example, two motor or motivational subsystems) can each
invest effort in avoidance actions. Acting builds the readiness to avoid, but it consumes a
shared motivational resource that recovers only slowly, and each module also benefits from the
other's readiness. The instance is therefore an episodic game between partially aligned modules
(partner benefit weight $w_i^4$) that must time their effort relative to the cue. The episode
concerns the expression of a learned avoidance response within one trial block; learning across
trials is not modeled.

Dopaminergic signaling links warning stimuli, avoidance actions, and delayed
reinforcement in acoustically signaled avoidance learning
\cite{Ilango2012DopamineAvoidance}. Work on biologically inspired adaptive control
emphasizes decentralized modules whose actions change the state experienced by
neighboring modules \cite{Schilling2019FluidAdaptivity,Schilling2020Hexapod}.
Neither line of work formulates a dynamic game; the model below is a proposed
continuous-time reduction.
\begin{figure}[p]
\centering
\includegraphics[width=\textwidth,height=\casefigheight,keepaspectratio]{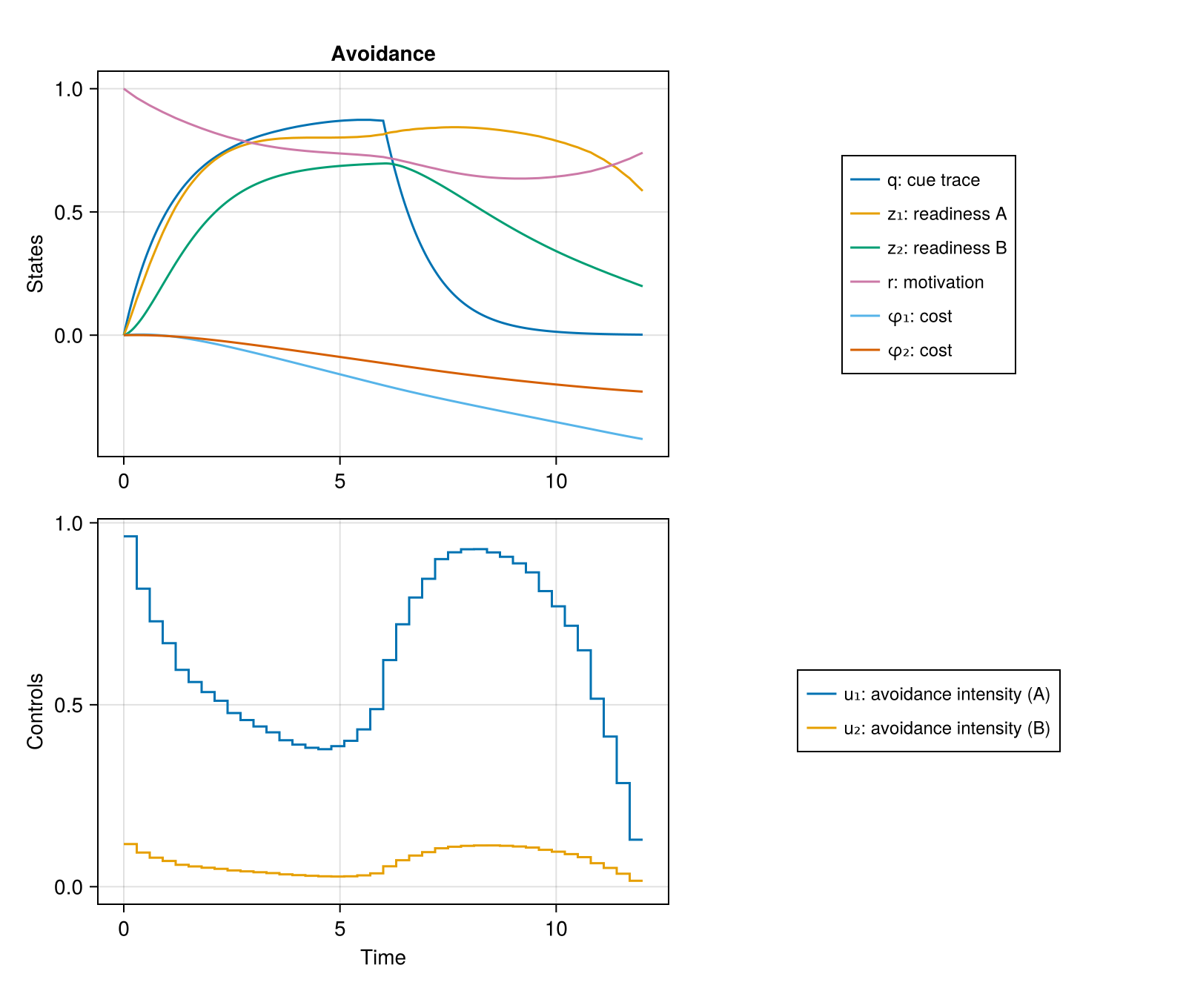}
\caption[Reference equilibrium candidate for \texttt{\caseidentifier}]{Reference equilibrium candidate for \texttt{\caseidentifier}. \textit{Interpretation.} While the warning cue is on ($t<6$), the cue trace $q$ rises quickly and saturates near
$0.87$, driving the avoidance readiness of both modules upward. Module A opens with almost
maximal avoidance intensity, lowers it to about $0.4$ once readiness has built up, and then
raises it again to above $0.9$ right after the cue is withdrawn, while its readiness $z_A$
stays high during the unsignaled second half; it finally releases effort toward the end of the
episode. Module B, which weights effort costs more strongly, mirrors this timing at a much lower
level (below $0.12$), so its readiness $z_B$ decays after cue offset and relies largely on the
cue itself. The shared motivation $r$ declines from its resting value during sustained
avoidance and partially recovers when A reduces effort at the end. No state or budget
restrictions are imposed; the contrast between the two modules follows solely from their
different effort weights and from the delayed benefit of acting after the cue has vanished.}\label{fig:\caseidentifier}
\end{figure}
\myheader{Literature and provenance.} The avoidance and decentralized-controller studies \citep{Ilango2012DopamineAvoidance,Schilling2019FluidAdaptivity,Schilling2020Hexapod} motivate the cue, readiness, and shared-resource mechanisms. The four-state ODE reduction, independent module objectives, control bounds, and all numerical coefficients are introduced for this benchmark.
\myheader{Dimensions.} $N=2$, $n_x=4$, $n_{\rm ODE}=6$, $(n_{u,i})_{i\in\mathcal P}=(2\times 1)$, $(n_{g,i})_{i\in\mathcal P}=(2\times 0)$.
\myheader{Categories.} \cat{organismal}, \cat{episodical}, \cat{regulation}.
\myheader{Game specification.} Let $x=(q,z_A,z_B,r)$ denote the cue trace $q$, the readiness $z_A,z_B$ of the two
modules, and the shared motivational state $r$. Two decision modules choose avoidance intensities
$u_A,u_B$:
\[
\begin{aligned}
\dot q&=a_q(s-q)-b_q(u_A+u_B)q,\\
\dot z_i&=a_z q(1-z_i)-d_z z_i+\gamma u_i r,\qquad i\in\{A,B\},\\
\dot r&=a_r(r_0-r)-c_r(u_A+u_B)-h_r(z_A+z_B).
\end{aligned}
\]
Here $s(t)$ is the warning cue, an external signal with $s(t)=1$ for $0\le t\le 6$ and $s(t)=0$ afterwards; the cue trace $q$ follows it at rate $a_q$ and is depleted by the avoidance intensities of both modules. The shared state $r$
represents a resource or motivational context; each module changes the context
in which the other module acts. 

These equations are a tractable benchmark inspired by the cited mechanisms,
not equations extracted from those publications.

Each module minimizes a weighted cost over $[0,T]$:
\[
J_i= w_i^1 J_{i1}+w_i^2 J_{i2}+w_i^3 J_{i3}+w_i^4 J_{i4},
\]
\[
\begin{aligned}
J_{i1}&=-\frac1T\int_0^T z_i(t)\,dt, &
J_{i2}&=\frac1T\int_0^T u_i(t)^2\,dt,\\
J_{i3}&=\frac1T\int_0^T \frac{K_r}{K_r+r(t)}\,dt, &
J_{i4}&=-\frac1T\int_0^T z_j(t)\,dt,\quad j\ne i.
\end{aligned}
\]
The first term rewards readiness to avoid, the second penalizes control effort,
the third penalizes an adverse motivational state, and the fourth represents
an optional partner benefit. Weights satisfy $w_i^k\geq0$ and
$\sum_k w_i^k=1$.

No state or budget bounds are enabled.
\myheader{Parameter instance.} %
\begin{lstlisting}
T = 12.0
aq = 0.8
bq = 0.25
az = 0.9
dz = 0.35
gamma = 0.45
ar = 0.6
r0 = 1.0
cr = 0.12
hr = 0.08
Kr = 0.35
x0 = (0.0, 0.0, 0.0, 1.0)
initial = (0.2, 0.3)
weights_A = (0.65, 0.15, 0.1, 0.1)
weights_B = (0.25, 0.45, 0.1, 0.2)
control_lower = [0.0, 0.0]
control_upper = [1.0, 1.0]
numerical_x0 = [0.0, 0.0, 0.0, 1.0, 0.0, 0.0]
\end{lstlisting}

\endgroup

\begingroup
\def\caseidentifier{batchFermentation}%
\subsection{batchFermentation: Division of labor in batch fermentation}\label{case:batchFermentation}
\myheader{Biological interpretation.} Industrial and natural fermentations often rely on a division of labor: an upstream strain
converts the primary substrate into an intermediate that downstream specialists turn into the
final product. In a closed batch, the substrate is finite, so the episode ends when a product
target has been reached. Each strain pays for its production activity with a growth burden,
and time is costly (for example, reactor occupancy or exposure to contaminants). The upstream
strain chooses the batch duration, so only its time weight is strategically relevant
(Section~\ref{sec:owned-duration}).
The game instance asks how the strains schedule their production when the timing of the intermediate
supply couples their incentives: three production controls trade batch time against production
effort and state-dependent costs.
\begin{figure}[p]
\centering
\includegraphics[width=\textwidth,height=\casefigheight,keepaspectratio]{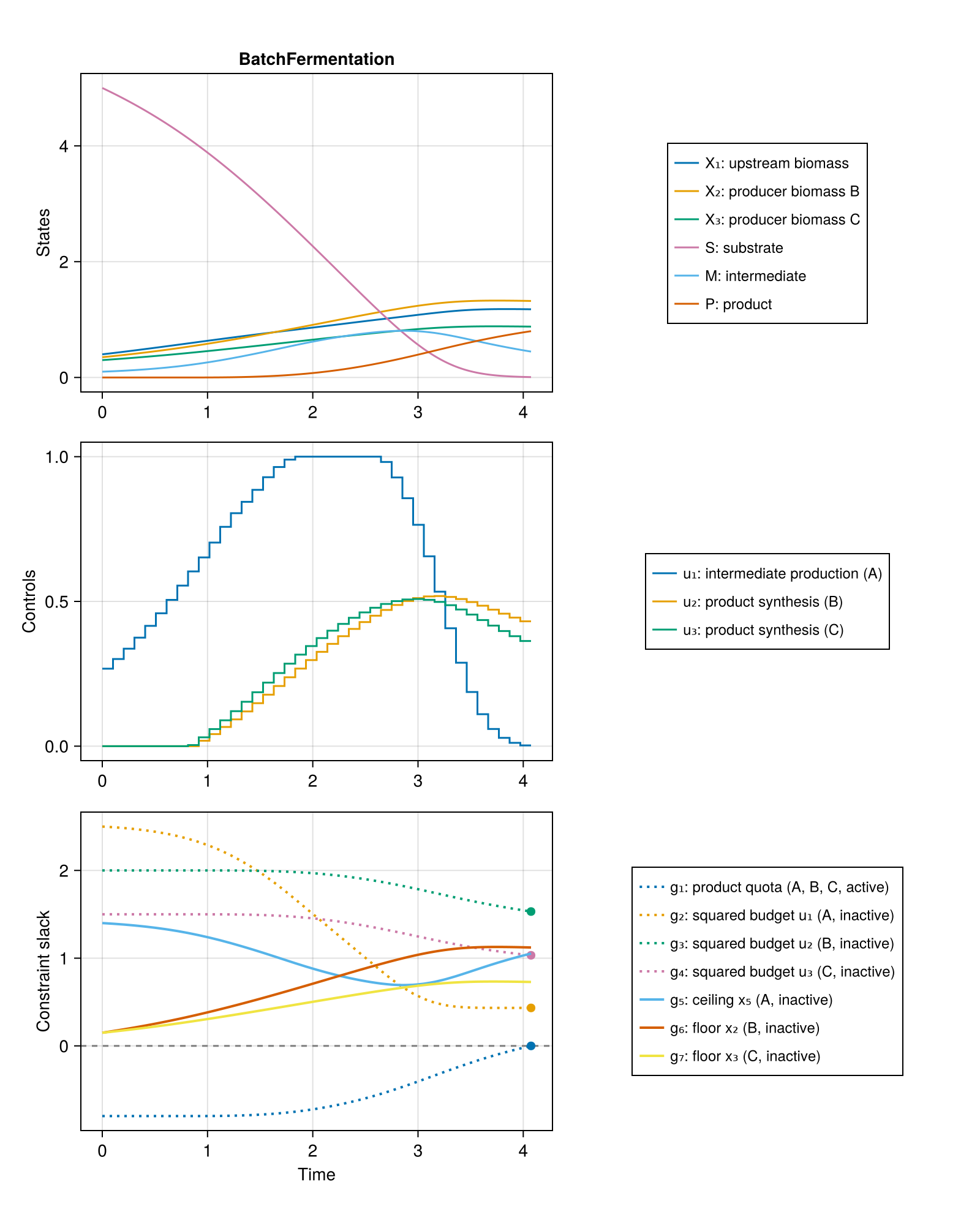}
\caption[Reference equilibrium candidate for \texttt{\caseidentifier}]{Reference equilibrium candidate for \texttt{\caseidentifier}. \textit{Interpretation.} The batch ends at $T\approx4.1$, exactly when the product reaches its quota (the terminal
equality $g_1$ is active for all three strains). The upstream strain A first ramps up
intermediate production from about $0.3$ to its maximum, holds it during the middle of the
batch, and switches it off as the substrate is exhausted near $t\approx3.8$. The producers B and
C wait until intermediate is available (their product synthesis is zero for $t<0.9$) and then
increase synthesis in parallel, peaking around $t\approx3$ and easing off slightly at the end.
Accordingly, the intermediate pool first accumulates, peaks near $t\approx2.9$, and is then drawn
down, while the product rises sigmoidally. All biomasses grow while substrate lasts and level
off afterwards. The effort budgets, the intermediate ceiling, and the producer floors remain
slack, so the sequential relay is shaped by the shared quota, the time cost, and the substrate
supply rather than by these limits.}\label{fig:\caseidentifier}
\end{figure}
\myheader{Literature and provenance.} A: upstream intermediate producer; B and C: product-synthesizing strains \citep{BenchmarkPMC5877992}. A also schedules batch completion.

$x=(X_A,X_B,X_C,S,M,P)$: strain biomasses, substrate, shared intermediate, and product.
Motivates division of labor and production burden. The three-strain ODEs and numerical defaults here are a reduced synthetic game, not a reproduction of the source model.
\myheader{Dimensions.} $N=3$, $n_x=6$, $n_{\rm ODE}=12$, $(n_{u,i})_{i\in\mathcal P}=(3\times 1)$, $(n_{g,i})_{i\in\mathcal P}=(3\times 3)$.
\myheader{Categories.} \cat{cellular}, \cat{episodical}, \cat{mutualism}, \cat{freeTime}.
\myheader{Game specification.} \[\begin{gathered}
h=S/(K+S),\quad r_A=p_Au_AX_Ah,\\
r_i=p_i u_iX_iM/(K+M)\ (i\in\{B,C\}),\\
\dot X_i=[g_i(1-0.5u_i)h-d]X_i,\\
\dot S=-\sum_{i\in\mathcal P} g_i(1-0.5u_i)hX_i/Y-r_A,\\
\dot M=r_A-r_B-r_C-d_MM,\quad \dot P=r_B+r_C.
\end{gathered}\]

$g=$ \texttt{growth}, $K=$ \texttt{saturation}, $d=$ \texttt{maintenance}, $p=$ \texttt{production}, $Y=$ \texttt{yield}, $d_M=$ \texttt{intermediate\_decay}. The allocation burden coefficient is the fixed dimensionless number 0.5.

$J_i=w_i^1T+w_i^2\int_0^T u_i^2\,dt+w_i^3\int_0^T L_i\,dt$.

$L_i=M^2-0.1\iota(i)X_i$. The duration term penalizes batch length, effort penalizes synthesis allocation, and the third candidate trades intermediate burden against own biomass. It may be negative.

$P(0)=0$, $0.8-P(T)\geq0$ and $P(T)-0.8\geq0$, shared by all players. Private effort budgets are $(2.5,2,1.5)$. A imposes $M\leq1.5$ (accumulated intermediate inhibits the upstream strain), B imposes $X_B\geq0.2$ and C imposes $X_C\geq0.15$ (each producer retains a viable population for subsequent use of the culture) along the trajectory. Other terminal states are free.

Enabled bounds:
\begin{center}\small\begin{tabular}{lll}\toprule Owners & Bound & Location\\\midrule
A, B, C & $y_{6}(T)=0.8$ & terminal \\
A & $y_{7}(T)\leq 2.5$ & terminal \\
B & $y_{8}(T)\leq 2$ & terminal \\
C & $y_{9}(T)\leq 1.5$ & terminal \\
A & $y_{5}(t)\leq 1.5$ & path \\
B & $0.2\leq y_{2}(t)$ & path \\
C & $0.15\leq y_{3}(t)$ & path \\
\bottomrule\end{tabular}\end{center}
\myheader{Parameter instance.} %
\begin{lstlisting}
x0 = (0.4, 0.35, 0.3, 5.0, 0.1, 0.0)
target = 0.8
growth = (0.7, 0.6, 0.5)
saturation = 0.5
maintenance = 0.03
production = (0.8, 0.7, 0.6)
yield = 0.7
intermediate_decay = 0.05
duration = 5.0
duration_bounds = (0.1, 20.0)
initial = (0.5, 0.45, 0.4)
weights = [0.75 0.15 0.1; 0.55 0.3 0.15; 0.45 0.25 0.3]
budgets = (2.5, 2.0, 1.5)
path_limits = (1.5, 0.2, 0.15)
control_lower = [0.0, 0.0, 0.0]
control_upper = [1.0, 1.0, 1.0]
numerical_x0 = [0.4, 0.35, 0.3, 5.0, 0.1, 0.0, 0.0, 0.0, 0.0, 0.0, 0.0, 0.0]
\end{lstlisting}

\endgroup

\begingroup
\def\caseidentifier{biofilm}%
\subsection{biofilm: Matrix investment and dispersal in a biofilm}\label{case:biofilm}
\myheader{Biological interpretation.} Bacterial biofilms develop in episodes: cells attach, build an extracellular polymeric matrix
(EPS) that protects them from loss, grow while substrate and oxygen diffuse in, and finally
disperse to colonize new surfaces. EPS is costly to the producer and, in this model, protects only
the producer's own cells; neighbors are affected indirectly through crowding and shared
resources \citep{BenchmarkPMC1783407}. Dispersal converts attached biomass into
offspring leaving the film. Strains sharing a surface therefore face a timing problem: invest
early in protection and growth, then switch to dispersal before resources are exhausted.
We model one colonization episode with three strains that differ in their preference for
attached biomass versus dispersal output.

Players are $N\geq2$ strains. $B_i$: attached biomass;
$E_i$: strain-associated extracellular polymeric substances (EPS);
$S,O$: shared substrate and oxygen.
$x=((B_i)_{i\in\mathcal P},(E_i)_{i\in\mathcal P},S,O)$.
Each strain controls polymer allocation $e_i$ and dispersal intensity
$d_i$, ordered $(e_i,d_i)_{i\in\mathcal P}$ in player order.
Dispersed biomass leaves this system; its departure is a reproductive-output proxy.
All controls depend on time alone.
\begin{figure}[p]
\centering
\includegraphics[width=\textwidth,height=\casefigheight,keepaspectratio]{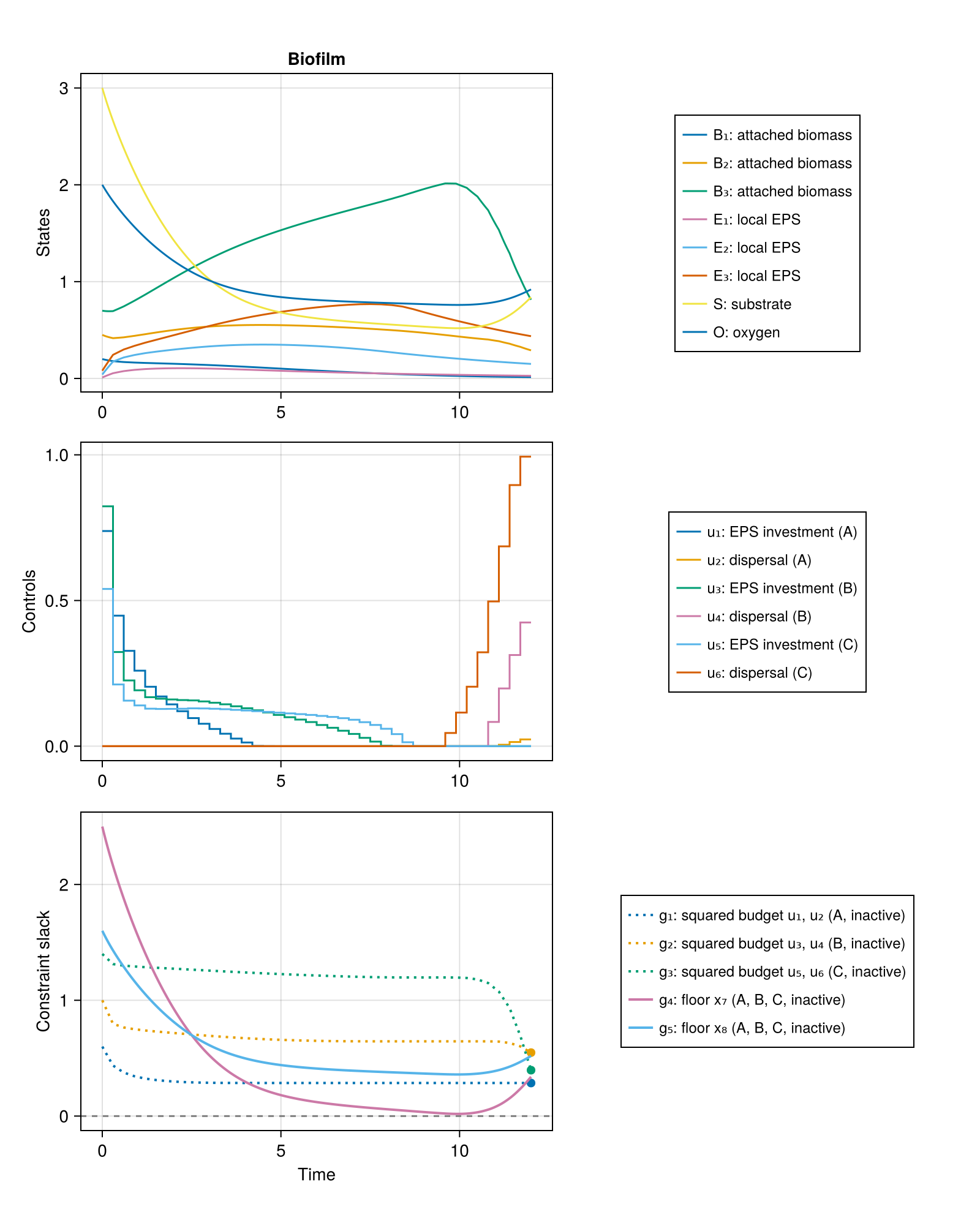}
\caption[Reference equilibrium candidate for \texttt{\caseidentifier}]{Reference equilibrium candidate for \texttt{\caseidentifier}. \textit{Interpretation.} All three strains invest in matrix production at the very beginning, consistent with protection acting
longest when built early, and then taper this investment: A stops after $t\approx4$, while B and C keep a low
level until $t\approx8$--$9$. Substrate and oxygen are drawn down from their inflow levels to
minima near $t\approx10$. Strain C, which values dispersal output most, builds the largest
attached population (from $0.7$ to about $2$ at $t\approx9.6$) and the largest EPS pool; the
small initial population of A declines throughout, and B stays nearly constant. In the final
phase C switches to strong dispersal, reaching maximal intensity at the end, B follows with a
smaller dispersal burst from $t\approx11$, and A disperses only marginally. Dispersal removes
much of C's attached biomass, so substrate and oxygen recover late in the episode while the EPS
pools decay. The effort budgets and the substrate and oxygen floors stay slack, although the
substrate floor is approached closely around $t\approx10$, exactly when dispersal begins.}\label{fig:\caseidentifier}
\end{figure}
\myheader{Literature and provenance.} The source studies competition between strains differing in extracellular polymer production using spatial individual-based simulations \citep{BenchmarkPMC1783407}. We extract resource allocation and matrix-mediated protection as mechanisms for an aggregate ODE game, adding controlled dispersal and nutrient restrictions. The local protection law is our proposed reduction.

All numerical defaults are synthetic and dimensionless; no fitted biological values or source equilibrium results are asserted.
\myheader{Dimensions.} $N=3$, $n_x=8$, $n_{\rm ODE}=14$, $(n_{u,i})_{i\in\mathcal P}=(3\times 2)$, $(n_{g,i})_{i\in\mathcal P}=(3\times 3)$.
\myheader{Categories.} \cat{cellular}, \cat{episodical}, \cat{competition}, \cat{variablePlayers}.
\myheader{Game specification.} \[
V_i=v_i^{\max}\frac{S}{h+S}\frac{O}{h+O}B_i,
\]
\[
\dot B_i=(1-e_i)Y_BV_i-
\left[\frac{m}{1+\sigma E_i}+d_i+\kappa\sum_{j\in\mathcal P}B_j\right]B_i,
\qquad \dot E_i=Y_E e_iV_i-\delta_EE_i,
\]
\[
\dot S=D(S_{\rm in}-S)-\sum_{i\in\mathcal P}V_i,\qquad
\dot O=D(O_{\rm in}-O)-c_O\sum_{i\in\mathcal P}V_i.
\]
$v_i^{\max},h,Y_B,Y_E,m,\sigma,\kappa,\delta_E$ map to
\texttt{vmax, half, biomass\_yield, eps\_yield, mortality, protection,
crowding, eps\_loss}; $D=\texttt{dilution}$ and $c_O=\texttt{oxygen\_cost}$.
The same uptake flux is split between biomass and EPS production.
No nutrient is consumed when $S=0$ or $O=0$.

Each player minimizes
\[
J_i=\frac1T\int_0^T\sum_{k=1}^3w_i^k\ell_{ik}\,dt,
\qquad w_i^k\geq0,\quad\sum_kw_i^k=1.
\]
There is no terminal objective. Each weight row is independently configurable.

\[
(\ell_{i1},\ell_{i2},\ell_{i3})=(-B_i,\ -d_iB_i,\ e_i^2+d_i^2).
\]
Candidates reward attached biomass and exported propagules, and penalize
strategy effort. Weight changes allow persistence versus dispersal preferences.
Export does not imply successful colonization elsewhere.

$\dot Z_i=e_i^2+d_i^2$, $Z_i(0)=0$,
$Z_i(T)\leq\mathrm{budgets}_i$.
Shared path floors $S\geq\mathrm{floors}_1$ and $O\geq\mathrm{floors}_2$
are enforced in every player's deviation problem.
The enabled restrictions are listed below.

For $q$ biological states, append $N$ cost integrals then $N$ budget integrals,
all initially zero; the cost integrals include the factor $1/T$, the ODE dimension is $q+2N$, and $J_i=x_{q+\iota(i)}(T)$.

Enabled bounds:
\begin{center}\small\begin{tabular}{lll}\toprule Owners & Bound & Location\\\midrule
A & $y_{12}(T)\leq 0.6$ & terminal \\
B & $y_{13}(T)\leq 1$ & terminal \\
C & $y_{14}(T)\leq 1.4$ & terminal \\
A, B, C & $0.5\leq y_{7}(t)$ & path \\
A, B, C & $0.4\leq y_{8}(t)$ & path \\
\bottomrule\end{tabular}\end{center}
\myheader{Parameter instance.} %
\begin{lstlisting}
N = 3
T = 12.0
vmax = [1.0, 1.2, 1.4]
half = 0.4
biomass_yield = 0.6
eps_yield = 1.5
mortality = 0.5
protection = 10.0
crowding = 0.08
eps_loss = 0.15
dilution = 0.5
substrate_in = 3.0
oxygen_in = 2.0
oxygen_cost = 0.5
x0 = [0.2, 0.45, 0.7, 0.01, 0.04, 0.08, 3.0, 2.0]
weights = [0.5 0.4 0.1; 0.35 0.5 0.15; 0.2 0.6 0.2]
budgets = [0.6, 1.0, 1.4]
floors = [0.5, 0.4]
initial = [0.2, 0.1, 0.2, 0.1, 0.2, 0.1]
control_lower = [0.0, 0.0, 0.0, 0.0, 0.0, 0.0]
control_upper = [1.0, 1.0, 1.0, 1.0, 1.0, 1.0]
numerical_x0 = [0.2, 0.45, 0.7, 0.01, 0.04, 0.08, 3.0, 2.0, 0.0, 0.0, 0.0, 0.0, 0.0, 0.0]
\end{lstlisting}

\endgroup

\begingroup
\def\caseidentifier{ciliate}%
\subsection{ciliate: Bacterium--ciliate interaction}\label{case:ciliate}
\myheader{Biological interpretation.} Bacteria and bacterivorous ciliates form a classical predator--prey system in aquatic microbial
food webs. Bacteria can express grazing defenses, for example filament formation, microcolonies, or
surface modifications, at a cost in growth, whereas ciliates adjust their feeding activity at a cost
in energy. Bacterial productivity in lakes and ponds follows the seasonal cycle of temperature and
nutrient supply. We model one full annual cycle that is assumed to repeat directly, after which both populations and both expressed
phenotypes must return to their initial values; the phenotypes follow the commanded defense and
feeding levels with a delay.

The physical state is $x=(x_1,x_2,x_3,x_4)$: bacterial density, ciliate density, expressed defense, and expressed feeding activity. Player A commands defense $u_A$ and player B commands feeding $u_B$.
\begin{figure}[p]
\centering
\includegraphics[width=\textwidth,height=\casefigheight,keepaspectratio]{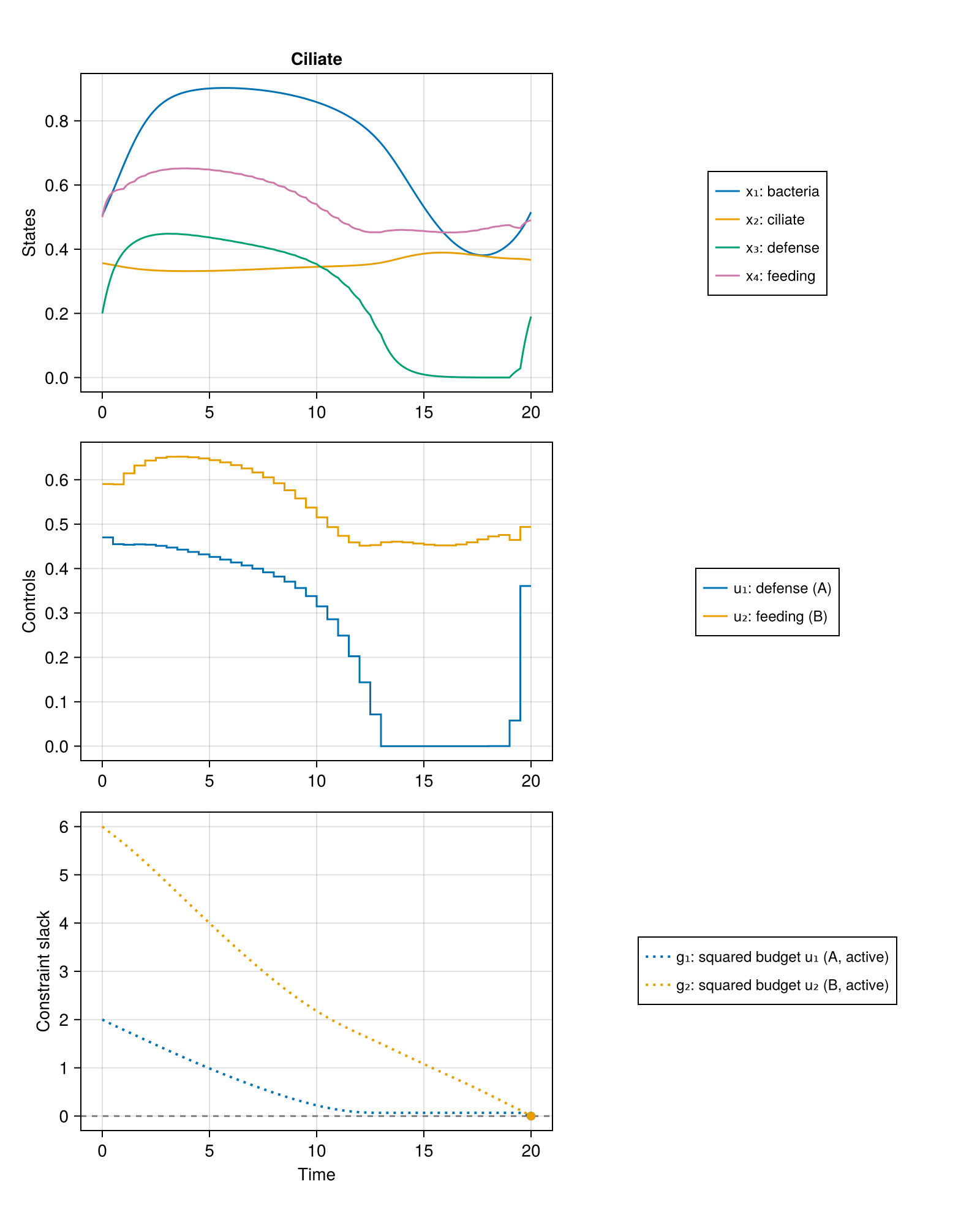}
\caption[Reference equilibrium candidate for \texttt{\caseidentifier}]{Reference equilibrium candidate for \texttt{\caseidentifier}. \textit{Interpretation.} During the productive half of the cycle, bacterial density rises from $0.5$ to about
$0.9$ and stays high until $t\approx10$; it then collapses to a minimum near $0.38$ at
$t\approx18$ when productivity is low, and recovers at the end of the cycle. The ciliate population
responds weakly and with a delay, varying only between $0.33$ and $0.39$ and peaking during the
bacterial decline. The bacteria invest in defense early, when their own density and the benefit of
protection are high, and let the defense command decline gradually; once productivity drops,
defense is switched off completely ($t\approx13$--$19$), consistent with the defense cost in the growth equation. The return condition then forces a short defense pulse at the end of the cycle that rebuilds
the expressed defense to its initial level. The ciliate raises feeding effort to about
$0.65$ while prey are abundant and lowers it to about $0.45$ during the bacterial trough. The
expressed phenotypes follow their commands with visible delays, especially the slower bacterial
defense. Both players exhaust their squared-effort budgets ($g_1,g_2$ active); the bacterial budget
is almost used up by $t\approx12$, and its remainder is spent on the return-induced final pulse.}\label{fig:\caseidentifier}
\end{figure}
\myheader{Literature and provenance.} \emph{Pseudomonas fluorescens} and \emph{Tetrahymena thermophila}
motivate an interaction between two cellular organisms.
\cite{Hiltunen2014Consumer} demonstrate the importance
of consumer co-evolution in this setting.

Rapid evolution of prey defense can reshape consumer--resource cycles
\citep{Yoshida2007Cryptic}. The following saturating feeding model is a proposed
coarse-grained game, rather than a fitted description of that experiment.

A representative bacterial lineage invests in defense that
reduces capture at a cost to growth.
A representative ciliate lineage invests in feeding effort
that increases food uptake at a metabolic cost.

Resource limitation, saturation, and effort costs create different objectives
for the two players.

The numerical values are synthetic hypotheses.
\myheader{Dimensions.} $N=2$, $n_x=4$, $n_{\rm ODE}=8$, $(n_{u,i})_{i\in\mathcal P}=(2\times 1)$, $(n_{g,i})_{i\in\mathcal P}=(2\times 9)$.
\myheader{Categories.} \cat{cellular}, \cat{periodical}, \cat{antagonism}.
\myheader{Game specification.} With $r(t)=r[1+a\sin(2\pi t/T_{\rm per})]$, $T_{\rm per}=T$, define the saturating feeding rate
\[
 F(x_1,x_3,x_4)=\frac{a_0x_4e^{-\alpha x_3}x_1}{1+ha_0x_4e^{-\alpha x_3}x_1}.
\]
The four-state dynamics are
\[
\begin{aligned}
 \dot x_1&=r(t)x_1(1-x_1/K)-c_Ax_3^2x_1-Fx_2,\\
 \dot x_2&=(eF-m_B-c_Bx_4^2)x_2,\\
 \dot x_3&=(u_A-x_3)/\tau_A,\qquad
 \dot x_4=(u_B-x_4)/\tau_B.
\end{aligned}
\]
The players minimize the cycle-average costs
\[
\begin{aligned}
 J_A&=\frac1T\int_0^T\left[-w_A^1\frac{x_1}{x_{1,0}}-w_A^2r(t)x_1+w_A^3u_A^2\right]\dt,\\
 J_B&=\frac1T\int_0^T\left[-w_B^1\frac{x_2}{x_{2,0}}-w_B^2eFx_2+w_B^3u_B^2\right]\dt.
\end{aligned}
\]
The first features reward mean relative abundance; the second reward gross production proxies. In the listed instance all feature scales equal one; otherwise each weight is divided by its positive \code{scales_A} or \code{scales_B} entry. Command effort remains distinct from expressed phenotype costs in the dynamics.

For A, $\int_0^T u_A^2\dt\leq U_A$; for B, $\int_0^T u_B^2\dt\leq U_B$. Optional density floors $x_1\geq\varepsilon_A$ and $x_2\geq\varepsilon_B$ are assigned to A and B respectively and are disabled here. Both players impose $|x_j(T)-x_j(0)|\leq\delta$, $j=1,\ldots,4$, with $\delta=10^{-2}$ in physical coordinates. Initial densities come from a constant-command reference cycle, and initial phenotypes equal these commands.

Numerical coordinates are $y=(\log x_1,\log x_2,x_3,x_4,z_A,z_B,b_A,b_B)$; $z_i$ accumulates the unaveraged running cost and $\dot b_i=u_i^2$. Quadratures start at zero. For the density components, return bounds become $\log(x_j(0)-\delta)\leq y_j(T)\leq\log(x_j(0)+\delta)$, rather than symmetric bands in log coordinates. The optional \code{phenotype=false} variant substitutes commands directly for expressed traits and removes their relaxation equations.

Enabled bounds:
\begin{center}\small\begin{tabular}{lll}\toprule Owners & Bound & Location\\\midrule
A & $y_{7}(T)\leq 2$ & terminal \\
B & $y_{8}(T)\leq 6$ & terminal \\
A, B & $|x_j(T)-x_j(0)|\leq 0.01$, $j=1,\ldots,4$ & terminal (return) \\
\bottomrule\end{tabular}\end{center}
\myheader{Parameter instance.} %
\begin{lstlisting}
forcing_period = 20.0
T = 20.0
forcing_amplitude = 0.8
periodicity_tolerance = 0.01
r = 1.0
K = 1.0
c_A = 0.15
a0 = 2.0
alpha = 2.0
h = 0.5
e = 0.5
m_B = 0.15
c_B = 0.1
x0 = (0.5055436644, 0.3564873025)
w_A = (0.9, 0.08, 0.02)
w_B = (0.4, 0.5, 0.1)
scales_A = (1.0, 1.0, 1.0)
scales_B = (1.0, 1.0, 1.0)
U_A = 2.0
U_B = 6.0
epsilon_A = nothing
epsilon_B = nothing
initial = (0.2, 0.5)
phenotype = true
tau_A = 0.75
tau_B = 0.25
phenotype0 = (0.2, 0.5)
control_lower = [0.0, 0.0]
control_upper = [1.0, 1.0]
numerical_x0 = [-0.6821208656, -1.031456657, 0.2, 0.5, 0.0, 0.0, 0.0, 0.0]
\end{lstlisting}

\endgroup

\begingroup
\def\caseidentifier{circadianCycle}%
\subsection{circadianCycle: PER/TIM circadian oscillator}\label{case:circadianCycle}
\myheader{Biological interpretation.} Circadian clocks are intrinsic biochemical oscillators whose period is close to, but not fixed
at, 24 hours. In \emph{Drosophila}, PER and TIM proteins are synthesized, phosphorylated,
combine into a complex, enter the nucleus, and repress their own transcription, closing a
delayed negative feedback loop. Light acts on this loop by accelerating TIM degradation, which
is how the clock is reset by the environment. We model one full cycle of the loop as a periodical
game between two regulatory influences: player A modulates \emph{per} transcription and sets the
cycle length, while player B, the light-input pathway (e.g., CRY-mediated TIM degradation)
modeled as a regulatory module, controls TIM degradation. Both prefer balanced PER/TIM levels
and low intervention effort. The duration cost represents a preference for a faster cycle
within the admissible range; it is a synthetic device to create a nontrivial duration decision,
not a biological fitness argument. Only A's duration weight is strategically relevant
(Section~\ref{sec:owned-duration}).

We use the ten-state PER/TIM model of Leloup and Goldbeter \citep{LeloupGoldbeter1998PerTim}
in the form given by \citet{Slaby2007-tr}, Appendix B, with its kinetic defaults (nM, hours).
\[
x=(M_P,P_0,P_1,P_2,M_T,T_0,T_1,T_2,C,C_N).
\]
$M_P,M_T$: per/tim mRNA; $P_j,T_j$: PER/TIM with $j$ phosphorylation
steps; $C,C_N$: cytosolic/nuclear complex. Nuclear complex represses transcription.

A controls \emph{per} transcription and duration. B controls TIM degradation,
the light-dependent intervention mechanism studied in the source.
The two-player game, centered control mapping, objectives, and restrictions are
new benchmark choices; they are not experimental parameter estimates.
\begin{figure}[p]
\centering
\includegraphics[width=\textwidth,height=\casefigheight,keepaspectratio]{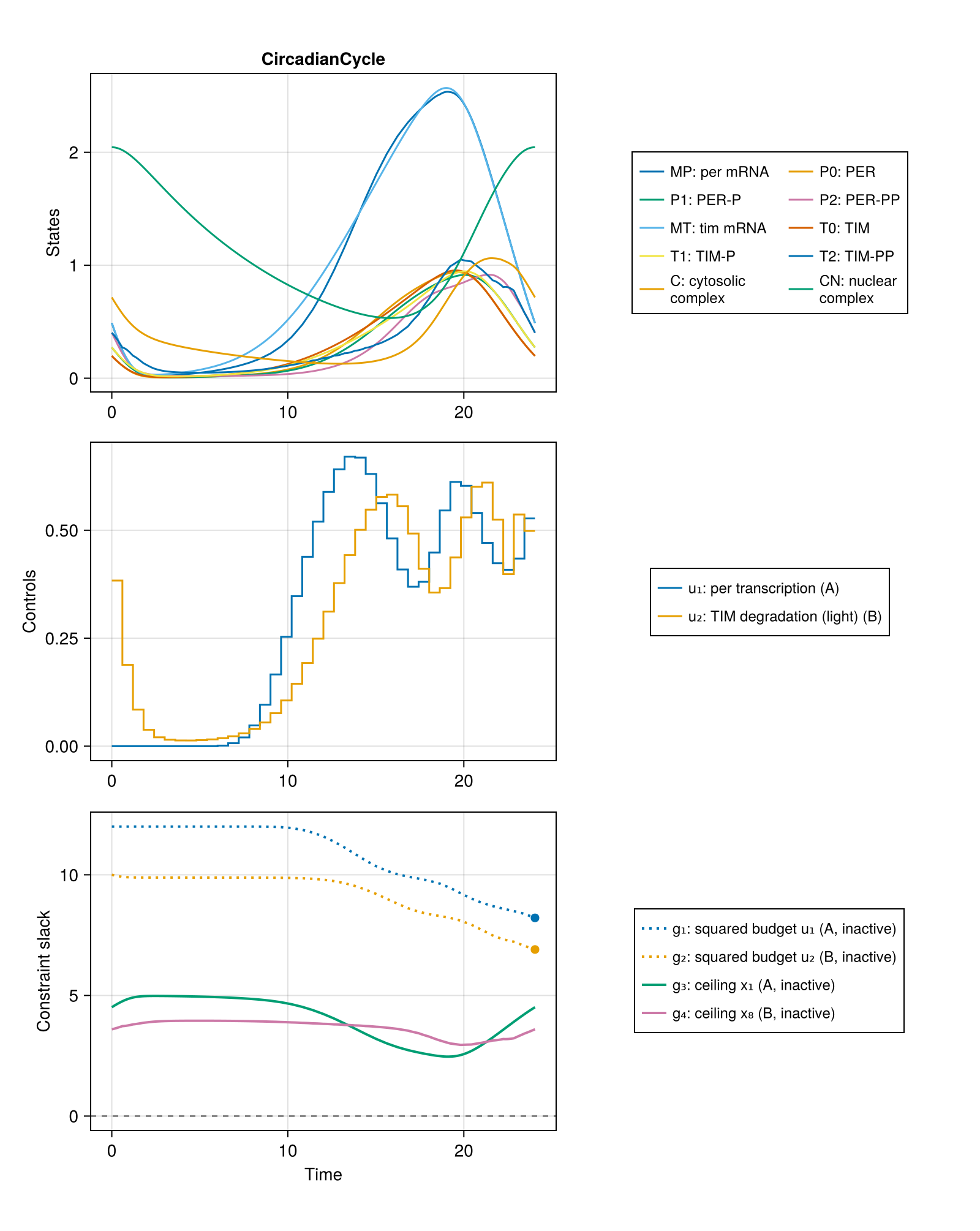}
\caption[Reference equilibrium candidate for \texttt{\caseidentifier}]{Reference equilibrium candidate for \texttt{\caseidentifier}. \textit{Interpretation.} Player A chooses a cycle length of about $24.04$ hours, slightly shorter than the uncontrolled
reference period of $24.13$ hours. During the first half of the cycle, when the nuclear complex
$C_N$ is high and represses transcription, both mRNAs and all PER and TIM forms fall to low
values; A keeps transcriptional intervention off and B quickly lowers TIM degradation after an
initial pulse. As repression wanes (with $C_N$ reaching its minimum near $t\approx16$), A raises
\emph{per} transcription to a maximum near $t\approx13$, followed with a delay of a few hours by
B's increase in TIM degradation, which keeps TIM from overshooting PER. The mRNAs peak around
$t\approx19$, the phosphorylated proteins and the cytosolic complex around $t\approx20$--$22$,
and the nuclear complex rebuilds to its initial level at the end of the cycle, closing the loop
within the return tolerance. The squared-effort budgets and the concentration ceilings remain
slack; the interventions are shaped by the phase structure of the feedback loop.}\label{fig:\caseidentifier}
\end{figure}
\myheader{Literature and provenance.} The ten concentration equations of \citet{LeloupGoldbeter1998PerTim} and the kinetic coefficients come from Appendix B of \citet{Slaby2007-tr}, with the notation correction stated below. Centered transcription and degradation controls, their separate objectives, player A's duration ownership, effort budgets, concentration ceilings, and approximate return conditions define the added game. The initial phase is constructed numerically from the source oscillator.
\myheader{Dimensions.} $N=2$, $n_x=10$, $n_{\rm ODE}=14$, $(n_{u,i})_{i\in\mathcal P}=(2\times 1)$, $(n_{g,i})_{i\in\mathcal P}=(2\times 22)$.
\myheader{Categories.} \cat{molecular}, \cat{periodical}, \cat{regulation}, \cat{freeTime}.
\myheader{Game specification.} For $X\in\{P,T\}$, define
\[\begin{gathered}
r_{1X}=V_{1X}\frac{X_0}{K_{1X}+X_0},\quad
r_{2X}=V_{2X}\frac{X_1}{K_{2X}+X_1},\\
r_{3X}=V_{3X}\frac{X_1}{K_{3X}+X_1},\quad
r_{4X}=V_{4X}\frac{X_2}{K_{4X}+X_2},\quad \Gamma=k_3P_2T_2-k_4C.
\end{gathered}\]
\[\begin{aligned}
\dot M_X&=\widetilde v_{sX}\frac{K_{IX}^{n}}{K_{IX}^{n}+C_N^{n}}
-v_{mX}\frac{M_X}{K_{mX}+M_X}-k_dM_X,\\
\dot X_0&=k_{sX}M_X-r_{1X}+r_{2X}-k_dX_0,\\
\dot X_1&=r_{1X}-r_{2X}-r_{3X}+r_{4X}-k_dX_1.
\end{aligned}\]

\[\begin{aligned}
\dot X_2&=r_{3X}-r_{4X}-\Gamma-\widetilde v_{dX}\frac{X_2}{K_{dX}+X_2}-k_dX_2,\\
\dot C&=\Gamma-k_1C+k_2C_N-k_{dC}C,\\
\dot C_N&=k_1C-k_2C_N-k_{dN}C_N.
\end{aligned}\]
\[\begin{gathered}
\widetilde v_{sP}=v_{sP}[1+0.6(u_A-0.5)],\quad \widetilde v_{sT}=v_{sT},\\
\widetilde v_{dP}=v_{dP},\quad \widetilde v_{dT}=v_{dT}+2(u_B-0.5).
\end{gathered}\]
At $(0.5,0.5)$, the source's Appendix B dynamics are recovered.
$\Gamma$ is the net complex-formation flux. On normalized time, multiply every physical derivative by the constant duration $T$.
The printed $K_1C$ in the source complex equation is interpreted as $k_1C$,
consistent with its nuclear equation, Figure 3, and Table B1.

\[
J_i=w_i^1T+w_i^2\int_0^T u_i^2\,dt+w_i^3\int_0^T L_i\,dt,
\]
\[
L_A=(M_P-M_T)^2,\qquad
L_B=(P_0+P_1+P_2-T_0-T_1-T_2)^2.
\]
These candidates penalize elapsed time, effort, and transcription/protein imbalance.

All ten concentrations satisfy $|x_j(T)-x_j(0)|\leq10^{-4}$ nM;
$12\leq T\leq36$ hours excludes
zero-time closure. Initial state is a fixed nonstationary reference-cycle point.
A has effort budget 12 and $M_P\leq5$; B has budget 10 and $T_2\leq4$.
This is an approximate cycle; zero return tolerance requests exact closure.
Cost quadratures are not periodic. These bounds do not prove exclusion of every
possible controlled multi-winding orbit.

Integrate Appendix B with controls $(0.5,0.5)$ from
$(1,.3,.5,.5,1,.3,.5,.5,.5,.5)$ to 1500 hours (Tsit5,
absolute/relative tolerances $10^{-12}/10^{-11}$).
Find late maxima of $C_N$ from $k_1C-(k_2+k_{dN})C_N=0$.

The stored reference period is $24.1346$ hours; successive reference
states differed by less than $6\cdot10^{-13}$ in that construction.
This is an uncontrolled periodic initializer, not an optimized game result.

Paired kinetic tuples below are ordered PER, TIM. The four $V$ and $K$ entries
are shared by both phosphorylation chains, as in Table B1.
Initial concentrations, unlike duration and controls, are fixed during optimization.

Enabled bounds:
\begin{center}\small\begin{tabular}{lll}\toprule Owners & Bound & Location\\\midrule
A & $y_{11}(T)\leq 12$ & terminal \\
B & $y_{12}(T)\leq 10$ & terminal \\
A & $y_{1}(t)\leq 5$ & path \\
B & $y_{8}(t)\leq 4$ & path \\
A, B & $|x_j(T)-x_j(0)|\leq 0.0001$, $j=1,\ldots,10$ & terminal (return) \\
\bottomrule\end{tabular}\end{center}
\myheader{Parameter instance.} %
\begin{lstlisting}
x0 = (0.4884645232, 0.1972717656, 0.2725715265, 0.4025780857, 0.4884645232, 0.1972717656, 0.2725715265, 0.4025780857, 0.7156565511, 2.044733003)
vs = (1.0, 1.0)
vm = (0.7, 0.7)
Km = (0.2, 0.2)
ks = (0.9, 0.9)
vd = (2.0, 2.0)
KI = (1.0, 1.0)
Kd = (0.2, 0.2)
V = (8.0, 1.0, 8.0, 1.0)
K = (2.0, 2.0, 2.0, 2.0)
k1 = 0.6
k2 = 0.2
k3 = 1.2
k4 = 0.6
kd = 0.01
kdC = 0.01
kdN = 0.01
hill = 4
transcription_gain = 0.6
light_gain = 2.0
duration = 24.13458674
duration_bounds = (12.0, 36.0)
periodicity_tolerance = 0.0001
initial = (0.5, 0.5)
weights = [0.65 0.2 0.15; 0.45 0.35 0.2]
budgets = (12.0, 10.0)
path_limits = (5.0, 4.0)
control_lower = [0.0, 0.0]
control_upper = [1.0, 1.0]
numerical_x0 = [0.4884645232, 0.1972717656, 0.2725715265, 0.4025780857, 0.4884645232, 0.1972717656, 0.2725715265, 0.4025780857, 0.7156565511, 2.044733003, 0.0, 0.0, 0.0, 0.0]
\end{lstlisting}

\endgroup

\begingroup
\def\caseidentifier{fishery}%
\subsection{fishery: Competing fleets harvest a renewable stock}\label{case:fishery}
\myheader{Biological interpretation.} Many fisheries are managed by rotation: a fishing ground is opened, fished down, and then
closed until the stock has recovered, after which the cycle repeats (pulse or rotational
harvesting). The length of such a rotation is a management decision with a clear biological
trade-off. Short rotations allow a per-opening catch allowance to be taken more often, but the
logistic stock needs time to rebuild between openings, and a deeply depleted stock is costly to
fish because catch per unit effort falls with biomass. We therefore model one recurring
rotation, which is periodic in the stock but not in the catch accounts.

\cite{Dockner1989Fishery}: a Gordon--Schaefer duopoly,
with Nash and Stackelberg solutions; an established fishery differential game.
\cite{Suri2008Harvesting}, p.~168: explicit finite-horizon
two-harvester dynamics and separate profit objectives with
linear--quadratic effort costs; we retain its stock and cost structure.

Players are fishing fleets, and the biological state is a renewable fish stock.
In rotational management the opening schedule is set by the licence holder or the management
authority; here it is assigned to fleet A, which also chooses the rotation length. Fleet B
fishes the same ground within each opening. The rotation structure, per-opening quotas,
finite effort limits, and numerical parameter values below define a synthetic benchmark,
not a calibrated fishery.
\begin{figure}[p]
\centering
\includegraphics[width=\textwidth,height=\casefigheight,keepaspectratio]{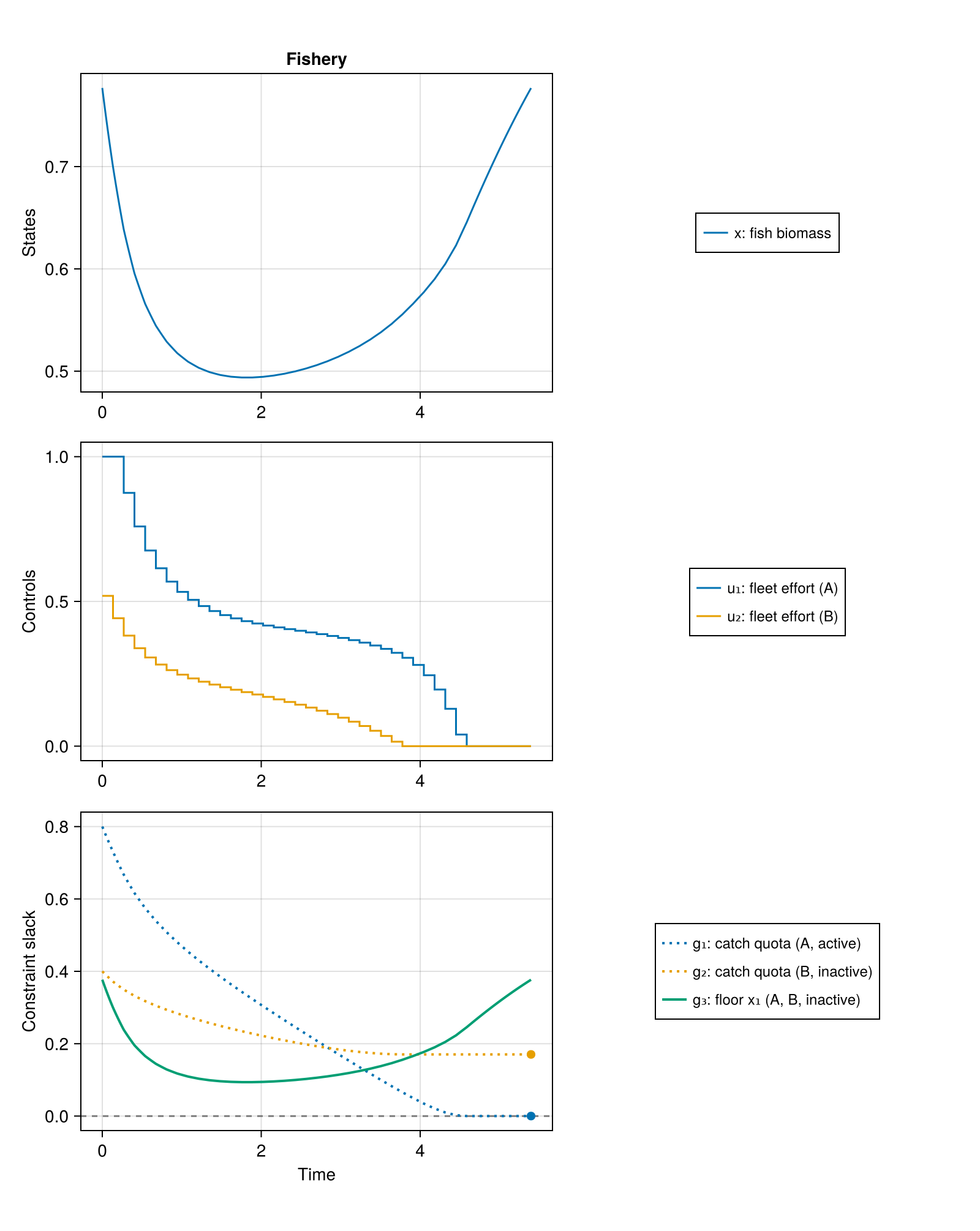}
\caption[Reference equilibrium candidate for \texttt{\caseidentifier}]{Reference equilibrium candidate for \texttt{\caseidentifier}. \textit{Interpretation.} Fleet A chooses a rotation of length $T\approx5.4$. Both fleets fish hardest when the ground
opens: A starts at full effort and B at about half effort, and both reduce effort as the stock
is fished down, consistent with catch per unit effort falling with biomass. Biomass therefore drops
steeply from its stationary reference level of about $0.78$ to a minimum near $0.49$ roughly
one third into the rotation. B stops fishing first, and A exhausts its per-rotation catch
quota ($g_1$ active) and closes the ground well before the end of the rotation. The final closed
phase lets the logistic stock recover to its initial level, which the return condition requires,
while B's quota and the shared stock floor remain slack. The chosen duration is consistent with a balance of two
effects: a shorter rotation would let A take its quota more often, but the stock could not
recover within the rotation without a much deeper and more costly depletion.}\label{fig:\caseidentifier}
\end{figure}
\myheader{Literature and provenance.} The stock balance and linear--quadratic fleet costs follow the finite-horizon harvesting formulation in \citet{Suri2008Harvesting}; \citet{Dockner1989Fishery} provides differential-game context. The duration owner, approximate stock return, catch quotas, finite effort limits, and normalized numerical instance are benchmark additions.
\myheader{Dimensions.} $N=2$, $n_x=1$, $n_{\rm ODE}=5$, $(n_{u,i})_{i\in\mathcal P}=(2\times 1)$, $(n_{g,i})_{i\in\mathcal P}=(2\times 4)$.
\myheader{Categories.} \cat{population}, \cat{periodical}, \cat{competition}, \cat{freeTime}.
\myheader{Game specification.} $x$: fish biomass; $u_i$: fleet $i$'s fishing effort, $i\in\{A,B\}$.
Fleet A additionally chooses the rotation duration $T\in[T_{\min},T_{\max}]$.
The dynamics are written on normalized time $s\in[0,1]$, with physical time $t=T s$;
the core performs no hidden time scaling. Catch rate is $h_i=q_iu_ix$, where $q_i$ is catchability, and
\[
\frac{dx}{ds}=T\left[r x(1-x/K)-q_Au_Ax-q_Bu_Bx\right],\qquad x(0)=x_0.
\]
$r$ is intrinsic growth; $K$ is carrying capacity. There is no external clock:
the default $x_0=K\bigl(1-(q_A\bar u_A+q_B\bar u_B)/r\bigr)$ is the stock that is stationary
under the reference efforts $\bar u_i$, so the reference policy returns for every duration.
Biomass and effort units are normalized; time, price, and cost coefficients use a
consistent arbitrary unit system.

Each fleet minimizes its average net cost rate over one rotation,
\[
J_i=\int_0^1\Bigl[-w_i^1p_iq_iu_ix+w_i^2\Bigl(a_iu_i+\frac{b_i}{2}u_i^2\Bigr)\Bigr]ds
=\frac1T\int_0^T\Bigl[\ldots\Bigr]dt,
\]
with $w_i^k\geq0$ and $w_i^1+w_i^2=1$; $p_i$ is fish price and $a_i,b_i$ are linear and
quadratic effort costs. Per-rotation catch accounts
\[
\frac{dH_i}{ds}=T q_iu_ix,\quad H_i(0)=0,\quad H_i(1)\leq Q_i,
\qquad x(s)\geq x_{\min},
\]
make short rotations attractive (more quota per unit time) but require a depletion--recovery
cycle within each rotation. Each quota is enforced for its fleet; the stock floor and the
return condition $|x(1)-x_0|\leq\delta$ are enforced in both players' deviation problems.

Julia uses the augmented order $(x,z_A,z_B,H_A,H_B)$, where $z_i$ accumulates the
integrand of $J_i$, $z_i(0)=0$, and $J_i=z_i(1)$; these four accumulators are
quadratures, not biological states. The duration is the scalar parameter owned by player A.

Enabled bounds:
\begin{center}\small\begin{tabular}{lll}\toprule Owners & Bound & Location\\\midrule
A & $y_{4}(T)\leq 0.8$ & terminal \\
B & $y_{5}(T)\leq 0.4$ & terminal \\
A, B & $0.4\leq y_{1}(t)$ & path \\
A, B & $|x_j(T)-x_j(0)|\leq 0.0001$, $j=1,\ldots,1$ & terminal (return) \\
\bottomrule\end{tabular}\end{center}
\myheader{Parameter instance.} %
\begin{lstlisting}
periodicity_tolerance = 0.0001
r = 0.8
K = 1.0
q = (0.7, 0.55)
price = (1.2, 1.0)
linear_cost = (0.06, 0.08)
quadratic_cost = (0.5, 0.65)
w_A = (0.8, 0.2)
w_B = (0.6, 0.4)
umax = (1.0, 1.0)
initial = (0.2, 0.07)
x0 = 0.776875
duration = 6.0
duration_bounds = (1.0, 12.0)
quotas = (0.8, 0.4)
floor = 0.4
control_lower = [0.0, 0.0]
control_upper = [1.0, 1.0]
numerical_x0 = [0.776875, 0.0, 0.0, 0.0, 0.0]
\end{lstlisting}

\endgroup

\begingroup
\def\caseidentifier{foraging}%
\subsection{foraging: Competitive foraging}\label{case:foraging}
\myheader{Biological interpretation.} Animals foraging on a shared, renewable food patch, for example a group of birds on a fruiting tree
or grazers on a pasture, convert food into energy reserves and part of it into a retained harvest
(such as caches or brood provisions). Intense foraging depletes the patch for everyone. A
management rule keeps the patch above $70\%$ of carrying capacity; with logistic regrowth this is
a conservative rule, not a collapse threshold. We model one foraging episode of three foragers
that differ in how they value their retained harvest relative to their foraging effort. Each player controls foraging effort subject to an effort budget.
\begin{figure}[p]
\centering
\includegraphics[width=\textwidth,height=\casefigheight,keepaspectratio]{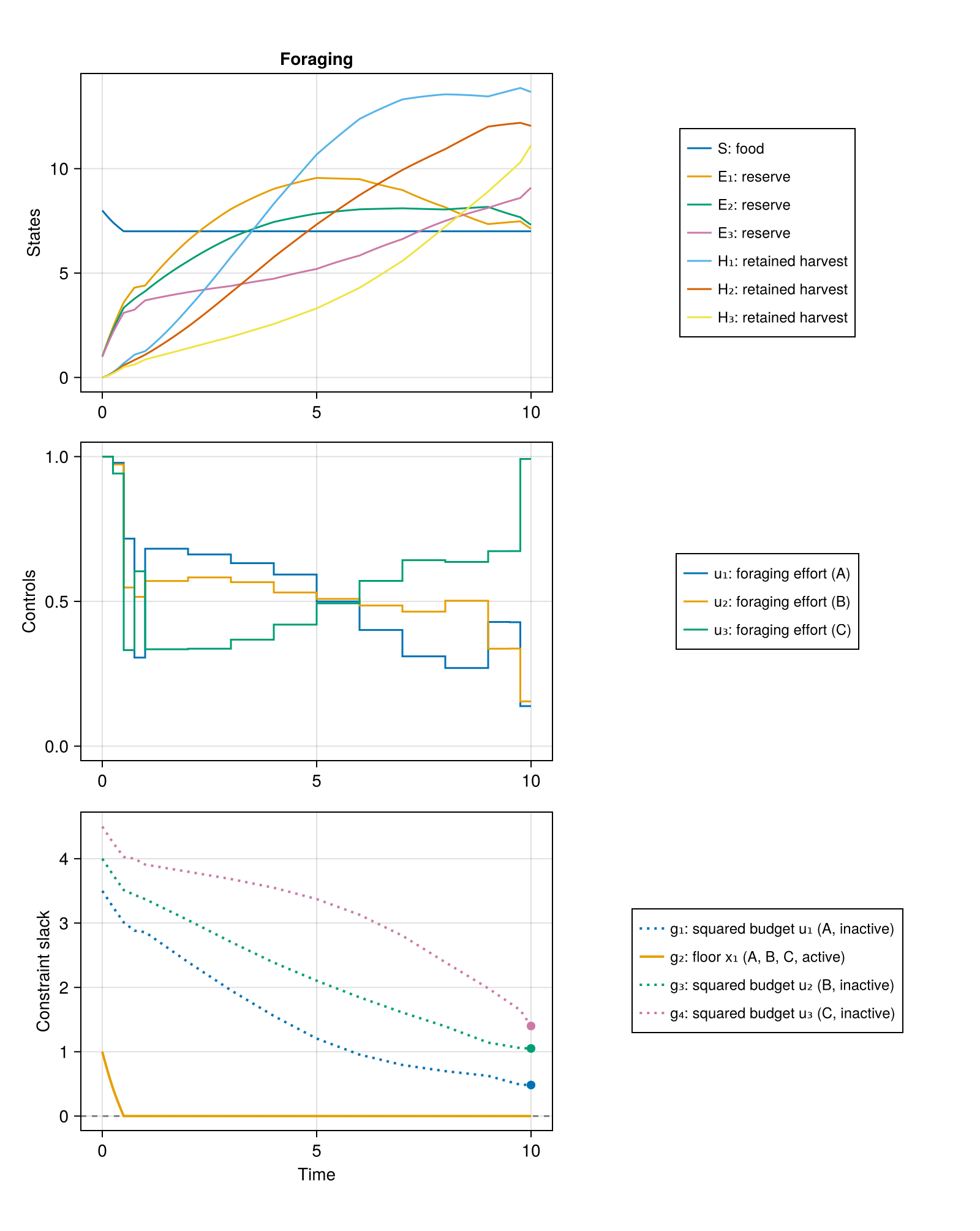}
\caption[Reference equilibrium candidate for \texttt{\caseidentifier}]{Reference equilibrium candidate for \texttt{\caseidentifier}. \textit{Interpretation.} All foragers start at maximal effort and draw the shared food stock down from $8$ to its
conservation floor of $7$ within about half a time unit; from then on the floor is active for the
rest of the episode ($g_2$), so the foragers share the patch's regrowth. Their schedules then diverge.
Forager A, which weights effort most, forages at about $0.68$ early and reduces effort to about $0.27$
near $t\approx8$ before a small final increase; B keeps an intermediate effort around $0.5$--$0.58$;
and C, which values retained harvest most, starts low ($0.33$) and increases its effort to about
$0.67$ toward the end. Because the floor fixes the total harvest rate, one forager's higher effort
lowers the others' intake: C's late increase coincides with A's decline. Energy reserves
reflect these schedules: A's reserve peaks at about $9.5$ near $t\approx5.5$ and then declines, B's
levels off near $8$, and C's keeps rising to about $9.1$. Retained harvest accumulates earliest for
A, which has the highest foraging and retention efficiencies, and reaches the largest total
($13.7$), whereas C's harvest accelerates late. None of the effort budgets is exhausted (used
$3.0$ of $3.5$, $2.9$ of $4$, and $3.1$ of $4.5$); the shared floor, not individual effort limits,
restrains the harvest.}\label{fig:\caseidentifier}
\end{figure}
\myheader{Literature and provenance.} Competitive foraging of mallards on two food patches was compared with the ideal free distribution \citep{Benchmarkforaging}. This supports competition for a shared renewable food supply, not the specific depletion and reserve dynamics. We extend the resource competition idea to three foragers with finite effort budgets; the particular state equations and budgets are synthetic.
All rates, budgets, weights, and initial values below are synthetic,
dimensionless benchmark choices, not fitted values from this source.
\myheader{Dimensions.} $N=3$, $n_x=7$, $n_{\rm ODE}=13$, $(n_{u,i})_{i\in\mathcal P}=(3\times 1)$, $(n_{g,i})_{i\in\mathcal P}=(3\times 2)$.
\myheader{Categories.} \cat{population}, \cat{episodical}, \cat{competition}, \cat{variablePlayers}.
\myheader{Game specification.} $S$: food; $E_i$: reserve; $H_i$: retained harvest (subject to loss). $H_i$ is a proxy fed by the reserves without depleting them. $u_i$: foraging effort. Players are individual foragers; $N=3$ by default.
\[\begin{gathered}
\dot S=rS(1-S/K)-S\sum_{j\in\mathcal P}c_ju_j,\\
\dot E_i=e_iu_iS-d_EE_i,\qquad \dot H_i=h_iu_iE_i-d_HH_i.
\end{gathered}\]
Controls are open-loop functions of time.
Initial-state order: $x=(S,(E_i)_{i\in\mathcal P},(H_i)_{i\in\mathcal P})$; values appear below.

Each player $i\in\mathcal P$ minimizes $J_i=T^{-1}\int_0^T\ell_i\,dt$, with
\[\ell_i=-w_i^1H_i+w_i^2u_i^2.\]
$w_i^k\geq0$, $\sum_kw_i^k=1$; weights are independently configurable.
\[
\dot Z_i=u_i^2,\quad Z_i(0)=0,\quad Z_i(T)\leq \mathcal B_i.
\]
Budgets are private to each player. The enabled path floor is $S\geq7$. No forager values the patch level itself: depletion results from competition for the shared regrowth, and the management floor restrains it.
When enabled, it is enforced in every player's deviation problem.

Enabled bounds:
\begin{center}\small\begin{tabular}{lll}\toprule Owners & Bound & Location\\\midrule
A & $y_{11}(T)\leq 3.5$ & terminal \\
A, B, C & $7\leq y_{1}(t)$ & path \\
B & $y_{12}(T)\leq 4$ & terminal \\
C & $y_{13}(T)\leq 4.5$ & terminal \\
\bottomrule\end{tabular}\end{center}
\myheader{Parameter instance.} %
\begin{lstlisting}
N = 3
T = 10.0
r = 0.8
K = 10.0
c = [0.12, 0.16, 0.2]
e = [0.8, 0.7272727273, 0.6666666667]
h = [0.6, 0.5454545455, 0.5]
de = 0.3
dh = 0.1
x0 = [8.0, 1.0, 1.0, 1.0, 0.0, 0.0, 0.0]
w_A = (0.3333333333, 0.6666666667)
w_B = (0.375, 0.625)
w_C = (0.8, 0.2)
budgets = [3.5, 4.0, 4.5]
floor = 7.0
initial = [0.2, 0.2, 0.2]
control_lower = [0.0, 0.0, 0.0]
control_upper = [1.0, 1.0, 1.0]
numerical_x0 = [8.0, 1.0, 1.0, 1.0, 0.0, 0.0, 0.0, 0.0, 0.0, 0.0, 0.0, 0.0, 0.0]
\end{lstlisting}

\endgroup

\begingroup
\def\caseidentifier{immuneClearance}%
\subsection{immuneClearance: Immune clearance and viral immune evasion}\label{case:immuneClearance}
\myheader{Biological interpretation.} An acute viral infection is a contest between the host immune response and the virus. Cytotoxic
T cells kill infected cells and antibodies neutralize free virus; both responses are costly, and
cytotoxic killing damages the infected tissue itself (immunopathology). Many viruses counter these
responses, for example by reducing antigen presentation on infected cells or by masking epitopes,
at a cost in replication. We model one infection episode as a game between the host (A), which
activates cytotoxic and antibody responses and ends the episode when the viral load has been reduced
to a prescribed low level, and the virus (B), which invests in immune evasion to maximize its
accumulated viral load, a proxy for transmission opportunity. Treating the viral population as a
player summarizes selection on viral gene expression within one infection; it is not a claim of
viral foresight. The host owns the episode duration, since clearance ends it.
\begin{figure}[p]
\centering
\includegraphics[width=\textwidth,height=\casefigheight,keepaspectratio]{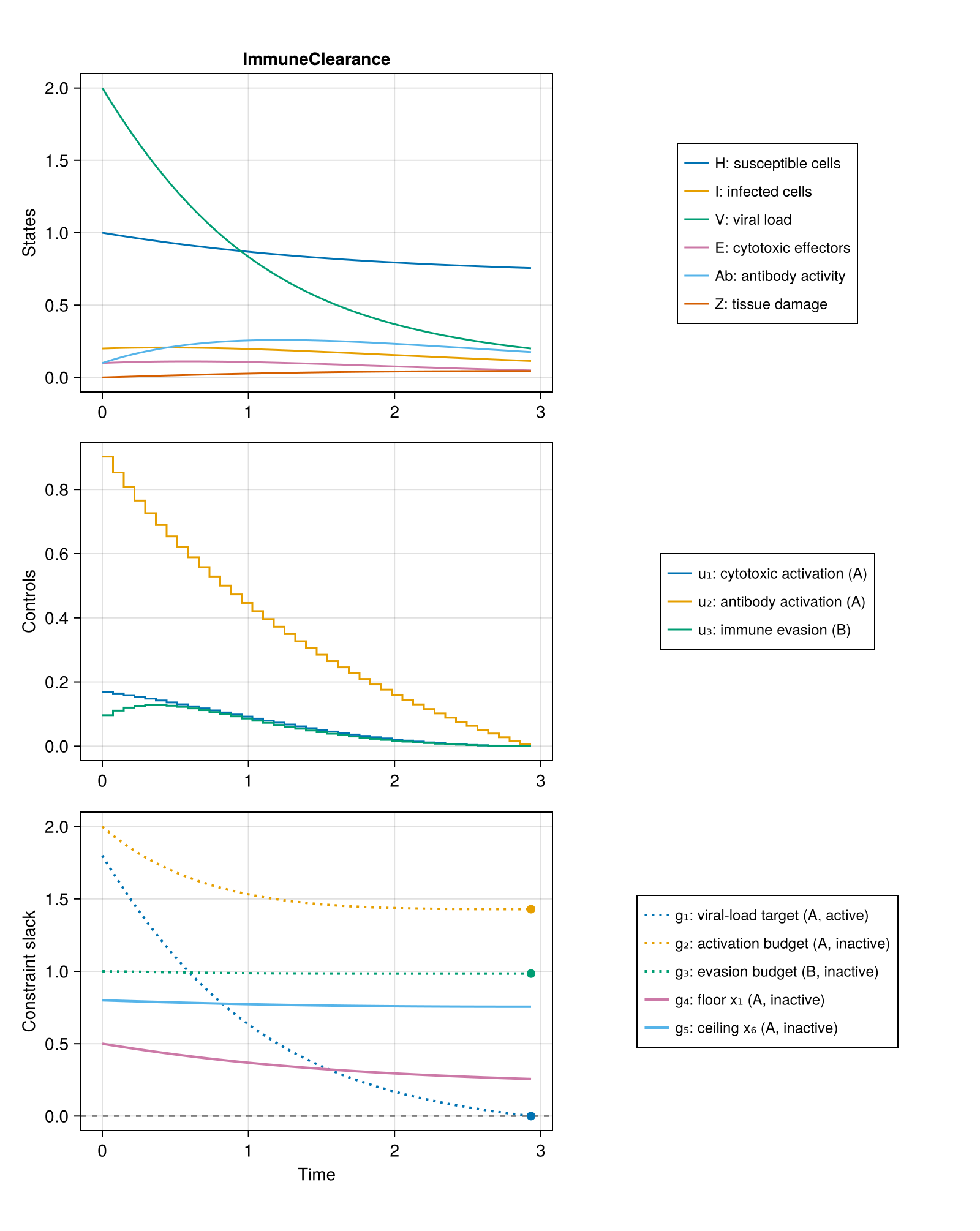}
\caption[Reference equilibrium candidate for \texttt{\caseidentifier}]{Reference equilibrium candidate for \texttt{\caseidentifier}. \textit{Interpretation.} The host ends the episode at $T\approx2.93$, when the viral load has fallen from $2$ to the
prescribed target $0.2$ (the host's terminal equality $g_1$ is active). The host relies mainly on
antibodies: antibody activation starts at $0.90$ and declines almost linearly to zero, whereas
cytotoxic activation stays below $0.17$ and declines in parallel. A mechanism consistent with this
split is that cytotoxic killing also causes tissue damage, which the host penalizes. The virus
invests in evasion early, rising from $0.10$ to a maximum of $0.13$ near $t\approx0.35$, and then
reduces it to zero as the viral load falls, well before antibody activity peaks; the virus uses
less than $2\%$ of its evasion budget. Antibody
activity rises from $0.1$ to a maximum of about $0.26$ near $t\approx1.2$ and then decays slowly,
cytotoxic effectors stay near $0.05$--$0.11$, infected cells decline from $0.2$ to about $0.11$, and
about a quarter of the susceptible cells become infected. Tissue damage accumulates to about $0.05$.
The host's activation budget, the evasion budget, the susceptible-cell floor, and the damage
ceiling remain slack, so the schedules are shaped by the host's time and burden costs and by the
virus's trade-off between evasion and production.}\label{fig:\caseidentifier}
\end{figure}
\myheader{Literature and provenance.} The target-cell infection and viral production/clearance kinetics follow the model of influenza A infection kinetics in humans \citep{BenchmarkimmuneClearance}. Saturating infection, the immune compartments, the damage law, the evasion mechanism, the host--virus game, and all coefficients are synthetic; this is not the fitted influenza model or a treatment recommendation.

$x=(H,I,V,E,Ab,Z)$: susceptible cells, infected cells, virus, cytotoxic effectors, antibody activity, and tissue damage.
\myheader{Dimensions.} $N=2$, $n_x=6$, $n_{\rm ODE}=10$, $(n_{u,i})_{i\in\mathcal P}=(2,1)$, $(n_{g,i})_{i\in\mathcal P}=(4,1)$.
\myheader{Categories.} \cat{cellular}, \cat{episodical}, \cat{antagonism}, \cat{freeTime}.
\myheader{Game specification.} Host A controls cytotoxic activation $u_{A1}$ and antibody activation $u_{A2}$; virus B controls
immune evasion $u_B$. With susceptible cells $H$, infected cells $I$, viral load $V$, cytotoxic
effectors $E$, antibody activity $Ab$, and tissue damage $Z$,
\[\begin{gathered}
b=\beta HV/(1+V),\quad \dot H=-b,\\
\dot I=b-[d_I+k_E(1-\sigma u_B)E]I,\\
\dot V=p(1-\kappa u_B)I-[c+(1-\sigma u_B)Ab]V,\\
\dot E=a_1u_{A1}V/(1+V)-d_1E,\qquad
\dot{Ab}=a_2u_{A2}V/(1+V)-d_2Ab,\\
\dot Z=h_1E^2+h_2I-r_ZZ.
\end{gathered}\]
Evasion reduces cytotoxic killing and antibody neutralization by the fraction $\sigma u_B$ and
costs the fraction $\kappa u_B$ of virion production.

$\beta=$ \texttt{beta}, $d_I=$ \texttt{infected\_decay}, $k_E=$ \texttt{killing}, $p=$ \texttt{production}, $c=$ \texttt{clearance}, $a=$ \texttt{activation}, $d=$ \texttt{effector\_decay}, $h=$ \texttt{damage}, $r_Z=$ \texttt{recovery}, $\sigma=$ \texttt{evasion}, $\kappa=$ \texttt{evasion\_cost}.

The players minimize
\[
\begin{aligned}
J_A&=w_A^1T+w_A^2\int_0^T(u_{A1}^2+u_{A2}^2)\dt+w_A^3\int_0^T(V^2+Z^2)\dt,\\
J_B&=-w_B^1\int_0^TV\dt+w_B^2\int_0^Tu_B^2\dt.
\end{aligned}
\]
The host trades a short infection against activation effort and against viral and pathological
burden; the virus trades accumulated viral load against the replication cost of evasion.

$V(0)=2$ and the clearance target $V(T)=0.2$ are restrictions of the host only; the virus is not
bound by them, so its deviations can make the host's target harder to reach. The host's squared
activation budget is $2$, and it imposes $H\geq0.5$ (a minimum of intact target tissue) and
$Z\leq0.8$ (a tolerable level of immunopathology) along the trajectory. The virus has a squared
evasion budget of $1$. Other terminal states are free.

Enabled bounds:
\begin{center}\small\begin{tabular}{lll}\toprule Owners & Bound & Location\\\midrule
A & $y_{3}(T)=0.2$ & terminal \\
A & $y_{7}(T)\leq 2$ & terminal \\
B & $y_{8}(T)\leq 1$ & terminal \\
A & $0.5\leq y_{1}(t)$ & path \\
A & $y_{6}(t)\leq 0.8$ & path \\
\bottomrule\end{tabular}\end{center}
\myheader{Parameter instance.} %
\begin{lstlisting}
x0 = (1.0, 0.2, 2.0, 0.1, 0.1, 0.0)
target = 0.2
beta = 0.25
infected_decay = 0.6
production = 0.7
clearance = 0.8
killing = 0.6
activation = (0.8, 0.6)
effector_decay = (0.5, 0.35)
damage = (0.2, 0.15)
recovery = 0.4
evasion = 0.8
evasion_cost = 0.3
duration = 5.0
duration_bounds = (0.1, 20.0)
initial = (0.5, 0.45, 0.3)
w_A = (0.5, 0.3, 0.2)
w_B = (0.6, 0.4)
budgets = (2.0, 1.0)
path_limits = (0.5, 0.8)
control_lower = [0.0, 0.0, 0.0]
control_upper = [1.0, 1.0, 1.0]
numerical_x0 = [1.0, 0.2, 2.0, 0.1, 0.1, 0.0, 0.0, 0.0, 0.0, 0.0]
\end{lstlisting}

\endgroup

\begingroup
\def\caseidentifier{metabolicPathway}%
\subsection{metabolicPathway: Enzyme allocation under budgets}\label{case:metabolicPathway}
\myheader{Biological interpretation.} Many metabolic pathways branch: a shared substrate is converted by parallel enzymes whose
products feed a common downstream pool. When the enzymes are expressed by different cells or
regulatory modules, each module pays for its own enzyme synthesis while all benefit from the
shared product, a classic public-good situation. Here the substrate is supplied with a recurring
(for example, feeding-related) oscillation, and we consider one full supply cycle after which all
pools return to their initial values. Each reaction also yields a module-specific byproduct that
accumulates and inhibits the module's own reaction (product inhibition), so enzyme investment must
be timed relative to the substrate supply and to the module's byproduct load. Because both modules value the same product, the game is a public-goods
game in which free-riding on the other module's expression is expected.
\begin{figure}[p]
\centering
\includegraphics[width=\textwidth,height=\casefigheight,keepaspectratio]{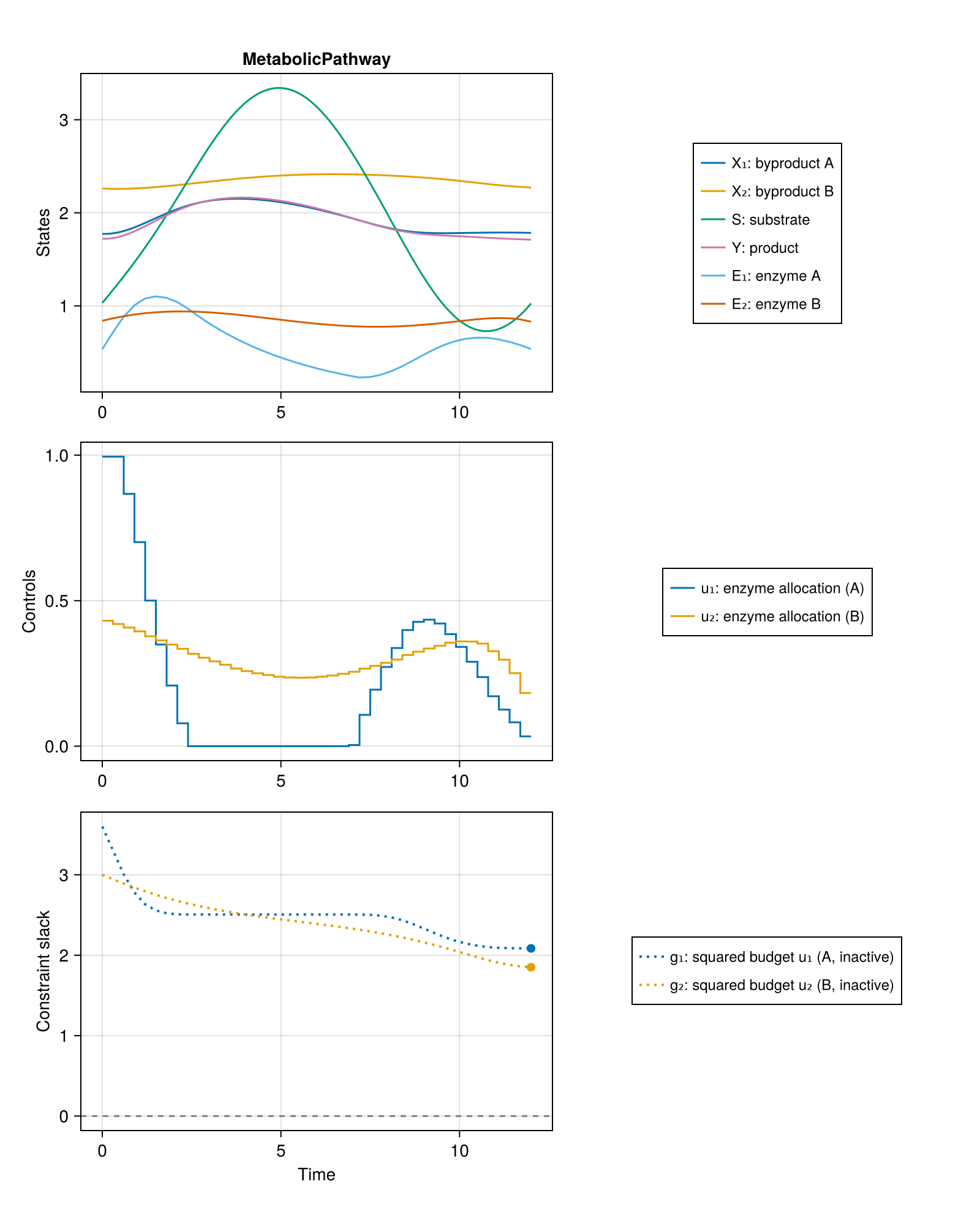}
\caption[Reference equilibrium candidate for \texttt{\caseidentifier}]{Reference equilibrium candidate for \texttt{\caseidentifier}. \textit{Interpretation.} Over one supply period ($T=12$), substrate rises from $1.0$ to a maximum of $3.3$ near
$t\approx5$ and falls to a minimum of $0.73$ near $t\approx10.8$. Module A, which values the shared
product most, expresses at the upper bound at the start, while its enzyme pool is low and substrate
is rising, builds its enzyme to $1.1$ by $t\approx1.5$, and then switches expression off completely
between $t\approx2.4$ and $t\approx7$, during the substrate peak. During this phase its
byproduct has accumulated to its maximum ($2.15$ near $t\approx3.8$) and inhibits its reaction. A
second, weaker expression phase (up to $0.43$) follows while substrate declines. Module B expresses
smoothly, from $0.43$ at the start down to $0.24$ near the substrate peak and up to $0.36$ near
$t\approx10$, so its enzyme pool varies only between $0.78$ and $0.94$. B's expression exceeds A's
during most of the cycle, so B does not free-ride on A's product. The shared product follows the
substrate with a lag and peaks at $2.16$ near $t\approx3.9$. Both squared-effort budgets remain
slack. Byproduct A, the product, and enzyme B end at the edge of their return bands
($|x_j(T)-x_j(0)|\approx10^{-2}$), so the return condition, rather than a budget, restricts the
end of the cycle.}\label{fig:\caseidentifier}
\end{figure}
\myheader{Literature and provenance.} Limited enzyme-production capacity is an established optimal-control modeling assumption \citep{Klipp2002-sg,Oyarzun2008PathwayArxiv,Oyarzun2009-ef}. We extend it to two competing expression modules with private budgets and unequal product/effort weights; the noncooperative split is a proposed benchmark.
All rates, budgets, weights, and initial values below are synthetic,
dimensionless benchmark choices, not fitted values from this source.
\myheader{Dimensions.} $N=2$, $n_x=6$, $n_{\rm ODE}=10$, $(n_{u,i})_{i\in\mathcal P}=(2\times 1)$, $(n_{g,i})_{i\in\mathcal P}=(2\times 13)$.
\myheader{Categories.} \cat{molecular}, \cat{periodical}, \cat{regulation}.
\myheader{Game specification.} $X_i$: byproduct pools; $S$: substrate; $Y$: product; $E_i$: active enzyme. $u_i$: enzyme-expression effort. $v_i$ denotes productive flux.
\[\begin{gathered}
v_i=\frac{V_iE_iS}{(K_i+S)(1+X_i/K_{X,i})},\quad \dot X_i=v_i-d_iX_i,\\
\dot S=D(S_{\rm in}-S)-\sum_{j\in\mathcal P}v_j,\quad \dot Y=\sum_{j\in\mathcal P}v_j-d_YY,\\
\dot E_i=\alpha_iu_i-\beta_iE_i.
\end{gathered}\]
Each reaction yields one unit of product and one unit of a module-specific byproduct $X_i$, which inhibits its own reaction with inhibition constant $K_{X,i}$ (\texttt{KX}) and is removed at rate $d_i$.
Controls are open-loop functions of time.
Initial-state order: $x=(X_A,X_B,S,Y,E_A,E_B)$; values appear below.

Each player $i\in\mathcal P$ minimizes $J_i=T^{-1}\int_0^T\ell_i\,dt$, with
\[\ell_i=-w_i^1Y+w_i^2u_i^2.\]
$w_i^k\geq0$, $\sum_kw_i^k=1$; weights are independently configurable.
\[
\dot Z_i=u_i^2,\quad Z_i(0)=0,\quad Z_i(T)\leq \mathcal B_i.
\]
Budgets are private to each player. The optional path floor is $S\geq\mathrm{floor}$.
When enabled, it is enforced in every player's deviation problem.

\emph{Periodical setting.} The substrate supply is modulated as
$S_{\rm in}(t)=S_{\rm in}[1+a\sin(2\pi t/T_{\rm per})]$ with forcing period $T_{\rm per}=T$. The fixed initial state lies on
the reference cycle obtained with constant expression efforts, and every player's problem
contains the return conditions $|x_j(T)-x_j(0)|\leq\delta$ for all six pools.

Enabled bounds:
\begin{center}\small\begin{tabular}{lll}\toprule Owners & Bound & Location\\\midrule
A & $y_{9}(T)\leq 3.6$ & terminal \\
B & $y_{10}(T)\leq 3$ & terminal \\
A, B & $|x_j(T)-x_j(0)|\leq 0.01$, $j=1,\ldots,6$ & terminal (return) \\
\bottomrule\end{tabular}\end{center}
\myheader{Parameter instance.} %
\begin{lstlisting}
initial = (0.2, 0.3)
forcing_period = 12.0
T = 12.0
forcing_amplitude = 0.8
periodicity_tolerance = 0.01
V = (1.2, 1.0)
K = (0.4, 0.6)
KX = (1.0, 1.0)
d = (0.1, 0.08)
D = 0.4
Sin = 3.0
dp = 0.2
al = (0.8, 0.7)
be = (0.3, 0.25)
x0 = (1.773752359, 2.261203403, 1.032172585, 1.721238584, 0.5333333333, 0.84)
w_A = (0.9, 0.1)
w_B = (0.7, 0.3)
budgets = (3.6, 3.0)
floor = nothing
control_lower = [0.0, 0.0]
control_upper = [1.0, 1.0]
numerical_x0 = [1.773752359, 2.261203403, 1.032172585, 1.721238584, 0.5333333333, 0.84, 0.0, 0.0, 0.0, 0.0]
\end{lstlisting}

\endgroup

\begingroup
\def\caseidentifier{microbial}%
\subsection{microbial: Microbial cross-feeding}\label{case:microbial}
\myheader{Biological interpretation.} Many microbial communities are held together by cross-feeding: strains that cannot synthesize
an essential metabolite (for example, an amino acid) obtain it from partners that secrete it,
while all strains compete for a common carbon source. Secretion is costly for the producer but
benefits its partner, and thereby indirectly the producer itself, because the partner supplies the
producer's own requirement. In a chemostat with a recurring feed with a period of $80$~h (for example, a
feed-and-starve regime of a maintained culture), substrate availability oscillates. We model one full feed
cycle of a maintained community of ten strains organized in reciprocal pairs; every strain and every
pool must return to its initial value, and each strain chooses its secretion allocation over the
cycle.

We extend the two-strain game of \citet{BangaVillaverde2026Hierarchical} to $N\geq2$ strains, each producing one metabolite
and choosing one allocation $u_i$. The state has $2N+1$ entries:
\[
x=((X_i)_{i\in\mathcal P},S,(M_i)_{i\in\mathcal P}).
\]
Let $a_{ij}\in\{0,1\}$ indicate that strain $i$ requires metabolite $M_j$,
with $a_{ii}=0$ and $d_i=\sum_{j\in\mathcal P} a_{ij}\geq1$.

Edges describe essential nutrients; multiple incoming nutrients
are complementary, not interchangeable.
Each metabolite is a shared extracellular pool. Its consumers compete
through depletion, while all strains compete for $S$.
This graph-based extension and the following $N=10$ setting are
proposals; the source treats only two strains.
\begin{figure}[p]
\centering
\includegraphics[width=\textwidth,height=\casefigheight,keepaspectratio]{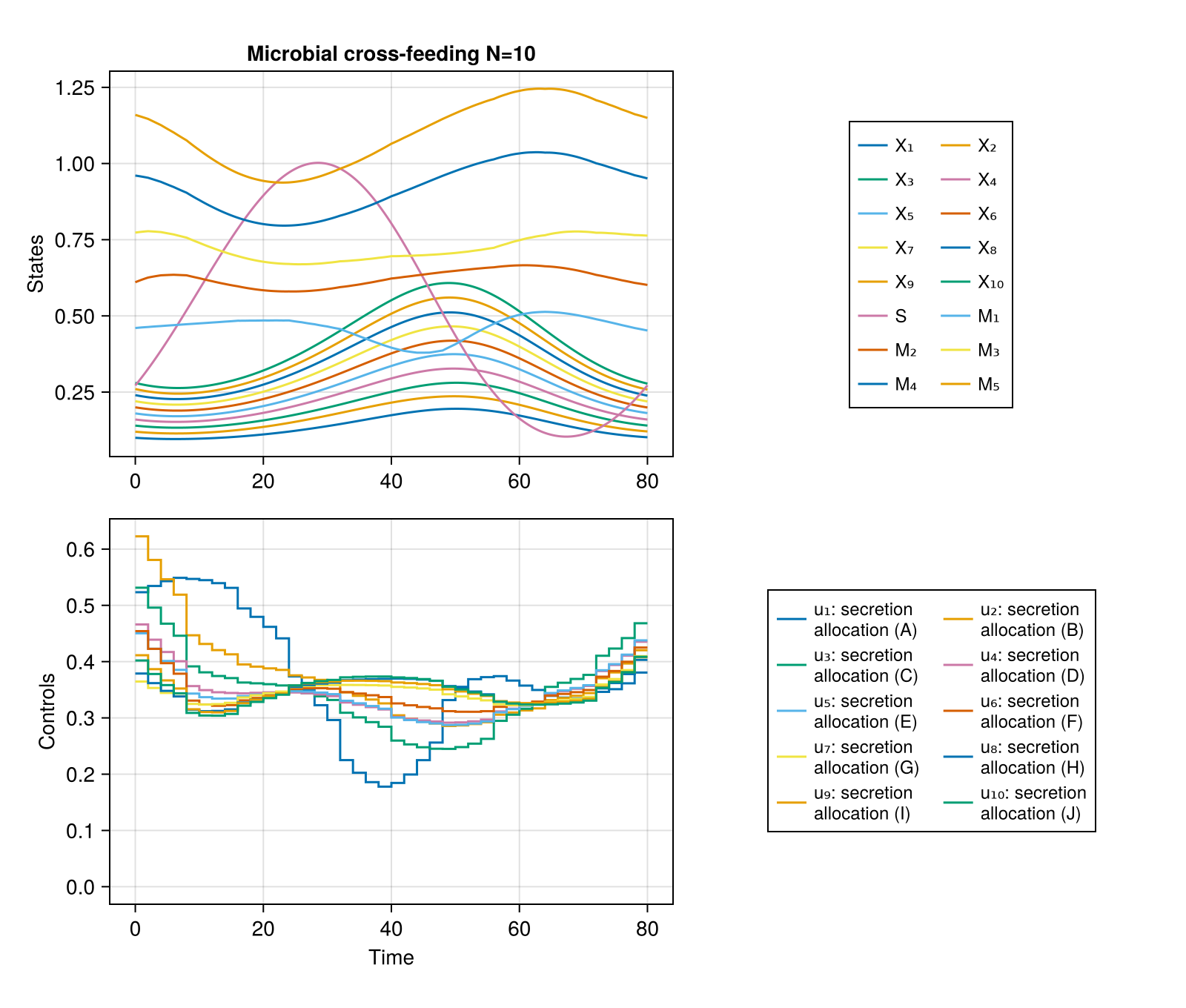}
\caption[Reference equilibrium candidate for \texttt{\caseidentifier}]{Reference equilibrium candidate for \texttt{\caseidentifier}. \textit{Interpretation.} The substrate (plotted scaled by $0.1$) rises from $2.7$ to about $10$ at $t\approx28$ during the
high-feed half of the cycle and falls to about $1$ near $t\approx67$. The strains grow with a delay
and roughly double their densities by $t\approx50$, then decline back to their initial values;
the metabolite pools vary more weakly and in a strain-specific way (only the pools $M_A$--$M_E$ are
drawn, since the plot shows at most $16$ states). All strains secrete strongly at the beginning of
the cycle, when their partners are still at low density and need the exchanged metabolites to
exploit the rising substrate, and reduce secretion as growth becomes substrate-limited. In the
low-feed phase between $t\approx30$ and $t\approx60$ several strains cut secretion markedly
(A down to $0.18$, C to $0.25$), because a partner's growth brings little benefit when substrate is
scarce. All strains raise secretion again toward the end of the cycle; this increase is forced by
the return condition on the metabolite pools rather than a biological prediction. The strains differ in their schedules according to their growth and
secretion traits and their partner-oriented weights. No budget or state bounds besides the return
conditions are enabled; the absolute return band of $10^{-2}$ gives the strains room to deviate
from the reference allocation.}\label{fig:\caseidentifier}
\end{figure}
\myheader{Literature and provenance.} The two-strain starting formulation is taken from \citet{BangaVillaverde2026Hierarchical}. The network extension, its ten-strain instance, and partner-averaged objective candidates are benchmark additions.
\myheader{Dimensions.} $N=10$, $n_x=21$, $n_{\rm ODE}=61$, $(n_{u,i})_{i\in\mathcal P}=(10\times 1)$, $(n_{g,i})_{i\in\mathcal P}=(10\times 42)$.
\myheader{Categories.} \cat{cellular}, \cat{periodical}, \cat{mutualism}, \cat{variablePlayers}.
\myheader{Game specification.} \[
\mu_i=\mu_{i,\max}(1-c_i u_i)\frac{S}{K_{S,i}+S}
       \prod_{j\in\mathcal P:a_{ij}=1}\frac{M_j}{K_{M,ij}+M_j}.
\]
\[
\begin{aligned}
\dot X_i&=(\mu_i-D)X_i,\qquad i\in\mathcal P,\\
\dot S&=D(S_{\mathrm{in}}-S)-\sum_{i\in\mathcal P}\frac{\mu_iX_i}{Y_{S,i}},\\
\dot M_j&=\alpha_j u_jX_j-DM_j
       -\sum_{i\in\mathcal P:a_{ij}=1}\frac{\mu_iX_i}{Y_{M,ij}},\qquad j\in\mathcal P.
\end{aligned}
\]
$Y_{M,ij}>0$ is the yield of consumer $i$ on metabolite $j$.
All consumers of a depleted $M_j$ have zero growth, preserving nonnegativity.

\textbf{Reduction check:} $N=2$, $a_{AB}=a_{BA}=1$,
$K_{M,ij}=K_{M,i}$ and $Y_{M,ij}=Y_{M,i}$ recover all five source equations.

With $\langle z\rangle=T^{-1}\int_0^Tz\,\mathrm dt$, the objectives are
\[
\begin{aligned}
J_{i1}&=-\langle X_i/B_0\rangle,&
J_{i2}&=\langle u_i^2+\epsilon_M M_i^2\rangle,\\
J_{i3}&=\left\langle\frac1{d_i}\sum_{j\in\mathcal P:a_{ij}=1}
                     \frac{K_{\mathrm{def}}}{K_{\mathrm{def}}+M_j}\right\rangle,&
J_{i4}&=-\left\langle\frac1{d_i}\sum_{j\in\mathcal P:a_{ij}=1}X_j/B_0\right\rangle.
\end{aligned}
\]
Player $i$ minimizes $J_i=\sum_{k=1}^4w_i^kJ_{ik}$,
$w_i\geq0$, $\sum_kw_i^k=1$.
Partner terms are averaged to avoid a purely degree-induced increase in scale.
No terminal cost; the only state restrictions are the return conditions below. Feature scales: $B_0=2$, $K_{\rm def}=0.35$, $\epsilon_M=10^{-5}$. The term $\epsilon_MM_i^2$ is a small regularization that penalizes accumulation of the own secreted metabolite.

Time is in h, $X_i$ in gDW L$^{-1}$, and $S,M_i$ in mmol L$^{-1}$.
$\mu_{i,\max},D$ have units h$^{-1}$; $K_{S,i},K_{M,ij},S_{\rm in}$
have concentration units. $\alpha_i$ is in mmol gDW$^{-1}$ h$^{-1}$;
$Y_{S,i},Y_{M,ij}$ are in gDW mmol$^{-1}$; $c_i,u_i$ are dimensionless.
$B_0$ has biomass units; $K_{\rm def}$ has metabolite units;
$\epsilon_M$ has inverse squared metabolite units.

For the default $N=10$, partners are reciprocal pairs
$(A,B),(C,D),\ldots,(I,J)$, coupled through shared substrate.
For odd $N>2$, the final strain additionally consumes its predecessor's
metabolite. Each player rewards partner biomass as well as its own.
All graph, kinetic, initial-state, and weight entries are configurable.
The $N=2$ parameters and graph remain unchanged.

\emph{Periodical setting.} The substrate feed is modulated as
$S_{\rm in}(t)=S_{\rm in}[1+a\sin(2\pi t/T_{\rm per})]$ with forcing period $T_{\rm per}=T$. Phase-zero strain densities are
fixed; the initial substrate and metabolite pools and the maximal growth rates $\mu_{i,\max}$ are
calibrated so that the reference cycle with constant allocations closes. Every player's problem
contains the return conditions $|x_j(T)-x_j(0)|\leq\delta$ for all $2N+1$ states. With a tight
band, these shared conditions pin every strain's allocation to its reference value, so the
archived instance uses the looser absolute band $\delta=10^{-2}$.

Enabled bounds:
\begin{center}\small\begin{tabular}{lll}\toprule Owners & Bound & Location\\\midrule
A, B, C, D, E, F, G, H, I, J & $|x_j(T)-x_j(0)|\leq 0.01$, $j=1,\ldots,21$ & terminal (return) \\
\bottomrule\end{tabular}\end{center}
\myheader{Parameter instance.} %
\begin{lstlisting}
forcing_period = 80.0
forcing_amplitude = 0.9
periodicity_tolerance = 0.01
N = 10
T = 80.0
B0 = 2.0
Kdef = 0.35
epsilon_M = 1.0e-5
mu_max = [0.1390691175, 0.1534971554, 0.1506718399, 0.1622049706, 0.1648060806, 0.1752300865, 0.1804094155, 0.190391922, 0.1970624911, 0.2069131292]
K_S = [5.7, 6.4, 7.1, 7.8, 8.5, 9.2, 9.9, 10.6, 11.3, 12.0]
K_M = [0.072 0.08 0.06 0.068 0.076 0.084 0.064 0.072 0.08 0.06; 0.076 0.084 0.064 0.072 0.08 0.06 0.068 0.076 0.084 0.064; 0.08 0.06 0.068 0.076 0.084 0.064 0.072 0.08 0.06 0.068; 0.084 0.064 0.072 0.08 0.06 0.068 0.076 0.084 0.064 0.072; 0.06 0.068 0.076 0.084 0.064 0.072 0.08 0.06 0.068 0.076; 0.064 0.072 0.08 0.06 0.068 0.076 0.084 0.064 0.072 0.08; 0.068 0.076 0.084 0.064 0.072 0.08 0.06 0.068 0.076 0.084; 0.072 0.08 0.06 0.068 0.076 0.084 0.064 0.072 0.08 0.06; 0.076 0.084 0.064 0.072 0.08 0.06 0.068 0.076 0.084 0.064; 0.08 0.06 0.068 0.076 0.084 0.064 0.072 0.08 0.06 0.068]
c = [0.17, 0.19, 0.21, 0.23, 0.25, 0.27, 0.29, 0.31, 0.33, 0.35]
alpha = [0.49, 0.53, 0.57, 0.61, 0.65, 0.69, 0.73, 0.77, 0.81, 0.85]
Y_S = [0.515, 0.53, 0.545, 0.56, 0.575, 0.59, 0.605, 0.62, 0.635, 0.65]
Y_M = [8.0 8.0 8.0 8.0 8.0 8.0 8.0 8.0 8.0 8.0; 8.0 8.0 8.0 8.0 8.0 8.0 8.0 8.0 8.0 8.0; 8.0 8.0 8.0 8.0 8.0 8.0 8.0 8.0 8.0 8.0; 8.0 8.0 8.0 8.0 8.0 8.0 8.0 8.0 8.0 8.0; 8.0 8.0 8.0 8.0 8.0 8.0 8.0 8.0 8.0 8.0; 8.0 8.0 8.0 8.0 8.0 8.0 8.0 8.0 8.0 8.0; 8.0 8.0 8.0 8.0 8.0 8.0 8.0 8.0 8.0 8.0; 8.0 8.0 8.0 8.0 8.0 8.0 8.0 8.0 8.0 8.0; 8.0 8.0 8.0 8.0 8.0 8.0 8.0 8.0 8.0 8.0; 8.0 8.0 8.0 8.0 8.0 8.0 8.0 8.0 8.0 8.0]
D = 0.05
S_in = 10.0
x0 = [0.1, 0.12, 0.14, 0.16, 0.18, 0.2, 0.22, 0.24, 0.26, 0.28, 2.716672258, 0.4602340172, 0.6105073918, 0.7733384874, 0.9605044771, 1.159541914, 1.383575536, 1.61888453, 1.879760125, 2.151413228, 2.449110336]
weights = [0.4 0.01 0.08 0.51; 0.4 0.01222222222 0.08 0.5077777778; 0.4 0.01444444444 0.08 0.5055555556; 0.4 0.01666666667 0.08 0.5033333333; 0.4 0.01888888889 0.08 0.5011111111; 0.4 0.02111111111 0.08 0.4988888889; 0.4 0.02333333333 0.08 0.4966666667; 0.4 0.02555555556 0.08 0.4944444444; 0.4 0.02777777778 0.08 0.4922222222; 0.4 0.03 0.08 0.49]
adjacency = [0 1 0 0 0 0 0 0 0 0; 1 0 0 0 0 0 0 0 0 0; 0 0 0 1 0 0 0 0 0 0; 0 0 1 0 0 0 0 0 0 0; 0 0 0 0 0 1 0 0 0 0; 0 0 0 0 1 0 0 0 0 0; 0 0 0 0 0 0 0 1 0 0; 0 0 0 0 0 0 1 0 0 0; 0 0 0 0 0 0 0 0 0 1; 0 0 0 0 0 0 0 0 1 0]
initial_control = [0.35, 0.35, 0.35, 0.35, 0.35, 0.35, 0.35, 0.35, 0.35, 0.35]
control_lower = [0.02, 0.02, 0.02, 0.02, 0.02, 0.02, 0.02, 0.02, 0.02, 0.02]
control_upper = [0.95, 0.95, 0.95, 0.95, 0.95, 0.95, 0.95, 0.95, 0.95, 0.95]
numerical_x0 = [0.1, 0.12, 0.14, 0.16, 0.18, 0.2, 0.22, 0.24, 0.26, 0.28, 2.716672258, 0.4602340172, 0.6105073918, 0.7733384874, 0.9605044771, 1.159541914, 1.383575536, 1.61888453, 1.879760125, 2.151413228, 2.449110336, 0.0, 0.0, 0.0, 0.0, 0.0, 0.0, 0.0, 0.0, 0.0, 0.0, 0.0, 0.0, 0.0, 0.0, 0.0, 0.0, 0.0, 0.0, 0.0, 0.0, 0.0, 0.0, 0.0, 0.0, 0.0, 0.0, 0.0, 0.0, 0.0, 0.0, 0.0, 0.0, 0.0, 0.0, 0.0, 0.0, 0.0, 0.0, 0.0, 0.0]
\end{lstlisting}

\endgroup

\begingroup
\def\caseidentifier{mitochondria}%
\subsection{mitochondria: Mitochondrial maintenance competition}\label{case:mitochondria}
\myheader{Biological interpretation.} Mitochondria continually accumulate damage and are maintained by repair, fusion, and biogenesis,
all of which consume energy. In tissues with recurring activity and rest, such as muscle under a
daily activity cycle, the cell's energetic capacity oscillates, and the maintenance of the
mitochondrial network has to be timed accordingly. Cells often carry more than one mitochondrial
DNA (mtDNA) variant (heteroplasmy), and mutant mtDNA can expand clonally within a cell because
mitochondria replicate and are removed independently of cell division; fusion and fission mix the
populations and thereby modulate this expansion \citep{BenchmarkPMC4440705}. We therefore interpret
the two players as a wild-type (A) and a variant (B) mtDNA subpopulation that draw on one shared
energetic pool: the variant contributes slightly less to the energy supply, converts energy into
biogenesis less efficiently, and is lost slightly faster, but has a larger replication allowance.
Each subpopulation invests in biogenesis to keep its own healthy mass high. The quadratic conversion
of healthy into damaged mass is a phenomenological density-dependent damage term, not a mechanism
taken from the source. We model one full activity cycle, after which all pools must return to their
initial values, with biogenesis effort limited by a budget.
\begin{figure}[p]
\centering
\includegraphics[width=\textwidth,height=\casefigheight,keepaspectratio]{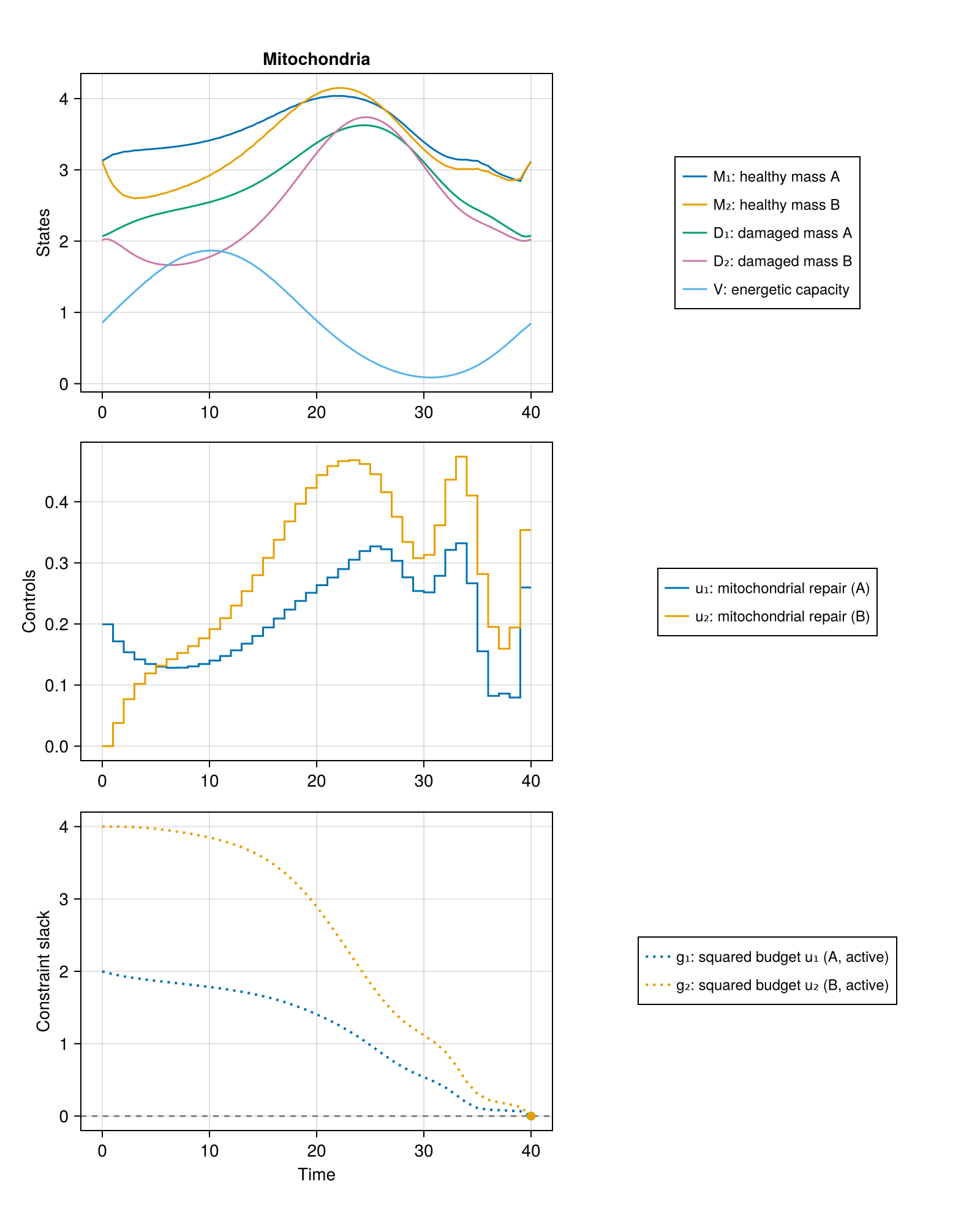}
\caption[Reference equilibrium candidate for \texttt{\caseidentifier}]{Reference equilibrium candidate for \texttt{\caseidentifier}. \textit{Interpretation.} The energetic capacity $V$ follows the activity cycle: it rises to about $1.9$ during the high-supply
first quarter ($t\approx10$), falls to a minimum near $0.09$ at $t\approx31$, and is restored at the
end of the cycle. Biogenesis driven by $V$ makes both healthy masses grow with a delay, peaking near
$t\approx22$ at about $4.0$--$4.2$; quadratic damage lets the damaged masses follow them and peak a
few time units later, and all four pools decline during the low-energy phase. Biogenesis effort
follows the healthy-mass peak: both populations increase effort steadily from
$t\approx5$ to $t\approx25$, reduce it while the energy pool is depleted, and make a second push
around $t\approx33$, when energy starts to recover. The variant subpopulation B, which starts without effort and
lets its healthy mass drop initially, invests more than A from $t\approx5$ onward. Shortly before
the end, both reduce effort sharply and then raise it again in the final interval; this
correction is forced by the return condition. Both effort budgets are exhausted ($g_1,g_2$ active).}\label{fig:\caseidentifier}
\end{figure}
\myheader{Literature and provenance.} \citet{BenchmarkPMC4440705} model the clonal expansion of mtDNA mutations and show that the effect of fusion and fission on this expansion depends on the cellular context. This motivates competing mtDNA subpopulations with a shared energetic environment as players. The source does not formulate a game; the controlled biogenesis, the individual objectives, the shared energy pool, and the periodic activity forcing are benchmark additions, and the wild-type/variant assignment of the parameter contrasts is an interpretation rather than a calibration.
All rates, budgets, weights, and initial values below are synthetic,
dimensionless benchmark choices, not fitted values from this source.
\myheader{Dimensions.} $N=2$, $n_x=5$, $n_{\rm ODE}=9$, $(n_{u,i})_{i\in\mathcal P}=(2\times 1)$, $(n_{g,i})_{i\in\mathcal P}=(2\times 11)$.
\myheader{Categories.} \cat{molecular}, \cat{periodical}, \cat{regulation}.
\myheader{Game specification.} $M_i$: healthy mass; $D_i$: damaged mass; $V$: energetic capacity. $u_i$: maintenance/biogenesis effort. Quadratic damage is a phenomenological crowding term.
\[\begin{gathered}
\dot M_i=s_i+a_iu_iV/(1+V)-d_iM_i-k_fM_i^2+k_rD_i,\\
\dot D_i=k_fM_i^2-(k_r+d_D)D_i,\\
\dot V=v_A+v_B-[d_V+\eta(D_A+D_B+u_A+u_B)]V.
\end{gathered}\]
Controls are open-loop functions of time.
Initial-state order: $x=(M_A,M_B,D_A,D_B,V)$; values appear below.

Each player $i\in\mathcal P$ minimizes $J_i=T^{-1}\int_0^T\ell_i\,dt$, with
\[\ell_i=-w_i^1M_i+w_i^2u_i^2+w_i^3D_i^2.\]
$w_i^k\geq0$, $\sum_kw_i^k=1$; weights are independently configurable.
\[
\dot Z_i=u_i^2,\quad Z_i(0)=0,\quad Z_i(T)\leq \mathcal B_i.
\]
Budgets are private to each player. The optional path floor is $M_A\geq\mathrm{floor}$.
When enabled, it is enforced in every player's deviation problem.

\emph{Periodical setting.} The energetic supply is modulated as
$(v_A+v_B)[1+a\sin(2\pi t/T_{\rm per})]$ with forcing period $T_{\rm per}=T$, representing recurring activity and rest.
The fixed initial state is a point on the reference cycle obtained with constant efforts.
Every player's problem contains the return conditions $|x_j(T)-x_j(0)|\leq\delta$ for all five
pools; the archived instance uses the absolute band $\delta=10^{-2}$, because a tight band lets
sequential best responses drift (Section~\ref{sec:methods}).

Enabled bounds:
\begin{center}\small\begin{tabular}{lll}\toprule Owners & Bound & Location\\\midrule
A & $y_{8}(T)\leq 2$ & terminal \\
B & $y_{9}(T)\leq 4$ & terminal \\
A, B & $|x_j(T)-x_j(0)|\leq 0.01$, $j=1,\ldots,5$ & terminal (return) \\
\bottomrule\end{tabular}\end{center}
\myheader{Parameter instance.} %
\begin{lstlisting}
initial = (0.2, 0.3)
forcing_period = 40.0
T = 40.0
forcing_amplitude = 0.9
periodicity_tolerance = 0.01
s = (0.1, 0.1)
amp = (5.0, 4.0)
d = (0.1, 0.12)
kf = 0.08
kr = 0.3
dd = 0.05
v = (0.5, 0.45)
dv = 0.5
eta = 0.1
x0 = (3.128867036, 3.124015553, 2.066994658, 2.016342084, 0.8547536323)
w_A = (0.7, 0.2, 0.1)
w_B = (0.45, 0.4, 0.15)
budgets = (2.0, 4.0)
floor = nothing
control_lower = [0.0, 0.0]
control_upper = [1.0, 1.0]
numerical_x0 = [3.128867036, 3.124015553, 2.066994658, 2.016342084, 0.8547536323, 0.0, 0.0, 0.0, 0.0]
\end{lstlisting}

\endgroup

\begingroup
\def\caseidentifier{motor}%
\subsection{motor: Physically coupled motor control}\label{case:motor}
\myheader{Biological interpretation.} Joint action, such as two people carrying an object or moving a shared handle, requires each
partner to reach their own goals while feeling the other's forces. In the experiments of
\citet{Chackochan2019JointAction}, two participants move handles that are connected by a virtual spring;
each has to pass through a private via-point on the way to a common target. Because the via-points
differ, the partners' preferences conflict, and the observed movements are well described by a
Nash equilibrium between two optimal controllers. We model one reaching episode of $0.5$ seconds
with arm dynamics driven by muscle activation and force states. The dynamics are linear and the
costs quadratic; nonlinearity enters only through the control bounds and the interior-time
via-point costs. This makes \texttt{motor} a candidate for an analytical cross-check of
numerical equilibria.

Two people move mechanically coupled handles toward a common target through
different via-points. The published model uses individual accuracy and effort
costs \cite{Chackochan2019JointAction}; related motor games are studied by
\cite{Braun2009Motor}.

Each player controls a two-dimensional muscle command $u_i(t)$.
Spring forces couple both movements and create conflicting preferences.
Our benchmark uses deterministic continuous-time dynamics and
open-loop policies. The paper also studies noise, partial information and
feedback policies; those are distinct extensions.

The published supplementary data contain per-trial indicators but not within-trial
trajectories, so the instance cannot be compared with measured trajectories directly.
\begin{figure}[p]
\centering
\includegraphics[width=\textwidth,height=\casefigheight,keepaspectratio]{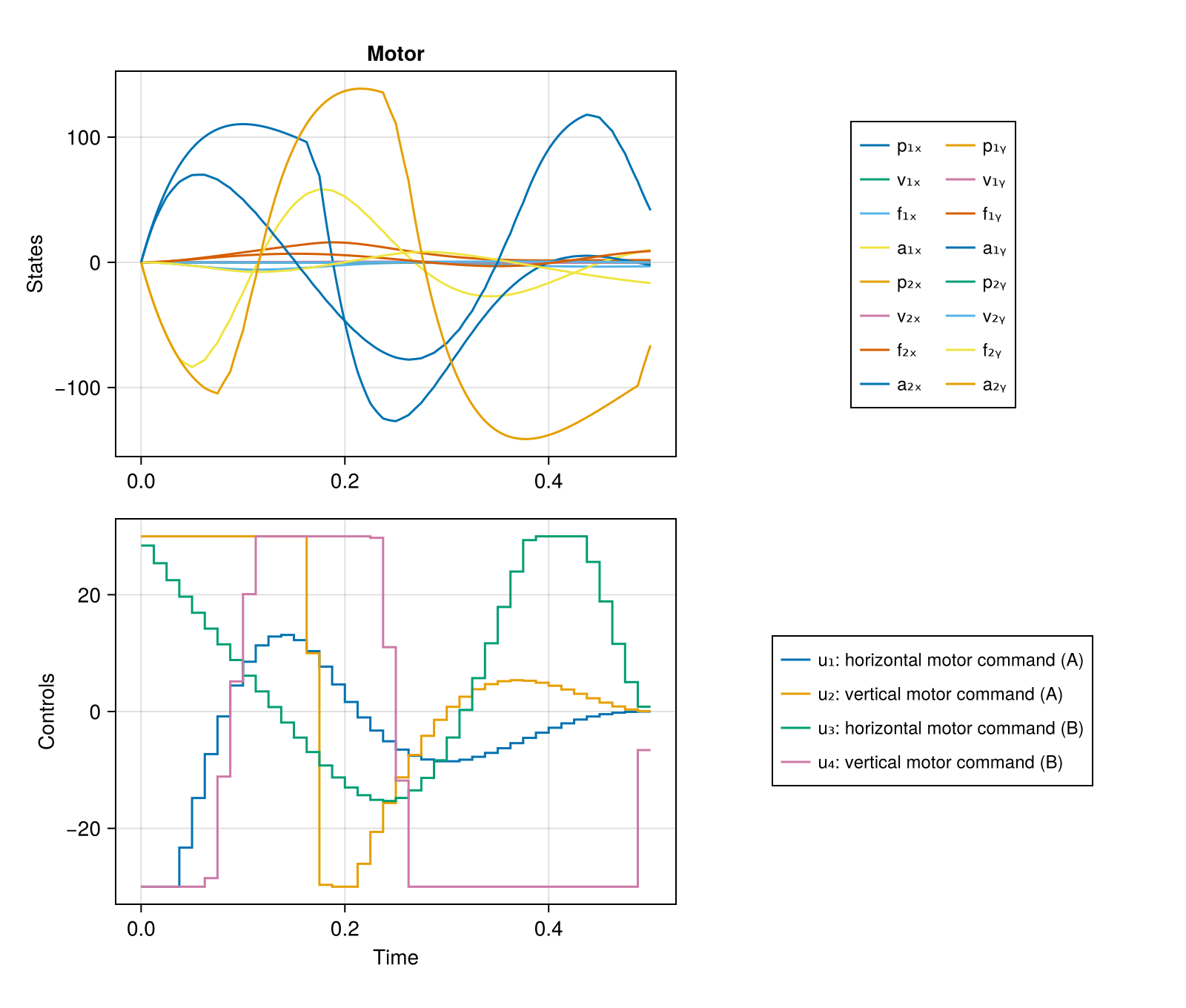}
\caption[Reference equilibrium candidate for \texttt{\caseidentifier}]{Reference equilibrium candidate for \texttt{\caseidentifier}. \textit{Interpretation.} The hand positions reveal the negotiated movement. Partner A must pass its via-point at $x=-0.02$ at
$t=0.2$, partner B its via-point at $x=+0.02$ at $t=0.4$, and both hands must end at the common
target $(0,0.1)$. Both partners use saturated bang-type commands at the bound $\pm30$ for large
parts of the movement; the control bounds are active, so the instance differs qualitatively from
the smooth minimum-effort solutions of the unconstrained source model: A first pulls to the left and pushes upward, reverses its vertical command
shortly before its via-point, and then relaxes; B first pushes to the right, reverses to the left
around $t\approx0.24$ while A is at its via-point, and pushes right again to reach its own
via-point at $t=0.4$, braking its vertical motion with a sustained negative command. A's hand
deviates leftward to $-0.025$ near $t\approx0.26$ and B's rightward to $0.035$ near $t\approx0.3$
before both converge on the target. The command reversals implement acceleration and braking and
compensate the spring force transmitted by the partner; no additional restrictions are imposed.
The state panel is dominated by the muscle force and activation variables, which reach magnitudes above $100$, whereas the hand positions move only a few centimeters.}\label{fig:\caseidentifier}
\end{figure}
\myheader{Literature and provenance.} The mechanical and muscle dynamics, physical coefficients, and accuracy/effort scales follow the supplementary model of \citet{Chackochan2019JointAction}. This benchmark fixes the via-point times and uses deterministic open-loop commands with explicit bounds; the published model also considers timing search, feedback, noise, and partial information. The parameter listing fixes the selected numerical instance.
\myheader{Dimensions.} $N=2$, $n_x=16$, $n_{\rm ODE}=18$, $(n_{u,i})_{i\in\mathcal P}=(2\times 2)$, $(n_{g,i})_{i\in\mathcal P}=(2\times 0)$.
\myheader{Categories.} \cat{organismal}, \cat{episodical}, \cat{regulation}.
\myheader{Game specification.} For $i\in\{A,B\}$ and partner $j$, let $p_i,v_i,f_i,a_i\in\mathbb R^2$.
In gravity-compensated coordinates, the source mechanics become
\[
\begin{aligned}
\dot p_i&=v_i,\\
\dot v_i&=\bigl(f_i+k(p_j-p_i)-b v_i\bigr)/m_i,\\
\dot f_i&=a_i,\\
\dot a_i&=(u_i-f_i-2\tau a_i)/\tau^2.
\end{aligned}
\]
There are 16 physical states and four scalar control components.
The spring force on A is the negative of that on B.
\[
m_A=m_B=2\ \mathrm{kg},\quad k=100\ \mathrm{N/m},\quad
b=10\ \mathrm{Ns/m},\quad \tau=0.10\ \mathrm{s}.
\]
These values and the second-order muscle dynamics are specified in the
paper's supplementary model note.

For fixed via-point times $t_i$ and common target $p_T$, minimize
$J_i=\sum_{k=1}^5w_i^kJ_{ik}$, $w_i\geq0$, $\sum_kw_i^k=1$, with
\[
\begin{aligned}
J_{i1}&=\|p_i(T)-p_T\|^2/s_T^2,&
J_{i2}&=\|v_i(T)\|^2/s_v^2,\\
J_{i3}&=\|p_i(t_i)-p_{V,i}\|^2/s_V^2,&
J_{i4}&=\frac1T\int_0^T\frac{\|p_i-p_j\|^2}{s_d^2}\,\mathrm dt,\\
J_{i5}&=\frac1T\int_0^T\frac{\|u_i\|^2}{s_u^2}\,\mathrm dt.
\end{aligned}
\]
Source scales: $s_T=0.005$ m, $s_v=0.005$ m/s, $s_V=0.0025$ m,
$s_d=0.025$ m, $s_u=15$ N. Equal normalized weights correspond to the
source's unit effort multiplier, up to an irrelevant common objective scale.

Fixing $t_i$ is our simplification; the source additionally searches over
crossing times. The implementation evaluates interior-time costs at control-grid points
$t_i$; both the coarse and refined grids must contain these times.

No state or budget bounds are enabled.
\myheader{Parameter instance.} %
\begin{lstlisting}
T = 0.5
m = (2.0, 2.0)
k = 100.0
b = 10.0
tau = 0.1
s_T = 0.005
s_v = 0.005
s_V = 0.0025
s_d = 0.025
s_u = 15.0
target = (0.0, 0.1)
via_A = (-0.02, 0.05)
via_B = (0.02, 0.05)
via_times = (0.2, 0.4)
w_A = (0.2, 0.2, 0.2, 0.2, 0.2)
w_B = (0.2, 0.2, 0.2, 0.2, 0.2)
umax = 30.0
x0 = [0.0, 0.0, 0.0, 0.0, 0.0, 0.0, 0.0, 0.0, 0.0, 0.0, 0.0, 0.0, 0.0, 0.0, 0.0, 0.0]
initial = (0.0, 0.0, 0.0, 0.0)
control_lower = [-30.0, -30.0, -30.0, -30.0]
control_upper = [30.0, 30.0, 30.0, 30.0]
numerical_x0 = [0.0, 0.0, 0.0, 0.0, 0.0, 0.0, 0.0, 0.0, 0.0, 0.0, 0.0, 0.0, 0.0, 0.0, 0.0, 0.0, 0.0, 0.0]
\end{lstlisting}

\endgroup

\begingroup
\def\caseidentifier{neural}%
\subsection{neural: Excitatory--inhibitory regulation}\label{case:neural}
\myheader{Biological interpretation.} Cortical circuits alternate between states of higher and lower excitability, and the balance
between excitation and inhibition must be maintained throughout. We interpret the Wilson--Cowan
populations as one local circuit subject to a slow modulation of background drive, e.g., by
neuromodulatory tone, over many population time constants; units are dimensionless. We
introduce two regulatory influences: one supports a higher target level of excitatory activity,
the other favors lower excitation and a moderate inhibitory level. The game is periodical: after
one full cycle of the background drive, both activities must return to their initial values.

Wilson--Cowan models describe coupled excitatory and inhibitory population
activity \cite{WilsonCowan1972}. We introduce two regulatory agents controlling
the drive to each population. Their distinct targets and effort costs form
the game extension used here.

Agent A chooses excitatory drive $u_A$ to support activity.
Agent B chooses inhibitory drive $u_B$ to regulate excitation.
The states $E,I$ are dimensionless activity variables.
\begin{figure}[p]
\centering
\includegraphics[width=\textwidth,height=\casefigheight,keepaspectratio]{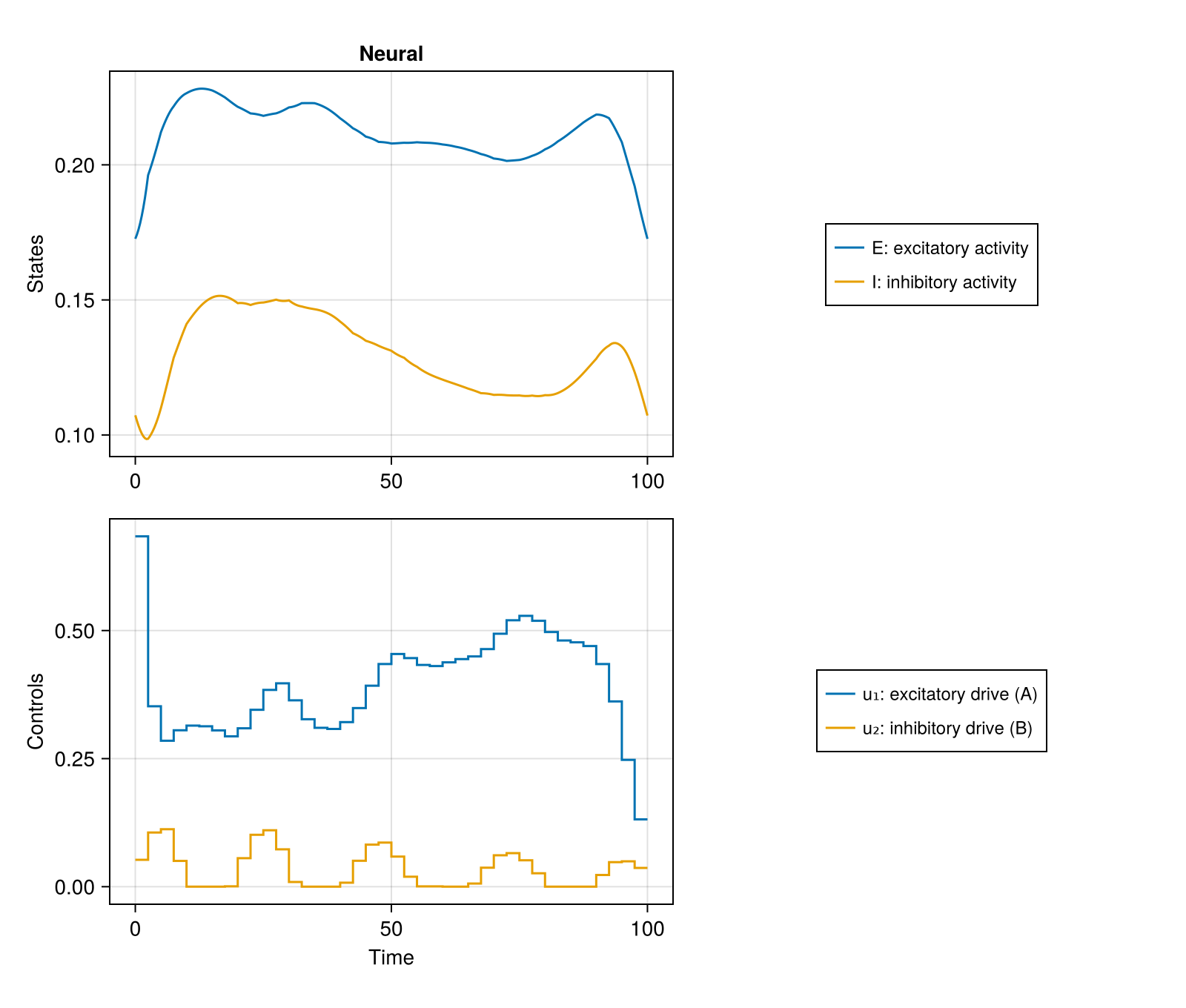}
\caption[Reference equilibrium candidate for \texttt{\caseidentifier}]{Reference equilibrium candidate for \texttt{\caseidentifier}. \textit{Interpretation.} Over the slow cycle ($T=100$), player A applies a strong excitatory drive at the beginning to
lift the excitatory activity $E$ from $0.17$ toward its target, then settles to a moderate drive
that follows the slow background oscillation (lower during the first half, when the background
drive is high, and higher in the second half), and finally withdraws drive, because the return condition requires $E$ to return
to its initial level at the end of the cycle. Player B, who prefers lower excitation, uses its
inhibitory drive sparingly in four short pulses of about $0.1$, which shape the inhibitory
activity $I$ without large cost. As a result $E$ stays between $0.20$ and $0.23$ for most of the
cycle, a compromise closer to B's target than to A's, with local maxima near $t\approx13$ and $t\approx90$, while $I$ rises to $0.15$ early and
relaxes toward $0.115$ in the second half. No effort budgets or additional path bounds are enabled; the return
conditions alone force the terminal descent of both activities.}\label{fig:\caseidentifier}
\end{figure}
\myheader{Literature and provenance.} The excitatory--inhibitory population structure follows \citet{WilsonCowan1972}. We use an unshifted sigmoid, add sinusoidal background drive, and assign the two control inputs to regulators with different tracking targets. The objective weights, control bounds, return tolerance, and tuned numerical instance are benchmark choices.
\myheader{Dimensions.} $N=2$, $n_x=2$, $n_{\rm ODE}=4$, $(n_{u,i})_{i\in\mathcal P}=(2\times 1)$, $(n_{g,i})_{i\in\mathcal P}=(2\times 4)$.
\myheader{Categories.} \cat{organismal}, \cat{periodical}, \cat{regulation}.
\myheader{Game specification.} Use a bounded, unshifted sigmoid variant of the population model:
\[
S_k(z)=\frac{1}{1+\exp[-a_k(z-b_k)]},\qquad k\in\{E,I\},
\]
\[
\begin{aligned}
\tau_E\dot E&=-E+(1-E)S_E(c_{EE}E-c_{EI}I+\Pi(t)+u_A),\\
\tau_I\dot I&=-I+(1-I)S_I(c_{IE}E-c_{II}I+Q+u_B).
\end{aligned}
\]
The unshifted sigmoid is an explicit modeling choice; parameter values are
starting guesses, not a reproduction of the 1972 experiments.
For $E(0),I(0)\in[0,1]$, the vector field points inward on the boundary.

Population-level input drives represent effective regulation.

Write $\langle z\rangle=T^{-1}\int_0^Tz\,\mathrm dt$ and define
\[
\begin{aligned}
J_A&=w_A^1\langle(E-E_A)^2\rangle
       +w_A^2\langle(u_A/\bar u)^2\rangle
       +w_A^3(E(T)-E_A)^2,\\
J_B&=w_B^1\langle(E-E_B)^2\rangle
       +w_B^2\langle(u_B/\bar u)^2\rangle
       +w_B^3\langle(I-I_B)^2\rangle.
\end{aligned}
\]
Both objectives are minimized, $w_i\geq0$, $\sum_kw_i^k=1$.
Choose $E_A>E_B$ to introduce conflicting excitation targets. Because the return band fixes
$E(T)\approx E(0)$, the terminal term $w_A^3(E(T)-E_A)^2$ is nearly constant and has no
strategic effect.
\[
w_A=(0.8,0.1,0.1),\quad w_B=(0.4,0.5,0.1).
\]
There are no effort budgets; the biological return bands below are imposed on both players.
The control bounds and all objective weights are synthetic choices.

\emph{Periodical setting.} The background excitatory input is modulated as
$\Pi(t)=\bar\Pi+a\sin(2\pi t/T_{\rm per})$ with forcing period $T_{\rm per}=T$, representing a slow excitability cycle. The fixed
initial activities lie on the reference cycle obtained with constant drives, and every player's
problem contains the return conditions $|x_j(T)-x_j(0)|\leq\delta$ for $E$ and $I$.

Enabled bounds:
\begin{center}\small\begin{tabular}{lll}\toprule Owners & Bound & Location\\\midrule
A, B & $|x_j(T)-x_j(0)|\leq 0.0001$, $j=1,\ldots,2$ & terminal (return) \\
\bottomrule\end{tabular}\end{center}
\myheader{Parameter instance.} %
\begin{lstlisting}
forcing_period = 100.0
T = 100.0
forcing_amplitude = 0.2
periodicity_tolerance = 0.0001
c_EE = 16.0
c_EI = 12.0
c_IE = 15.0
c_II = 3.0
a_E = 1.3
a_I = 2.0
b_E = 4.0
b_I = 3.7
tau_E = 8.0
tau_I = 8.0
P = 1.0
Q = 0.0
E_A = 0.35
E_B = 0.15
I_B = 0.2
umax = 2.0
x0 = (0.1726288025, 0.1072703456)
initial = (0.5, 0.5)
w_A = (0.8, 0.1, 0.1)
w_B = (0.4, 0.5, 0.1)
control_lower = [0.0, 0.0]
control_upper = [2.0, 2.0]
numerical_x0 = [0.1726288025, 0.1072703456, 0.0, 0.0]
\end{lstlisting}

\endgroup

\begingroup
\def\caseidentifier{phage}%
\subsection{phage: Bacterium--phage interaction}\label{case:phage}
\myheader{Biological interpretation.} Bacteria and their viruses (phages) coevolve in an arms race: bacteria express surface
modifications, restriction systems, or CRISPR immunity that reduce adsorption, and phages
respond with receptor-binding variants that increase infectivity, often at a cost in burst size.
Within a single growth episode, for example one batch culture or one bloom, these traits are
regulated phenotypically rather than by mutation, and their expression follows the commanded
level only with a delay. We model one such episode: the bacterial population invests in
resistance against the cost of reduced growth, while the phage population adjusts the
composition of its released progeny, e.g., via phase-variable or host-range-variant
subpopulations, and thereby their average infectivity, to maximize its own net expansion.

Use five states to distinguish population densities, the latent period of infection,
and expressed phenotypes:
\[
x=(x_1,x_2,x_3,x_4,x_5)^\top
=\begin{pmatrix}\text{susceptible bacteria}\\\text{infected bacteria}\\\text{free-phage density}\\
\text{expressed resistance}\\\text{expressed infectivity}\end{pmatrix}.
\]
A chooses target resistance $u_A$; B chooses target infectivity $u_B$:
\[
\tau_A\dot x_4=u_A-x_4,\quad \tau_B\dot x_5=u_B-x_5.
\]
Positive $\tau_A,\tau_B$ describe response/turnover times, not mutation rates.
If $x_4(0),x_5(0)\in[0,1]$, these phenotype bounds are invariant.

The bacterial target represents regulation; the phage target represents progeny composition. Free phage cannot
instantaneously alter an infection apparatus. Heritable changes to target
policies belong to a separate evolutionary timescale.
\begin{figure}[p]
\centering
\includegraphics[width=\textwidth,height=\casefigheight,keepaspectratio]{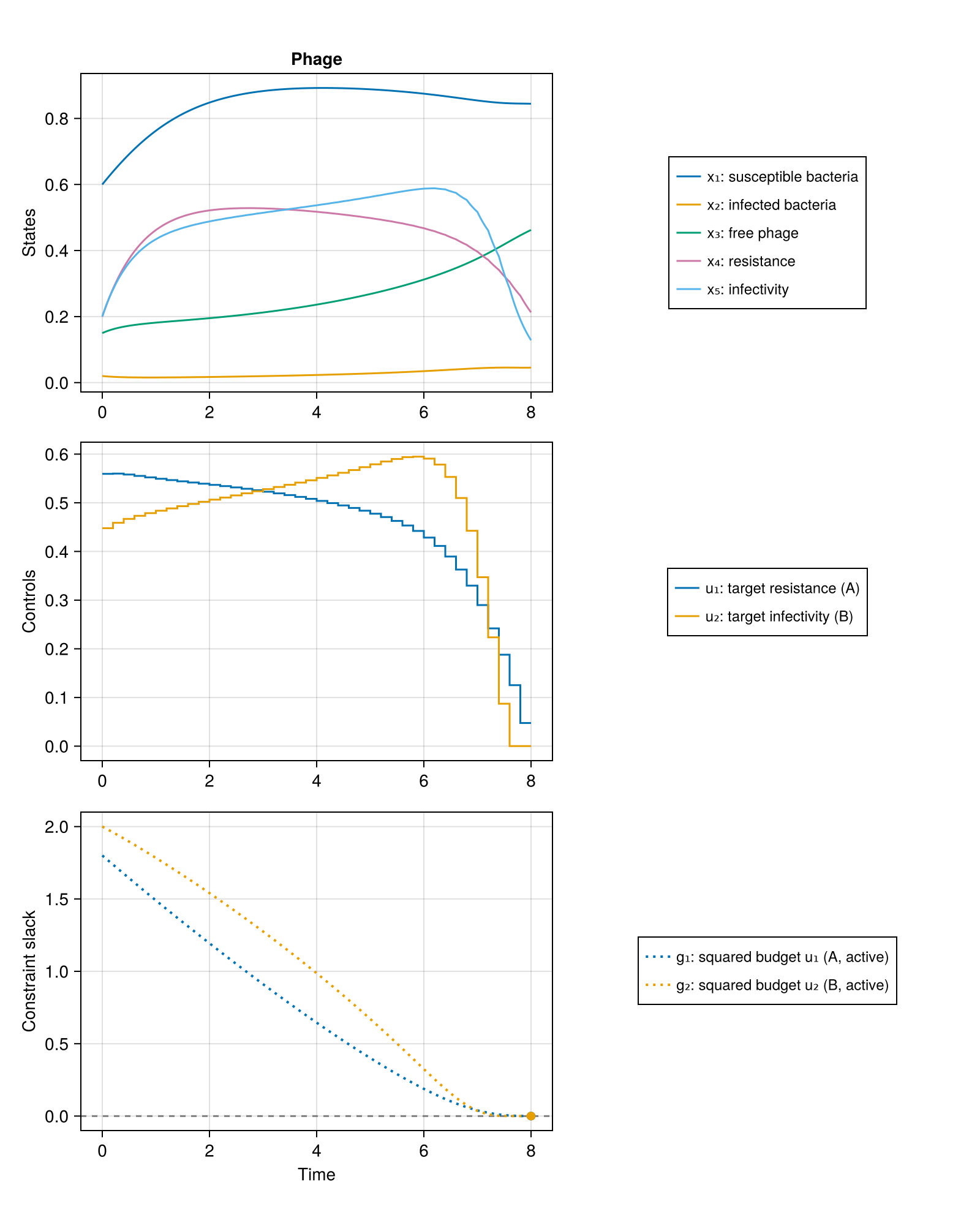}
\caption[Reference equilibrium candidate for \texttt{\caseidentifier}]{Reference equilibrium candidate for \texttt{\caseidentifier}. \textit{Interpretation.} The bacterium starts with a resistance target of $0.56$ and lowers it slowly during most of the
episode, then faster after $t\approx6$ (to $0.05$ at the end). The phage raises its infectivity
target from $0.45$ to $0.60$ at $t\approx5.8$ and then cuts it to zero during the last $0.4$ time
units. The expressed phenotypes follow with delays: expressed resistance peaks at $0.53$ near
$t\approx2.8$, and expressed infectivity at $0.59$ near $t\approx6.2$. Susceptible bacteria grow
from $0.60$ to $0.89$ at $t\approx4.1$ and then decline to $0.84$ as free phage accumulate
roughly threefold (from $0.15$ to $0.46$), accelerating toward the end. The latent compartment
stays small (between $0.016$ and $0.046$) but delays the phage return on infectivity. Both players
exhaust their squared-effort budgets ($g_1,g_2$ active). A mechanism consistent with the late
relaxation of both targets is that effort spent near $T$ changes the expressed phenotype, and
through the latent period the released progeny, only after the episode has ended; this
interpretation was not tested by a separate comparison.}\label{fig:\caseidentifier}
\end{figure}
\myheader{Literature and provenance.} \emph{Pseudomonas fluorescens} SBW25 and phage $\Phi2$ motivate this game.
\cite{Paterson2010Coevolution} compared phage evolution
with a co-evolving host and with a fixed host genotype.

Experiments with \emph{P. aeruginosa} and different phages distinguish
arms races from fluctuating selection \cite{Betts2014Contrasted}.
These studies motivate reciprocal adaptation; the following equations and
continuous controls are our modeling hypotheses.

A representative bacterial lineage invests in resistance.
A representative phage lineage invests in infectivity.

Resistance reduces infection but costs growth. Infectivity increases adsorption
but reduces progeny yield. Their objectives form a general-sum game.

The model represents the latent period by one exponentially distributed infected-cell
stage and omits explicit resources, spatial structure, and genotype diversity.
\myheader{Dimensions.} $N=2$, $n_x=5$, $n_{\rm ODE}=9$, $(n_{u,i})_{i\in\mathcal P}=(2\times 1)$, $(n_{g,i})_{i\in\mathcal P}=(2\times 1)$.
\myheader{Categories.} \cat{cellular}, \cat{episodical}, \cat{antagonism}.
\myheader{Game specification.} Let adsorption and burst size depend on the expressed phenotypes:
\[
a(x_4,x_5;p)=a_0e^{\alpha_Bx_5-\alpha_Ax_4},\qquad
\beta(x_5;p)=1+(\beta_0-1)e^{-c_Bx_5^2},\quad\beta_0>1.
\]
With $a=a(x_4,x_5;p)$ and $\beta=\beta(x_5;p)$,
\[
f(x,u_A,u_B;p)=\begin{pmatrix}
r x_1[1-(x_1+x_2)/K]-c_Ax_4^2x_1-a\,x_1x_3\\
a\,x_1x_3-x_2/\tau_L\\
\beta x_2/\tau_L-a(x_1+x_2)x_3-m_Bx_3\\
(u_A-x_4)/\tau_A\\(u_B-x_5)/\tau_B
\end{pmatrix}.
\]
\[
p=(r,K,c_A,a_0,\alpha_A,\alpha_B,\beta_0,c_B,m_B,\tau_L,\tau_A,\tau_B).
\]
$r,K$: resource-limited growth, with infected cells occupying capacity; $c_A$: resistance cost;
$m_B$: phage decay; $\tau_L$: mean latent period between adsorption and lysis.
$\alpha_A,\alpha_B$: adsorption sensitivities; $c_B$: infectivity--yield trade-off.
Free phage are lost by adsorption to susceptible and infected cells; the burst size is
evaluated with the infectivity expressed at lysis. The latent compartment delays the phage
return on infectivity investment, which is the mechanism that distinguishes this game from
the ciliate consumer--resource game.
Positive initial densities remain positive over a finite horizon.

Minimize $J_A=\sum_{k=1}^3w_A^kJ_{Ak}$ and
$J_B=\sum_{k=1}^3w_B^kJ_{Bk}$ with nonnegative, normalized weights.
All features use fixed unit scales in the scaled benchmark.

\renewcommand{\arraystretch}{1.5}
{\setlength{\tabcolsep}{4pt}\begin{tabular}{@{}lll@{}}
Objective & Weighted term & Interpretation\\\hline
$J_A$ & $-w_A^1\log[x_1(T)/x_{1,0}]$ & Maximize net expansion of uninfected cells\\
& $-w_A^2\int_0^T r x_1\,\mathrm dt$ & Maximize gross birth proxy\\
& $+w_A^3\int_0^T u_A^2\,\mathrm dt$ & Minimize target effort\\\hline
$J_B$ & $-w_B^1\log[x_3(T)/x_{3,0}]$ & Maximize net lineage expansion\\
& $-w_B^2\int_0^T\beta(x_5;p)x_2/\tau_L\,\mathrm dt$ & Maximize gross phage release\\
& $+w_B^3\int_0^T u_B^2\,\mathrm dt$ & Minimize target effort\\
\end{tabular}}

The gross-production hypotheses can favor turnover rather than persistence.
The bacterial $r x_1$ feature treats density limitation as crowding loss;
this decomposition is an additional biological hypothesis.

Use $\theta=(U_A,U_B,\varepsilon_A,\varepsilon_B,
\sigma_{U_A},\sigma_{U_B},\sigma_{\varepsilon_A},\sigma_{\varepsilon_B})$.
A switch $\sigma\in\{0,1\}$ determines whether its restriction is enabled.
\[
\begin{array}{ll}
g_A[x,u_A,u_B;p,\theta]\geq0:&U_A-\int_0^T u_A^2\,\mathrm dt\geq0\quad\text{if }\sigma_{U_A}=1,\\
&x_1(t)-\varepsilon_A\geq0\quad\text{if }\sigma_{\varepsilon_A}=1;\\[2mm]
g_B[x,u_A,u_B;p,\theta]\geq0:&U_B-\int_0^T u_B^2\,\mathrm dt\geq0\quad\text{if }\sigma_{U_B}=1,\\
&x_3(t)-\varepsilon_B\geq0\quad\text{if }\sigma_{\varepsilon_B}=1.
\end{array}
\]
Effort budgets concern target investment; expressed-phenotype costs already
enter $f$. Optional density floors encode persistence hypotheses.

Enabled bounds:
\begin{center}\small\begin{tabular}{lll}\toprule Owners & Bound & Location\\\midrule
A & $y_{8}(T)\leq 1.8$ & terminal \\
B & $y_{9}(T)\leq 2$ & terminal \\
\bottomrule\end{tabular}\end{center}
\myheader{Parameter instance.} %
\begin{lstlisting}
T = 8.0
r = 1.0
K = 1.0
c_A = 0.1
a0 = 0.3
alpha_A = 2.0
alpha_B = 1.5
beta0 = 4.0
c_B = 0.3
m_B = 0.4
latent = 0.5
tau_A = 0.75
tau_B = 0.5
x0 = (0.6, 0.02, 0.15, 0.2, 0.2)
w_A = (0.6, 0.3, 0.1)
w_B = (0.9, 0.0, 0.1)
scales_A = (1.0, 1.0, 1.0)
scales_B = (1.0, 1.0, 1.0)
U_A = 1.8
U_B = 2.0
epsilon_A = nothing
epsilon_B = nothing
umax = (1.0, 1.0)
initial = (0.2, 0.2)
control_lower = [0.0, 0.0]
control_upper = [1.0, 1.0]
numerical_x0 = [-0.5108256238, -3.912023005, -1.897119985, 0.2, 0.2, 0.0, 0.0, 0.0, 0.0]
\end{lstlisting}

\endgroup

\begingroup
\def\caseidentifier{plantNitrogen}%
\subsection{plantNitrogen: Competing plants and nitrogen}\label{case:plantNitrogen}
\myheader{Biological interpretation.} In temperate grasslands and crops, neighboring plants compete for a pulse of soil nitrogen in
early spring, when mineralization has built up a pool that roots can exploit before the growing
season proceeds. Nitrogen taken up early builds photosynthetic capacity, which later drives carbon gain, but the
plants also shade one another, and root uptake capacity is costly to build and maintain. Roots
cannot deplete soil mineral nitrogen below a minimal concentration at which net uptake ceases;
this threshold is represented as a shared floor on $S$. We model one growing season; its end is not required
to reproduce the initial state, since senescence and winter processes connect consecutive seasons.
\begin{figure}[p]
\centering
\includegraphics[width=\textwidth,height=\casefigheight,keepaspectratio]{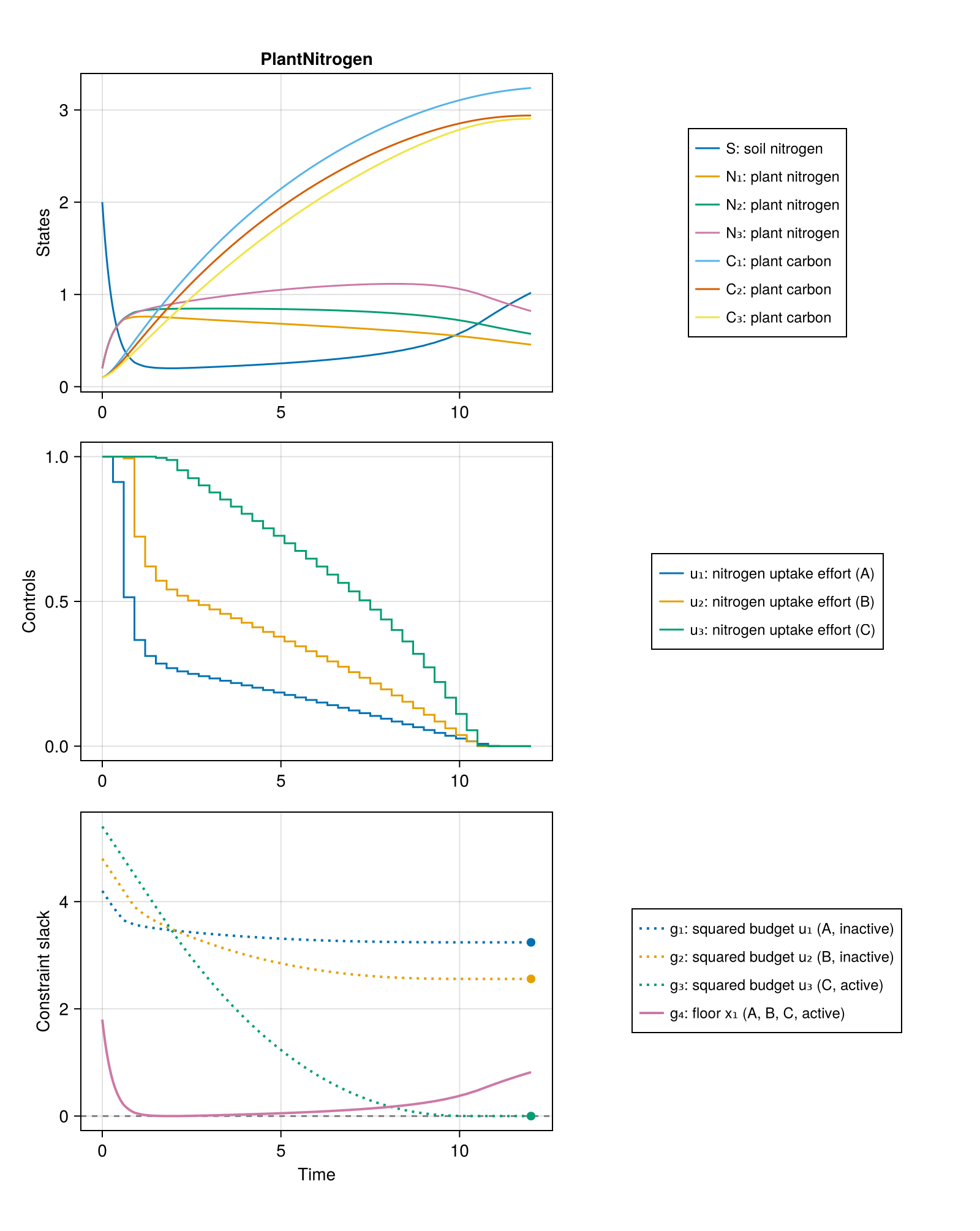}
\caption[Reference equilibrium candidate for \texttt{\caseidentifier}]{Reference equilibrium candidate for \texttt{\caseidentifier}. \textit{Interpretation.} All three plants start the season at maximal uptake and draw the soil nitrogen pool down from
$2$ to the uptake threshold of $0.2$ within less than two time units ($g_4$ active); afterwards the
floor limits the joint uptake, and the plants share the slowly replenished pool. They then
reduce uptake at different rates: A, which weighs effort costs most, drops to about $0.3$
almost immediately and tapers off; B keeps an intermediate schedule; and C, which values carbon
most, maintains high uptake longest and exhausts its effort budget ($g_3$ active). Plant nitrogen
therefore peaks early and moderately for A ($t\approx1.2$), later for B ($t\approx3$), and latest
and highest for C ($t\approx8$), after which it declines through turnover ($\rho_i$).
Carbon rises in a sigmoidal fashion for all plants and saturates toward the end of the season as
mutual shading increases; A reaches the largest carbon pool despite its smaller nitrogen store.
Once uptake has ceased around $t\approx11$, soil nitrogen recovers to about $1$.}\label{fig:\caseidentifier}
\end{figure}
\myheader{Literature and provenance.} Competition for nitrogen and asymmetric light acquisition motivate independent plant players \citep{BenchmarkplantNitrogen}. The players are separate plants rather than organs of one plant, so that independent objectives are not attributed to coordinated organs. Competition between separate plants for a shared soil resource has been formulated as a game in which each plant overinvests in roots relative to growing alone, a prediction confirmed for soybeans \citep{Gersani2001RootCommons}; here the corresponding investment is a time-dependent uptake effort.
All rates, budgets, weights, and initial values below are synthetic,
dimensionless benchmark choices, not fitted values from this source.
\myheader{Dimensions.} $N=3$, $n_x=7$, $n_{\rm ODE}=13$, $(n_{u,i})_{i\in\mathcal P}=(3\times 1)$, $(n_{g,i})_{i\in\mathcal P}=(3\times 2)$.
\myheader{Categories.} \cat{population}, \cat{seasonal}, \cat{competition}, \cat{variablePlayers}.
\myheader{Game specification.} $S$: soil nitrogen; $N_i$: plant nitrogen; $C_i$: carbon pool. $u_i$: uptake investment. $L$ is algebraic available light, not an ODE state.
\[\begin{gathered}
L=\frac{\mathrm{light}}{1+\mathrm{shade}\sum_{j\in\mathcal P}C_j},\quad
\dot S=D(S_{\rm in}-S)-S\sum_{j\in\mathcal P}u_j,\\
\dot N_i=u_iS-\rho_iN_i,\qquad
\dot C_i=V_iN_i\frac{L}{K_i+L}-d_iC_i.
\end{gathered}\]
Controls are open-loop functions of time.
Initial-state order: $x=(S,(N_i)_{i\in\mathcal P},(C_i)_{i\in\mathcal P})$; values appear below.

Each player $i\in\mathcal P$ minimizes $J_i=T^{-1}\int_0^T\ell_i\,dt$, with
\[\ell_i=-w_i^1C_i+w_i^2u_i^2+w_i^3N_i^2.\]
The penalty $N_i^2$ represents a cost of storing nitrogen, e.g., nitrate accumulation.
$w_i^k\geq0$, $\sum_kw_i^k=1$; weights are independently configurable.
\[
\dot Z_i=u_i^2,\quad Z_i(0)=0,\quad Z_i(T)\leq \mathcal B_i.
\]
Budgets are private to each player. The enabled path floor is $S\geq0.2$.
When enabled, it is enforced in every player's deviation problem.

Enabled bounds:
\begin{center}\small\begin{tabular}{lll}\toprule Owners & Bound & Location\\\midrule
A & $y_{11}(T)\leq 4.2$ & terminal \\
B & $y_{12}(T)\leq 4.8$ & terminal \\
C & $y_{13}(T)\leq 5.4$ & terminal \\
A, B, C & $0.2\leq y_{1}(t)$ & path \\
\bottomrule\end{tabular}\end{center}
\myheader{Parameter instance.} %
\begin{lstlisting}
N = 3
T = 12.0
rho = [0.1, 0.12, 0.14]
V = [1.0, 0.9090909091, 0.8333333333]
K = [0.3, 0.4, 0.5]
d = [0.08, 0.1, 0.12]
D = 0.2
Sin = 2.0
light = 1.0
shade = 0.1
x0 = [2.0, 0.2, 0.2, 0.2, 0.1, 0.1, 0.1]
w_A = (0.15, 0.8, 0.05)
w_B = (0.35, 0.55, 0.1)
w_C = (0.65, 0.25, 0.1)
budgets = [4.2, 4.8, 5.4]
floor = 0.2
initial = [0.2, 0.2, 0.2]
control_lower = [0.0, 0.0, 0.0]
control_upper = [1.0, 1.0, 1.0]
numerical_x0 = [2.0, 0.2, 0.2, 0.2, 0.1, 0.1, 0.1, 0.0, 0.0, 0.0, 0.0, 0.0, 0.0]
\end{lstlisting}

\endgroup

\begingroup
\def\caseidentifier{plantWater}%
\subsection{plantWater: Root investment across shared soil layers}\label{case:plantWater}
\myheader{Biological interpretation.} In seasonally dry ecosystems, plants must build their root systems during the wet part of the
season in order to extract water that is stored in the soil when rainfall declines. Soil layers
differ in how much water they hold and how quickly it is replenished, and species differ in the
depth at which their roots are most effective. Roots are costly to build and to maintain, and
water taken from a layer is no longer available to neighbors rooting in the same layer. We
model one growing season with a wet first half and a dry second half; four plants allocate root
investment across the layers and are limited by a seasonal water-uptake allowance. A and D
specialize on the top layer, B on the middle and C on the deep layer (rows of \code{uptake}). The
seasonal allowance is a synthetic proxy for a cumulative transpiration limit imposed by
hydraulic safety; it is not a physiological quota.

Players are $N\geq2$ plants in $L\geq2$ shared soil layers.
$W_k$: soil water; $R_{ik}$: plant $i$'s roots in layer $k$; $B_i$: biomass.
Each plant owns $L$ root-investment controls $u_{ik}$.
State order is soil water by layer, roots by player then layer, and plant biomass by player;
controls also list all layers of one plant before the next plant.
The archived dimensions are given below.
All controls depend on time alone.
\begin{figure}[p]
\centering
\includegraphics[width=\textwidth,height=\casefigheight,keepaspectratio]{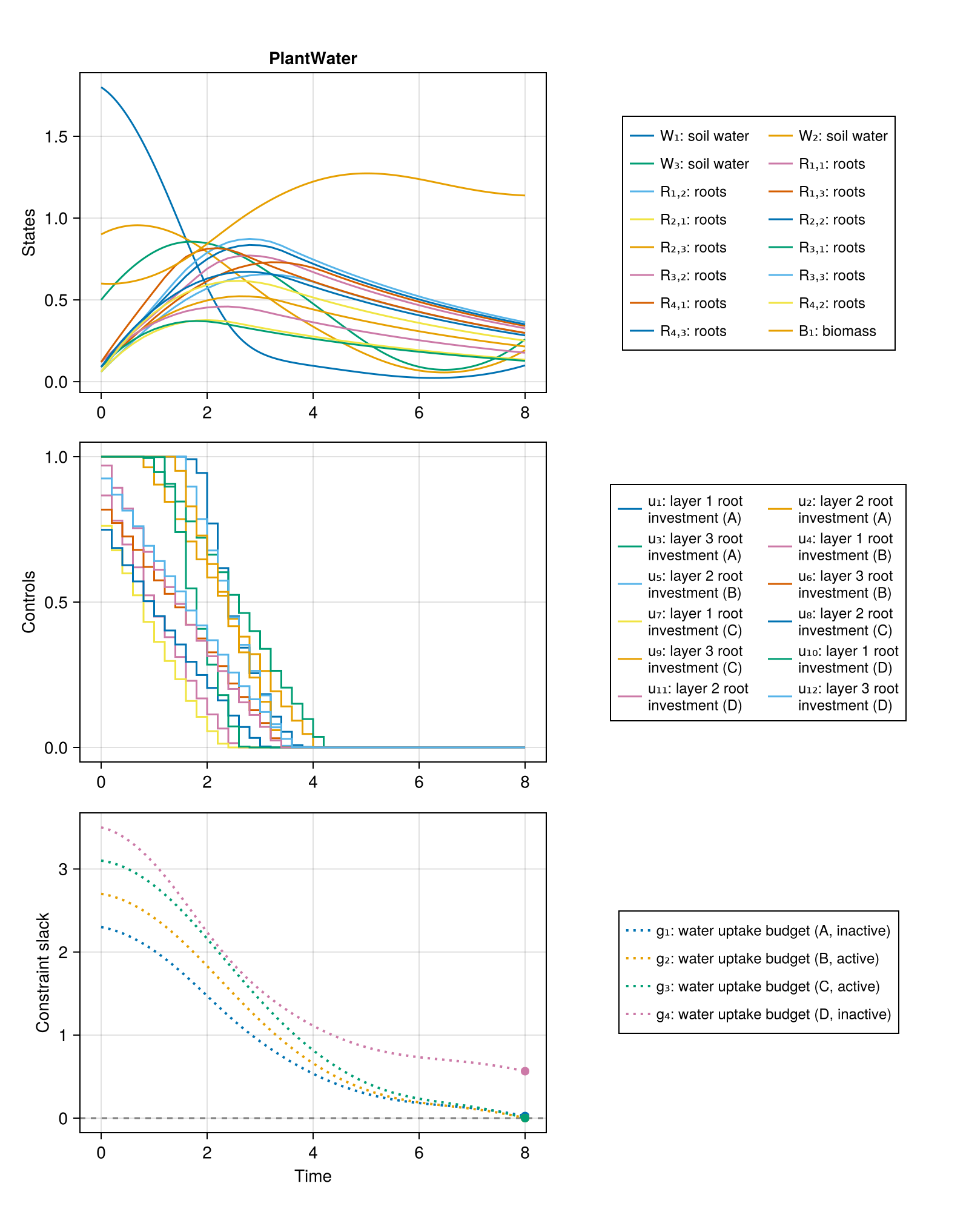}
\caption[Reference equilibrium candidate for \texttt{\caseidentifier}]{Reference equilibrium candidate for \texttt{\caseidentifier}. \textit{Interpretation.} All plants invest heavily in roots at the start of the season, mostly at the upper bound, and
shut off new investment layer by layer between $t\approx2$ and $t\approx4$, well before the dry
second half. Root stocks accordingly peak between $t\approx2$ and $t\approx3$ and then decay
slowly, but the established roots keep drawing water: the top layer, initially the wettest, is
depleted almost completely by $t\approx6$, and the middle and deep layers first gain water from
the wet-season rain and are then drawn down as rainfall fades. Each plant builds the largest root
stock in its preferred layer, and the two top-layer specialists split that layer. Biomass grows
by roughly a factor of two until $t\approx4$--$5$ and declines slightly during the dry end of the
season, when uptake no longer covers maintenance. Plants B and C exhaust their seasonal water
allowances ($g_2,g_3$ active), whereas A and D remain below theirs, so the allowance, rather
than root cost alone, limits the middle- and deep-layer specialists.}\label{fig:\caseidentifier}
\end{figure}
\myheader{Literature and provenance.} The source treats root proliferation as an evolutionary game in shared resource patches \citep{BenchmarkPMC2141669}. We construct a finite-horizon water-allocation game with explicit root and biomass dynamics, seasonal replenishment, and water-use quotas. Soil layers replace the source spatial geometry; this is not its calibrated model. A game-theoretic model with a glasshouse test predicts and observes root overproliferation and reduced reproductive yield when soybean plants share soil, a tragedy of the commons \citep{Gersani2001RootCommons}; the present game allows the analogous overinvestment in roots under a shared water resource.

All numerical defaults are synthetic and dimensionless; no fitted biological values or source equilibrium results are asserted.
\myheader{Dimensions.} $N=4$, $n_x=19$, $n_{\rm ODE}=27$, $(n_{u,i})_{i\in\mathcal P}=(4\times 3)$, $(n_{g,i})_{i\in\mathcal P}=(4\times 1)$.
\myheader{Categories.} \cat{population}, \cat{seasonal}, \cat{competition}, \cat{variablePlayers}.
\myheader{Game specification.} Water uptake and seasonally replenished soil obey
\[
V_{ik}=c_{ik}R_{ik}\frac{W_k}{h+W_k},\qquad
\dot W_k=I_k[1+s\sin(2\pi t/T)]-\delta W_k-\sum_{i\in\mathcal P}V_{ik}.
\]
\[
\dot R_{ik}=g u_{ik}\frac{B_i}{1+B_i}-d_RR_{ik},
\]
\[
\dot B_i=B_i\left[\frac{e\sum_kV_{ik}}{1+B_i}
-d_B-\kappa B_i-c_u\sum_ku_{ik}^2\right].
\]
$c_{ik}=\texttt{uptake}$, $h=\texttt{half}$, $I_k=\texttt{rain}$,
$s=\texttt{seasonality}$, $\delta=\texttt{evaporation}$;
$g,d_R,e,d_B,\kappa,c_u$ correspond to \texttt{root\_growth, root\_loss,
assimilation, mortality, crowding, root\_cost}.
$0\leq s\leq1$ preserves nonnegative inflows.

Each player minimizes
\[
J_i=\frac1T\int_0^T\sum_{k=1}^3w_i^k\ell_{ik}\,dt,
\qquad w_i^k\geq0,\quad\sum_kw_i^k=1.
\]
There is no terminal objective. Each weight row is independently configurable.

\[
(\ell_{i1},\ell_{i2},\ell_{i3})=
\left(-B_i,\ \sum_ku_{ik}^2,\ \sum_kR_{ik}^2\right).
\]
Candidates represent biomass benefit, investment effort, and root maintenance.
The growth equation also accounts for a physiological investment cost;
its coefficient and the preference weights have distinct roles.

Water allowance: $\dot Z_i=\sum_kV_{ik}$, $Z_i(0)=0$,
$Z_i(T)\leq\mathrm{budgets}_i$ (water amount, not effort).
For each plant's deviation enforce $W_k(t)\geq\mathrm{floors}_k$ for every enabled layer.
Default: four allowances; soil-water floors are disabled. Enabling all three
layer floors adds twelve player/layer entries (three distinct shared floors). Water is removed once in the soil balance, irrespective of
how many players enforce a floor.

For $q$ biological states, append $N$ cost integrals then $N$ budget integrals,
all initially zero; the cost integrals include the factor $1/T$, the ODE dimension is $q+2N$, and $J_i=x_{q+\iota(i)}(T)$.

Enabled bounds:
\begin{center}\small\begin{tabular}{lll}\toprule Owners & Bound & Location\\\midrule
A & $y_{24}(T)\leq 2.3$ & terminal \\
B & $y_{25}(T)\leq 2.7$ & terminal \\
C & $y_{26}(T)\leq 3.1$ & terminal \\
D & $y_{27}(T)\leq 3.5$ & terminal \\
\bottomrule\end{tabular}\end{center}
\myheader{Parameter instance.} %
\begin{lstlisting}
N = 4
layers = 3
T = 8.0
rain = [0.25, 0.4, 0.55]
seasonality = 0.8
evaporation = 0.12
uptake = [0.9 0.3 0.3; 0.324 0.972 0.324; 0.348 0.348 1.044; 1.116 0.372 0.372]
half = 0.4
root_growth = 0.9
root_loss = 0.18
assimilation = 1.2
mortality = 0.05
crowding = 0.03
root_cost = 0.02
x0 = [1.8, 0.9, 0.5, 0.12, 0.06, 0.09, 0.06, 0.09, 0.12, 0.09, 0.12, 0.06, 0.12, 0.06, 0.09, 0.6, 0.8, 1.0, 1.2]
budgets = [2.3, 2.7, 3.1, 3.5]
weights = [0.85 0.06 0.09; 0.8 0.08666666667 0.1133333333; 0.75 0.1133333333 0.1366666667; 0.7 0.14 0.16]
floors = [nothing, nothing, nothing]
initial = [0.2, 0.2, 0.2, 0.2, 0.2, 0.2, 0.2, 0.2, 0.2, 0.2, 0.2, 0.2]
control_lower = [0.0, 0.0, 0.0, 0.0, 0.0, 0.0, 0.0, 0.0, 0.0, 0.0, 0.0, 0.0]
control_upper = [1.0, 1.0, 1.0, 1.0, 1.0, 1.0, 1.0, 1.0, 1.0, 1.0, 1.0, 1.0]
numerical_x0 = [1.8, 0.9, 0.5, 0.12, 0.06, 0.09, 0.06, 0.09, 0.12, 0.09, 0.12, 0.06, 0.12, 0.06, 0.09, 0.6, 0.8, 1.0, 1.2, 0.0, 0.0, 0.0, 0.0, 0.0, 0.0, 0.0, 0.0]
\end{lstlisting}

\endgroup

\begingroup
\def\caseidentifier{pollination}%
\subsection{pollination: Nectar investment and pollinator visitation}\label{case:pollination}
\myheader{Biological interpretation.} Flowering plants and their pollinators interact during a limited flowering season. Plants invest
in nectar, which is collected by visitors that transfer pollen, while
pollinators allocate foraging time among plant species according to the nectar reward and their
own handling efficiency. Nectar is costly for the plant, visitation is costly for the pollinator,
and pollinator populations depend on the nectar they collect. We model one flowering season of two
plant and two pollinator species; the plant state $P_i$ is reproductive biomass (flowers and
developing seeds), which pollinator visits increase within the season; each pollinator is more efficient on one plant species, and all
players are limited by effort budgets and population floors. The end of the season is linked to the
next one through seed and overwintering stages, not through periodic return.

There are $n_p\geq2$ plant and $n_q\geq2$ pollinator players ($N=n_p+n_q$).
$P_i$: reproductive biomass, $\nu_i$: nectar, $A_j$: pollinators;
The state order is plant abundances, their nectar pools, then pollinator abundances, alphabetically within each group.
Plant $i$ chooses nectar investment $s_i$. Pollinator $j$ chooses $n_p$
visitation intensities $a_{ji}$, with actual allocation
$\pi_{ji}=a_{ji}/(1+\sum_{k\in\mathcal P_{\rm plant}}a_{jk})$; residual time is not spent visiting.
Control order: all $s_i$, then the $n_p$ intensities of each pollinator.
All controls depend on time alone. The player set $\mathcal P_{\rm plant}$ contains the first $n_p$ players, and $\mathcal P_{\rm poll}=\mathcal P\setminus\mathcal P_{\rm plant}$. In the archived instance, these are $\{A,B\}$ and $\{C,D\}$.
\begin{figure}[p]
\centering
\includegraphics[width=\textwidth,height=\casefigheight,keepaspectratio]{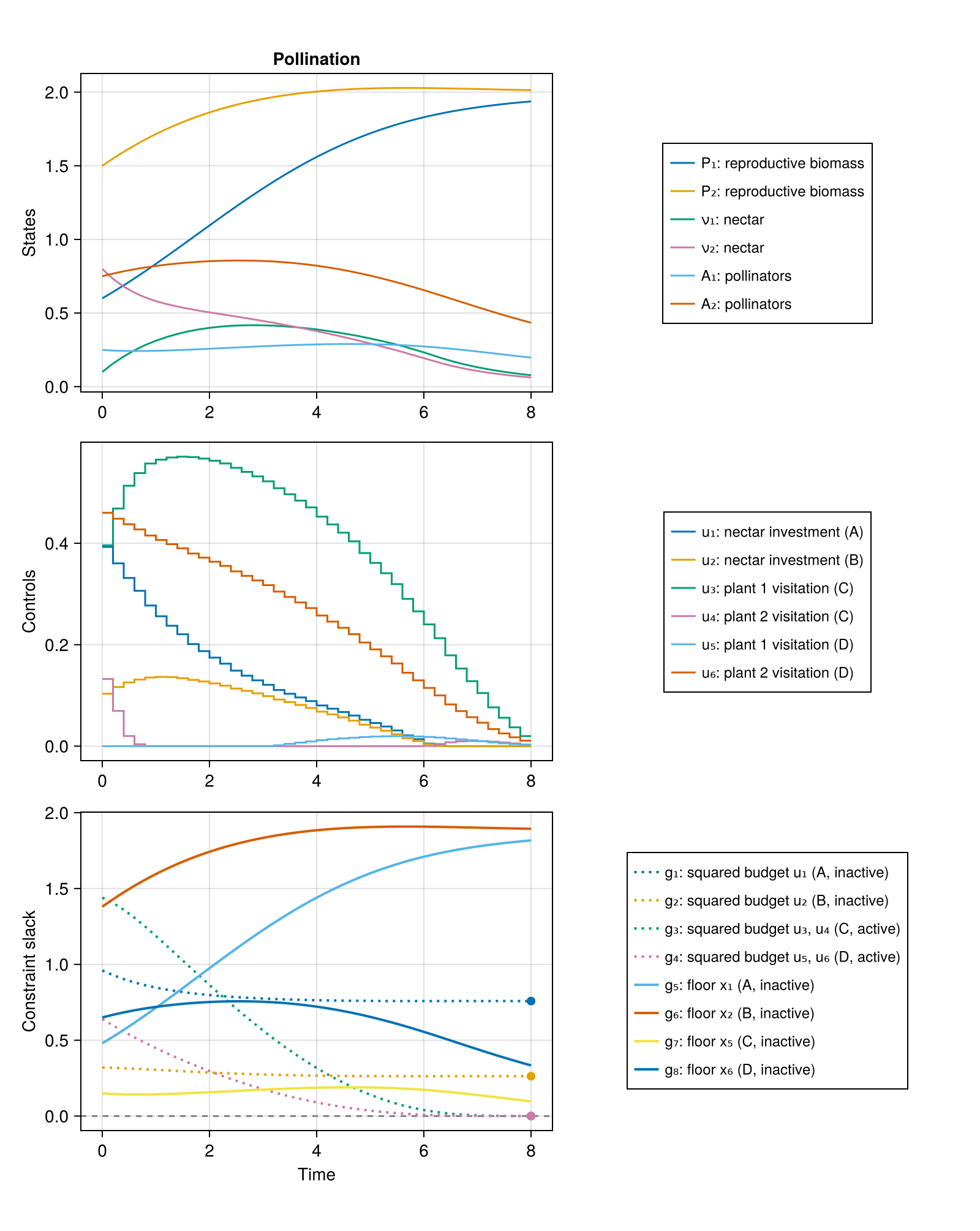}
\caption[Reference equilibrium candidate for \texttt{\caseidentifier}]{Reference equilibrium candidate for \texttt{\caseidentifier}. \textit{Interpretation.} Each pollinator specializes on the plant it handles most efficiently: C concentrates its visits on
plant A after briefly visiting plant B at the very start, and D visits plant B almost exclusively,
with only marginal visits to A late in the season. C's visitation rises to a maximum of about
$0.57$ near $t\approx1.5$ and declines afterwards, while D's decreases steadily from $0.46$; both
pollinators exhaust their seasonal effort budgets ($g_3,g_4$ active). Plant A, which starts small,
invests strongly in nectar early (from $0.39$) and then reduces investment, consistent with the
visit-driven gain in its reproductive biomass, while plant B, already large, invests little
throughout; both stop nectar production late in the season, consistent with the finite horizon:
nectar produced late can hardly raise the time-averaged reproductive biomass. Nectar of plant A accumulates to a maximum near
$t\approx2.8$ before being depleted by visits, whereas plant B's initial nectar pool is drawn down
from the start. Plant reproductive biomasses grow logistically, A more than threefold, and approach the
carrying capacity. Pollinator D's population first increases and then declines to $0.43$ as nectar
runs out, and C's population varies only moderately. The plant budgets and all population floors
remain slack.}\label{fig:\caseidentifier}
\end{figure}
\myheader{Literature and provenance.} The source links pollinator foraging preferences to competing plant population dynamics \citep{BenchmarkPMC4978411}. We introduce independently controlled nectar production and visitation in a multi-plant, multi-pollinator game with explicit nectar pools. Source adaptive-preference dynamics are replaced by open-loop controls; equivalence is not claimed. Plant--pollinator reward interactions have also been analyzed as a signalling game, in which the stability of honest nectar advertisement depends on signal costs and pollinator learning \citep{Sun2018PollinatorSignalling}; the present game contains no signal and treats nectar as a direct reward.

All numerical defaults are synthetic and dimensionless; no fitted biological values or source equilibrium results are asserted.
\myheader{Dimensions.} $N=4$, $n_x=6$, $n_{\rm ODE}=14$, $(n_{u,i})_{i\in\mathcal P}=(1,1,2,2)$, $(n_{g,i})_{i\in\mathcal P}=(4\times 2)$.
\myheader{Categories.} \cat{population}, \cat{seasonal}, \cat{mutualism}, \cat{variablePlayers}.
\myheader{Game specification.} Visits remove nectar at rate
\[
V_{ji}=\kappa_{ji}\pi_{ji}A_j\frac{\nu_i}{h+\nu_i}.
\]
\[
\dot P_i=P_i\left[r(1-P_i/K)+\frac{b\sum_{j\in\mathcal P_{\rm poll}}V_{ji}}{h+P_i}-c_\nu s_i^2\right],
\qquad \dot \nu_i=\alpha s_iP_i-\delta_\nu\nu_i-\sum_{j\in\mathcal P_{\rm poll}}V_{ji},
\]
\[
\dot A_j=e\sum_{i\in\mathcal P_{\rm plant}}V_{ji}
-\left[m+c_AA_j+c_u\sum_{i\in\mathcal P_{\rm plant}} a_{ji}^2\right]A_j.
\]
$\kappa=\texttt{visit}$ (pollinator rows); $h=\texttt{half}$;
$b,c_\nu,\alpha,\delta_\nu,e,m,c_A,c_u$ map to \texttt{benefit, nectar\_cost,
nectar\_rate, nectar\_loss, conversion, mortality, crowding, effort\_cost}.
Nectar uptake vanishes at a depleted nectar pool.

Each player minimizes
\[
J_i=\frac1T\int_0^T\sum_{k=1}^3w_i^k\ell_{ik}\,dt,
\qquad w_i^k\geq0,\quad\sum_kw_i^k=1.
\]
There is no terminal objective. Each weight row is independently configurable.
For plant $i$:
\[
(\ell_{i1},\ell_{i2},\ell_{i3})=(-P_i,\ s_i^2,\ \nu_i^2).
\]
For pollinator $j$, with $I_j=\sum_{i\in\mathcal P_{\rm plant}}V_{ji}$:
\[
(\ell_{j,1},\ell_{j,2},\ell_{j,3})=
\left(-A_j,\ \sum_{i\in\mathcal P_{\rm plant}} a_{ji}^2,\ \frac1{1+I_j}\right).
\]
Plants value reproductive biomass and penalize investment and surplus nectar;
pollinators value abundance and penalize effort and food shortfall.

Use $\dot Z_i=s_i^2$ for plants and
$\dot Z_{j}=\sum_{i\in\mathcal P_{\rm plant}} a_{ji}^2$ for pollinators; $Z(0)=0$ and
$Z_i(T)\leq\mathrm{budgets}_i$.
Population path floors apply to the corresponding plant or pollinator
(the nectar states have no floor). Default: four budgets and four path
population restrictions.

For $q$ biological states, append $N$ cost integrals then $N$ budget integrals,
all initially zero; the cost integrals include the factor $1/T$, the ODE dimension is $q+2N$, and $J_i=x_{q+\iota(i)}(T)$.

Enabled bounds:
\begin{center}\small\begin{tabular}{lll}\toprule Owners & Bound & Location\\\midrule
A & $y_{11}(T)\leq 0.96$ & terminal \\
B & $y_{12}(T)\leq 0.32$ & terminal \\
C & $y_{13}(T)\leq 1.44$ & terminal \\
D & $y_{14}(T)\leq 0.64$ & terminal \\
A & $0.12\leq y_{1}(t)$ & path \\
B & $0.12\leq y_{2}(t)$ & path \\
C & $0.1\leq y_{5}(t)$ & path \\
D & $0.1\leq y_{6}(t)$ & path \\
\bottomrule\end{tabular}\end{center}
\myheader{Parameter instance.} %
\begin{lstlisting}
plants = 2
pollinators = 2
T = 8.0
r = 0.5
K = 2.0
benefit = 0.4
nectar_cost = 0.04
nectar_rate = 1.5
nectar_loss = 0.5
half = 0.3
visit = [0.95 0.2; 0.23 1.0925]
conversion = 1.5
mortality = 0.2
crowding = 0.06
effort_cost = 0.04
x0 = [0.6, 1.5, 0.1, 0.8, 0.25, 0.75]
weights = [0.88 0.08 0.04; 0.65 0.3 0.05; 0.8 0.06 0.14; 0.5 0.3 0.2]
budgets = [0.96, 0.32, 1.44, 0.64]
floors = [0.12, 0.12, 0.1, 0.1]
initial = [0.2, 0.2, 0.2, 0.2, 0.2, 0.2]
control_lower = [0.0, 0.0, 0.0, 0.0, 0.0, 0.0]
control_upper = [1.0, 1.0, 1.0, 1.0, 1.0, 1.0]
numerical_x0 = [0.6, 1.5, 0.1, 0.8, 0.25, 0.75, 0.0, 0.0, 0.0, 0.0, 0.0, 0.0, 0.0, 0.0]
\end{lstlisting}

\endgroup

\begingroup
\def\caseidentifier{predatorPrey}%
\subsection{predatorPrey: Multiple prey/predator game}\label{case:predatorPrey}
\myheader{Biological interpretation.} Predator--prey communities in enriched environments, such as plankton in eutrophic lakes or
rodents and their predators, often oscillate on their own: prey increase until predators become
abundant, predators then deplete the prey and decline, and the cycle repeats (the paradox of
enrichment \citep{rosenzweig1971paradox}). Rapid prey evolution can shape such consumer--resource
cycles \citep{Yoshida2007Cryptic}. Within such cycles, prey can invest in defense or faster growth, and predators in
hunting effort or in reducing their own mortality. The cycle length is not imposed by an external
clock; it emerges from the interactions and can be modified by the players' strategies. We model
one full cycle of a community with three prey and two predator species. Prey species A is
assigned the choice of the cycle duration; this uses the same duration-ownership convention as \caseref{circadianCycle}.

There are $n_{\rm prey}\geq2$ prey and $n_{\rm pred}\geq2$ predators, each a player ($N=n_{\rm prey}+n_{\rm pred}$).
$X_i$: prey density; $Y_j$: predator density. State order is
$x=((X_i)_{i\in\mathcal P_{\rm prey}},(Y_j)_{j\in\mathcal P_{\rm pred}})$.
Prey controls $(d_i,g_i)$ are defense and growth investment;
predator controls $(h_j,v_j)$ are hunting and maintenance.
Controls are ordered by player, with two components each.
All controls depend on time alone. The player set $\mathcal P_{\rm prey}$ contains the first $n_{\rm prey}$ players, and $\mathcal P_{\rm pred}=\mathcal P\setminus\mathcal P_{\rm prey}$. In the archived instance, these are $\{A,B,C\}$ and $\{D,E\}$.
\begin{figure}[p]
\centering
\includegraphics[width=\textwidth,height=\casefigheight,keepaspectratio]{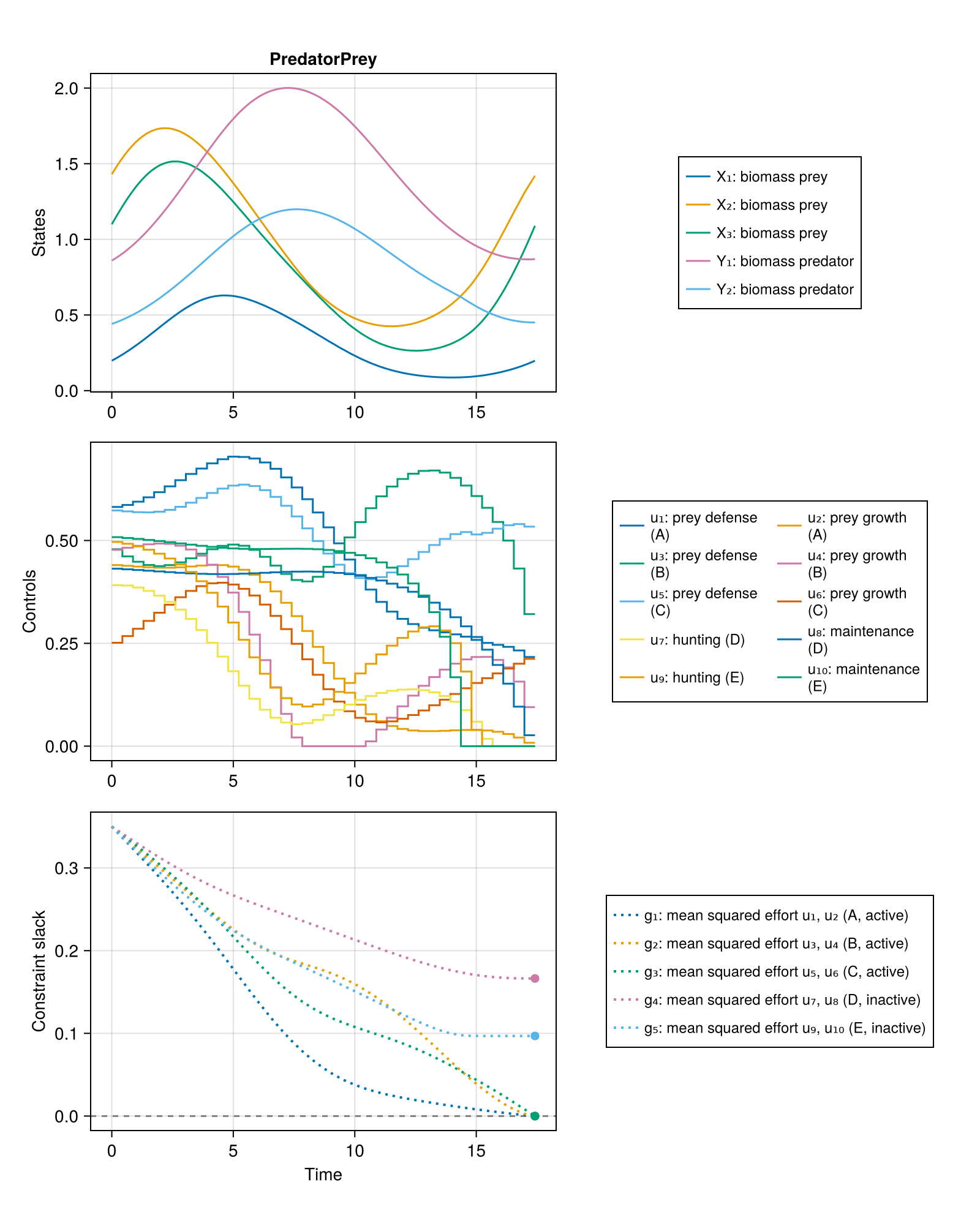}
\caption[Reference equilibrium candidate for \texttt{\caseidentifier}]{Reference equilibrium candidate for \texttt{\caseidentifier}. \textit{Interpretation.} Prey A lengthens the cycle to $T\approx17.4$, compared with the constant-control reference period
of $16.6$. The states show the classical phase shift of a trophic cycle: the prey peak early
(B and C near $t\approx2.2$--$2.6$, A near $t\approx4.6$), the predators follow and peak near
$t\approx7.2$--$7.6$, the prey collapse to minima between $t\approx11$ and $t\approx14$, and all
populations return to their initial levels at the end of the cycle. The prey keep substantial
defense throughout, raising it while predators are abundant (A up to $0.70$, B up to $0.67$
near $t\approx13$), and use growth investment mainly early, when predation pressure is low; B
stops growth investment in the middle of the cycle and resumes it slightly at the end. The
predators hunt and invest in maintenance while prey are available and reduce both efforts to zero
in the final third, when prey are scarce, consistent with the low return on hunting effort. The three prey exhaust their
mean-effort budgets ($g_1$--$g_3$ active), whereas the predators' budgets remain slack.}\label{fig:\caseidentifier}
\end{figure}
\myheader{Literature and provenance.} The source studies how predator-induced defenses alter predator--prey dynamics \citep{BenchmarkpredatorPrey}. We use controlled, saturating predator--prey interactions with logistic prey growth, and assign separate lineage objectives. The multispecies game and its control costs are our extension; source bifurcation results do not automatically carry over.

All numerical defaults are synthetic and dimensionless; no fitted biological values or source equilibrium results are asserted.
\myheader{Dimensions.} $N=5$, $n_x=5$, $n_{\rm ODE}=15$, $(n_{u,i})_{i\in\mathcal P}=(5\times 2)$, $(n_{g,i})_{i\in\mathcal P}=(5\times 11)$.
\myheader{Categories.} \cat{population}, \cat{periodical}, \cat{antagonism}, \cat{variablePlayers}, \cat{freeTime}.
\myheader{Game specification.} Set $L_{ji}=a_{ji}(1+h_j)/(1+d_i)$, the controlled attack coefficient, and the
Holling type-II saturation $S_j=1+\eta\sum_{k\in\mathcal P_{\rm prey}}L_{jk}X_k$ with handling time $\eta$.
Player A (prey 1) additionally chooses the cycle duration $T$; the dynamics are written on
normalized time $s\in[0,1]$ with physical time $t=T s$:
\[
\frac{dX_i}{ds}=T X_i\left[r_i(1-X_i/K_i)+b g_i
-c(d_i^2+g_i^2)-\sum_{j\in\mathcal P_{\rm pred}}\frac{L_{ji}Y_j}{S_j}\right],
\]
\[
\frac{dY_j}{ds}=T Y_j\left[\frac{\gamma}{S_j}\sum_{i\in\mathcal P_{\rm prey}}L_{ji}X_i-
\frac{\mu_j}{1+v_j}-kY_j-c(h_j^2+v_j^2)\right].
\]
$a_{ji}$ is the attack matrix (predator row, prey column); $\gamma$ is conversion,
$\mu_j$ mortality, and $k$ predator self-limitation;
$b=\texttt{growth\_benefit}$ and $c=\texttt{strategy\_cost}$.
With enriched prey (large $K_i$) and saturating predation, the constant reference strategies
produce a stable limit cycle (paradox of enrichment \citep{rosenzweig1971paradox}); there is no external clock.
The fixed $x_0$ is the point of this cycle where prey 1 crosses its cycle mean upward,
and the reference period $T_{\rm ref}$ initializes $T\in[T_{\rm ref}/2,2T_{\rm ref}]$.
At a zero population boundary its derivative vanishes.

Each player minimizes a cycle average,
\[
J_i=\int_0^1\sum_{k=1}^3w_i^k\ell_{ik}\,ds,
\qquad w_i^k\geq0,\quad\sum_kw_i^k=1,
\]
without terminal objective. Let $n_i$ denote player $i$'s population, $\bar n_i$ its mean over
the reference cycle, and $e_i$ the sum of its two squared controls:
\[
(\ell_{i1},\ell_{i2},\ell_{i3})=
\left(-n_i/\bar n_i,\ e_i,\ [n_i/\bar n_i-1]^2\right).
\]
These candidates reward abundance, penalize strategy effort, and penalize departures from the
reference abundance. Effort budgets bound the mean squared effort per unit time,
$\dot Z_i=e_i$, $Z_i(0)=0$, $Z_i(1)\leq\mathrm{budgets}_i$; optional path floors
$n_i(s)\geq\mathrm{floors}_i$ are disabled by default. Every player's problem contains the
return conditions $|n_j(1)-n_j(0)|\leq\delta$ for all populations.

For $q$ biological states, append $N$ cost integrals then $N$ budget integrals,
all initially zero; the ODE dimension is $q+2N$ and $J_i=x_{q+\iota(i)}(1)$.

Enabled bounds:
\begin{center}\small\begin{tabular}{lll}\toprule Owners & Bound & Location\\\midrule
A & $y_{11}(T)\leq 0.35$ & terminal \\
B & $y_{12}(T)\leq 0.35$ & terminal \\
C & $y_{13}(T)\leq 0.35$ & terminal \\
D & $y_{14}(T)\leq 0.35$ & terminal \\
E & $y_{15}(T)\leq 0.35$ & terminal \\
A, B, C, D, E & $|x_j(T)-x_j(0)|\leq 0.01$, $j=1,\ldots,5$ & terminal (return) \\
\bottomrule\end{tabular}\end{center}
\myheader{Parameter instance.} %
\begin{lstlisting}
prey = 3
predators = 2
periodicity_tolerance = 0.01
r = [0.9, 1.0, 1.1]
K = [3.0, 3.0, 3.0]
attack = [0.78 0.54 0.66; 0.54 0.66 0.78]
handling = 0.5
conversion = 0.45
mortality = [0.3, 0.35]
crowding = 0.03
growth_benefit = 0.25
strategy_cost = 0.12
default_size = true
x0 = [0.1981184964, 1.43031613, 1.099886694, 0.8595951089, 0.4409164492]
scale = [0.1988327166, 0.9306449029, 0.7068407768, 1.259519777, 0.6566602541]
duration = 16.58214832
duration_bounds = (8.291074161, 33.16429664)
weights = [0.65 0.25 0.1; 0.65 0.25 0.1; 0.65 0.25 0.1; 0.75 0.2 0.05; 0.75 0.2 0.05]
budgets = [0.35, 0.35, 0.35, 0.35, 0.35]
floors = [nothing, nothing, nothing, nothing, nothing]
initial = [0.2, 0.2, 0.2, 0.2, 0.2, 0.2, 0.2, 0.2, 0.2, 0.2]
control_lower = [0.0, 0.0, 0.0, 0.0, 0.0, 0.0, 0.0, 0.0, 0.0, 0.0]
control_upper = [1.0, 1.0, 1.0, 1.0, 1.0, 1.0, 1.0, 1.0, 1.0, 1.0]
numerical_x0 = [0.1981184964, 1.43031613, 1.099886694, 0.8595951089, 0.4409164492, 0.0, 0.0, 0.0, 0.0, 0.0, 0.0, 0.0, 0.0, 0.0, 0.0]
\end{lstlisting}

\endgroup

\begingroup
\def\caseidentifier{proteinFolding}%
\subsection{proteinFolding: Chaperone assistance with finite effort}\label{case:proteinFolding}
\myheader{Biological interpretation.} When cells experience proteotoxic stress, for example a heat shock, many proteins unfold at once.
Molecular chaperones bind unfolded proteins and folding intermediates and assist their folding to
the native state, but chaperone capacity is itself limited: it is consumed while chaperones are
engaged and restored only slowly by synthesis. Two chaperone systems with different efficiencies
(for example, an ATP-dependent chaperonin and a cheaper holdase) share this capacity. We model the
recovery episode after a stress pulse, in which each system chooses its assistance effort under a
finite effort budget and aims to restore the native protein pool. Because both systems value the
same native pool, the game isolates free-riding on assistance effort.
\begin{figure}[p]
\centering
\includegraphics[width=\textwidth,height=\casefigheight,keepaspectratio]{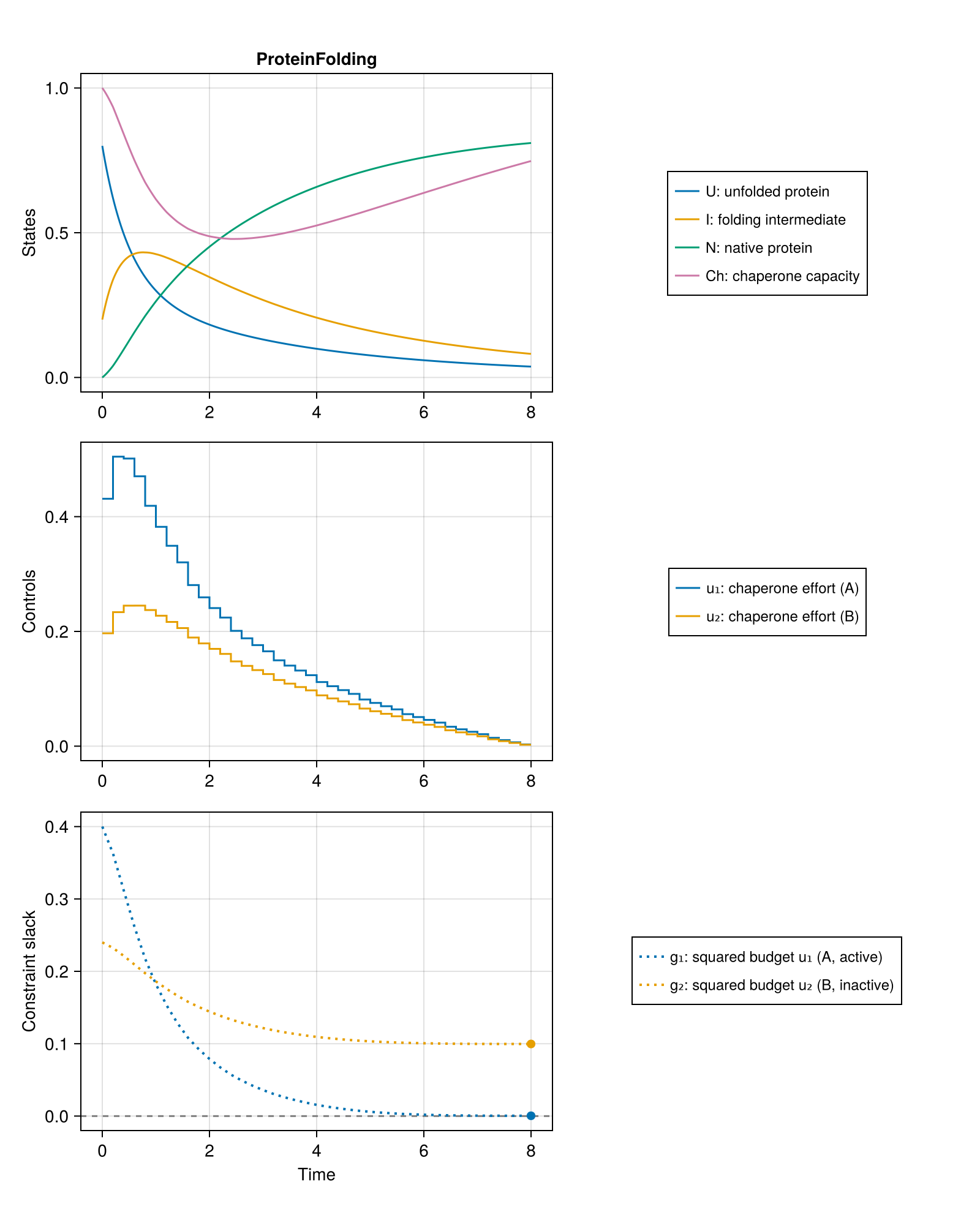}
\caption[Reference equilibrium candidate for \texttt{\caseidentifier}]{Reference equilibrium candidate for \texttt{\caseidentifier}. \textit{Interpretation.} After the stress pulse, most protein is unfolded ($0.8$) and chaperone capacity is full. Both
chaperone systems respond immediately: the efficient system A starts near $0.43$, peaks at about
$0.5$, and then reduces its effort steadily, whereas the less efficient system B contributes about
half as much and tapers in parallel. Assisted folding consumes chaperone capacity faster than it is
resynthesized, so capacity falls to about half its initial value near $t\approx2.4$ and recovers
slowly afterwards as the need for assistance declines. Unfolded protein is converted quickly into
folding intermediates, which peak near $t\approx0.75$ at more than twice their initial level, and
the native pool rises sigmoidally to about $0.81$ by the end of the episode. System A exhausts its
effort budget ($g_1$ active), while B stays below its budget; the transient depletion of the shared
chaperone capacity limits how much the two systems can accelerate folding at the peak of demand.}\label{fig:\caseidentifier}
\end{figure}
\myheader{Literature and provenance.} A kinetic model of bacterial proteostasis shows that chaperone-assisted folding trades energy consumption against chaperone utilization \citep{BenchmarkPMC5380058}; this supports the finite, consumable chaperone capacity and the effort costs, not the division of assistance into competing players. Two chaperone modules have different effort preferences; the game is a hypothetical decomposition of cellular control, not a published Nash model.
All rates, budgets, weights, and initial values below are synthetic,
dimensionless benchmark choices, not fitted values from this source.
\myheader{Dimensions.} $N=2$, $n_x=4$, $n_{\rm ODE}=8$, $(n_{u,i})_{i\in\mathcal P}=(2\times 1)$, $(n_{g,i})_{i\in\mathcal P}=(2\times 1)$.
\myheader{Categories.} \cat{molecular}, \cat{episodical}, \cat{regulation}.
\myheader{Game specification.} $U,I,N$: unfolded, intermediate, native protein; $Ch$: chaperone capacity. $u_i$: module effort. $v=(c_1u_A+c_2u_B)\,Ch\,I$ is assisted conversion, with efficiencies $(c_1,c_2)=$ \code{c}; since $c_1>c_2$, system A is the more efficient one.
\[\begin{gathered}
\dot U=-k_fU+k_uI-d_UU,\quad
\dot I=k_fU-(k_u+k_n)I-v,\\
\dot N=k_nI-d_NN+v,\quad
\dot{Ch}=s_0-c_a(u_A+u_B)\,Ch\,I-d_a\,Ch.
\end{gathered}\]
Controls are open-loop functions of time.
Initial-state order: $x=(U,I,N,Ch)$; values appear below. The supply rate is $s_0=$ \code{A0}.

Each player $i\in\mathcal P$ minimizes $J_i=T^{-1}\int_0^T\ell_i\,dt$, with
\[\ell_i=w_i^1(N-1)^2+w_i^2u_i^2.\]
$w_i^k\geq0$, $\sum_kw_i^k=1$; weights are independently configurable.
\[
\dot Z_i=u_i^2,\quad Z_i(0)=0,\quad Z_i(T)\leq \mathcal B_i.
\]
Budgets are private to each player. The optional path floor is $N\geq\mathrm{floor}$.
When enabled, it is enforced in every player's deviation problem.

Enabled bounds:
\begin{center}\small\begin{tabular}{lll}\toprule Owners & Bound & Location\\\midrule
A & $y_{7}(T)\leq 0.4$ & terminal \\
B & $y_{8}(T)\leq 0.24$ & terminal \\
\bottomrule\end{tabular}\end{center}
\myheader{Parameter instance.} %
\begin{lstlisting}
initial = (0.2, 0.3)
T = 8.0
kf = 1.5
ku = 0.6
kn = 0.3
du = 0.02
dn = 0.01
c = (0.8, 0.5)
ca = 2.0
da = 0.2
A0 = 0.2
x0 = (0.8, 0.2, 0.0, 1.0)
w_A = (0.85, 0.15)
w_B = (0.75, 0.25)
budgets = (0.4, 0.24)
floor = nothing
control_lower = [0.0, 0.0]
control_upper = [1.0, 1.0]
numerical_x0 = [0.8, 0.2, 0.0, 1.0, 0.0, 0.0, 0.0, 0.0]
\end{lstlisting}

\endgroup

\begingroup
\def\caseidentifier{seedGermination}%
\subsection{seedGermination: Seed-cohort germination under shared water}\label{case:seedGermination}
\myheader{Biological interpretation.} Seeds in the same soil patch compete for water during imbibition and germination in the
germination window after the first autumn rains. Rapid activation of the embryo lets a seed
emerge first and secure light and space, but it consumes water that neighboring seeds also need,
and activation itself is metabolically costly. Once the first seed has completed radicle
emergence, the competitive situation changes fundamentally. We model this seasonal window for
three seed cohorts sharing one water pool; it ends when the leading cohort reaches radicle
emergence. Cohort A owns the duration, so only its time weight is strategically relevant
(Section~\ref{sec:owned-duration}).
\begin{figure}[p]
\centering
\includegraphics[width=\textwidth,height=\casefigheight,keepaspectratio]{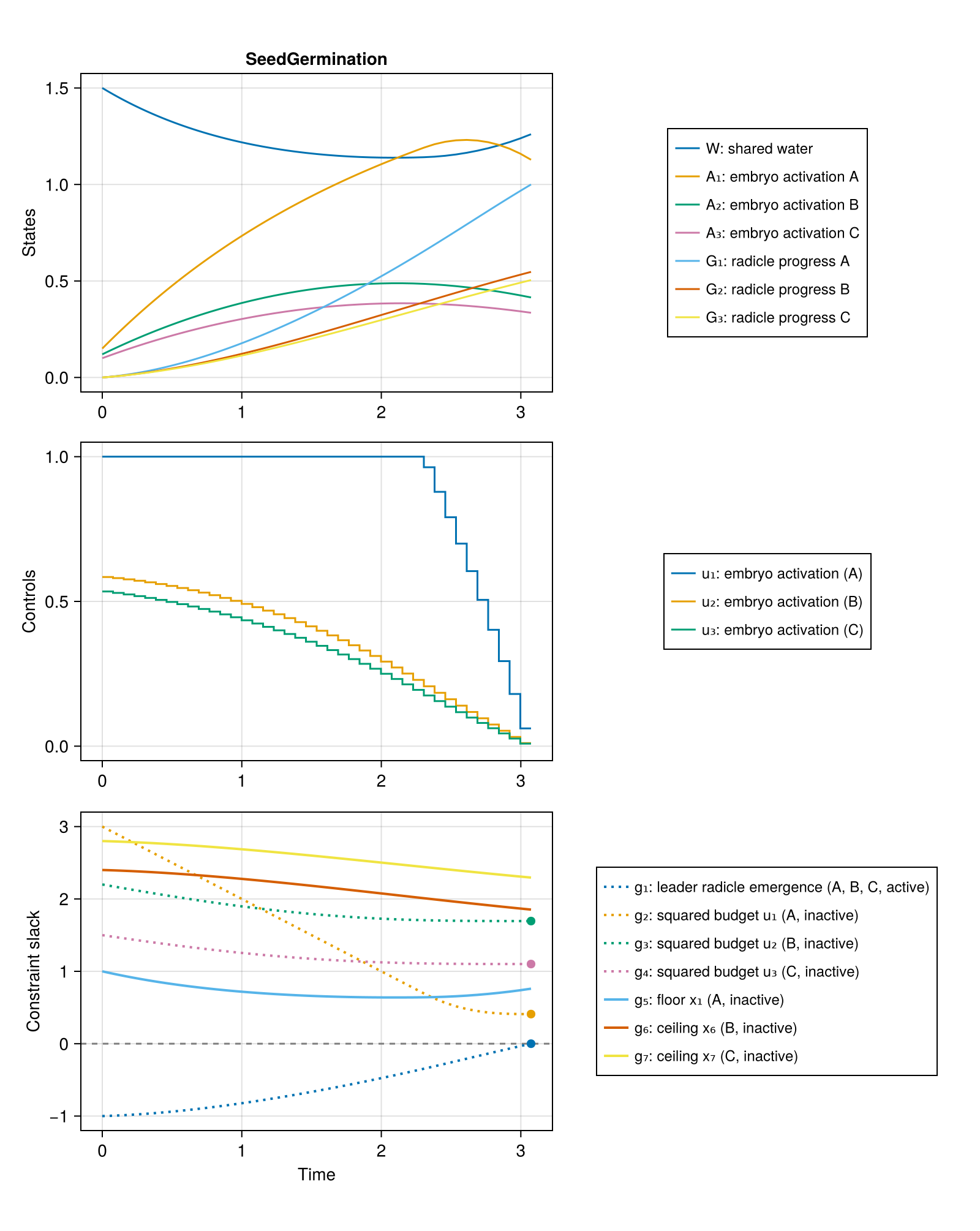}
\caption[Reference equilibrium candidate for \texttt{\caseidentifier}]{Reference equilibrium candidate for \texttt{\caseidentifier}. \textit{Interpretation.} The leading cohort A activates its embryos at the maximal rate for three quarters of the episode
and relaxes only shortly before its radicle progress reaches the emergence threshold at
$T\approx3.1$ (the terminal equality $g_1$ is active for all players). Cohorts B and C, which
weigh time less and effort and water burden more, start with moderate activation ($0.58$ and
$0.53$) and reduce it steadily. Embryo activation of A rises to about $1.2$ and saturates, while
B and C peak near $0.49$ and $0.39$ at $t\approx2.2$ and decline afterwards. Radicle progress
accumulates in a convex fashion and reaches $1$ for A but only about $0.55$ and $0.5$ for B and
C, so A emerges well ahead of its competitors. The shared water pool falls from $1.5$ to about
$1.14$ while activation is high and recovers once the cohorts relax. The effort budgets, the
water floor, and the radicle-progress ceilings of B and C remain slack.}\label{fig:\caseidentifier}
\end{figure}
\myheader{Literature and provenance.} A, B, C: three seed cohorts controlling activation effort \citep{BenchmarkseedGermination}. Cohort A controls the emergence schedule.

$x=(W,A_A,A_B,A_C,G_A,G_B,G_C)$: shared water, activation levels, and cumulative radicle development.
The hydrotime analysis of \citet{BenchmarkseedGermination} supports only the mechanism that germination rate increases with water availability above a threshold. Germination timing under competition among individuals has been analyzed as an evolutionary game, with competition acting through growth, mortality, or recruitment \citep{BesshoIwasa2012Germination}; our game resolves one germination window and does not predict evolutionarily stable dates. The activation states and interacting-cohort game are synthetic and are not the hydrotime equations of \citet{BenchmarkseedGermination}.
\myheader{Dimensions.} $N=3$, $n_x=7$, $n_{\rm ODE}=13$, $(n_{u,i})_{i\in\mathcal P}=(3\times 1)$, $(n_{g,i})_{i\in\mathcal P}=(3\times 3)$.
\myheader{Categories.} \cat{population}, \cat{seasonal}, \cat{competition}, \cat{freeTime}.
\myheader{Game specification.} \[\begin{gathered}
\dot W=q(W_{\max}-W)-cW\sum_{j\in\mathcal P} u_j,\\
\dot A_i=a_i u_i W/(K+W)-d_iA_i,\\
\dot G_i=g_i A_iW/(K+W),\qquad i\in\mathcal P.
\end{gathered}\]

$q=$ \texttt{water\_supply}, $W_{\max}=$ \texttt{water\_capacity}, $c=$ \texttt{water\_use}, $a=$ \texttt{activation}, $d=$ \texttt{decay}, $g=$ \texttt{growth}, $K=$ \texttt{saturation}.

$J_i=w_i^1T+w_i^2\int_0^T u_i^2\,dt+w_i^3\int_0^T L_i\,dt$.

$L_i=(G_i-\bar G)^2$, $\bar G=\sum_{j\in\mathcal P}G_j/N$. The duration term penalizes episode length, effort penalizes activation, and development mismatch penalizes loss of cohort synchrony. The synchrony term represents benefits of emerging with the cohort, such as predator satiation; it conflicts with the advantage of emerging first.

$G_A(0)=0$, $1-G_A(T)\geq0$ and $G_A(T)-1\geq0$, shared by all players. Private effort budgets are $(3,2.2,1.5)$. A imposes $W\geq0.5$ (a minimum soil moisture for embryo survival), and B and C impose $G_B\leq2.4$ and $G_C\leq2.8$ (caps on radicle elongation); the ceilings cannot bind before $G_A=1$ ends the episode. This is leader-cohort emergence, not simultaneous emergence.

Enabled bounds:
\begin{center}\small\begin{tabular}{lll}\toprule Owners & Bound & Location\\\midrule
A, B, C & $y_{5}(T)=1$ & terminal \\
A & $y_{8}(T)\leq 3$ & terminal \\
B & $y_{9}(T)\leq 2.2$ & terminal \\
C & $y_{10}(T)\leq 1.5$ & terminal \\
A & $0.5\leq y_{1}(t)$ & path \\
B & $y_{6}(t)\leq 2.4$ & path \\
C & $y_{7}(t)\leq 2.8$ & path \\
\bottomrule\end{tabular}\end{center}
\myheader{Parameter instance.} %
\begin{lstlisting}
x0 = (1.5, 0.15, 0.12, 0.1, 0.0, 0.0, 0.0)
target = 1.0
water_supply = 0.4
water_capacity = 2.0
water_use = 0.2
activation = (1.0, 0.85, 0.7)
decay = (0.4, 0.35, 0.3)
growth = (0.5, 0.6, 0.7)
saturation = 0.4
duration = 5.0
duration_bounds = (0.1, 20.0)
initial = (0.5, 0.45, 0.4)
weights = [0.75 0.15 0.1; 0.55 0.3 0.15; 0.45 0.25 0.3]
budgets = (3.0, 2.2, 1.5)
path_limits = (0.5, 2.4, 2.8)
control_lower = [0.0, 0.0, 0.0]
control_upper = [1.0, 1.0, 1.0]
numerical_x0 = [1.5, 0.15, 0.12, 0.1, 0.0, 0.0, 0.0, 0.0, 0.0, 0.0, 0.0, 0.0, 0.0]
\end{lstlisting}

\endgroup

\begingroup
\def\caseidentifier{sexRatio}%
\subsection{sexRatio: Sex allocation among foundresses}\label{case:sexRatio}
\myheader{Biological interpretation.} When several foundresses lay eggs in the same patch, for example parasitoid wasps in one host
cluster or fig wasps in one fig, their sons compete for mates with each other before the
daughters disperse (local mate competition). Hamilton's classical argument
\citep{Hamilton1967SexRatios} predicts
female-biased sex ratios that become less biased as more foundresses share the patch. Here the
foundresses provision their brood from a shared resource that is replenished seasonally, and
producing sons is more resource-intensive. Each foundress chooses the fraction of its offspring
that are male over one reproductive season; offspring success depends on the numbers of sons and
daughters of all foundresses, and the shared brood resource must not fall below a floor. This
floor is a synthetic persistence restriction; noncooperative foundresses would not impose it on
themselves. The controls represent male-allocation fractions rather than discrete offspring
choices.
\begin{figure}[p]
\centering
\includegraphics[width=\textwidth,height=\casefigheight,keepaspectratio]{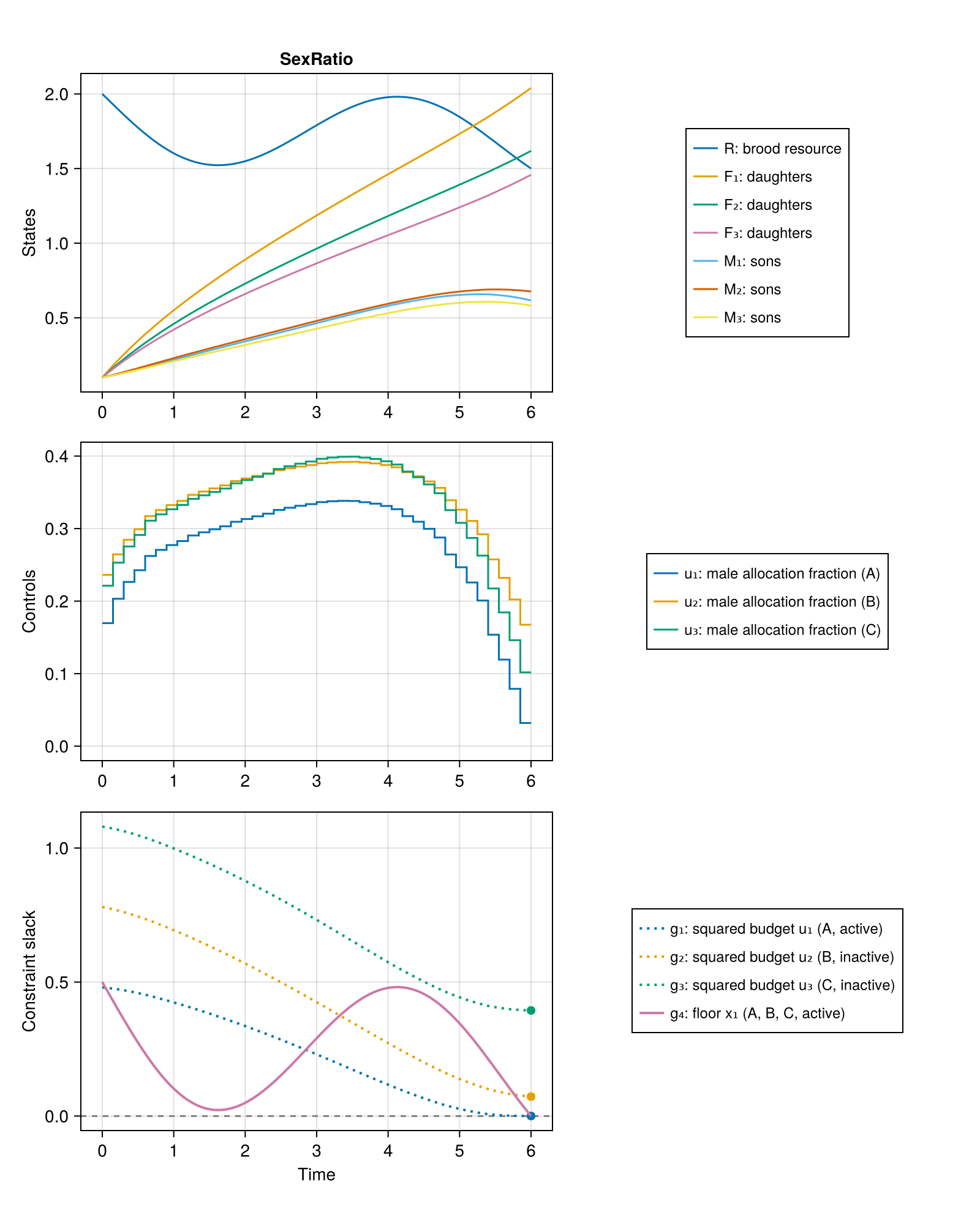}
\caption[Reference equilibrium candidate for \texttt{\caseidentifier}]{Reference equilibrium candidate for \texttt{\caseidentifier}. \textit{Interpretation.} All foundresses produce female-biased broods throughout the season: male allocation rises from
about $0.2$ early to a maximum near $0.34$--$0.40$ in the middle of the season and declines again
toward the end, when newly produced sons have little time left to mate. Foundress A, whose squared-allocation
budget is exhausted ($g_1$ active) and which also weighs mating success most, is the most
female-biased, especially late in the season (final allocation about $0.03$). The shared brood
resource follows the seasonal supply: it falls from $2$ to a minimum near $1.5$ at
$t\approx1.6$, recovers to about $2$ at $t\approx4$ when supply peaks, and is drawn down again to
the protective floor at the end of the season ($g_4$ active). Daughters accumulate almost linearly
(to $2.0$, $1.6$, and $1.5$ for A, B, and C), whereas sons level off at about $0.6$--$0.7$ as
male production declines. The effort budgets of B and C remain slack.
For comparison, Hamilton's unbeatable sex ratio under local mate competition with three
foundresses is $(N-1)/(2N)=1/3$ \citep{Hamilton1967SexRatios}. The seasonal allocations
($\approx0.2$--$0.4$) bracket this value, but the present objective is a within-season
mating proxy, the male-cost term and the squared-allocation budgets also bias allocation
toward daughters, so agreement is not a test of the classical result.}\label{fig:\caseidentifier}
\end{figure}
\myheader{Literature and provenance.} Local mate competition and mate limitation support multiple foundresses and frequency-dependent mating returns \citep{BenchmarkPMC5830725}. The following finite-horizon proxy is not a reproduction of the published evolutionary stability analysis.
All rates, budgets, weights, and initial values below are synthetic,
dimensionless benchmark choices, not fitted values from this source.
\myheader{Dimensions.} $N=3$, $n_x=7$, $n_{\rm ODE}=13$, $(n_{u,i})_{i\in\mathcal P}=(3\times 1)$, $(n_{g,i})_{i\in\mathcal P}=(3\times 2)$.
\myheader{Categories.} \cat{population}, \cat{seasonal}, \cat{competition}, \cat{variablePlayers}.
\myheader{Game specification.} $R$: brood resource; $F_i,M_i$: surviving daughters and sons. $u_i$: male fraction chosen by foundress $i$. Default $N=3$; $F_\Sigma=\sum_{j\in\mathcal P}F_j$, $M_\Sigma=\sum_{j\in\mathcal P}M_j$.
\[\begin{gathered}
b_i^*=b_iR/(K+R),\quad
\dot R=s(t)-\mathrm{decay}R-\mathrm{cost}\sum_{j\in\mathcal P}b_j^*(1+c_m u_j),\\
s(t)=\mathrm{supply}[1-a_s\cos(2\pi t/T)],\quad
a_s=\mathrm{seasonality},\ c_m=\mathrm{male\_cost},\\
\dot F_i=b_i^*(1-u_i)-d_fF_i,\qquad\dot M_i=b_i^*u_i-d_mM_i,\\
W_i=\frac{F_iM_\Sigma+M_iF_\Sigma}{\mathrm{mate\_half}+M_\Sigma}.
\end{gathered}\]
Controls are open-loop functions of time.
Initial-state order: $x=(R,(F_i)_{i\in\mathcal P},(M_i)_{i\in\mathcal P})$; values appear below.

Each player $i\in\mathcal P$ minimizes $J_i=T^{-1}\int_0^T\ell_i\,dt$, with
\[\ell_i=-w_i^1W_i+w_i^2(u_i-0.5)^2+w_i^3u_i^2.\]
The penalty $(u_i-0.5)^2$, the term $u_i^2$, and the budget on $\int_0^Tu_i^2\,dt$ are synthetic:
the squared terms represent increasing marginal costs of producing sons (e.g., larger
provisioning per son), in addition to the male cost $c_m$ in $\dot R$. Fitness is usually
measured at the end of the season, as mated daughters at dispersal; here it is the time
average of the mating proxy $W_i$.
$w_i^k\geq0$, $\sum_kw_i^k=1$; weights are independently configurable.
\[
\dot Z_i=u_i^2,\quad Z_i(0)=0,\quad Z_i(T)\leq \mathcal B_i.
\]
Budgets are private to each player. The default path floor is $R\geq1.5$.
Seasonal replenishment and sex-dependent provisioning allow drawdown and
recovery. The shared floor restricts feasible allocations jointly.

Enabled bounds:
\begin{center}\small\begin{tabular}{lll}\toprule Owners & Bound & Location\\\midrule
A & $y_{11}(T)\leq 0.48$ & terminal \\
B & $y_{12}(T)\leq 0.78$ & terminal \\
C & $y_{13}(T)\leq 1.08$ & terminal \\
A, B, C & $1.5\leq y_{1}(t)$ & path \\
\bottomrule\end{tabular}\end{center}
\myheader{Parameter instance.} %
\begin{lstlisting}
N = 3
T = 6.0
b = [1.0, 0.8695652174, 0.7692307692]
df = 0.12
dm = 0.2
supply = 0.5
seasonality = 0.8
male_cost = 2.0
decay = 0.1
cost = 0.15
K = 1.0
mate_half = 0.5
x0 = [2.0, 0.1, 0.1, 0.1, 0.1, 0.1, 0.1]
w_A = (0.8, 0.1, 0.1)
w_B = (0.65, 0.2, 0.15)
w_C = (0.7, 0.05, 0.25)
budgets = [0.48, 0.78, 1.08]
floor = 1.5
initial = [0.2, 0.2, 0.2]
control_lower = [0.0, 0.0, 0.0]
control_upper = [1.0, 1.0, 1.0]
numerical_x0 = [2.0, 0.1, 0.1, 0.1, 0.1, 0.1, 0.1, 0.0, 0.0, 0.0, 0.0, 0.0, 0.0]
\end{lstlisting}

\endgroup

\begingroup
\def\caseidentifier{tupelo}%
\subsection{tupelo: Tupelo--leafminer case study}\label{case:tupelo}
\myheader{Biological interpretation.} Leaf-mining insects lay their eggs on young leaves in spring; the larvae then feed inside the
leaf blade until they mature or the leaf's resources run out. For the host tree, chemical
defense reduces egg and larval survival but diverts carbon from growth, and the value of leaf
tissue declines over the season as leaves age. For the herbivore, the timing of egg hatching
determines when larvae meet the host's defense and how much fresh leaf tissue remains. This
"cold war" over timing within one season is a classical plant--herbivore scenario; the end of the
season is linked to the next one through overwintering, not through periodic return.

The host tree \emph{Nyssa sylvatica} and the leafminer
\emph{Antispila nysaefoliella} motivate a seasonal game. 
\cite{Low2013ColdWar} already studied optimal plant defense and herbivore
hatching in this system, including a co-evolutionary calculation.

The plant (tupelo) chooses a defense investment schedule $u_A(t)$. In \citet{Low2013ColdWar} defense is induced; here it is an open-loop schedule, i.e., constitutive within the season.
The herbivore (leafminer) chooses hatching hazard $u_B(t)$.

The biological states are usable leaf area $L$, unhatched eggs $E$,
and feeding larvae $H$ per representative leaf, $x=(L,E,H)$.
The hatching flux is $u_BE$.

As in \citet{Low2013ColdWar}, larvae consume leaf tissue at a constant per-capita rate and leave the
cohort through pupation or death. The alternative of Section~\ref{sec:tupelo-variants} instead uses a
saturating tissue availability $\rho(L)=L/(K+L)$ in consumption and maturation.

Defense increases egg and larval mortality and thus reduces the cohort and its future feeding pressure. 
Hatching transfers eggs to larvae, it determines when larvae encounter defense and declining seasonal resources.
\begin{figure}[p]
\centering
\includegraphics[width=\textwidth,height=\casefigheight,keepaspectratio]{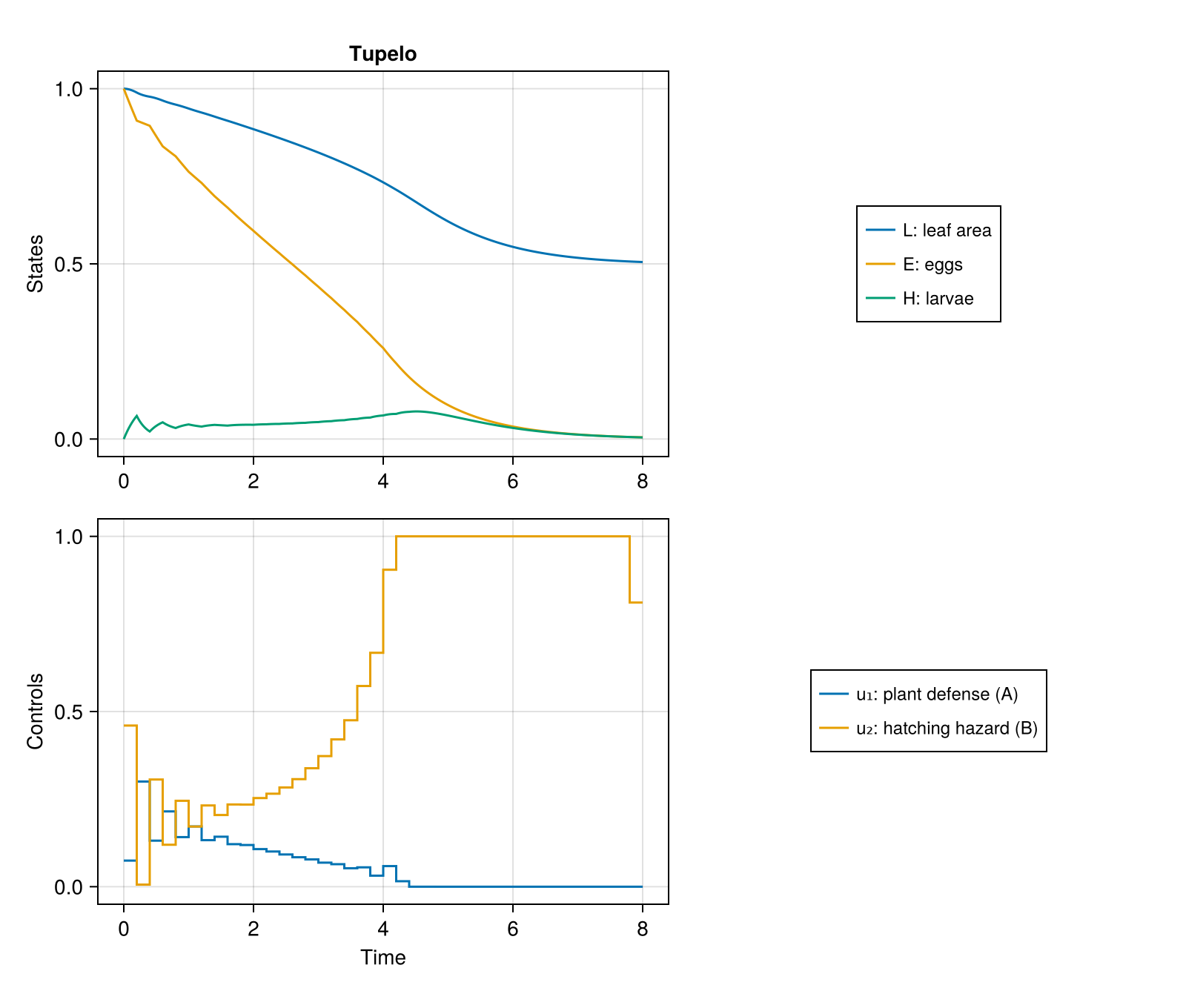}
\caption[Reference equilibrium candidate for \texttt{\caseidentifier}]{Reference equilibrium candidate for \texttt{\caseidentifier}. \textit{Interpretation.} The tree defends early with at most $0.30$, mostly between $0.1$ and $0.2$, and stops defending at
$t\approx4.4$, before the leaf value declines. The leafminer releases a small hatching pulse at the start,
hatches at a reduced hazard of about $0.2$--$0.3$ while defense lasts, and hatches at maximal hazard from
$t\approx4.2$ until shortly before $T$. Eggs decline to $0.26$ by $t=4$ and are nearly exhausted by $t=6$;
larvae stay below $0.08$ and peak at $t\approx4.5$. Usable leaf area decreases to about $0.50$, and
46\% of the loss occurs after $t=4$. The pattern agrees qualitatively with the delayed mass hatching of
\citet{Low2013ColdWar}; unlike in their feedback model, the leafminer does not avoid defense completely
(Section~\ref{sec:tupelo-reference}).}\label{fig:\caseidentifier}
\end{figure}
\myheader{Literature and provenance.} \citet{Low2013ColdWar} supply the model structure (their Eqs.~1--7) and the rates of their coevolution scenario~11 (their Table~2), mapped here to one time unit of 25 days. Open-loop schedules replace their induced defense, a Nash candidate replaces their evolutionarily stable strategies, and small quadratic cost terms regularize the responses. The photoperiod scale $t^*$ and the season length are our assumptions, since the source does not state them.
\myheader{Dimensions.} $N=2$, $n_x=3$, $n_{\rm ODE}=8$, $(n_{u,i})_{i\in\mathcal P}=(2\times 1)$, $(n_{g,i})_{i\in\mathcal P}=(2\times 0)$.
\myheader{Categories.} \cat{population}, \cat{seasonal}, \cat{antagonism}.
\myheader{Game specification.} We use
\[
f(t,x,u_A,u_B;p)=\begin{pmatrix}
-b\rho(L)H\\
-(u_B+\mu_E+d_Eu_A)E\\
u_B E-(\mu_H+d_Hu_A+\kappa\rho(L))H
\end{pmatrix},\qquad x(0)=\begin{pmatrix}L_0\\E_0\\H_0\end{pmatrix},
\]
with $p=(b,K,\mu_E,\mu_H,d_E,d_H,\kappa,\delta_A,\delta_B,t^*,z)$ and $\rho(L)=1$ in the reference instance
(\code{saturation=false}); \code{saturation=true} gives $\rho(L)=L/(K+L)$.

$b$ is feeding intensity and $K$ is tissue half-saturation.
$\mu_E,\mu_H$ are baseline mortality rates and $d_E,d_H$ are defense-induced mortality coefficients.
The maturation flux $\kappa\rho(L)H$ contributes to player B's objective.
The rates map scenario~11 of \citet{Low2013ColdWar} with one time unit of 25 days (Section~\ref{sec:tupelo-setup}).
In the parameter listing, $\delta_A=$ \code{p_decay} and $\delta_B=$ \code{q_decay}.
The reference weights are $w_A=(0.49,0.49,0.02,0)$ and $w_B=(0.98,0.02,0)$.

Implementation comment: nonnegativity, decreasing leaf area, and decreasing cohort size are checked
by independent integration; they are not extra NLP constraints.

\renewcommand{\arraystretch}{1.35}
\begin{tabular}{@{}lll@{}}
Objective & Weighted term & Interpretation\\\hline
$J_A$ & $-w_A^1\int_0^T p(t)L\,\mathrm dt$ & Maximize value-weighted leaf area\\
& $+w_A^2\int_0^T u_AL\,\mathrm dt$ & Minimize linear defense cost\\
& $+w_A^3\int_0^T u_A^2L\,\mathrm dt$ & Minimize quadratic defense cost\\
& $-w_A^4L(T)$ & Maximize terminal leaf tissue\\\hline
$J_B$ & $-w_B^1\int_0^T e^{-\delta_B t}\kappa\rho(L)H\,\mathrm dt$ & Maximize discounted maturation\\
& $+w_B^2\int_0^T u_B^2E\,\mathrm dt$ & Minimize quadratic hatching cost\\
& $-w_B^3H(T)$ & Reward larvae maturing after $T$\\
\end{tabular}

The leaf value is the photoperiod curve $p(t)=e^{-(t/t^*)^z}$ of \citet{Low2013ColdWar} (\code{photoperiod=true}),
with $t^*=6$ as our assumption; \code{photoperiod=false} gives $p(t)=e^{-\delta_At}$.
The terminal larva term counts larvae that complete development after the modeled window; its weight is zero in the reference instance; Section~\ref{sec:tupelo-variants} mentions a positive weight as an example of a preference variant.

Collect restriction switches (a zero switch omits a particular restriction) and bounds in
\[
\theta=(\sigma_{U_A},\sigma_{U_B},\sigma_{C_A},U_A,U_B,C_A),
\qquad \sigma_{U_A},\sigma_{U_B},\sigma_{C_A}\in\{0,1\}.
\]
For $g_A[x,u_A,u_B;p,\theta]\geq0$, impose
\[
\begin{cases}
U_A-\int_0^T u_A\,\mathrm dt\geq0 &\text{if }\sigma_{U_A}=1,\\
C_A-\int_0^T u_AL\,\mathrm dt\geq0 &\text{if }\sigma_{C_A}=1.
\end{cases}
\]
Note that the carbon allowance $C_A-\int_0^T u_AL\,\mathrm dt\geq0$ depends on $u_B$ and $p$ through $L$.

For $g_B[x,u_A,u_B;p,\theta]\geq0$, impose
\[
U_B-\int_0^T u_B\,\mathrm dt\geq0\quad\text{if }\sigma_{U_B}=1.
\]
All switches are zero in the reference instance, as in \citet{Low2013ColdWar}; the alternative of Section~\ref{sec:tupelo-variants} sets them to one.

No state or budget bounds are enabled.
\myheader{Parameter instance.} %
\begin{lstlisting}
T = 8.0
L0 = 1.0
E0 = 1.0
H0 = 0.0
b = 1.5
K = 0.15
muE = 0.0
muH = 1.0
dE = 0.25
dH = 12.5
kappa = 1.0
p_decay = 0.2
q_decay = 0.125
t_star = 6.0
z = 9.0
umax = 1.0
vmax = 1.0
w_A = (0.49, 0.49, 0.02, 0.0)
w_B = (0.98, 0.02, 0.0)
scales_P = (1.0, 1.0, 1.0, 1.0)
scales_H = (1.0, 1.0, 1.0)
U_A = nothing
U_B = nothing
C_A = nothing
saturation = false
photoperiod = true
control_lower = [0.0, 0.0]
control_upper = [1.0, 1.0]
numerical_x0 = [1.0, 1.0, 0.0, 0.0, 0.0, 0.0, 0.0, 0.0]
\end{lstlisting}

\endgroup

\begingroup
\def\caseidentifier{woundRepair}%
\subsection{woundRepair: Wound repair and remodeling}\label{case:woundRepair}
\myheader{Biological interpretation.} Cutaneous wound healing proceeds through overlapping phases: an inflammatory phase in which
neutrophils and macrophages clear debris and release growth factors, a proliferative phase in
which fibroblasts migrate in and deposit collagen, and a remodeling phase. Fibroblast activity
depends on signals from inflammatory cells and on oxygen supply, while prolonged inflammation
damages tissue and delays closure. We model one healing episode that ends when the wound
deficit has fallen to a prescribed residual size, and we assign inflammatory recruitment,
fibroblast recruitment, and oxygen support through angiogenesis to three players who weigh
effort, excess inflammation, and collagen quality differently. Inflammatory recruitment (A) owns
the completion time; as for other regulatory modules, this is a synthetic allocation of decision
authority (Section~\ref{sec:assumptions}), and only A's time weight is strategically relevant.
\begin{figure}[p]
\centering
\includegraphics[width=\textwidth,height=\casefigheight,keepaspectratio]{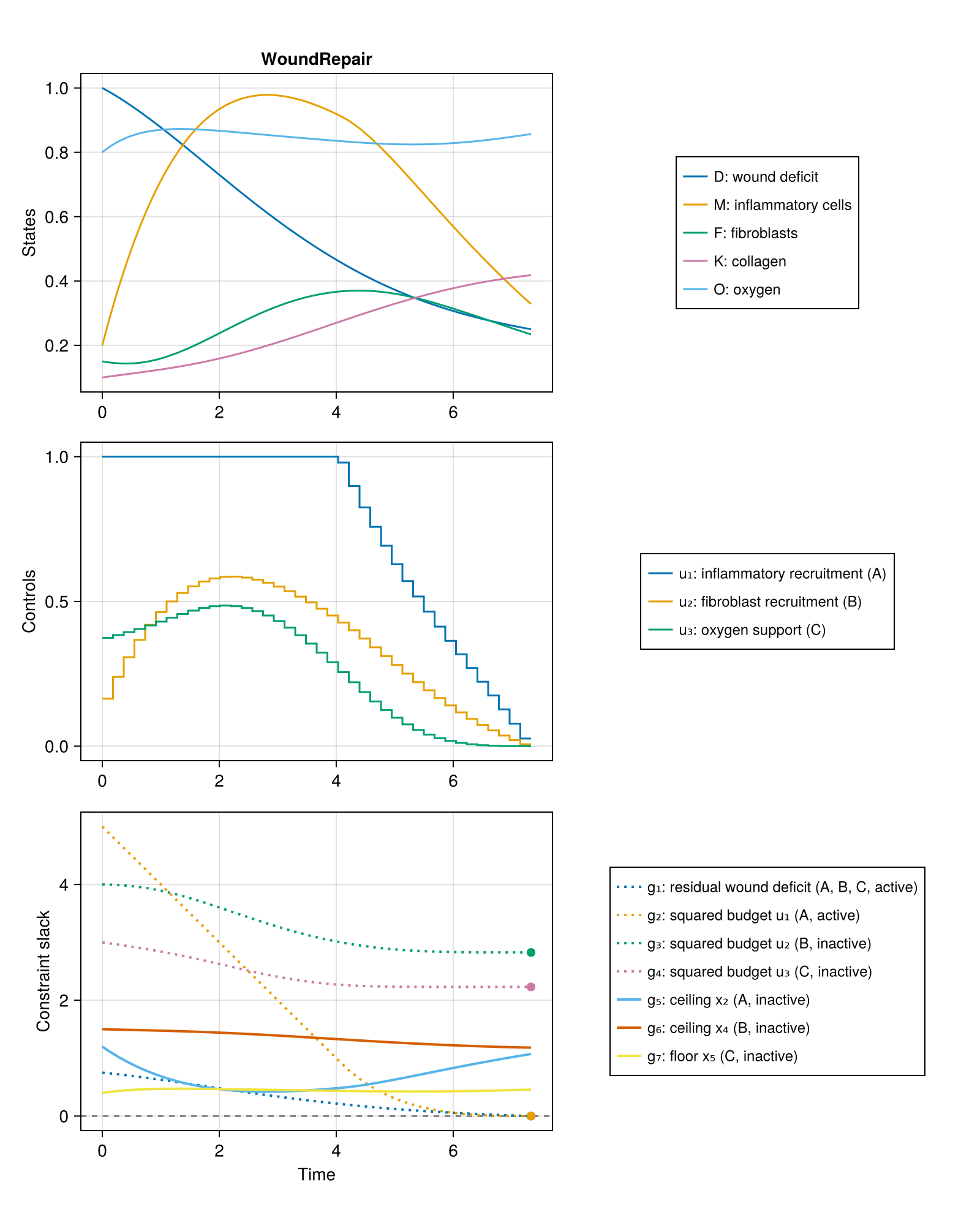}
\caption[Reference equilibrium candidate for \texttt{\caseidentifier}]{Reference equilibrium candidate for \texttt{\caseidentifier}. \textit{Interpretation.} Healing ends at $T\approx7.3$, when the wound deficit has fallen from $1$ to the prescribed
residual $0.25$ (terminal equality $g_1$ active for all players). Inflammatory recruitment starts
at its maximum and stays there until $t\approx4$, then decreases to zero, using up its entire
effort budget ($g_2$ active). Fibroblast recruitment and oxygen support rise during the first two
time units to about $0.59$ and $0.49$ and decline afterwards, oxygen support earliest. The state
trajectories reproduce the classical sequence of healing phases: inflammatory cells rise from
$0.2$ to a peak near $1$ at $t\approx2.8$ and then resolve, fibroblasts follow with a delay and peak
near $t\approx4.4$, and collagen accumulates steadily to about $0.42$ while the deficit decreases.
Oxygen increases slightly while support is given and stays above its floor. The budgets of B
and C and the inflammation, collagen, and oxygen bounds remain slack; the transition from
inflammation to matrix deposition is driven by the shared closure target and the time cost.}\label{fig:\caseidentifier}
\end{figure}
\myheader{Literature and provenance.} Macrophage--fibroblast--collagen interactions motivate wound-healing models \citep{BenchmarkPMC2839370}. Here a reduced five-state system assigns inflammatory recruitment, fibroblast recruitment, and oxygen support to separate players. Duration ownership, objective weights, terminal wound target, and other restrictions are introduced with synthetic coefficients.
\myheader{Dimensions.} $N=3$, $n_x=5$, $n_{\rm ODE}=11$, $(n_{u,i})_{i\in\mathcal P}=(3\times 1)$, $(n_{g,i})_{i\in\mathcal P}=(3\times 3)$.
\myheader{Categories.} \cat{cellular}, \cat{episodical}, \cat{regulation}, \cat{freeTime}.
\myheader{Game specification.} \[\begin{gathered}
\dot D=-(r_1+r_2M+r_3F)D,\quad \dot M=b_1u_A D-d_1M,\\
\dot F=b_2u_B MO/(1+O)-d_2F,\\
\dot K=k_KFO/(1+O)-d_3K,\\
\dot O=k_O(1+u_C)(1-O)-c_O(M+F)O.
\end{gathered}\]

$r=$ \texttt{repair}, $b=$ \texttt{recruitment}, $d=$ \texttt{decay}, $k_K=$ \texttt{collagen\_rate}, $k_O=$ \texttt{oxygen\_supply}, $c_O=$ \texttt{oxygen\_use}, $K_\star=$ \texttt{collagen\_target}.

$J_i=w_i^1T+w_i^2\int_0^T u_i^2\,dt+w_i^3\int_0^T L_i\,dt$.

$L_i=(K-K_\star)^2+(0.2+0.1\iota(i))M^2$. The duration term penalizes healing time, effort penalizes intervention, and $L_i$ balances matrix restoration against inflammation.

$D(0)=1$, $0.25-D(T)\geq0$ and $D(T)-0.25\geq0$, shared by all players. Private effort budgets are $(5,4,3)$. A imposes $M\leq1.4$ (limiting collateral damage by excess inflammation), B imposes $K\leq1.6$ (avoiding fibrotic overproduction of matrix), and C imposes $O\geq0.4$ (avoiding hypoxia) along the trajectory. Other terminal states are free.

Enabled bounds:
\begin{center}\small\begin{tabular}{lll}\toprule Owners & Bound & Location\\\midrule
A, B, C & $y_{1}(T)=0.25$ & terminal \\
A & $y_{6}(T)\leq 5$ & terminal \\
B & $y_{7}(T)\leq 4$ & terminal \\
C & $y_{8}(T)\leq 3$ & terminal \\
A & $y_{2}(t)\leq 1.4$ & path \\
B & $y_{4}(t)\leq 1.6$ & path \\
C & $0.4\leq y_{5}(t)$ & path \\
\bottomrule\end{tabular}\end{center}
\myheader{Parameter instance.} %
\begin{lstlisting}
x0 = (1.0, 0.2, 0.15, 0.1, 0.8)
target = 0.25
repair = (0.04, 0.12, 0.22)
recruitment = (0.8, 0.6)
decay = (0.5, 0.25, 0.08)
oxygen_supply = 0.7
oxygen_use = 0.15
collagen_rate = 0.5
collagen_target = 0.6
duration = 5.0
duration_bounds = (0.1, 20.0)
initial = (0.5, 0.45, 0.4)
weights = [0.75 0.15 0.1; 0.55 0.3 0.15; 0.45 0.25 0.3]
budgets = (5.0, 4.0, 3.0)
path_limits = (1.4, 1.6, 0.4)
control_lower = [0.0, 0.0, 0.0]
control_upper = [1.0, 1.0, 1.0]
numerical_x0 = [1.0, 0.2, 0.15, 0.1, 0.8, 0.0, 0.0, 0.0, 0.0, 0.0, 0.0]
\end{lstlisting}

\endgroup

\clearpage
\section{Numerical settings and solver statistics}\label{app:numerics}
This appendix collects the numerical details of Section~\ref{sec:methods}.

\paragraph*{Transcription and nonlinear programs.}
Controls are piecewise constant on $M$ equal intervals. Path restrictions are imposed at four equidistant points per control interval (eight in \caseref{antColonies}). Single shooting eliminates the state variables by integration with the explicit Tsitouras 5/4 Runge--Kutta method \code{Tsit5} at absolute and relative tolerance $10^{-9}$ ($10^{-10}$ in four instances); forward-mode automatic differentiation through the integrator supplies derivatives. The nonlinear programs are solved by Ipopt with tolerance $10^{-8}$ and a limited-memory Hessian approximation; three instances (\caseref{batchFermentation}, \caseref{circadianCycle}, \caseref{immuneClearance}) use exact Hessians. Multiple shooting and the backends BlockSQP2 and Uno are available but not used for the reference records.

\paragraph*{Iteration settings.}
The damping ranges from $\alpha=0.15$ (\caseref{neural}) to $\alpha=1$, the undamped sweep used in nine instances. Primal warm starts are enabled in five instances, and the round limit $k_{\max}$ ranges from 20 to 300. Warm starts reuse the player's previous primal controls on the same grid. Refinement transfers the coarse controls to the fine grid by interval averaging.

\paragraph*{Audits.}
The interface allows instance-specific reporting scales $\gamma_i$ in \eqref{eq:gain}, separate from the feature scales in the objectives. All reported audits use a single start per player and round, namely the current profile or, when primal warm starts are enabled and the grid layout is compatible, the player's previous response. The implementation supports additional constant-control starts for audits; they were not used for the reference records. For each assigned state restriction $c$ in the bound form \eqref{eq:bound-form}, let $\mathcal S_c$ contain either the terminal time or all audit sample times, i.e., eight equidistant points per control interval; the audit thus samples path restrictions at least as densely as the optimization. With $[a]_+=\max\{0,a\}$, the implemented violation is
\begin{equation}\label{eq:violation}
 \begin{aligned} V_i={}&\sum_{j:\,\operatorname{owner}(j)=i}\sum_{k=1}^{M}
 \left([\underline u_j-u_{jk}]_+ +[u_{jk}-\overline u_j]_+\right)
 \\ &+\sum_{c:\,\operatorname{owner}(c)=i}\max_{t\in\mathcal S_c}
 \left([\ell_c-y_{q_c}(t)]_+ +[y_{q_c}(t)-r_c]_+\right),
\end{aligned}
\end{equation}
where $y$ is the numerical state of Section~\ref{sec:formulation}. An absent finite bound contributes zero. Thus, violations are summed across restrictions after taking the maximum over time for each path restriction. The violation is not normalized across physical units or state transformations; the saved records call it the constraint norm (\code{norm_g}). The completed profile is reintegrated with \code{Tsit5} at absolute and relative tolerance $10^{-10}$.

\paragraph*{Independent checks of responses.}
A reintegration of every response at tolerance $10^{-10}$ must satisfy the player's restrictions and reproduce the shooting objective within $10^{-5}\gamma_i$; multiple shooting additionally checks its matching residuals. A response that fails one of these checks stops the computation, so every saved round passed them. The implementation logs unsuccessful nonlinear-programming return codes and uses the returned iterate if the checks pass; a successful return code is not itself a requirement. After an unsuccessful return code, a fallback solver can be started from the returned iterate. Damping can violate nonlinear endpoint restrictions. Such intermediate profiles remain in the iteration history, while acceptance requires the completed profile to pass the feasibility audit.

\paragraph*{Reproduction, timing, and solver statistics.}
Reintegrating every saved final control profile with the released model definitions reproduces the saved terminal states, including the objective and budget quadratures, to double precision; the maximum sampled violation of these reintegrations is $9.1\cdot10^{-7}$, below $\varepsilon_g=10^{-6}$. Instance-specific solver options are distributed with the models, and the archived release fixes the Julia environment used to produce the records. The CPU times in Table~\ref{tab:reference-results} are medians over all response and audit calls of one grid; the median excludes the effect of compilation in the first call of each grid. They were measured with Julia~1.11 on a machine with eight logical CPU cores and 8\,GB memory, with one instance running at a time. The per-response solver logs distinguish outer convergence from inner solver termination. All final deviation-audit calls on both grids report \code{Success} for every instance. In \caseref{motor}, nine of the 140 intermediate and final response calls report \code{MaxIters}; the other 131 report \code{Success}. All recorded calls in the remaining 24 instances report \code{Success}. The motor iterates were accepted under the independent feasibility and objective-consistency checks described above.

\bibliography{ioc-bibliography}


\begin{thebibliography}{70}
\ifx \bisbn   \undefined \def \bisbn  #1{ISBN #1}\fi
\ifx \binits  \undefined \def \binits#1{#1}\fi
\ifx \bauthor  \undefined \def \bauthor#1{#1}\fi
\ifx \batitle  \undefined \def \batitle#1{#1}\fi
\ifx \bjtitle  \undefined \def \bjtitle#1{#1}\fi
\ifx \bvolume  \undefined \def \bvolume#1{\textbf{#1}}\fi
\ifx \byear  \undefined \def \byear#1{#1}\fi
\ifx \bissue  \undefined \def \bissue#1{#1}\fi
\ifx \bfpage  \undefined \def \bfpage#1{#1}\fi
\ifx \blpage  \undefined \def \blpage #1{#1}\fi
\ifx \burl  \undefined \def \burl#1{\textsf{#1}}\fi
\ifx \doiurl  \undefined \def \doiurl#1{\url{https://doi.org/#1}}\fi
\ifx \betal  \undefined \def \betal{\textit{et al.}}\fi
\ifx \binstitute  \undefined \def \binstitute#1{#1}\fi
\ifx \binstitutionaled  \undefined \def \binstitutionaled#1{#1}\fi
\ifx \bctitle  \undefined \def \bctitle#1{#1}\fi
\ifx \beditor  \undefined \def \beditor#1{#1}\fi
\ifx \bpublisher  \undefined \def \bpublisher#1{#1}\fi
\ifx \bbtitle  \undefined \def \bbtitle#1{#1}\fi
\ifx \bedition  \undefined \def \bedition#1{#1}\fi
\ifx \bseriesno  \undefined \def \bseriesno#1{#1}\fi
\ifx \blocation  \undefined \def \blocation#1{#1}\fi
\ifx \bsertitle  \undefined \def \bsertitle#1{#1}\fi
\ifx \bsnm \undefined \def \bsnm#1{#1}\fi
\ifx \bsuffix \undefined \def \bsuffix#1{#1}\fi
\ifx \bparticle \undefined \def \bparticle#1{#1}\fi
\ifx \barticle \undefined \def \barticle#1{#1}\fi
\bibcommenthead
\ifx \bconfdate \undefined \def \bconfdate #1{#1}\fi
\ifx \botherref \undefined \def \botherref #1{#1}\fi
\ifx \url \undefined \def \url#1{\textsf{#1}}\fi
\ifx \bchapter \undefined \def \bchapter#1{#1}\fi
\ifx \bbook \undefined \def \bbook#1{#1}\fi
\ifx \bcomment \undefined \def \bcomment#1{#1}\fi
\ifx \oauthor \undefined \def \oauthor#1{#1}\fi
\ifx \citeauthoryear \undefined \def \citeauthoryear#1{#1}\fi
\ifx \endbibitem  \undefined \def \endbibitem {}\fi
\ifx \bconflocation  \undefined \def \bconflocation#1{#1}\fi
\ifx \arxivurl  \undefined \def \arxivurl#1{\textsf{#1}}\fi
\csname PreBibitemsHook\endcsname

\bibitem[\protect\citeauthoryear{Adler and Gordon}{2003}]{AdlerGordon2003Ants}
\begin{barticle}
\bauthor{\bsnm{Adler}, \binits{F.R.}},
\bauthor{\bsnm{Gordon}, \binits{D.M.}}:
\batitle{Optimization, conflict, and nonoverlapping foraging ranges in ants}.
\bjtitle{The American Naturalist}
\bvolume{162}(\bissue{5}),
\bfpage{529}--\blpage{543}
(\byear{2003})
\doiurl{10.1086/378856}
\end{barticle}
\endbibitem

\bibitem[\protect\citeauthoryear{Andersson et~al.}{2018}]{Andersson2018}
\begin{barticle}
\bauthor{\bsnm{Andersson}, \binits{J.A.E.}},
\bauthor{\bsnm{Gillis}, \binits{J.}},
\bauthor{\bsnm{Horn}, \binits{G.}},
\bauthor{\bsnm{Rawlings}, \binits{J.B.}},
\bauthor{\bsnm{Diehl}, \binits{M.}}:
\batitle{{CasADi} -- {A} software framework for nonlinear optimization and
  optimal control}.
\bjtitle{Mathematical Programming Computation}
\bvolume{11}(\bissue{1}),
\bfpage{1}--\blpage{36}
(\byear{2018})
\doiurl{10.1007/s12532-018-0139-4}
\end{barticle}
\endbibitem

\bibitem[\protect\citeauthoryear{Allen et~al.}{2023}]{Allen2023SpaceGym}
\begin{bchapter}
\bauthor{\bsnm{Allen}, \binits{R.E.}},
\bauthor{\bsnm{Rachlin}, \binits{Y.}},
\bauthor{\bsnm{Ruprecht}, \binits{J.}},
\bauthor{\bsnm{Loughran}, \binits{S.}},
\bauthor{\bsnm{Varey}, \binits{J.}},
\bauthor{\bsnm{Viggh}, \binits{H.}}:
\bctitle{{SpaceGym}: Discrete and differential games in non-cooperative space
  operations}.
In: \bbtitle{2023 IEEE Aerospace Conference},
pp. \bfpage{1}--\blpage{12}.
\bpublisher{IEEE},
\blocation{Big Sky, MT, USA}
(\byear{2023}).
\doiurl{10.1109/AERO55745.2023.10115968}
\end{bchapter}
\endbibitem

\bibitem[\protect\citeauthoryear{Baccam
  et~al.}{2006}]{BenchmarkimmuneClearance}
\begin{barticle}
\bauthor{\bsnm{Baccam}, \binits{P.}},
\bauthor{\bsnm{Beauchemin}, \binits{C.}},
\bauthor{\bsnm{Macken}, \binits{C.A.}},
\bauthor{\bsnm{Hayden}, \binits{F.G.}},
\bauthor{\bsnm{Perelson}, \binits{A.S.}}:
\batitle{Kinetics of influenza a virus infection in humans}.
\bjtitle{Journal of virology}
\bvolume{80}(\bissue{15}),
\bfpage{7590}--\blpage{7599}
(\byear{2006})
\doiurl{10.1128/jvi.01623-05}
\end{barticle}
\endbibitem

\bibitem[\protect\citeauthoryear{Bessho and
  Iwasa}{2012}]{BesshoIwasa2012Germination}
\begin{barticle}
\bauthor{\bsnm{Bessho}, \binits{K.}},
\bauthor{\bsnm{Iwasa}, \binits{Y.}}:
\batitle{Variability in the evolutionarily stable seasonal timing of
  germination and maturation of annuals and the mode of competition}.
\bjtitle{Journal of Theoretical Biology}
\bvolume{304},
\bfpage{66}--\blpage{80}
(\byear{2012})
\doiurl{10.1016/j.jtbi.2012.03.020}
\end{barticle}
\endbibitem

\bibitem[\protect\citeauthoryear{Betts et~al.}{2014}]{Betts2014Contrasted}
\begin{barticle}
\bauthor{\bsnm{Betts}, \binits{A.}},
\bauthor{\bsnm{Kaltz}, \binits{O.}},
\bauthor{\bsnm{Hochberg}, \binits{M.E.}}:
\batitle{Contrasted coevolutionary dynamics between a bacterial pathogen and
  its bacteriophages}.
\bjtitle{Proceedings of the National Academy of Sciences of the United States
  of America}
\bvolume{111}(\bissue{30}),
\bfpage{11109}--\blpage{11114}
(\byear{2014})
\doiurl{10.1073/pnas.1406763111}
\end{barticle}
\endbibitem

\bibitem[\protect\citeauthoryear{Başar and Olsder}{1982}]{BasarOlsder1998}
\begin{bbook}
\bauthor{\bsnm{Başar}, \binits{T.}},
\bauthor{\bsnm{Olsder}, \binits{G.}}:
\bbtitle{Dynamic Noncooperative Game Theory},
\bedition{2}nd edn.
\bpublisher{Society for Industrial and Applied Mathematics},
\blocation{Philadelphia}
(\byear{1982}).
\doiurl{10.1016/s0076-5392(08)x6169-7}
\end{bbook}
\endbibitem

\bibitem[\protect\citeauthoryear{Braun et~al.}{2009}]{Braun2009Motor}
\begin{barticle}
\bauthor{\bsnm{Braun}, \binits{D.A.}},
\bauthor{\bsnm{Ortega}, \binits{P.A.}},
\bauthor{\bsnm{Wolpert}, \binits{D.M.}}:
\batitle{Nash equilibria in multi-agent motor interactions}.
\bjtitle{PLoS Computational Biology}
\bvolume{5}(\bissue{8}),
\bfpage{1000468}
(\byear{2009})
\doiurl{10.1371/journal.pcbi.1000468}
\end{barticle}
\endbibitem

\bibitem[\protect\citeauthoryear{Bradford}{1990}]{BenchmarkseedGermination}
\begin{barticle}
\bauthor{\bsnm{Bradford}, \binits{K.J.}}:
\batitle{A water relations analysis of seed germination rates}.
\bjtitle{Plant Physiology}
\bvolume{94}(\bissue{2}),
\bfpage{840}--\blpage{849}
(\byear{1990})
\doiurl{10.1104/pp.94.2.840}
\end{barticle}
\endbibitem

\bibitem[\protect\citeauthoryear{Banga and Sager}{2025}]{Banga2025}
\begin{barticle}
\bauthor{\bsnm{Banga}, \binits{J.R.}},
\bauthor{\bsnm{Sager}, \binits{S.}}:
\batitle{Generalized inverse optimal control and its application in biology}.
\bjtitle{Annual Reviews in Control}
\bvolume{60},
\bfpage{101029}
(\byear{2025})
\doiurl{10.1016/j.arcontrol.2025.101029}
\end{barticle}
\endbibitem

\bibitem[\protect\citeauthoryear{Banga and
  Villaverde}{2026}]{BangaVillaverde2026Hierarchical}
\begin{botherref}
\oauthor{\bsnm{Banga}, \binits{J.R.}},
\oauthor{\bsnm{Villaverde}, \binits{A.F.}}:
A Hierarchical Framework for Inverse Problems in Biological Cybernetics:
  Opportunities and Limitations.
Consejo Superior de Investigaciones Científicas
(2026).
\doiurl{10.20350/digitalCSIC/29135}
\end{botherref}
\endbibitem

\bibitem[\protect\citeauthoryear{Caillau
  et~al.}{2025}]{OptimalControlProblemsRepository}
\begin{botherref}
\oauthor{\bsnm{Caillau}, \binits{J.-B.}},
\oauthor{\bsnm{Cots}, \binits{O.}},
\oauthor{\bsnm{Gergaud}, \binits{J.}},
\oauthor{\bsnm{Martinon}, \binits{P.}}:
{OptimalControlProblems.jl}: A collection of optimal control problems with
  {ODE}s in {Julia}.
Software, version 0.2.2.
Accessed September 23, 2026
(2025).
\doiurl{10.5281/zenodo.17013180} .
\url{https://github.com/control-toolbox/OptimalControlProblems.jl}
\end{botherref}
\endbibitem

\bibitem[\protect\citeauthoryear{Cumming et~al.}{2010}]{BenchmarkPMC2839370}
\begin{barticle}
\bauthor{\bsnm{Cumming}, \binits{B.D.}},
\bauthor{\bsnm{McElwain}, \binits{D.L.S.}},
\bauthor{\bsnm{Upton}, \binits{Z.}}:
\batitle{A mathematical model of wound healing and subsequent scarring}.
\bjtitle{Journal of the Royal Society, Interface}
\bvolume{7}(\bissue{42}),
\bfpage{19}--\blpage{34}
(\byear{2010})
\doiurl{10.1098/rsif.2008.0536}
\end{barticle}
\endbibitem

\bibitem[\protect\citeauthoryear{Chackochan and
  Sanguineti}{2019}]{Chackochan2019JointAction}
\begin{barticle}
\bauthor{\bsnm{Chackochan}, \binits{V.T.}},
\bauthor{\bsnm{Sanguineti}, \binits{V.}}:
\batitle{Incomplete information about the partner affects the development of
  collaborative strategies in joint action}.
\bjtitle{PLOS Computational Biology}
\bvolume{15}(\bissue{12}),
\bfpage{1006385}
(\byear{2019})
\doiurl{10.1371/journal.pcbi.1006385}
\end{barticle}
\endbibitem

\bibitem[\protect\citeauthoryear{Dockner et~al.}{1989}]{Dockner1989Fishery}
\begin{barticle}
\bauthor{\bsnm{Dockner}, \binits{E.}},
\bauthor{\bsnm{Feichtinger}, \binits{G.}},
\bauthor{\bsnm{Mehlmann}, \binits{A.}}:
\batitle{Noncooperative solutions for a differential game model of fishery}.
\bjtitle{Journal of Economic Dynamics and Control}
\bvolume{13}(\bissue{1}),
\bfpage{1}--\blpage{20}
(\byear{1989})
\doiurl{10.1016/0165-1889(89)90008-0}
\end{barticle}
\endbibitem

\bibitem[\protect\citeauthoryear{Da-Silva et~al.}{2012}]{BenchmarkantWater}
\begin{barticle}
\bauthor{\bsnm{Da-Silva}, \binits{A.C.}},
\bauthor{\bsnm{Navas}, \binits{C.A.}},
\bauthor{\bsnm{Ribeiro}, \binits{P.L.}}:
\batitle{Dealing with water deficit in atta ant colonies: large ants scout for
  water while small ants transport it}.
\bjtitle{Biology open}
\bvolume{1}(\bissue{9}),
\bfpage{827}--\blpage{830}
(\byear{2012})
\doiurl{10.1242/bio.2012703}
\end{barticle}
\endbibitem

\bibitem[\protect\citeauthoryear{Facchinei and
  Kanzow}{2010}]{FacchineiKanzow2010}
\begin{barticle}
\bauthor{\bsnm{Facchinei}, \binits{F.}},
\bauthor{\bsnm{Kanzow}, \binits{C.}}:
\batitle{Generalized nash equilibrium problems}.
\bjtitle{Annals of Operations Research}
\bvolume{175}(\bissue{1}),
\bfpage{177}--\blpage{211}
(\byear{2010})
\doiurl{10.1007/s10479-009-0653-x}
\end{barticle}
\endbibitem

\bibitem[\protect\citeauthoryear{Fridovich-Keil
  et~al.}{2020}]{FridovichKeil2020}
\begin{bchapter}
\bauthor{\bsnm{Fridovich-Keil}, \binits{D.}},
\bauthor{\bsnm{Ratner}, \binits{E.}},
\bauthor{\bsnm{Peters}, \binits{L.}},
\bauthor{\bsnm{Dragan}, \binits{A.D.}},
\bauthor{\bsnm{Tomlin}, \binits{C.J.}}:
\bctitle{Efficient iterative linear-quadratic approximations for nonlinear
  multi-player general-sum differential games}.
In: \bbtitle{2020 IEEE International Conference on Robotics and Automation
  (ICRA)},
pp. \bfpage{1475}--\blpage{1481}.
\bpublisher{IEEE},
\blocation{Paris, France}
(\byear{2020}).
\doiurl{10.1109/ICRA40945.2020.9197129} .
\burl{https://arxiv.org/abs/1909.04694}
\end{bchapter}
\endbibitem

\bibitem[\protect\citeauthoryear{Greenwald et~al.}{2018}]{BenchmarkantResource}
\begin{barticle}
\bauthor{\bsnm{Greenwald}, \binits{E.E.}},
\bauthor{\bsnm{Baltiansky}, \binits{L.}},
\bauthor{\bsnm{Feinerman}, \binits{O.}}:
\batitle{Individual crop loads provide local control for collective food intake
  in ant colonies}.
\bjtitle{eLife}
\bvolume{7},
\bfpage{31730}
(\byear{2018})
\doiurl{10.7554/elife.31730}
\end{barticle}
\endbibitem

\bibitem[\protect\citeauthoryear{Gersani et~al.}{2001}]{Gersani2001RootCommons}
\begin{barticle}
\bauthor{\bsnm{Gersani}, \binits{M.}},
\bauthor{\bsnm{Brown}, \binits{J.S.}},
\bauthor{\bsnm{O'Brien}, \binits{E.E.}},
\bauthor{\bsnm{Maina}, \binits{G.M.}},
\bauthor{\bsnm{Abramsky}, \binits{Z.}}:
\batitle{Tragedy of the commons as a result of root competition}.
\bjtitle{Journal of Ecology}
\bvolume{89}(\bissue{4}),
\bfpage{660}--\blpage{669}
(\byear{2001})
\doiurl{10.1046/j.0022-0477.2001.00609.x}
\end{barticle}
\endbibitem

\bibitem[\protect\citeauthoryear{Goldberg
  et~al.}{2017}]{BenchmarkplantNitrogen}
\begin{barticle}
\bauthor{\bsnm{Goldberg}, \binits{D.E.}},
\bauthor{\bsnm{Martina}, \binits{J.P.}},
\bauthor{\bsnm{Elgersma}, \binits{K.J.}},
\bauthor{\bsnm{Currie}, \binits{W.S.}}:
\batitle{Plant size and competitive dynamics along nutrient gradients}.
\bjtitle{The American naturalist}
\bvolume{190}(\bissue{2}),
\bfpage{229}--\blpage{243}
(\byear{2017})
\doiurl{10.1086/692438}
\end{barticle}
\endbibitem

\bibitem[\protect\citeauthoryear{Gonz{\'a}lez-Olivares
  et~al.}{2017}]{BenchmarkpredatorPrey}
\begin{barticle}
\bauthor{\bsnm{Gonz{\'a}lez-Olivares}, \binits{E.}},
\bauthor{\bsnm{Gonz{\'a}lez-Ya{\~n}ez}, \binits{B.}},
\bauthor{\bsnm{Becerra-Klix}, \binits{R.}},
\bauthor{\bsnm{Ramos-Jiliberto}, \binits{R.}}:
\batitle{Multiple stable states in a model based on predator-induced defenses}.
\bjtitle{Ecological Complexity}
\bvolume{32},
\bfpage{111}--\blpage{120}
(\byear{2017})
\doiurl{10.1016/j.ecocom.2017.10.004}
\end{barticle}
\endbibitem

\bibitem[\protect\citeauthoryear{Hamilton}{1967}]{Hamilton1967SexRatios}
\begin{barticle}
\bauthor{\bsnm{Hamilton}, \binits{W.D.}}:
\batitle{Extraordinary sex ratios}.
\bjtitle{Science}
\bvolume{156}(\bissue{3774}),
\bfpage{477}--\blpage{488}
(\byear{1967})
\doiurl{10.1126/science.156.3774.477}
\end{barticle}
\endbibitem

\bibitem[\protect\citeauthoryear{Harper}{1982}]{Benchmarkforaging}
\begin{barticle}
\bauthor{\bsnm{Harper}, \binits{D.G.C.}}:
\batitle{Competitive foraging in mallards: ``ideal free' ducks}.
\bjtitle{Animal Behaviour}
\bvolume{30}(\bissue{2}),
\bfpage{575}--\blpage{584}
(\byear{1982})
\doiurl{10.1016/s0003-3472(82)80071-7}
\end{barticle}
\endbibitem

\bibitem[\protect\citeauthoryear{Hiltunen and
  Becks}{2014}]{Hiltunen2014Consumer}
\begin{barticle}
\bauthor{\bsnm{Hiltunen}, \binits{T.}},
\bauthor{\bsnm{Becks}, \binits{L.}}:
\batitle{Consumer co-evolution as an important component of the
  eco-evolutionary feedback}.
\bjtitle{Nature Communications}
\bvolume{5}(\bissue{1}),
\bfpage{5226}
(\byear{2014})
\doiurl{10.1038/ncomms6226}
\end{barticle}
\endbibitem

\bibitem[\protect\citeauthoryear{Hummert
  et~al.}{2014}]{hummert2014evolutionary}
\begin{barticle}
\bauthor{\bsnm{Hummert}, \binits{S.}},
\bauthor{\bsnm{Bohl}, \binits{K.}},
\bauthor{\bsnm{Basanta}, \binits{D.}},
\bauthor{\bsnm{Deutsch}, \binits{A.}},
\bauthor{\bsnm{Werner}, \binits{S.}},
\bauthor{\bsnm{Theißen}, \binits{G.}},
\bauthor{\bsnm{Schroeter}, \binits{A.}},
\bauthor{\bsnm{Schuster}, \binits{S.}}:
\batitle{Evolutionary game theory: cells as players}.
\bjtitle{Molecular BioSystems}
\bvolume{10}(\bissue{12}),
\bfpage{3044}--\blpage{3065}
(\byear{2014})
\doiurl{10.1039/c3mb70602h}
\end{barticle}
\endbibitem

\bibitem[\protect\citeauthoryear{Hamelin and
  Lewis}{2010}]{HamelinLewis2010Territoriality}
\begin{barticle}
\bauthor{\bsnm{Hamelin}, \binits{F.M.}},
\bauthor{\bsnm{Lewis}, \binits{M.A.}}:
\batitle{A differential game theoretical analysis of mechanistic models for
  territoriality}.
\bjtitle{Journal of Mathematical Biology}
\bvolume{61}(\bissue{5}),
\bfpage{665}--\blpage{694}
(\byear{2010})
\doiurl{10.1007/s00285-009-0316-1}
\end{barticle}
\endbibitem

\bibitem[\protect\citeauthoryear{Ilango
  et~al.}{2012}]{Ilango2012DopamineAvoidance}
\begin{barticle}
\bauthor{\bsnm{Ilango}, \binits{A.}},
\bauthor{\bsnm{Shumake}, \binits{J.}},
\bauthor{\bsnm{Wetzel}, \binits{W.}},
\bauthor{\bsnm{Scheich}, \binits{H.}},
\bauthor{\bsnm{Ohl}, \binits{F.W.}}:
\batitle{The role of dopamine in the context of aversive stimuli with
  particular reference to acoustically signaled avoidance learning}.
\bjtitle{Frontiers in Neuroscience}
\bvolume{6},
\bfpage{132}
(\byear{2012})
\doiurl{10.3389/fnins.2012.00132}
\end{barticle}
\endbibitem

\bibitem[\protect\citeauthoryear{Klipp et~al.}{2002}]{Klipp2002-sg}
\begin{barticle}
\bauthor{\bsnm{Klipp}, \binits{E.}},
\bauthor{\bsnm{Heinrich}, \binits{R.}},
\bauthor{\bsnm{Holzhütter}, \binits{H.-G.}}:
\batitle{Prediction of temporal gene expression. metabolic opimization by
  re-distribution of enzyme activities}.
\bjtitle{European Journal of Biochemistry}
\bvolume{269}(\bissue{22}),
\bfpage{5406}--\blpage{5413}
(\byear{2002})
\doiurl{10.1046/J.1432-1033.2002.03223.X}
\end{barticle}
\endbibitem

\bibitem[\protect\citeauthoryear{Kanzow and Steck}{2016}]{KanzowSteck2016}
\begin{barticle}
\bauthor{\bsnm{Kanzow}, \binits{C.}},
\bauthor{\bsnm{Steck}, \binits{D.}}:
\batitle{Augmented lagrangian methods for the solution of generalized nash
  equilibrium problems}.
\bjtitle{SIAM Journal on Optimization}
\bvolume{26}(\bissue{4}),
\bfpage{2034}--\blpage{2058}
(\byear{2016})
\doiurl{10.1137/16M1068256}
\end{barticle}
\endbibitem

\bibitem[\protect\citeauthoryear{Le~Cleac'h et~al.}{2021}]{LeCleach2022ALGAMES}
\begin{barticle}
\bauthor{\bsnm{Le~Cleac'h}, \binits{S.}},
\bauthor{\bsnm{Schwager}, \binits{M.}},
\bauthor{\bsnm{Manchester}, \binits{Z.}}:
\batitle{{ALGAMES}: A fast augmented lagrangian solver for constrained dynamic
  games}.
\bjtitle{Autonomous Robots}
\bvolume{46}(\bissue{1}),
\bfpage{201}--\blpage{215}
(\byear{2021})
\doiurl{10.1007/s10514-021-10024-7}
\end{barticle}
\endbibitem

\bibitem[\protect\citeauthoryear{Low et~al.}{2013}]{Low2013ColdWar}
\begin{barticle}
\bauthor{\bsnm{Low}, \binits{C.}},
\bauthor{\bsnm{Ellner}, \binits{S.P.}},
\bauthor{\bsnm{Holden}, \binits{M.H.}}:
\batitle{Optimal control and cold war dynamics between plant and herbivore}.
\bjtitle{The American Naturalist}
\bvolume{182}(\bissue{2}),
\bfpage{25}--\blpage{39}
(\byear{2013})
\doiurl{10.1086/670810}
\end{barticle}
\endbibitem

\bibitem[\protect\citeauthoryear{Leloup and
  Goldbeter}{1998}]{LeloupGoldbeter1998PerTim}
\begin{barticle}
\bauthor{\bsnm{Leloup}, \binits{J.-C.}},
\bauthor{\bsnm{Goldbeter}, \binits{A.}}:
\batitle{A model for circadian rhythms in {Drosophila} incorporating the
  formation of a complex between the {PER} and {TIM} proteins}.
\bjtitle{Journal of Biological Rhythms}
\bvolume{13}(\bissue{1}),
\bfpage{70}--\blpage{87}
(\byear{1998})
\doiurl{10.1177/074873098128999934}
\end{barticle}
\endbibitem

\bibitem[\protect\citeauthoryear{Lanctot et~al.}{2019}]{Lanctot2019}
\begin{botherref}
\oauthor{\bsnm{Lanctot}, \binits{M.}},
\oauthor{\bsnm{Lockhart}, \binits{E.}},
\oauthor{\bsnm{Lespiau}, \binits{J.-B.}},
\oauthor{\bsnm{Zambaldi}, \binits{V.}},
\oauthor{\bsnm{Upadhyay}, \binits{S.}},
\oauthor{\bsnm{Pérolat}, \binits{J.}},
\oauthor{\bsnm{Srinivasan}, \binits{S.}},
\oauthor{\bsnm{Timbers}, \binits{F.}},
\oauthor{\bsnm{Tuyls}, \binits{K.}},
\oauthor{\bsnm{Omidshafiei}, \binits{S.}}, et al.:
{OpenSpiel}: A Framework for Reinforcement Learning in Games.
arXiv:1908.09453
(2019).
\url{https://arxiv.org/abs/1908.09453}
\end{botherref}
\endbibitem

\bibitem[\protect\citeauthoryear{Lehtonen and
  Schwanz}{2018}]{BenchmarkPMC5830725}
\begin{barticle}
\bauthor{\bsnm{Lehtonen}, \binits{J.}},
\bauthor{\bsnm{Schwanz}, \binits{L.E.}}:
\batitle{Mate limitation and sex ratio evolution}.
\bjtitle{Royal Society open science}
\bvolume{5}(\bissue{2}),
\bfpage{171135}
(\byear{2018})
\doiurl{10.1098/rsos.171135}
\end{barticle}
\endbibitem

\bibitem[\protect\citeauthoryear{Molloy et~al.}{2022}]{molloy2022inverse}
\begin{bbook}
\bauthor{\bsnm{Molloy}, \binits{T.L.}},
\bauthor{\bsnm{Charaja}, \binits{J.I.}},
\bauthor{\bsnm{Hohmann}, \binits{S.}},
\bauthor{\bsnm{Perez}, \binits{T.}}:
\bbtitle{Inverse Optimal Control and Inverse Noncooperative Dynamic Game
  Theory}.
\bpublisher{Springer},
\blocation{Cham}
(\byear{2022}).
\doiurl{10.1007/978-3-030-93317-3}
\end{bbook}
\endbibitem

\bibitem[\protect\citeauthoryear{McNamara
  et~al.}{2001}]{mcnamara2001optimality}
\begin{barticle}
\bauthor{\bsnm{McNamara}, \binits{J.}},
\bauthor{\bsnm{Houston}, \binits{A.}},
\bauthor{\bsnm{Collins}, \binits{E.J.}}:
\batitle{Optimality models in behavioral biology}.
\bjtitle{Siam Review}
\bvolume{43}(\bissue{3}),
\bfpage{413}--\blpage{466}
(\byear{2001})
\doiurl{10.1137/S0036144500385263}
\end{barticle}
\endbibitem

\bibitem[\protect\citeauthoryear{Molloy et~al.}{2020}]{Molloy2020-vg}
\begin{barticle}
\bauthor{\bsnm{Molloy}, \binits{T.L.}},
\bauthor{\bsnm{Inga}, \binits{J.}},
\bauthor{\bsnm{Flad}, \binits{M.}},
\bauthor{\bsnm{Ford}, \binits{J.J.}},
\bauthor{\bsnm{Perez}, \binits{T.}},
\bauthor{\bsnm{Hohmann}, \binits{S.}}:
\batitle{Inverse open-loop noncooperative differential games and inverse
  optimal control}.
\bjtitle{IEEE Transactions on Automatic Control}
\bvolume{65}(\bissue{2}),
\bfpage{897}--\blpage{904}
(\byear{2020})
\doiurl{10.1109/TAC.2019.2921835}
\end{barticle}
\endbibitem

\bibitem[\protect\citeauthoryear{Maynard~Smith}{1982}]{Smith1982}
\begin{bbook}
\bauthor{\bsnm{Maynard~Smith}, \binits{J.}}:
\bbtitle{Evolution and the Theory of Games}.
\bpublisher{Cambridge University Press},
\blocation{Cambridge}
(\byear{1982}).
\doiurl{10.1017/cbo9780511806292}
\end{bbook}
\endbibitem

\bibitem[\protect\citeauthoryear{Maynard~Smith and
  Price}{1973}]{MaynardSmithPrice1973}
\begin{barticle}
\bauthor{\bsnm{Maynard~Smith}, \binits{J.}},
\bauthor{\bsnm{Price}, \binits{G.R.}}:
\batitle{The logic of animal conflict}.
\bjtitle{Nature}
\bvolume{246}(\bissue{5427}),
\bfpage{15}--\blpage{18}
(\byear{1973})
\doiurl{10.1038/246015a0}
\end{barticle}
\endbibitem

\bibitem[\protect\citeauthoryear{Nudelman et~al.}{2004}]{Nudelman2004}
\begin{bchapter}
\bauthor{\bsnm{Nudelman}, \binits{E.}},
\bauthor{\bsnm{Wortman}, \binits{J.}},
\bauthor{\bsnm{Shoham}, \binits{Y.}},
\bauthor{\bsnm{Leyton-Brown}, \binits{K.}}:
\bctitle{Run the {GAMUT}: A comprehensive approach to evaluating game-theoretic
  algorithms}.
In: \bbtitle{Proceedings of the Third International Joint Conference on
  Autonomous Agents and Multiagent Systems, 2004. AAMAS 2004.},
pp. \bfpage{880}--\blpage{887}.
\bpublisher{IEEE Computer Society},
\blocation{New York, NY, USA}
(\byear{2004}).
\doiurl{10.1109/AAMAS.2004.238}
\end{bchapter}
\endbibitem

\bibitem[\protect\citeauthoryear{O'Brien et~al.}{2007}]{BenchmarkPMC2141669}
\begin{barticle}
\bauthor{\bsnm{O'Brien}, \binits{E.E.}},
\bauthor{\bsnm{Brown}, \binits{J.S.}},
\bauthor{\bsnm{Moll}, \binits{J.D.}}:
\batitle{Roots in space: a spatially explicit model for below-ground
  competition in plants}.
\bjtitle{Proceedings of the Royal Society B: Biological Sciences}
\bvolume{274}(\bissue{1612}),
\bfpage{929}--\blpage{935}
(\byear{2007})
\doiurl{10.1098/rspb.2006.0113}
\end{barticle}
\endbibitem

\bibitem[\protect\citeauthoryear{Oyarz{\'u}n
  et~al.}{2008}]{Oyarzun2008PathwayArxiv}
\begin{botherref}
\oauthor{\bsnm{Oyarz{\'u}n}, \binits{D.}},
\oauthor{\bsnm{Ingalls}, \binits{B.}},
\oauthor{\bsnm{Middleton}, \binits{R.}},
\oauthor{\bsnm{Kalamatianos}, \binits{D.}}:
Optimal metabolic pathway activation.
Elsevier BV
(2008).
\doiurl{10.3182/20080706-5-kr-1001.02130} .
\url{https://arxiv.org/abs/0801.2447}
\end{botherref}
\endbibitem

\bibitem[\protect\citeauthoryear{Oyarzún et~al.}{2009}]{Oyarzun2009-ef}
\begin{barticle}
\bauthor{\bsnm{Oyarzún}, \binits{D.}},
\bauthor{\bsnm{Ingalls}, \binits{B.}},
\bauthor{\bsnm{Middleton}, \binits{R.}},
\bauthor{\bsnm{Kalamatianos}, \binits{D.}}:
\batitle{Sequential activation of metabolic pathways: a dynamic optimization
  approach}.
\bjtitle{Bulletin of Mathematical Biology}
\bvolume{71}(\bissue{8}),
\bfpage{1851}--\blpage{1872}
(\byear{2009})
\doiurl{10.1007/s11538-009-9427-5}
\end{barticle}
\endbibitem

\bibitem[\protect\citeauthoryear{Parker and
  Maynard~Smith}{1990}]{Parker1990-zj}
\begin{barticle}
\bauthor{\bsnm{Parker}, \binits{G.A.}},
\bauthor{\bsnm{Maynard~Smith}, \binits{J.}}:
\batitle{Optimality theory in evolutionary biology}.
\bjtitle{Nature}
\bvolume{348}(\bissue{6296}),
\bfpage{27}--\blpage{33}
(\byear{1990})
\doiurl{10.1038/348027a0}
\end{barticle}
\endbibitem

\bibitem[\protect\citeauthoryear{Patterson and Rao}{2014}]{PattersonRao2014}
\begin{barticle}
\bauthor{\bsnm{Patterson}, \binits{M.A.}},
\bauthor{\bsnm{Rao}, \binits{A.V.}}:
\batitle{{GPOPS-II}: A {MATLAB} software for solving multiple-phase optimal
  control problems using hp-adaptive {Gaussian} quadrature collocation methods
  and sparse nonlinear programming}.
\bjtitle{ACM Transactions on Mathematical Software}
\bvolume{41}(\bissue{1}),
\bfpage{1}--\blpage{1137}
(\byear{2014})
\doiurl{10.1145/2558904}
\end{barticle}
\endbibitem

\bibitem[\protect\citeauthoryear{Paterson
  et~al.}{2010}]{Paterson2010Coevolution}
\begin{barticle}
\bauthor{\bsnm{Paterson}, \binits{S.}},
\bauthor{\bsnm{Vogwill}, \binits{T.}},
\bauthor{\bsnm{Buckling}, \binits{A.}},
\bauthor{\bsnm{Benmayor}, \binits{R.}},
\bauthor{\bsnm{Spiers}, \binits{A.}},
\bauthor{\bsnm{Thomson}, \binits{N.}},
\bauthor{\bsnm{Quail}, \binits{M.}},
\bauthor{\bsnm{Smith}, \binits{F.}},
\bauthor{\bsnm{Walker}, \binits{D.}},
\bauthor{\bsnm{Libberton}, \binits{B.}}, \betal:
\batitle{Antagonistic coevolution accelerates molecular evolution}.
\bjtitle{Nature}
\bvolume{464}(\bissue{7286}),
\bfpage{275}--\blpage{278}
(\byear{2010})
\doiurl{10.1038/nature08798}
\end{barticle}
\endbibitem

\bibitem[\protect\citeauthoryear{Revilla and K{\v
  r}ivan}{2016}]{BenchmarkPMC4978411}
\begin{barticle}
\bauthor{\bsnm{Revilla}, \binits{T.A.}},
\bauthor{\bsnm{K{\v r}ivan}, \binits{V.}}:
\batitle{Pollinator foraging adaptation and coexistence of competing plants}.
\bjtitle{PloS one}
\bvolume{11}(\bissue{8}),
\bfpage{0160076}
(\byear{2016})
\doiurl{10.1371/journal.pone.0160076}
\end{barticle}
\endbibitem

\bibitem[\protect\citeauthoryear{Rao et~al.}{2020}]{Rao2020NutritionalAnts}
\begin{botherref}
\oauthor{\bsnm{Rao}, \binits{F.}},
\oauthor{\bsnm{Messan}, \binits{M.R.}},
\oauthor{\bsnm{Marquez}, \binits{A.}},
\oauthor{\bsnm{Smith}, \binits{N.}},
\oauthor{\bsnm{Kang}, \binits{Y.}}:
Nutritional Regulation Influencing Colony Dynamics and Task Allocations in
  Social Insect Colonies.
Informa UK Limited
(2020).
\doiurl{10.1080/17513758.2020.1786859} .
\url{https://arxiv.org/abs/2005.03432}
\end{botherref}
\endbibitem

\bibitem[\protect\citeauthoryear{Rosenzweig}{1971}]{rosenzweig1971paradox}
\begin{barticle}
\bauthor{\bsnm{Rosenzweig}, \binits{M.L.}}:
\batitle{Paradox of enrichment: destabilization of exploitation ecosystems in
  ecological time}.
\bjtitle{Science}
\bvolume{171}(\bissue{3969}),
\bfpage{385}--\blpage{387}
(\byear{1971})
\doiurl{10.1126/science.171.3969.385}
\end{barticle}
\endbibitem

\bibitem[\protect\citeauthoryear{Sager}{2012}]{Sager2012-ts}
\begin{bchapter}
\bauthor{\bsnm{Sager}, \binits{S.}}:
\bctitle{{A} benchmark library of mixed-integer optimal control problems}.
In: \beditor{\bsnm{Lee}, \binits{J.}},
\beditor{\bsnm{Leyffer}, \binits{S.}} (eds.)
\bbtitle{{M}ixed {I}nteger {N}onlinear {P}rogramming},
pp. \bfpage{631}--\blpage{670}.
\bpublisher{Springer},
\blocation{New York, NY}
(\byear{2012}).
\doiurl{10.1007/978-1-4614-1927-3_22}
\end{bchapter}
\endbibitem

\bibitem[\protect\citeauthoryear{{SciML}}{n.d.}]{CorleoneRepository}
\begin{botherref}
\oauthor{\bsnm{{SciML}}}:
{Corleone.jl}: Optimal Control meets Scientific Machine Learning.
Software repository, accessed September 23, 2026
(n.d.).
\url{https://github.com/SciML/Corleone.jl}
\end{botherref}
\endbibitem

\bibitem[\protect\citeauthoryear{Santra et~al.}{2017}]{BenchmarkPMC5380058}
\begin{barticle}
\bauthor{\bsnm{Santra}, \binits{M.}},
\bauthor{\bsnm{Farrell}, \binits{D.W.}},
\bauthor{\bsnm{Dill}, \binits{K.A.}}:
\batitle{Bacterial proteostasis balances energy and chaperone utilization
  efficiently}.
\bjtitle{Proceedings of the National Academy of Sciences of the United States
  of America}
\bvolume{114}(\bissue{13}),
\bfpage{2654}--\blpage{2661}
(\byear{2017})
\doiurl{10.1073/pnas.1620646114}
\end{barticle}
\endbibitem

\bibitem[\protect\citeauthoryear{Stephens and Krebs}{1986}]{StephensKrebs1986}
\begin{bbook}
\bauthor{\bsnm{Stephens}, \binits{D.W.}},
\bauthor{\bsnm{Krebs}, \binits{J.R.}}:
\bbtitle{Foraging Theory}.
\bpublisher{Princeton University Press},
\blocation{Princeton, NJ}
(\byear{1986}).
\doiurl{10.2307/j.ctvs32s6b}
\end{bbook}
\endbibitem

\bibitem[\protect\citeauthoryear{Schilling et~al.}{2020}]{Schilling2020Hexapod}
\begin{bchapter}
\bauthor{\bsnm{Schilling}, \binits{M.}},
\bauthor{\bsnm{Konen}, \binits{K.}},
\bauthor{\bsnm{Ohl}, \binits{F.W.}},
\bauthor{\bsnm{Korthals}, \binits{T.}}:
\bctitle{Decentralized deep reinforcement learning for a distributed and
  adaptive locomotion controller of a hexapod robot}.
In: \bbtitle{2020 IEEE/RSJ International Conference on Intelligent Robots and
  Systems},
pp. \bfpage{5335}--\blpage{5342}.
\bpublisher{IEEE},
\blocation{Las Vegas, NV, USA}
(\byear{2020}).
\doiurl{10.1109/IROS45743.2020.9341754} .
\burl{https://arxiv.org/abs/2005.11164}
\end{bchapter}
\endbibitem

\bibitem[\protect\citeauthoryear{Schuster et~al.}{2008}]{Schuster2008}
\begin{barticle}
\bauthor{\bsnm{Schuster}, \binits{S.}},
\bauthor{\bsnm{Kreft}, \binits{J.-U.}},
\bauthor{\bsnm{Schroeter}, \binits{A.}},
\bauthor{\bsnm{Pfeiffer}, \binits{T.}}:
\batitle{Use of game-theoretical methods in biochemistry and biophysics}.
\bjtitle{Journal of Biological Physics}
\bvolume{34}(\bissue{1–2}),
\bfpage{1}--\blpage{17}
(\byear{2008})
\doiurl{10.1007/s10867-008-9101-4}
\end{barticle}
\endbibitem

\bibitem[\protect\citeauthoryear{Sun
  et~al.}{2018}]{Sun2018PollinatorSignalling}
\begin{barticle}
\bauthor{\bsnm{Sun}, \binits{S.}},
\bauthor{\bsnm{Leshowitz}, \binits{M.I.}},
\bauthor{\bsnm{Rycht{\'a}{\v r}}, \binits{J.}}:
\batitle{The signalling game between plants and pollinators}.
\bjtitle{Scientific Reports}
\bvolume{8}(\bissue{1}),
\bfpage{6686}
(\byear{2018})
\doiurl{10.1038/s41598-018-24779-0}
\end{barticle}
\endbibitem

\bibitem[\protect\citeauthoryear{Schilling
  et~al.}{2019}]{Schilling2019FluidAdaptivity}
\begin{bchapter}
\bauthor{\bsnm{Schilling}, \binits{M.}},
\bauthor{\bsnm{Ritter}, \binits{H.}},
\bauthor{\bsnm{Ohl}, \binits{F.W.}}:
\bctitle{From crystallized adaptivity to fluid adaptivity in deep reinforcement
  learning -- insights from biological systems on adaptive flexibility}.
In: \bbtitle{2019 IEEE International Conference on Systems, Man and
  Cybernetics},
pp. \bfpage{1472}--\blpage{1478}.
\bpublisher{IEEE},
\blocation{Bari, Italy}
(\byear{2019}).
\doiurl{10.1109/SMC.2019.8914302} .
\burl{https://arxiv.org/abs/1908.05348}
\end{bchapter}
\endbibitem

\bibitem[\protect\citeauthoryear{Slaby et~al.}{2007}]{Slaby2007-tr}
\begin{barticle}
\bauthor{\bsnm{Slaby}, \binits{O.}},
\bauthor{\bsnm{Sager}, \binits{S.}},
\bauthor{\bsnm{Shaik}, \binits{O.S.}},
\bauthor{\bsnm{Kummer}, \binits{U.}},
\bauthor{\bsnm{Lebiedz}, \binits{D.}}:
\batitle{Optimal control of self-organized dynamics in cellular signal
  transduction}.
\bjtitle{Mathematical and Computer Modelling of Dynamical Systems}
\bvolume{13}(\bissue{5}),
\bfpage{487}--\blpage{502}
(\byear{2007})
\doiurl{10.1080/13873950701243969}
\end{barticle}
\endbibitem

\bibitem[\protect\citeauthoryear{Suri}{2008}]{Suri2008Harvesting}
\begin{botherref}
\oauthor{\bsnm{Suri}, \binits{R.}}:
Optimal harvesting strategies for fisheries: A differential equations approach.
PhD thesis,
Massey University
(2008).
Two-player differential game formulated on p. 168.
\url{https://mro.massey.ac.nz/server/api/core/bitstreams/6dc7ad99-1d4f-4b56-a43e-b4965827292f/content}
\end{botherref}
\endbibitem

\bibitem[\protect\citeauthoryear{Tam et~al.}{2015}]{BenchmarkPMC4440705}
\begin{barticle}
\bauthor{\bsnm{Tam}, \binits{Z.Y.}},
\bauthor{\bsnm{Gruber}, \binits{J.}},
\bauthor{\bsnm{Halliwell}, \binits{B.}},
\bauthor{\bsnm{Gunawan}, \binits{R.}}:
\batitle{Context-dependent role of mitochondrial fusion-fission in clonal
  expansion of mtdna mutations}.
\bjtitle{PLoS computational biology}
\bvolume{11}(\bissue{5}),
\bfpage{1004183}
(\byear{2015})
\doiurl{10.1371/journal.pcbi.1004183}
\end{barticle}
\endbibitem

\bibitem[\protect\citeauthoryear{Tu et~al.}{2005}]{Tu2005MetabolicCycle}
\begin{barticle}
\bauthor{\bsnm{Tu}, \binits{B.P.}},
\bauthor{\bsnm{Kudlicki}, \binits{A.}},
\bauthor{\bsnm{Rowicka}, \binits{M.}},
\bauthor{\bsnm{McKnight}, \binits{S.L.}}:
\batitle{Logic of the yeast metabolic cycle: temporal compartmentalization of
  cellular processes}.
\bjtitle{Science}
\bvolume{310}(\bissue{5751}),
\bfpage{1152}--\blpage{1158}
(\byear{2005})
\doiurl{10.1126/science.1120499}
\end{barticle}
\endbibitem

\bibitem[\protect\citeauthoryear{Tsoi et~al.}{2018}]{BenchmarkPMC5877992}
\begin{barticle}
\bauthor{\bsnm{Tsoi}, \binits{R.}},
\bauthor{\bsnm{Wu}, \binits{F.}},
\bauthor{\bsnm{Zhang}, \binits{C.}},
\bauthor{\bsnm{Bewick}, \binits{S.}},
\bauthor{\bsnm{Karig}, \binits{D.}},
\bauthor{\bsnm{You}, \binits{L.}}:
\batitle{Metabolic division of labor in microbial systems}.
\bjtitle{Proceedings of the National Academy of Sciences of the United States
  of America}
\bvolume{115}(\bissue{10}),
\bfpage{2526}--\blpage{2531}
(\byear{2018})
\doiurl{10.1073/pnas.1716888115}
\end{barticle}
\endbibitem

\bibitem[\protect\citeauthoryear{Verschueren et~al.}{2022}]{Verschueren2022}
\begin{barticle}
\bauthor{\bsnm{Verschueren}, \binits{R.}},
\bauthor{\bsnm{Frison}, \binits{G.}},
\bauthor{\bsnm{Kouzoupis}, \binits{D.}},
\bauthor{\bsnm{Frey}, \binits{J.}},
\bauthor{\bsnm{Duijkeren}, \binits{N.}},
\bauthor{\bsnm{Zanelli}, \binits{A.}},
\bauthor{\bsnm{Novoselnik}, \binits{B.}},
\bauthor{\bsnm{Albin}, \binits{T.}},
\bauthor{\bsnm{Quirynen}, \binits{R.}},
\bauthor{\bsnm{Diehl}, \binits{M.}}:
\batitle{{acados}---a modular open-source framework for fast embedded optimal
  control}.
\bjtitle{Mathematical Programming Computation}
\bvolume{14}(\bissue{1}),
\bfpage{147}--\blpage{183}
(\byear{2022})
\doiurl{10.1007/s12532-021-00208-8}
\end{barticle}
\endbibitem

\bibitem[\protect\citeauthoryear{Wilson and Cowan}{1972}]{WilsonCowan1972}
\begin{barticle}
\bauthor{\bsnm{Wilson}, \binits{H.R.}},
\bauthor{\bsnm{Cowan}, \binits{J.D.}}:
\batitle{Excitatory and inhibitory interactions in localized populations of
  model neurons}.
\bjtitle{Biophysical Journal}
\bvolume{12}(\bissue{1}),
\bfpage{1}--\blpage{24}
(\byear{1972})
\doiurl{10.1016/S0006-3495(72)86068-5}
\end{barticle}
\endbibitem

\bibitem[\protect\citeauthoryear{Witzig et~al.}{2020}]{MitochondrialClock2020}
\begin{barticle}
\bauthor{\bsnm{Witzig}, \binits{M.}},
\bauthor{\bsnm{Grimm}, \binits{A.}},
\bauthor{\bsnm{Schmitt}, \binits{K.}},
\bauthor{\bsnm{Lejri}, \binits{I.}},
\bauthor{\bsnm{Frank}, \binits{S.}},
\bauthor{\bsnm{Brown}, \binits{S.A.}},
\bauthor{\bsnm{Eckert}, \binits{A.}}:
\batitle{Clock-controlled mitochondrial dynamics correlates with cyclic
  pregnenolone synthesis}.
\bjtitle{Cells}
\bvolume{9}(\bissue{10}),
\bfpage{2323}
(\byear{2020})
\doiurl{10.3390/cells9102323}
\end{barticle}
\endbibitem

\bibitem[\protect\citeauthoryear{Xavier and Foster}{2007}]{BenchmarkPMC1783407}
\begin{barticle}
\bauthor{\bsnm{Xavier}, \binits{J.B.}},
\bauthor{\bsnm{Foster}, \binits{K.R.}}:
\batitle{Cooperation and conflict in microbial biofilms}.
\bjtitle{Proceedings of the National Academy of Sciences of the United States
  of America}
\bvolume{104}(\bissue{3}),
\bfpage{876}--\blpage{881}
(\byear{2007})
\doiurl{10.1073/pnas.0607651104}
\end{barticle}
\endbibitem

\bibitem[\protect\citeauthoryear{Yoshida et~al.}{2007}]{Yoshida2007Cryptic}
\begin{barticle}
\bauthor{\bsnm{Yoshida}, \binits{T.}},
\bauthor{\bsnm{Ellner}, \binits{S.P.}},
\bauthor{\bsnm{Jones}, \binits{L.E.}},
\bauthor{\bsnm{Bohannan}, \binits{B.J.M.}},
\bauthor{\bsnm{Lenski}, \binits{R.E.}},
\bauthor{\bsnm{Hairston}, \binits{J.} \bsuffix{Nelson~G.}}:
\batitle{Cryptic population dynamics: Rapid evolution masks trophic
  interactions}.
\bjtitle{PLoS Biology}
\bvolume{5}(\bissue{9}),
\bfpage{235}
(\byear{2007})
\doiurl{10.1371/journal.pbio.0050235}
\end{barticle}
\endbibitem

\bibitem[\protect\citeauthoryear{Zhu and Borrelli}{2023}]{ZhuBorrelli2023}
\begin{bchapter}
\bauthor{\bsnm{Zhu}, \binits{E.L.}},
\bauthor{\bsnm{Borrelli}, \binits{F.}}:
\bctitle{A sequential quadratic programming approach to the solution of
  open-loop generalized nash equilibria}.
In: \bbtitle{2023 IEEE International Conference on Robotics and Automation},
pp. \bfpage{3211}--\blpage{3217}.
\bpublisher{IEEE},
\blocation{London, United Kingdom}
(\byear{2023}).
\doiurl{10.1109/ICRA48891.2023.10160799} .
\burl{https://arxiv.org/abs/2203.16478}
\end{bchapter}
\endbibitem

\bibitem[\protect\citeauthoryear{Zhou et~al.}{2022}]{Zhou2022Algae}
\begin{barticle}
\bauthor{\bsnm{Zhou}, \binits{Q.S.}},
\bauthor{\bsnm{Gao}, \binits{Y.}},
\bauthor{\bsnm{Hou}, \binits{J.M.}},
\bauthor{\bsnm{Wang}, \binits{T.}},
\bauthor{\bsnm{Tang}, \binits{L.}}:
\batitle{Preference of carbon absorption determines the competitive ability of
  algae along atmospheric {CO2} concentration}.
\bjtitle{Ecology and Evolution}
\bvolume{12}(\bissue{7}),
\bfpage{9079}
(\byear{2022})
\doiurl{10.1002/ece3.9079}
\end{barticle}
\endbibitem

\end{thebibliography}
\end{document}